\documentclass[12pt]{amsart}
\usepackage{amscd}
\usepackage[dvips]{graphicx}
\usepackage{amscd}
\usepackage{latexsym}
\usepackage{amssymb}

\newtheorem{theorem}{Theorem}[section]
\newtheorem{proposition}[theorem]{Proposition}
\newtheorem{prop}[theorem]{Proposition}
\newtheorem{lemma}[theorem]{Lemma}
\newtheorem{corollary}[theorem]{Corollary}

\newtheorem{remark}[theorem]{Remark}
\newtheorem{definition}[theorem]{Definition}
\newtheorem{notation}[theorem]{Notation}

\newtheorem{example}[theorem]{Example}

\newcommand{\prf}{{\it Proof.} }

\newcommand{\qedd}{\hspace*{\fill}\hbox{$\Box$}}

\newcommand{\gl}[2]{\mbox{${\rm GL}_{#1}(#2)$}}

\newcommand{\tr}{\mbox{${\rm tr}$}}

\def\Spec{\mathop{\rm Spec}}

\newcommand{\cal}{\mathcal}
\newcommand{\mat}[2]{\mbox{${\rm M}_{#1}(#2)$}} 
\begin{document}
\begin{center}{\large\bf The moduli of representations 
of degree $2$} \\
\end{center}

\begin{center}{Kazunori NAKAMOTO\footnote{The author  
was partially supported by Grant-in-Aid in Scientific Research (C) (No.~23540044 
and No.~15K04814) 
from JSPS.} \\ \medskip 
Center for Medical Education and Sciences, Faculty of Medicine \\ 
University of Yamanashi \\  nakamoto@yamanashi.ac.jp} 
\end{center}

\bigskip 

\begin{quote}\scriptsize 
{\bf Abstract.} 
There are $6$ types of $2$-dimensional representations in general. 
For any groups and any monoids,  
we can construct the moduli of $2$-dimensional representations for each type: 
the moduli of absolutely irreducible representations, 
representations with Borel mold, representations with semi-simple mold, 
representations with unipotent mold, 
representations with unipotent mold over ${\Bbb F}_2$, 
and representations with scalar mold.  We can also construct them  
for any associative algebras. 
\end{quote} 

\begin{quote}\scriptsize 
\textup{2010} {\it Mathematics Subject Classification}. 
Primary 14D22; Secondary 14D20, 20M30, 20C99, 16G99.
\end{quote}

\begin{quote}\scriptsize 
\end{quote}

%
%

\section{Introduction}

In this paper we deal with the moduli of representations of degree $2$. 
We can classify $2$-dimensional representations into 
$6$ types in general. 
For any groups and any monoids,  
we can construct the moduli of $2$-dimensional representations for each type. 
For any associative algebras, we can also construct them for each type. 

In \cite{Borel} we have introduced the notion of mold. 
A mold is, so to say, a subalgebra of the full matrix ring. 
More precisely, a subsheaf of ${\mathcal O}_X$-algebras 
${\mathcal A} \subseteq {\rm M}_n({\mathcal O}_X)$ on a scheme $X$ 
is called a {\it mold} if ${\mathcal A}$ is a subbundle of 
${\rm M}_n({\mathcal O}_X)$. 
Let $\Gamma$ be a group or a monoid. 
By a homomorphism $\rho : \Gamma \to {\rm M}_n(\Gamma(X, {\mathcal O}_X))$, 
we understand an {\it $n$-dimensional representation 
of $\Gamma$ on a scheme $X$.}  
We say that a representation $\rho$ {\it has a mold ${\mathcal A}$} if 
the subsheaf of ${\mathcal O}_X$-algebras 
${\mathcal O}_X[\rho(\Gamma)]$ of ${\rm M}_n({\mathcal O}_X)$ 
generated by $\rho(\Gamma)$ coincides with ${\mathcal A}$. 
It is effective to classify representations with respect to 
molds for constructing the moduli of equivalence classes of representations. 
If we try to construct the moduli of equivalence classes of all representations 
without classifying representations with respect to molds, then 
two representations which have the same composition factors coincide 
as points of the moduli even if they are not equivalent.  
For separating such representations in the moduli, 
we need to collect only representations which have the same mold. 
For example, we have constructed the moduli  
of equivalence classes of 
absolutely irreducible representations denoted by 
${\rm Ch}_n(\Gamma)_{\rm air}$ in \cite{Nkmt00}, where $\rho : \Gamma 
\to {\rm M}_n(\Gamma(X, {\mathcal O}_X))$ is 
{\it absolutely irreducible} if ${\mathcal O}_X[\rho(\Gamma)]$ coincides with 
the full matrix ring ${\rm M}_n({\mathcal O}_X)$. 
We also have constructed the 
moduli of equivalence classes of 
representations with Borel mold denoted by 
${\rm Ch}_n(\Gamma)_{B}$ in \cite{Borel}, where  
$\rho : \Gamma 
\to {\rm M}_n(\Gamma(X, {\mathcal O}_X))$ is a {\it representation with 
Borel mold} if for each $x \in X$ there exists $P \in {\rm GL}_n({\mathcal O}_X(U))$ 
on a neighbourhood $U$ of $x$ such that 
$P \cdot {\mathcal O}_U[\rho(\Gamma)] \cdot P^{-1}$ coincides with  
the subsheaf of ${\mathcal O}_U$-algebras of ${\rm M}_n({\mathcal O}_U)$  
consisting of upper triangular matrices.  
The author calls the plan to construct the moduli of equivalence classes of 
representations for any 
suitable molds ``mold program''. 
In this article, we will complete the mold program of degree $2$.

For $k$-subalgebras $A$ and $B$ of the full matrix ring 
${\rm M}_2(k)$ over an algebraically 
closed field $k$, we say that $A$ and $B$ 
are {\it equivalent} if there exists $P \in {\rm GL}_2(k)$ such that 
$P^{-1}AP = B$.  
There are $5$ equivalence classes of $k$-subalgebras $A$ of 
${\rm M}_2(k)$:  
(1) $A = {\rm M}_2(k)$, 
(2) $A = \left\{ 
\left( \begin{array}{cc}
\ast & \ast \\
0 & \ast 
\end{array} 
\right) \right\}$, 
(3) $A = \left\{ 
\left( \begin{array}{cc}
\ast & 0 \\
0 & \ast 
\end{array} 
\right) \right\}$,
(4) $A = \left\{ \begin{array}{c|c}  
\left( \begin{array}{cc}
a & b \\
0 & a
\end{array} 
\right) & a, b \in k  
\end{array} \right\}$, 
(5) $A = \left\{ \begin{array}{c|c} 
\left( \begin{array}{cc}
a & 0 \\
0 & a
\end{array} 
\right) & a \in k 
\end{array} \right\}$. 
Let $\rho : \Gamma \to 
{\rm M}_2(k)$ be a $2$-dimensional representations 
of a group or a monoid $\Gamma$.   
By equivalence classes of the subalgebra $k[\rho(\Gamma)]$ 
of ${\rm M}_2(k)$, we classify $2$-dimensional representations into 
$6$-types (not $5$-types!). 
For each cases (1)--(5), we say that $\rho$ is 
(1) an {\it absolutely irreducible representation},  
(2) a {\it representation with Borel mold}, 
(3) a {\it representation with semi-simple mold}, 
(4) a {\it representation with unipotent mold},   
(5) a {\it representation with scalar mold}, 
respectively. 
In the case (4), we need to divide representations with unipotent mold into 
$2$ types: (4-a) when ${\rm ch} k \neq 2$, we say $\rho$ is a 
{\it representation with unipotent mold}, and 
(4-b) when ${\rm ch} k = 2$, we say $\rho$ is a 
{\it representation with unipotent mold 
over ${\Bbb F}_2$}. 
It is natural to divide the case (4) into $2$ types for 
constructing the``good'' moduli of representations with unipotent mold. 
Here, by constructing the ``good'' moduli of representations, we understand 
constructing smooth moduli schemes of representations at least for 
free monoids (more precisely, see the beginning of \S 5). 
Hence there are $6$ types of $2$-dimensional representations in general. 

In \S 3, we introduce the notions of (1), (2), (3), (4-a), (4-b), (5) on $2$-dimensional 
representations on arbitrary schemes $X$ (Definitions~\ref{def:airandborel},  \ref{def:ssmold}, \ref{def:umold}, \ref{def:uf2mold}, and \ref{def:repwith}). 
For $2$-dimensional representations $\rho_1, \rho_2$ on $X$, 
we say that $\rho_1$, $\rho_2$ are {\it equivalent} 
(or $\rho_1 \sim \rho_2$) if  
there exists a $\Gamma(X, {\mathcal O}_X)$-algebra isomorphism 
$\sigma : {\rm M}_2(\Gamma(X, {\mathcal O}_X)) \to 
{\rm M}_2(\Gamma(X, {\mathcal O}_X))$ such that 
$\sigma(\rho_1(\gamma)) = \rho_2(\gamma)$ 
for any $\gamma \in \Gamma$.
If $\rho_1 \sim \rho_2$, then 
for each $x \in X$ there exists $P \in {\rm GL}_2({\mathcal O}_X(U))$ 
on a neighbourhood $U$ of $x$ such that $P^{-1}\rho_1(\gamma) P = \rho_2(\gamma)$ 
on $U$ for any $\gamma \in \Gamma$. We have constructed the moduli of 
equivalence classes of representations in the cases (1) absolutely irreducible 
representations and (2) representations with Borel mold in 
\cite{Nkmt00} and \cite{Borel}, respectively. 
In the case (5) representations with scalar mold, we can easily construct the moduli  (Theorem~\ref{th:moduliscalar}).  
In the cases (3) representations with semi-simple mold, (4-a) representations 
with unipotent mold, and (4-b) representations with unipotent mold over ${\Bbb F}_2$, 
we have the following theorems: 

\begin{theorem}[Theorem~\ref{th:moduliss}]  
There exists a fine moduli scheme 
${\rm Ch}_2(\Gamma)_{\rm s.s.}$ 
associated to the sheafification 
${\mathcal E}q\mathcal{SS}_2(\Gamma)$ of the functor
\[
\begin{array}{ccl}
 ({\bf Sch})^{op} & \to & ({\bf Sets}) \\
 X & \mapsto & 
\left\{ 
\begin{array}{l}  
\mbox{ $2$-dimensional representations } \\ 
\mbox{ with semi-simple mold of $\Gamma$ on } X  
\end{array} 
\right\} 
\Big/\sim
\end{array}
\] with respect to Zariski topology 
for arbitrary group or monoid $\Gamma$. 
The moduli ${\rm Ch}_2(\Gamma)_{\rm s.s.}$ is separated over ${\Bbb Z}$;
if $\Gamma$ is a finitely generated group or monoid, then 
${\rm Ch}_2(\Gamma)_{\rm s.s.}$ is of finite type over ${\Bbb Z}$. 
\end{theorem}

\begin{theorem}[Theorem~\ref{th:moduliu}]  
There exists a fine moduli scheme ${\rm Ch}_2(\Gamma)_{u}$    
associated to the sheafification ${\mathcal E}q \:\mathcal{U}_2(\Gamma)$ of the 
functor 
\[
\begin{array}{ccl}
 ({\bf Sch}/{\Bbb Z}[1/2])^{op} & \to & ({\bf Sets}) \\
 X & \mapsto & 
\left\{ 
\begin{array}{l}  
\mbox{ $2$-dimensional representations } \\ 
\mbox{ with unipotent mold of $\Gamma$ on } X  
\end{array} 
\right\} 
\Big/\sim
\end{array}
\] with respect to Zariski topology 
for arbitrary group or monoid $\Gamma$.   
The moduli ${\rm Ch}_2(\Gamma)_{u}$ is separated over ${\Bbb Z}[1/2]$; 
if $\Gamma$ is a finitely generated group or monoid, then 
${\rm Ch}_2(\Gamma)_{u}$ is of finite type over ${\Bbb Z}[1/2]$. 
\end{theorem}

\begin{theorem}[Theorem~\ref{th:moduliuf2}]   
There exists a fine moduli scheme ${\rm Ch}_2(\Gamma)_{u/{\Bbb F}_2}$ 
associated to the sheafification ${\mathcal E}q \:\mathcal{U}_2(\Gamma)_{{\Bbb F}_2}$ 
of the functor   
\[
\begin{array}{ccl}
 ({\bf Sch}/{\Bbb F}_2)^{op} & \to & ({\bf Sets}) \\
 X & \mapsto & 
\left\{ 
\begin{array}{l}  
\mbox{ $2$-dimensional representations with } \\ 
\mbox{ unipotent mold over ${\Bbb F}_2$ of $\Gamma$ on } X  
\end{array} 
\right\} 
\Big/\sim
\end{array}
\] 
with respect to Zariski topology 
for arbitrary group or monoid $\Gamma$. 
The moduli ${\rm Ch}_2(\Gamma)_{u/{\Bbb F}_2}$ is separated over ${\Bbb F}_2$; 
if $\Gamma$ is a finitely generated group or monoid, then 
${\rm Ch}_2(\Gamma)_{u/{\Bbb F}_2}$ is of finite type over ${\Bbb F}_2$. 
\end{theorem}

For any associative algebra $A$ over any commutative ring $R$, 
we also obtain the same theorems on $2$-dimensional representations of $A$ over 
$R$:  There exist fine moduli schemes 
${\rm Ch}_2(A)_{\rm s.s.}$, ${\rm Ch}_2(A)_{u}$, 
and ${\rm Ch}_2(A)_{u/{\Bbb F}_2}$ separated over $R$. 
If $A$ is finitely generated associative algebra over $R$, then 
the fine moduli schemes ${\rm Ch}_2(A)_{\rm s.s.}$, ${\rm Ch}_2(A)_{u}$, 
and ${\rm Ch}_2(A)_{u/{\Bbb F}_2}$ are of finite type over $R$ (Remarks~\ref{remark:algss}, \ref{remark:algu1}, \ref{remark:constructu}, \ref{remark:alguf2}, and \ref{remark:constructuf2}). 
These theorems are main results of this article. 

\bigskip 

As a continuation of this article, we can deal with 
the absolutely irreducible representations parts of the representation 
variety and the character variety: ${\rm Rep}_2(\Gamma)_{\rm air}$ and 
${\rm Ch}_2(\Gamma)_{\rm air}$ in \cite{Nkmt00}. 
For a group or a monoid $\Gamma$, the 
representation variety ${\rm Rep}_2(\Gamma)$ is 
the affine scheme representing the contravariant functor 
which maps each scheme $X$ to the set of $2$-dimensional representations 
of $\Gamma$ on $X$. 
For $\ast = {\rm air}$, $B$, ${\rm s.s.}$, $u$, $u/{\Bbb F}_2$, or ${\rm scalar}$, 
${\rm Rep}_2(\Gamma)_{\ast}$ denotes the 
subscheme of ${\rm Rep}_2(\Gamma)$ consisting of 
$2$-dimensional representations with the mold corresponding 
to $\ast$. 
For a field $k$, the set of $k$-rational points of   
the representation variety ${\rm Rep}_2(\Gamma)$ is the disjoint union of 
the sets of $k$-rational points of 
${\rm Rep}_2(\Gamma)_{\rm air}$, ${\rm Rep}_2(\Gamma)_{B}$, 
${\rm Rep}_2(\Gamma)_{\rm s.s.}$, ${\rm Rep}_2(\Gamma)_{u}$ (or 
${\rm Rep}_2(\Gamma)_{u/{\Bbb F}_2}$), and ${\rm Rep}_2(\Gamma)_{\rm scalar}$. 
Hence for a finitely generated group or monoid $\Gamma$ 
and for the finite field ${\Bbb F}_q$, the number of ${\Bbb F}_q$-rational 
points of ${\rm Rep}_2(\Gamma)_{\rm air}$ can be 
calculated from those of ${\rm Rep}_2(\Gamma)$ and 
the others ${\rm Rep}_2(\Gamma)_{\ast}$.  
Since ${\rm Rep}_2(\Gamma)_{\rm air} \to {\rm Ch}_2(\Gamma)_{\rm air}$ 
is a ${\rm PGL}_2$-principal fibre bundle, 
the number of ${\Bbb F}_q$-rational points of ${\rm Ch}_2(\Gamma)_{\rm air}$ 
can be also calculated from the result of ${\rm Rep}_2(\Gamma)_{\rm air}$.   
Similarly, the virtual Hodge polynomials of ${\rm Rep}_2(\Gamma)_{\rm air}$ 
and ${\rm Ch}_2(\Gamma)_{\rm air}$ over ${\Bbb C}$ 
can be calculated from those of ${\rm Rep}_2(\Gamma)$ and 
the others ${\rm Rep}_2(\Gamma)_{\ast}$ over ${\Bbb C}$. 
The existence of such geometric objects 
as the moduli of representations with several molds 
helps us to understand relations between the numbers of 
equivalence classes of representations of $\Gamma$ over ${\Bbb F}_q$ and virtual Hodge polynomials of the moduli ({\it cf.} \cite{topos2}).

In \cite{topos2}, the authors deal with 
the case that $\Gamma$ is the free monoid $\Upsilon_{m}$ of rank $m$. 
Since ${\rm Rep}_2(\Upsilon_m)$ is isomorphic to 
${\rm M}_2 \times \cdots \times {\rm M}_2$ ($m$ times),  
the numbers of the ${\Bbb F}_q$-rational points of 
${\rm Rep}_2(\Upsilon_m)_{\rm air}$ and 
${\rm Ch}_2(\Upsilon_m)_{\rm air}$ have been calculated explicitly.  
(The author needs to mention that 
our strategy to calculate the numbers of the 
${\Bbb F}_q$-rational points is essentially same as \cite{SL2-ch} and \cite{SL3-ch}.  
Moreover, the method of \cite{Reineke} is much easier than 
our strategy.)  
We have also calculated the virtual Hodge polynomials of 
${\rm Rep}_2(\Upsilon_m)_{\rm air}$ and 
${\rm Ch}_2(\Upsilon_m)_{\rm air}$. 
We see that the Hasse-Weil zeta functions of 
${\rm Rep}_2(\Upsilon_m)_{\rm air}$ and 
${\rm Ch}_2(\Upsilon_m)_{\rm air}$ satisfy functional equations.

\bigskip 

The organization of this article is as follows: 
in \S 2, we review representations and molds on schemes. 
We also review the moduli of absolutely irreducible representations 
and the moduli of representations with Borel mold. 
In \S 3, we introduce several molds of degree $2$: 
semi-simple mold, unipotent mold, unipotent mold over ${\Bbb F}_2$, and 
scalar mold. We also introduce the moduli of representations with 
scalar mold. In \S 4, we construct the moduli of 
equivalence classes of representations with semi-simple mold. 
In \S 5, we construct the moduli of 
equivalence classes of representations with unipotent mold over 
${\Bbb Z}[1/2]$.  
In \S 6, we construct the moduli of 
equivalence classes of representations with unipotent mold over ${\Bbb F}_2$. 
In \S 7, 
we deal with different approach from \S 6. 
The approach in \S 7 gives us another construction of the moduli of 
equivalence classes of representations with unipotent mold over ${\Bbb F}_2$ 
by using derivations as in \S 5. 
In \S 8, 
we reformulate the moduli functors   
by using the notion of representations generating 
sheaves of ${\mathcal O}_X$-algebras which define molds of rank $2$.  
In \S 9, we deal with discriminants which describe the absolutely 
irreducible representation part ${\rm Rep}_2(\Gamma)_{\rm air}$ 
in the representation variety ${\rm Rep}_2(\Gamma)$ 
as an appendix. 

\bigskip 

The author would like to thank Takeshi Torii for his essential ideas and 
important suggestions on the moduli of representations. 
Although his name does not appear in the list of the authors, his 
contribution to this paper is not ignorable. 
This article has been inspired by his 
descriptions of the moduli of representations and 
approaches from viewpoints of algebraic topology, and so on,     
which will be written in \cite{topos2}.  

The author would like to express his gratitude to the referee for  
suggesting several important points. 
Example \ref{ex:describe-chuf2m}, \S 7, and \S 8 have been inspired by the referee. 
The author also wants to thank Michiaki Inaba. 
He suggested the proof of Lemma \ref{lemma:univgeomdescent}, which states 
the ``descent'' of universal geometric quotients.

\section{Preliminaries} 

In this section, we review representations and molds on schemes. 
(For details, see \cite{Nkmt00} and \cite{Borel}.)

\begin{definition}[\cite{Nkmt00}]\rm
Let $\Gamma$ be a group or a monoid. 
By a {\it representation} of $\Gamma$ on a scheme $X$, we understand 
a group homomorphism (or a monoid homomorphism) 
$\rho : \Gamma \to {\rm M}_n(\Gamma(X, {\mathcal O}_X))$. 
For two representations $\rho$ and $\rho'$, we say that 
$\rho$ and $\rho'$ are {\it equivalent} to each other (or 
$\rho \sim \rho'$) if 
there exists an ${\mathcal O}_X$-algebra isomorphism $\sigma : 
{\rm M}_n(\Gamma(X, {\mathcal O}_X)) \to {\rm M}_n(\Gamma(X, {\mathcal O}_X))$
such that $\sigma(\rho(\gamma)) = \rho'(\gamma)$ for each $\gamma \in \Gamma$. 
\end{definition}

\begin{remark}[\cite{Nkmt00}]\rm 
Let $\rho$ and $\rho'$ be 
$n$-dimensional representations 
of $\Gamma$ on $X$.  
If $\rho \sim \rho'$, then 
for each $x \in X$ there exists $P \in {\rm GL}_n({\mathcal O}_X(U))$ 
on a neighbourhood $U$ of $x$ such that $P^{-1}\rho(\gamma) P = \rho'(\gamma)$ 
on $U$ for any $\gamma \in \Gamma$. 
Indeed, the group scheme ${\rm PGL}_n$ over ${\Bbb Z}$ 
represents the functor 
\[ \begin{array}{ccl}
({\bf Sch})^{op} & \to & ({\bf Sets}) \\
X & \mapsto & {\rm Aut}_{{\mathcal O}_X\mbox{-alg}} ({\rm M}_n({\mathcal O}_X)) . 
\end{array} 
\] 
For details, see \cite[Definition~6.1 and Theorem~6.2]{Nkmt00}. 
\end{remark} 

\begin{definition}[\cite{Nkmt00}]\rm
Let $\Gamma$ be a group or a monoid. The following contravariant functor 
is representable by an affine scheme: 
\[
\begin{array}{ccccl}
{\rm Rep}_n(\Gamma) & : & ({\bf Sch})^{op} & \to & ({\bf Sets}) \\
 & & X & \mapsto & \{ \rho : \text{ rep. of deg $n$ for $\Gamma$ on $X$ } \}. 
\end{array}
\]
We call the affine scheme ${\rm Rep}_n(\Gamma)$ the {\it representation variety}  
of degree $n$ for $\Gamma$. The group scheme ${\rm PGL}_n$ over ${\mathbb Z}$ acts on 
${\rm Rep}_n(\Gamma)$ by $\rho \mapsto P^{-1}\rho P$. 
Each ${\rm PGL}_n$-orbit forms an equivalence class of representations. 
If $\Gamma$ is a finitely generated group (or monoid), then 
${\rm Rep}_n(\Gamma)$ is of finite type over ${\mathbb Z}$. 
\end{definition}

\begin{definition}[\cite{Borel}]\rm 
Let ${\mathcal A}$ be a subsheaf of ${\rm M}_n({\mathcal O}_X)$ of 
${\mathcal O}_X$-algebras on a scheme $X$. 
We say that ${\mathcal A}$ is a {\it mold} on $X$ if 
${\rm M}_n({\mathcal O}_X)/{\mathcal A}$ is locally free. 
Let ${\rm rank}{\mathcal A}$ denote the rank of a mold ${\mathcal A}$ 
as a locally free sheaf. 
For two molds ${\mathcal A}, {\mathcal B} \subseteq {\rm M}_n({\mathcal O}_X)$ on $X$, 
we say that ${\mathcal A}$ and ${\mathcal B}$ are {\it locally equivalent} 
if there exist an open covering $X = \cup_{i \in I} U_i$ 
and $P_i \in {\rm GL}_n({\mathcal O}_X(U_i))$ such that 
$P_i ({\mathcal A}\mid_{U_i}) P_i^{-1} = {\mathcal B}\mid_{U_i}$ for 
each $i \in I$. 
\end{definition}

Here let us introduce an example of molds. 

\begin{example}[\cite{Borel}]\rm 
We define the mold ${\mathcal B}_n$ on ${\rm Spec}\;{\Bbb Z}$ by 
\[ {\mathcal B}_n := \{ (a_{ij}) \in {\rm M}_n({\Bbb Z}) \mid a_{ij} = 0 
\mbox{ for each $i>j$ } \}. \]
For a mold ${\mathcal A} \subseteq {\rm M}_n({\mathcal O}_X)$ 
on a scheme $X$, we say that ${\mathcal A}$ 
is a {\it Borel mold} if ${\mathcal A}$ and 
${\mathcal B}_n\otimes_{\Bbb Z}{\mathcal O}_X$ are locally equivalent 
to each other. 
\end{example}

\begin{definition}[\cite{Borel}]\rm
Let ${\mathcal A}$ be a mold on a scheme $X$. 
For a representation $\rho$ of $\Gamma$ on $X$, we say that 
$\rho$ has {\it mold type} ${\mathcal A}$ if 
the image $\rho(\Gamma)$ generates ${\mathcal A}$ as an 
${\mathcal O}_X$-algebra. 
\end{definition}

\begin{definition}[\cite{Borel}]\label{def:airandborel}\rm
Let $\rho$ be an $n$-dimensional representation of $\Gamma$ 
on a scheme $X$. We say that $\rho$ is an {\it absolutely irreducible representation}  
(or {\it air}) if $\rho$ has mold type ${\rm M}_n({\mathcal O}_X)$. 
We also say that $\rho$ is a {\it representation with Borel mold} if 
$\rho$ has a Borel mold type. 
\end{definition}

\begin{proposition}[\cite{Nkmt00}, \cite{Borel}]
The contravariant functor 
\[
\begin{array}{ccccl}
{\rm Rep}_n(\Gamma)_{\rm air} & : & ({\bf Sch})^{op} & \to & ({\bf Sets}) \\
 & & X & \mapsto & \{ \text{ air of degree $n$ for $\Gamma$ on $X$ } \} 
\end{array}
\]
is representable by an open subscheme ${\rm Rep}_n(\Gamma)_{\rm air}$ of 
${\rm Rep}_n(\Gamma)$. 
The contravariant functor 
\[
\begin{array}{ccccl}
{\rm Rep}_n(\Gamma)_{B} & : & ({\bf Sch})^{op} & \to & ({\bf Sets}) \\
 & & X & \mapsto & 
\left\{ 
\begin{array}{l}
 \text{ rep. with Borel mold of degree $n$ } \\ 
 \text{ for $\Gamma$ on $X$  } 
\end{array} 
\right\} 
\end{array}
\]
is representable by a subscheme ${\rm Rep}_n(\Gamma)_{B}$ of 
${\rm Rep}_n(\Gamma)$. 
The action of ${\rm PGL}_n$ on ${\rm Rep}_n(\Gamma)$ induces 
the ones of ${\rm PGL}_n$ on ${\rm Rep}_n(\Gamma)_{\rm air}$ and ${\rm Rep}_n(\Gamma)_{B}$. 
\end{proposition}

\medskip 

For absolutely irreducible representations, there exists a coarse moduli scheme. 

\begin{theorem}[\cite{Nkmt00}]
There exists a coarse moduli scheme ${\rm Ch}_n(\Gamma)_{\rm air}$ 
separated over ${\Bbb Z}$ associated to 
the following functor:
\[
\begin{array}{ccccl}
\mathcal{E}q\mathcal{AIR}_n(\Gamma) & : & ({\bf Sch})^{op} & \to & ({\bf Sets}) \\
  & & X & \mapsto & \{ \rho : \mbox{ air of degree $n$ for $\Gamma$ on } X \}/\sim.
\end{array}
\]
Furthermore, the canonical morphism ${\rm Rep}_n(\Gamma)_{\rm air} 
\to {\rm Ch}_n(\Gamma)_{\rm air}$ gives a universal geometric quotient 
of ${\rm Rep}_n(\Gamma)_{\rm air}$ by ${\rm PGL}_n$.  
If $\Gamma$ is a finitely generated group (or monoid), then 
the moduli ${\rm Ch}_n(\Gamma)_{\rm air}$ is of finite type over ${\Bbb Z}$. 
\end{theorem}

\medskip 

For representations with Borel mold, there exists a fine moduli scheme. 

\begin{theorem}[\cite{Borel}]
There exists a fine moduli scheme ${\rm Ch}_n(\Gamma)_{B}$ 
separated over ${\Bbb Z}$ associated to the 
sheafification $\mathcal{E}q\mathcal{B}_n(\Gamma)$  
of the following functor with respect to Zariski topology:   
\[
\begin{array}{ccl}
 ({\bf Sch})^{op} & \to & ({\bf Sets}) \\
   X & \mapsto & 
\left\{ 
\begin{array}{l} 
\text{ rep. with Borel mold } \\
\text{ of degree $n$ for $\Gamma$  }  
\end{array}
\right\} 
\Big/\sim.
\end{array}
\]  
Furthermore, the canonical morphism ${\rm Rep}_n(\Gamma)_{B} 
\to {\rm Ch}_n(\Gamma)_{B}$ gives a universal geometric quotient 
of ${\rm Rep}_n(\Gamma)_{B}$ by ${\rm PGL}_n$.  
If $\Gamma$ is a finitely generated group (or monoid), then 
the moduli is of finite type over ${\Bbb Z}$. 
\end{theorem}

%
%

\section{The degree $2$ case}

From now on, we deal mainly with the degree $2$ case. 

\bigskip 

Let $A_2(\Gamma)$ be the coordinate ring of the 
representation variety of degree $2$ for a group or a monoid 
$\Gamma$. Let $\sigma_{\Gamma} : \Gamma \to {\rm M}_2(A_2(\Gamma))$ be 
the universal representation of degree $2$ for $\Gamma$. 

\begin{definition}\rm
Let $A_2(\Gamma)^{\rm Ch}$ be the subalgebra of 
$A_2(\Gamma)$ generated by 
$\{ {\rm tr}(\sigma_{\Gamma}(\gamma)), 
\det(\sigma_{\Gamma}(\gamma)) \mid \gamma \in \Gamma \}$ 
over ${\mathbb Z}$. 
We denote ${\rm Spec} A_2(\Gamma)^{\rm Ch}$ by ${\rm Ch}_2(\Gamma)$. 
\end{definition} 

\bigskip 

In {\cite[Example 1.3]{Borel}} we investigated the moduli of molds:  
\begin{eqnarray*}
  {\rm Mold}_{2, 1}&  = & \Spec{{\Bbb Z}},  \\
  {\rm Mold}_{2, 2}&  = & {\Bbb P}^2_{\Bbb Z},  \\
  {\rm Mold}_{2, 3}&  = & {\Bbb P}^1_{\Bbb Z},  \\
  {\rm Mold}_{2, 4}&  = & \Spec{\Bbb Z}.  
\end{eqnarray*}
Let $k$ be an algebraically closed field, for simplicity.  
Let us classify $k$-subalgebras $A$ of ${\rm M}_2(k)$ up to 
inner automorphisms of ${\rm M}_2(k)$ for explaining molds of degree $2$. 
In the case $\dim A = 4$, $A$ is equal to ${\rm M}_2(k)$. 
For any subalgebra $A$ of dimension $3$, 
there exists $P \in {\rm GL}_2(k)$ such that
$P^{-1}AP = {\mathcal B}_2(k) := \{ (a_{ij}) \in {\rm M}_2(k) \mid a_{21}=0 \}$.
In the case $\dim A = 2$, 
there exists $X \in A$ such that $A = k I_2 + k X$. 
For $[X] \in {\rm M}_2(k)/kI_2$, we can define a mold $A = k I_2 + k X$,  
which is independent from 
choosing a representative $X \in {\rm M}_2(k)$ of $[X]$. 
This is the reason why ${\rm Mold}_{2, 2}(k)  \cong {\mathbb P}_{\ast}({\rm M}_2(k)/k I_2) = {\Bbb P}^2_{k}$.  
There exist two types of molds of rank $2$. 
The one is a semi-simple algebra, and the other is a non-semi-simple algebra. 
In other words, the former is 
$\left\{ \left. \left( 
\begin{array}{cc}
a & 0 \\
0 & b 
\end{array}
\right) \right|  a, b \in k \right\}$, and the 
latter is 
$\left\{ \left. \left( 
\begin{array}{cc}
a & b \\
0 & a 
\end{array}
\right) \right|  a, b \in k \right\}$ up to inner automorphisms. 
Of course, a subalgebra $A$ of dimension $1$ is equal to $k I_2$. 

\bigskip 

By using the classification of $k$-subalgebras of ${\rm M}_2(k)$, 
we introduce several molds of degree $2$. 
For the case of rank $4$, we consider the full matrix ring mold ${\rm M}_2({\mathcal O}_X)$.  
For the case of rank $3$, we introduced Borel molds. 

Here we introduce several types of molds of rank $2$. 
There are two types of molds of rank $2$: 
the semi-simple subalgebra case and the non-semi-simple subalgebra case. 
Moreover we can divide the non-semi-simple $2$-dimensional subalgebra case  
into two types: the ${\rm ch} \neq 2$ type and the ${\rm ch} =2$ type. 

\begin{notation}\rm
Let $R$ be a commutative ring. 
For $X \in {\rm M}_2(R)$, we denote  
${\rm tr}(X)^2 - 4\det(X)$ by $m(X)$. 
Remark that $m(X) = 2{\rm tr}(X^2) - ({\rm tr}(X))^2$. 
\end{notation}

\begin{remark}\label{remark:m}\rm
For $X \in {\rm M}_2(R)$, $m(X)$ is the discriminant of the characteristic polynomial 
of $X$. If $R$ is a field, then $m(X) \neq 0$ if and only if $X$ is semi-simple and not scalar. 
\end{remark}

\begin{definition}\label{def:ssmold}\rm 
Let $X$ be a scheme. Let 
${\mathcal A} \subseteq {\rm M}_2({\mathcal O}_X)$ be a rank $2$ mold
on $X$. 
We say that ${\mathcal A}$ is {\it semi-simple} 
if there exists $P_x \in {\mathcal A}_x$ such that 
$m(P_x) \not\equiv 0$ in the residue field $k(x)$ for each $x \in X$. 
\end{definition}

\begin{definition}\label{def:umold}\rm
Let $X$ be a scheme over ${\mathbb Z}[1/2]$. Let 
${\mathcal A} \subseteq {\rm M}_2({\mathcal O}_X)$ be a rank $2$ mold
on $X$. 
We say that ${\mathcal A}$ is {\it unipotent} if 
$m(A) = 0$ for 
each $A \in {\mathcal A}(U)$ and for each open set $U \subseteq X$.
\end{definition} 

\begin{definition}\label{def:uf2mold}\rm
Let $X$ be a scheme over ${\mathbb F}_2$. Let 
${\mathcal A} \subseteq {\rm M}_2({\mathcal O}_X)$ be a rank $2$ mold
on $X$. 
We say that ${\mathcal A}$ is {\it unipotent over ${\mathbb F}_2$} if 
${\rm tr}(A) = 0$ for 
each $A \in {\mathcal A}(U)$ and for each open set $U \subseteq X$.
\end{definition} 

\begin{remark}\rm
The name "unipotent" seems to be strange. 
However, the author calls non-semi-simple molds of rank $2$ 
unipotent molds because  
each unipotent mold over an algebraically closed field $k$ is generated by a unipotent matrix 
of ${\rm M}_2(k)$. 
\end{remark} 

\bigskip 

For each type of molds of rank $2$, we introduce representations with a given mold. 

\begin{definition}\label{def:repwith}\rm
For a $2$-dimensional representation $\rho$ for a group or a monoid $\Gamma$ on a scheme $X$, 
we say that $\rho$ is a {\it representation with semi-simple mold} if 
${\mathcal O}_X[\rho(\Gamma)]$ is a semi-simple mold on $X$.  
When $X$ is a scheme over ${\mathbb Z}[1/2]$ (or over ${\mathbb F}_2$), 
we say that $\rho$ is a {\it representation with unipotent mold} 
(or {\it unipotent mold over ${\mathbb F}_2$}) if 
${\mathcal O}_X[\rho(\Gamma)]$ is a unipotent mold (or a unipotent mold over 
${\mathbb F}_2$, respectively) 
on $X$. 
\end{definition} 

\bigskip 

For each case of molds of rank $2$, we construct the moduli of 
representations in \S 4-\S 6. 

\bigskip 

Finally, we consider molds of rank $1$. 
This case is trivial. Indeed, any mold of rank $1$ 
is the mold consisting of 
scalar matrices. 
Let us introduce the following definition for any degree. 

\begin{definition}\rm
Let $X$ be a scheme. 
We say that ${\mathcal A} \subseteq {\rm M}_n({\mathcal O}_X)$ 
is a {\it scalar mold} if ${\mathcal A}$ is a rank $1$ mold
on $X$. 
In other words, ${\mathcal A}$ is a scalar mold if and only if 
${\mathcal A} = {\mathcal O}_X\cdot I_n$. 
\end{definition}

\begin{definition}\rm
For an $n$-dimensional representation $\rho$ for a group or a monoid $\Gamma$ on a scheme $X$, 
we say that $\rho$ is a {\it representation with scalar mold} if 
${\mathcal O}_X[\rho(\Gamma)]$ is a scalar mold on $X$.  
\end{definition} 

\begin{proposition}\label{prop:repscalar} 
The contravariant functor 
\[
\begin{array}{ccccl}
{\rm Rep}_n(\Gamma)_{\rm scalar} & : & ({\bf Sch})^{op} & \to & ({\bf Sets}) \\
 & & X & \mapsto & 
\left\{ 
\begin{array}{l}
 \text{ rep. with scalar mold } \\ 
 \text{ of degree $n$ for $\Gamma$ on $X$  } 
\end{array} 
\right\} 
\end{array}
\]
is representable by a closed subscheme ${\rm Rep}_n(\Gamma)_{\rm scalar}$ of 
${\rm Rep}_n(\Gamma)$. 
The induced action of ${\rm PGL}_n$ on ${\rm Rep}_n(\Gamma)_{\rm scalar}$ is trivial.  
\end{proposition} 

\prf
Let $A_n(\Gamma)$ be the coordinate ring of 
the representation variety ${\rm Rep}_n(\Gamma)$.   
Let $\sigma_{\Gamma} : \Gamma \to 
{\rm M}_n(A_n(\Gamma))$ be the universal representation of 
degree $n$ for $\Gamma$. 
We denote by $I$ the ideal of $A_n(\Gamma)$ generated by 
$\{ \sigma_{\Gamma}(\gamma)_{ij} \mid 1 \le i \neq j \le n, \gamma \in \Gamma \} 
\cup \{ \sigma_{\Gamma}(\gamma)_{ii} - \sigma_{\Gamma}(\gamma)_{jj} 
\mid  1 \le i < j \le n, \gamma \in \Gamma \}$. 
Then it is easy to check that  
${\rm Rep}_n(\Gamma)_{\rm scalar}$ is 
representable by the affine scheme 
${\rm Spec} A_n(\Gamma)/I$. 
Since $I$ is ${\rm PGL}_n$-invariant and 
the action of ${\rm PGL}_n$ on $A_n(\Gamma)/I$ is trivial, 
the induced action of ${\rm PGL}_n$ on 
${\rm Rep}_n(\Gamma)_{\rm scalar}$ is trivial. 
\qed 

\begin{theorem}\label{th:moduliscalar} 
There exists a fine moduli scheme ${\rm Ch}_n(\Gamma)_{\rm scalar}$ 
separated over ${\Bbb Z}$ associated to 
the following contravariant functor: 
\[
\begin{array}{ccccl}
\mathcal{E}q\mathcal{S}_n(\Gamma) & : & ({\bf Sch})^{op} & \to & ({\bf Sets}) \\
  & & X & \mapsto & 
\left\{ 
\begin{array}{l} 
\mbox{ rep. with scalar mold  } \\
\mbox{ of degree $n$ for $\Gamma$ on $X$ } \\ 
\end{array}
\right\}
\Big/\sim.
\end{array}
\] 
The moduli ${\rm Ch}_n(\Gamma)_{\rm scalar}$ is isomorphic to 
${\rm Rep}_n(\Gamma)_{\rm scalar}$. 
Moreover, they are isomorphic to 
${\rm Rep}_1(\Gamma) \cong {\rm Ch}_1(\Gamma) := 
{\rm Rep}_1(\Gamma)/{\rm PGL}_1$. 
In particular, if $\Gamma$ is a finitely generated group (or monoid), then 
the moduli is of finite type over ${\Bbb Z}$. 
\end{theorem}

\prf
Since the action of ${\rm PGL}_n$ on ${\rm Rep}_n(\Gamma)_{\rm scalar}$ 
is trivial, the affine scheme 
${\rm Rep}_n(\Gamma)_{\rm scalar}$ also represents the functor  
$\mathcal{E}q\mathcal{S}_n(\Gamma)$. 
We easily see that 
$A_n(\Gamma)/I \cong A_1(\Gamma)$, where $I$ is defined in 
the proof of Proposition \ref{prop:repscalar}. 
The action of ${\rm PGL}_1 \cong {\rm Spec} \: {\mathbb Z}$ 
on ${\rm Rep}_1(\Gamma)$ is trivial. 
Hence we see that ${\rm Rep}_n(\Gamma)_{\rm scalar} 
\cong {\rm Rep}_1(\Gamma) \cong {\rm Ch}_1(\Gamma)$. 
If $\Gamma$ is finitely generated, then 
${\rm Rep}_n(\Gamma)$ is of finite type over ${\mathbb Z}$, 
and therefore so is ${\rm Rep}_n(\Gamma)_{\rm scalar}$. 
\qed

%
%
%

\section{Semi-simple mold}

In \S 4-\S 6, we only deal with rank $2$ molds of degree $2$. 
In this section, we investigate the semi-simple mold case.

\begin{definition}\label{def:reprkle2}\rm
  Let $\sigma_{\Gamma} : \Gamma \to {\rm M}_2(A_2(\Gamma))$ be 
the universal representation of degree $2$ for a group or a monoid $\Gamma$. 
For $\alpha, \beta, \gamma \in \Gamma$, 
we define the matrix $M(\alpha, \beta, \gamma)$ by 
\[
M(\alpha, \beta, \gamma) := 
\left(
\begin{array}{ccc}
\sigma_{\Gamma}(\alpha)_{11} & \sigma_{\Gamma}(\beta)_{11} & 
\sigma_{\Gamma}(\gamma)_{11} \\
\sigma_{\Gamma}(\alpha)_{12} & \sigma_{\Gamma}(\beta)_{12} & 
\sigma_{\Gamma}(\gamma)_{12} \\
\sigma_{\Gamma}(\alpha)_{21} & \sigma_{\Gamma}(\beta)_{21} & 
\sigma_{\Gamma}(\gamma)_{21} \\
\sigma_{\Gamma}(\alpha)_{22} & \sigma_{\Gamma}(\beta)_{22} & 
\sigma_{\Gamma}(\gamma)_{22} \\
\end{array}
\right).
\]
We define the closed subscheme ${\rm Rep}_2(\Gamma)_{{\rm rk} \le 2}$ 
of ${\rm Rep}_2(\Gamma)$ by 
\[
{\rm Rep}_2(\Gamma)_{{\rm rk}\le 2} := 
\left\{ \rho \in {\rm Rep}_2(\Gamma)\:
\begin{array}{|c}
\mbox{ all } (3 \times 3) \mbox{ minor determinants of } \\
M(\alpha, \beta, \gamma) \mbox{ are } 0 \mbox{ for each } 
\alpha, \beta, \gamma \in \Gamma
\end{array}
\right\}.
\]
We also define the open subscheme ${\rm Rep}_2(\Gamma)_{{\rm rk} 2}$ 
of the affine scheme ${\rm Rep}_2(\Gamma)_{{\rm rk}\le 2}$
by 
\[
{\rm Rep}_2(\Gamma)_{{\rm rk} 2} := 
\{
\rho \in {\rm Rep}_2(\Gamma) \mid {\mathcal O}_{X}[\rho(\Gamma)] 
\mbox{ is a rank $2$ mold }
\}.
\]
\end{definition}

\bigskip 

\begin{definition}\rm
We define the {\it representation variety with semi-simple mold} 
of degree $2$ for a group or a monoid $\Gamma$ by 
\[
\begin{array}{ccccl}
{\rm Rep}_2(\Gamma)_{\rm s.s.}  & : & ({\bf Sch})^{op} & \to &
({\bf Sets}) \\
 & & X & \mapsto & \{ \rho \in {\rm Rep}_2(\Gamma) 
\mid \rho \mbox{ has a semi-simple mold } \}.
\end{array}
\]
We easily see that  ${\rm Rep}_2(\Gamma)_{\rm s.s.}$ is  
an open subscheme of ${\rm Rep}_2(\Gamma)_{{\rm rk} 2}$. 
\end{definition}

\begin{remark}\rm
The scheme ${\rm Rep}_2(\Gamma)_{\rm s.s.}$ is an open subscheme 
of the affine scheme ${\rm Rep}_2(\Gamma)_{{\rm rk}\le 2}$ where 
$m(\sigma_{\Gamma}(\gamma))$ does not vanish for some $\gamma \in \Gamma$ 
by Remark~\ref{remark:m}. 
Here recall that $m(\sigma_{\Gamma}(\gamma)) = {\rm tr}(\sigma_{\Gamma}(\gamma))^2 
-4\det(\sigma_{\Gamma}(\gamma))$. 
\end{remark}

\bigskip

Let us denote by $A_2(\Gamma)_{{\rm rk} \le 2}$ the coordinate ring
of the affine scheme ${\rm Rep}_2(\Gamma)_{{\rm rk}\le 2}$.
We define $A_2(\Gamma)_{{\rm rk} \le 2}^{\rm Ch}$ 
as the subring 
of $A_2(\Gamma)_{{\rm rk} \le 2}$ generated by 
$\{ {\rm tr}(\sigma_{\Gamma}(\gamma)), \det(\sigma_{\Gamma}(\gamma)) \mid
\gamma \in \Gamma \}$ over ${\Bbb Z}$. 
We also denote by ${\rm Ch}_2(\Gamma)_{{\rm rk}\le 2}$ 
the spectrum of $A_2(\Gamma)_{{\rm rk} \le 2}^{\rm Ch}$. 
We define the open subscheme ${\rm Ch}_2(\Gamma)_{\rm s.s.}$ 
of ${\rm Ch}_2(\Gamma)_{{\rm rk}\le 2}$ by 
\[
{\rm Ch}_2(\Gamma)_{\rm s.s.} :=  \bigcup_{\gamma \in \Gamma} 
{\rm Spec} \; (A_2(\Gamma)_{{\rm rk} \le 2}^{\rm Ch})_{m
(\sigma_{\Gamma}(\gamma))}.
\]
Then we have the canonical morphism 
\[
\pi_{\Gamma, \rm s.s.} : {\rm Rep}_2(\Gamma)_{\rm s.s.} 
\to {\rm Ch}_2(\Gamma)_{\rm s.s.}.
\]
For $\gamma \in \Gamma$ we define 
\[
\begin{array}{ccl}
{\rm Rep}_2(\Gamma)_{\rm s.s., \gamma} & :=  & 
\{ x \in {\rm Rep}_2(\Gamma)_{\rm s.s.} \mid 
m(\sigma_{\Gamma}(\gamma)) \not\equiv 0 \text{ in } k(x) \} \\  
& = & {\rm Spec} \; (A_2(\Gamma)_{{\rm rk} \le 2})_{m
(\sigma_{\Gamma}(\gamma))} \\
\end{array} 
\]
and 
\[
\begin{array}{ccl}  
{\rm Ch}_2(\Gamma)_{\rm s.s., \gamma} & := &  
\{ x \in {\rm Ch}_2(\Gamma)_{\rm s.s.} \mid 
m(\sigma_{\Gamma}(\gamma)) \not\equiv 0 \text{ in } k(x) \} \\ 
& = &  {\rm Spec} \; (A_2(\Gamma)_{{\rm rk} \le 2}^{\rm Ch})_{m
(\sigma_{\Gamma}(\gamma))}.   \\
\end{array} 
\]
Then we have the canonical morphism 
\[
\pi_{\Gamma, \rm s.s., \gamma} : {\rm Rep}_2(\Gamma)_{\rm s.s., \gamma} 
\to {\rm Ch}_2(\Gamma)_{\rm s.s.,\gamma}.
\]

For a group or a monoid $\Gamma$, we have the following diagram for each $\gamma \in \Gamma$: 
\[
\begin{array}{ccccc}
{\rm Rep}_2(\Gamma)_{{\rm rk} \le 2} & \supseteq & {\rm Rep}_2(\Gamma)_{\rm s.s.} & \supseteq 
& {\rm Rep}_2(\Gamma)_{\rm s.s., \gamma} \\ 
\downarrow & & \downarrow & & \downarrow \\ 
{\rm Ch}_2(\Gamma)_{{\rm rk} \le 2} & \supseteq & {\rm Ch}_2(\Gamma)_{\rm s.s.} & \supseteq 
& {\rm Ch}_2(\Gamma)_{\rm s.s., \gamma} . 
\end{array}
\]

\bigskip

\begin{proposition}\label{prop:finitetype}
If $\Gamma$ is a finitely generated group or monoid, then 
${\rm Rep}_2(\Gamma)_{{\rm rk} \le 2}$ and 
${\rm Ch}_2(\Gamma)_{{\rm rk} \le 2}$ are of finite 
type over ${\Bbb Z}$. 
\end{proposition} 

\prf
Let $S = \{ \alpha_1, \ldots, \alpha_n \}$ be a set of generators of $\Gamma$. 
We may assume that $\alpha_{i}^{-1}$ is also an element of $S$ for each $1 \le i \le n$ 
if $\Gamma$ is a group. The coordinate ring $A_2(\Gamma)_{{\rm rk} \le 2}$  
is generated by all entries of $\sigma_{\Gamma}(\alpha_i)$ 
for $1 \le i \le n$ over ${\Bbb Z}$. Hence  ${\rm Rep}_2(\Gamma)_{{\rm rk} \le 2}$ 
is of finite type over ${\Bbb Z}$. 
Let $A_i := \sigma_{\Gamma}(\alpha_i)$ for $1 \le i \le n$. 
Then the coordinate ring $A_2(\Gamma)_{{\rm rk} \le 2}^{\rm Ch}$ is 
generated by $\{ \det(A_i) \mid 1 \le i \le n \}$ and $\{ {\rm tr}(A_{i_1}A_{i_2}\cdots A_{i_k}) \mid 
1 \le i_1 < i_2 < \cdots < i_k \le n \}$ over ${\Bbb Z}$. 
Indeed, we can verify it by using the following equalities: 
\begin{eqnarray*}
{\rm tr}(X^2Y) & = & {\rm tr}(X){\rm tr}(XY) - \det(X){\rm tr}(Y) \\
{\rm tr}(XYZ) & = & -{\rm tr}(XZY)  + {\rm tr}(X){\rm tr}(YZ) + {\rm tr}(Y){\rm tr}(ZX)  \\ 
& & + {\rm tr}(Z){\rm tr}(YX) - {\rm tr}(X){\rm tr}(Y){\rm tr}(Z) 
\end{eqnarray*}  
for $2 \times 2$ matrices $X, Y, Z$. 
These equalities have been well known (For proofs see \cite{kj-Saito93} or \cite[Appendix]{Nkmt02}). 
Therefore ${\rm Ch}_2(\Gamma)_{{\rm rk} \le 2}$ is of finite 
type over ${\Bbb Z}$. 
\qed

\bigskip 

\begin{definition}\rm 
  Let $\Upsilon_1 = \langle \alpha \rangle$ be the free monoid of 
rank $1$. We call the morphism 
$\pi_{{\Upsilon}_1, \rm s.s., \alpha} : {\rm Rep}_2({\Upsilon}_1)_{\rm s.s., \alpha} 
\to {\rm Ch}_2({\Upsilon}_1)_{\rm s.s., \alpha}$ 
the {\it prototype} with semi-simple mold of degree $2$.  

  Let ${\rm F}_1 = \langle \alpha \rangle$ be the free group of 
rank $1$. We call the morphism 
$\pi_{{\rm F}_1, \rm s.s., \alpha} : {\rm Rep}_2({\rm F}_1)_{\rm s.s., \alpha} 
\to {\rm Ch}_2({\rm F}_1)_{\rm s.s., \alpha}$ 
the {\it prototype for group representations} 
with semi-simple mold of degree $2$.  
\end{definition}

\bigskip

Let $\sigma_{{\Upsilon}_1}$
be the universal representation of degree $2$ for 
${\Upsilon}_1$. 
Put $\sigma_{{\Upsilon}_1}(\alpha) = 
\left(
  \begin{array}{cc}
    a & b \\
    c & d 
  \end{array}
\right)$. 
Then we see that the coordinate ring $A_2({\Upsilon}_1)$ 
of ${\rm Rep}_2({\Upsilon}_1)$ is isomorphic to the polynomial ring 
${\Bbb Z}[a, b, c, d]$. 
Note that ${\rm Rep}_2({\Upsilon}_1) = {\rm Rep}_2({\Upsilon}_1)_{{\rm rk} \le 2}$ 
and that ${\rm Rep}_2(\Upsilon_1)_{{\rm rk} 2} = 
D(a-d) \cup D(b) \cup D(c) \subseteq {\rm Rep}_2({\Upsilon}_1) = 
\Spec{{\Bbb Z}[a, b, c, d]}$. 

Put $D := ad-bc$ and $T := a+d$. 
Let $A_2(\Upsilon_1)^{\rm Ch}$ be the subalgebra of 
$A_2({\Upsilon}_1)$ generated by $\{ {\rm tr}(\sigma_{{\Upsilon}_1}(\gamma)), 
\det (\sigma_{{\Upsilon}_1}(\gamma)) \mid \gamma \in \Upsilon_{1} \}$ 
over ${\mathbb Z}$. 
Then $A_2(\Upsilon_1)^{\rm Ch}$ is isomorphic to the polynomial ring 
${\Bbb Z}[T, D]$.
Set ${\rm Ch}_2({\Upsilon}_1) := {\rm Spec} A_2(\Upsilon_1)^{\rm Ch}$. 
Then ${\rm Ch}_2({\Upsilon}_1) = {\rm Ch}_2({\Upsilon}_1)_{{\rm rk} \le 2}$. 

\begin{proposition}\label{prop:mformula}  
Let $R$ be a commutative ring. Let $A \in {\rm M}_2(R)$. 
For each $n \in {\mathbb N}$, 
there exists $f(x, y) \in {\mathbb Z}[x, y]$ such that 
$m(A^n) = m(A) f({\rm tr}A, \det A)$.  
\end{proposition}

\prf
Let us claim that 
\[ 
\begin{array}{l} 
m(A^n) =  \\ 
\left\{ 
\begin{array}{ll} 
\displaystyle 
m(A)\cdot\left[ \sum_{k=0}^{(n-3)/2} \det(A)^k {\rm tr}(A^{n-2k-1}) + \det(A)^{(n-1)/2} 
\right]^2 & (n : \text{ odd } )  \\
\displaystyle  m(A)\cdot\left[ \sum_{k=0}^{(n-2)/2} 
\det(A)^k {\rm tr}(A^{n-2k-1}) \right]^2 & (n : \text{ even } ).  \\
\end{array}  
\right. 
\end{array} 
\]
Since ${\rm tr}(A^k)$ can be expressed by a polynomial in ${\mathbb Z}[ {\rm tr}(A), \det (A) ]$ 
for each $k \in {\mathbb N}$, the statement follows from this claim.  
It only suffices to prove that this claim holds for 
$A = \sigma_{{\Upsilon}_1}(\alpha) \in {\rm M}_2(A_2(\Upsilon_1))$. 
 
For $A \in {\rm M}_2(k)$ with an algebraically closed field $k$, 
let $\lambda, \mu$ be eigenvalues of $A$. 
Note that $m(A)=(\lambda - \mu)^2$ and $m(A^n) = (\lambda^n - \mu^n)^2$. 
Then 
\begin{eqnarray*}
m(A^n) & = & (\lambda - \mu)^2 (\lambda^{n-1} + \lambda^{n-2} \mu + 
\cdots + \lambda \mu^{n-2} + \mu^{n-1} )^2 \\
 & = & m(A) \{ {\rm tr}(A^{n-1}) + \det(A) {\rm tr}(A^{n-3}) + 
\det(A)^2 {\rm tr}(A^{n-5}) + \cdots \}^2. 
\end{eqnarray*} 
Hence the claim holds for $A \in {\rm M}_2(k)$ with $k = \overline{k}$.   
Because the claim holds for an algebraic closure 
$k$ of the quotient field $Q(A_2(\Upsilon_1))$ of 
$A_2(\Upsilon_1)$, it also holds for $Q(A_2(\Upsilon_1))$ and for $A_2(\Upsilon_1)$. 
This completes the proof.   
\qed 

\begin{remark}\rm 
Using Proposition \ref{prop:mformula},  
we easily see that ${\rm Rep}_2({\Upsilon}_1)_{\rm s.s.} = 
{\rm Rep}_2({\Upsilon}_1)_{\rm s.s., \alpha}$ and that 
${\rm Ch}_2({\Upsilon}_1)_{\rm s.s.} = 
{\rm Ch}_2({\Upsilon}_1)_{\rm s.s., \alpha}$.  
Note that $m(A^{-1}) = m(A) (\det (A))^{-2}$ for $A \in {\rm GL}_2(R)$ 
with a commutative ring $R$. 
Hence we also see that 
${\rm Rep}_2({\rm F}_1)_{\rm s.s.} = 
{\rm Rep}_2({\rm F}_1)_{\rm s.s., \alpha}$ and that 
${\rm Ch}_2({\rm F}_1)_{\rm s.s.} = 
{\rm Ch}_2({\rm F}_1)_{\rm s.s., \alpha}$.
\end{remark}

\bigskip 

We have the following diagram for the free monoid $\Upsilon_1 = \langle \alpha \rangle$:  
\[
\begin{array}{ccccccc}
{\rm Rep}_2(\Upsilon_1) & & & & & & \\ 
\parallel & & & & & & \\
{\rm Rep}_2(\Upsilon_1)_{{\rm rk} \le 2} & \supset & 
{\rm Rep}_2(\Upsilon_1)_{{\rm rk} 2} & \supset & {\rm Rep}_2(\Upsilon_1)_{\rm s.s.} & =  
& {\rm Rep}_2(\Upsilon_1)_{\rm s.s., \alpha} \\ 
\downarrow & & \downarrow & & \downarrow & & \downarrow \\ 
{\rm Ch}_2(\Upsilon_1) & = & {\rm Ch}_2(\Upsilon_1)_{{\rm rk} \le 2} 
& \supset & {\rm Ch}_2(\Upsilon_1)_{\rm s.s.} & =  
& {\rm Ch}_2(\Upsilon_1)_{\rm s.s., \alpha} . 
\end{array}
\]

\bigskip 

Put $m := T^2 -4D$.  
The morphism $\pi : {\rm Rep}_2(\Upsilon_1)_{{\rm rk} 2} 
\to {\rm Ch}_2(\Upsilon_1)$ is given by 
$D(a-d) \cup D(b) \cup D(c) \to 
{\rm Spec} {\Bbb Z}[T, D]$. 
The prototype $\pi_{{\Upsilon}_1, \rm s.s.} 
: {\rm Rep}_2({\Upsilon}_1)_{\rm s.s.} 
\to {\rm Ch}_2({\Upsilon}_1)_{\rm s.s.}$
is induced by the ring homomorphism 
${\Bbb Z}[T, D]_{m} \to 
{\Bbb Z}[a, b, c, d]_{m}$.

\bigskip 

\begin{lemma}\label{lemma-norm-nonsca}
  Let $(R, m)$ be a local ring. 
Let $A \in {\rm M}_2(R)$.
Suppose that 
$(A \mod m)$
 is not a scalar matrix of 
${\rm M}_2(R/m)$.
Then there exists $P \in {\rm GL}_2(R)$ 
such that 
$$P^{-1}AP = 
\left(
  \begin{array}{cc}
    0 & -\det(A) \\
    1 & {\rm tr}(A) 
  \end{array}
\right).$$
If $Q \in {\rm M}_2(R)$ satisfies 
$AQ = QA$, then 
$Q = \lambda I_2 + \mu A$ for some $\lambda, \mu \in R$. 
\end{lemma}

\prf
Put $A = 
\left(
  \begin{array}{cc}
    a & b \\
    c & d 
  \end{array}
\right)$.
From the assumption, at least one of $a-d$, $b$, $c$ is 
contained in $R^{\times}$.
Assume that $b \in R^{\times}$.
Then the vectors $e_2 := {}^{t}(0, 1)$ and $A e_2 \in R^{2}$ 
form a basis of $R^{2}$. With respect to the basis 
$\{ e_2,  A e_2 \}$, the linear map $A : R^2 \to R^2$ can be 
expressed as 
$\left(
  \begin{array}{cc}
    0 & -\det(A) \\
    1 & {\rm tr}(A) 
  \end{array}
\right)$.
In the case $c \in R^{\times}$, we can choose $\{ e_1, A e_1 \}$ 
as a basis of $R^2$, where $e_1 := {}^t (1, 0)$.  
Then we can change $A$ into the form which we want.
If $a-d \in R^{\times}$ and $b, c \not\in R^{\times}$, then 
the vectors $e_1+e_2 = {}^{t}(1, 1)$ and $A (e_1+e_2)$ 
form a basis of $R^2$. 
Similarly we can change $A$ into the desired form.

To prove the latter part of the statement, 
we may assume that 
$A = \left(
  \begin{array}{cc}
    0 & -\det(A) \\
    1 & {\rm tr}(A) 
  \end{array}
\right)$.
By direct calculation, we see that 
$AQ = QA$ implies $Q = \lambda I_2 + \mu A$ for some 
$\lambda, \mu \in R$. 
\qed

\begin{proposition}\label{lemmassff}
The morphism 
${\rm Rep}_2(\Upsilon_1)_{{\rm rk} 2} \to 
{\rm Ch}_2(\Upsilon_1)$ is smooth and surjective. 
In particular, it is faithfully flat.   
\end{proposition}

\prf 
Let $I$ be an ideal of a local ring $R$ with $I^2=0$. 
For a given commutative diagram 
\[
\begin{array}{ccc}
{\rm Rep}_2(\Upsilon_1)_{{\rm rk} 2} & \to & 
{\rm Ch}_2(\Upsilon_1) \\ 
\uparrow & & \uparrow  \\
{\rm Spec} R/I & \to & {\rm Spec} R, \\ 
\end{array} 
\]
we obtain $(T, D) \in R^2$ and $\overline{A} \in {\rm Rep}_2(\Upsilon_1)_{{\rm rk}2}(R/I) 
\subset {\rm M}_2(R/I)$ 
such that ${\rm tr}(\overline{A}) \equiv T$ and  
$\det (\overline{A}) \equiv D \; ({\rm mod} \; I)$.  
By Lemma \ref{lemma-norm-nonsca}, there exists $\overline{P} \in {\rm GL}_2(R/I)$ 
such that  
\[
\overline{P}^{\:-1}\:\overline{A}\:\overline{P} \equiv B :=  
\left(
  \begin{array}{cc}
    0 & -D \\
    1 & T 
  \end{array}
\right) \; ({\rm mod}\; I). 
\] 
Let us take $P \in {\rm GL}_2(R)$ such that $P \equiv \overline{P} \; ({\rm mod}\; I)$. 
Put $A = P B P^{-1}$. 
Then $A \in {\rm Rep}_2(\Upsilon_1)_{{\rm rk} 2}(R)$ such that 
${\rm tr}(A) = T$ and  
$\det (A) = D$.  Hence we obtain a morphism 
${\rm Spec} R \to {\rm Rep}_2(\Upsilon_1)_{{\rm rk} 2}$ satisfying 
the commutativity. 
This implies that the morphism ${\rm Rep}_2(\Upsilon_1)_{{\rm rk} 2} \to 
{\rm Ch}_2(\Upsilon_1)$ is smooth. 
Surjectivity follows from that we can take such matrix as $B$ above for 
a given $k$-valued point $(T, D) \in {\rm Ch}_2(\Upsilon_1)(k)$ 
with a field $k$.  Since smoothness implies flatness, 
it is faithfully flat.  
\qed 

\begin{lemma}\label{lemma:detformula} 
Let $R$ be a commutative ring. 
For $X, Y \in {\rm M}_2(R)$ and $a, b \in R$, we have 
\[
\det(aX + bY) = a^2 \det(X) + b^2 \det(Y) + ab({\rm tr}(X){\rm tr}(Y) 
- {\rm tr}(XY)). 
\]
\end{lemma} 

\prf
By direct calculation, we can check the formula above. 
\qed

\begin{lemma}\label{lemma:eqder} 
Let $(R, m, k)$ be an artinian local ring, and  
$I$ be an ideal of $R$ with $mI = 0$. 
For $A \in {\rm M}_2(R)$, let us define the 
$k$-linear map $[A, -] : {\rm M}_2(I) \to {\rm M}_2(I)$ 
by $X \mapsto AX-XA$.  If  $(A \mod m)$ is not a scalar matrix 
of ${\rm M}_2(k)$, then   
\[
{\rm Im} [A, -] = \{ Y \in {\rm M}_2(I) \mid {\rm tr}(Y)={\rm tr}(AY)=0 \}. 
\]    
\end{lemma}

\prf 
Since $mI =0$, we can regard $I$ as a vector space 
over $R/m = k$. Put $d := \dim_{k} I < \infty$. 
Set $N := \{ Y \in {\rm M}_2(I) \mid {\rm tr}(Y)={\rm tr}(AY)=0 \}$. 
If $Y = [A, X] \in {\rm Im}[A, -]$, then 
${\rm tr}(Y) = {\rm tr}(AX) -{\rm tr}(XA) = 0$ and 
${\rm tr}(AY) = {\rm tr}(AAX) - {\rm tr}(AXA) = 0$. 
Hence ${\rm Im}[A, -] \subseteq N$. 
For showing that ${\rm Im} [A, -] = N$, 
we prove that the dimensions of the both sides coincide. 
In order to calculate the dimensions, we may change $A$ into 
\[
P^{-1}AP = 
\left(
\begin{array}{cc}
0 & -D \\
1 & T \\
\end{array}
\right)  
\]
for suitable $P \in {\rm GL}_2(R)$ 
by considering the automorphism ${\rm Ad}(P) : 
{\rm M}_2(I) \to {\rm M}_2(I)$ by $X \mapsto P^{-1}XP$.   

If $X \in {\rm Ker}[A, -]$, then 
$X = \lambda I_2 + \mu A$ for some $\lambda, \mu \in R$ 
by Lemma \ref{lemma-norm-nonsca}. 
Since $X \in {\rm M}_2(I)$, we get  
$\lambda, \mu \in I$.  
Hence we see that $\dim_{k} {\rm Ker}[A, -] = 
\dim_{k} (I\cdot I_2 + I\cdot A) = 2d$ 
and that $\dim_{k} {\rm Im}[A, -] = 
\dim_{k} {\rm M}_2(I) - \dim_{k} {\rm Ker}[A, -] = 2d$.  
On the other hand, if $X \in N$, then 
${\rm tr}(X)={\rm tr}(AX)=0$. By direct calculation, 
we have $\dim_{k} N = 2d$. Thus 
we have proved that $\dim_{k} {\rm Im}[A, -] 
= \dim_{k} N = 2d$ and that ${\rm Im} [A, -] = N$. 
\qed 

\bigskip 

Let $s : {\rm Ch}_2(\Upsilon_1) \to {\rm Rep}_2(\Upsilon_1)_{{\rm rk} 2}$ by 
$(T, D) \mapsto 
\left(
\begin{array}{cc}
0 & -D \\
1 & T \\
\end{array}
\right)
$.
Then $\pi \circ s = 1_{{\rm Ch}_2(\Upsilon_1)}$. 
We have the following proposition: 

\begin{proposition}\label{prop:sectionff} 
The composition of the morphisms  
\[
\begin{array}{ccccl}
  {\rm Ch}_2(\Upsilon_1) \times {\rm PGL}_2 & \stackrel{(s, id)}{\to} & 
{\rm Rep}_2(\Upsilon_1)_{{\rm rk} 2}\times 
{\rm PGL}_2  & \stackrel{\sigma}{\to} & {\rm Rep}_2(\Upsilon_1)_{{\rm rk} 2} \\
((T, D), P) & \mapsto & 
\left( 
\left(
\begin{array}{cc}
0 & -D \\
1 & T \\
\end{array}
\right)
, P
\right) & \mapsto & P^{-1}
\left(
\begin{array}{cc}
0 & -D \\
1 & T \\
\end{array}
\right) 
P 
\end{array}
\]  
is smooth and surjective. In particular, it is faithfully flat. 
\end{proposition} 

\prf
Let $(R, m, k)$ be an artinian local ring, and let 
$I$ be an ideal of $R$ with $mI = 0$. 
For a given commutative diagram 
\[
\begin{array}{ccc}
  {\rm Ch}_2(\Upsilon_1) \times {\rm PGL}_2 & \stackrel{\sigma \circ (s, id) }{\to} & 
{\rm Rep}_2(\Upsilon_1)_{{\rm rk} 2} \\
\uparrow  & & \uparrow \\
{\rm Spec} \: R/I & \to & {\rm Spec} \: R, \\
\end{array}
\]  
we obtain $A \in {\rm Rep}_2(\Upsilon_1)_{{\rm rk} 2}(R)$, 
$(\overline{T}, \overline{D}) 
\in {\rm Ch}_2(\Upsilon_1)(R/I)$, and  
$\overline{P} \in {\rm GL}_2(R/I)$ such that 
$\overline{P}^{\;-1}  
\left(
\begin{array}{cc}
0 & - \overline{D} \\ 
1 & \overline{T} \\
\end{array}
\right) 
 \overline{P} \equiv A$ 
in ${\rm M}_2(R/I)$. 
Take $P \in {\rm GL}_2(R)$ such that $(P \mod I) = \overline{P}$.  
Put $D := \det A$ and $T := {\rm tr}(A)$. 
Note that $(D \mod I) = \overline{D}$ and $(T \mod I)  
= \overline{T}$. 

Set $Y := P^{\;-1}  
\left(
\begin{array}{cc}
0 & - D \\ 
1 & T \\
\end{array}
\right) 
P - A \in {\rm M}_2(I)$. 
Let us show that ${\rm tr}(AY) = 0$. 
Remark that $\det Y = 0$ by $I^2 = 0$ and that ${\rm tr}(Y) = 0$. 
Using Lemma \ref{lemma:detformula}, 
we have 
\begin{eqnarray*}
D & = &  \det (P^{-1} \left(
\begin{array}{cc}
0 & - D \\ 
1 & T \\
\end{array}
\right) 
P)  \\ 
 & = & \det(A+Y) \\
 & = & \det A + \det Y +{\rm tr}(A){\rm tr}(Y)-{\rm tr}(AY) \\
 & = & D - {\rm tr}(AY). 
\end{eqnarray*}
Hence we have proved ${\rm tr}(AY)=0$.

By Lemma \ref{lemma:eqder}, we have 
$Y \in {\rm Im}[A, -]$. There exists $X \in {\rm M}_2(I)$ 
such that $[A, X] = Y$. 
Put $P':=P(I_2 - X) \in {\rm GL}_2(R)$. 
Then $P'^{-1} 
\left(
\begin{array}{cc}
0 & - D \\ 
1 & T \\
\end{array}
\right) 
P' = (I_2 + X)(A+Y)(I_2 - X) = A + Y - [A, X]
= A$. 

Now let us define the morphism  
${\rm Spec} R \to {\rm Ch}_2(\Upsilon_1) 
\times {\rm PGL}_2$ corresponding to 
$((T, D), P')$. Verifying  
the commutativity, we see that the morphism is smooth. 
By Lemma \ref{lemma-norm-nonsca} we see that 
the morphism is surjective. 
Hence it is faithfully flat. 
\qed

\begin{proposition}\label{propssff}
The morphism 
\[
\begin{array}{ccc}
  {\rm Rep}_2(\Upsilon_1)_{{\rm rk} 2}\times{\rm PGL}_2 & \to & 
{\rm Rep}_2(\Upsilon_1)_{{\rm rk} 2}\times_{{\rm Ch}_2(\Upsilon_1)}
{\rm Rep}_2(\Upsilon_1)_{{\rm rk} 2} \\
(\rho, P) & \mapsto & (\rho, P^{-1}\rho P)
\end{array}
\]  
is smooth and surjective. In particular, it is faithfully flat. 
\end{proposition}

\prf
Let $(R, m, k)$ be an artinian local ring, and let 
$I$ be an ideal of $R$ with $mI = 0$. 
For a given commutative diagram 
\[
\begin{array}{ccc}
  {\rm Rep}_2(\Upsilon_1)_{{\rm rk} 2}\times{\rm PGL}_2 & \to & 
{\rm Rep}_2(\Upsilon_1)_{{\rm rk} 2}\times_{{\rm Ch}_2(\Upsilon_1)}
{\rm Rep}_2(\Upsilon_1)_{{\rm rk} 2} \\
\uparrow  & & \uparrow \\
{\rm Spec} \: R/I & \to & {\rm Spec} \: R, \\
\end{array}
\]  
we obtain $(A, B) \in {\rm Rep}_2(\Upsilon_1)_{{\rm rk} 2}(R)
\times_{{\rm Ch}_2(\Upsilon_1)(R)}
{\rm Rep}_2(\Upsilon_1)_{{\rm rk} 2}(R)$  
and $\overline{P} \in {\rm GL}_2(R/I)$ such that 
$\overline{P}^{\;-1} A \overline{P} \equiv B$ 
in ${\rm M}_2(R/I)$. 
For proving that the morphism is smooth, we define a 
morphism ${\rm Spec} R \to {\rm Rep}_2(\Upsilon_1)_{{\rm rk} 2}
\times{\rm PGL}_2$ satisfying the commutativity. 
Put $T = {\rm tr}(A)={\rm tr}(B)$ and $D = \det A = \det B$. 
Let us take $P \in {\rm GL}_2(R)$ such that 
$(P \mod I) = \overline{P}$. 
Then set $C := P^{-1}AP - B \in {\rm M}_2(I)$. 

Let us show that ${\rm tr}(C) = {\rm tr}(BC)=0$. 
Indeed, ${\rm tr}(C) = {\rm tr}(P^{-1}AP)-{\rm tr}(B) = T-T = 0$. 
Note that $\det C = 0$ by $I^2=0$. 
Using Lemma \ref{lemma:detformula}, 
we have   
\begin{eqnarray*}
D = \det (P^{-1}AP) & = & \det(B+C) \\
 & = & \det B + \det C +{\rm tr}(B){\rm tr}(C)-{\rm tr}(BC) \\
 & = & D - {\rm tr}(BC). 
\end{eqnarray*}
Hence we have verified ${\rm tr}(BC) = 0$. 

By Lemma \ref{lemma:eqder}, we have 
$C \in {\rm Im}[B, -]$. There exists $X \in {\rm M}_2(I)$ 
such that $[B, X] = C$. 
Put $P':=P(I_2 - X) \in {\rm GL}_2(R)$. 
Then $P'^{-1}AP' = (I_2 + X)(B+C)(I_2 - X) = B + C - [B, X]
= B$. 
Now let us define the morphism  
${\rm Spec} R \to {\rm Rep}_2(\Upsilon_1)_{{\rm rk} 2}
\times{\rm PGL}_2$ corresponding to 
$(A, P')$. Since $P'^{-1}AP'=B$, we can verify 
the commutativity. Therefore the morphism is smooth. 

By Lemma \ref{lemma-norm-nonsca} we see that 
the morphism is surjective. 
Hence it is faithfully flat. 
\qed

\bigskip 

Let us introduce the following two lemmas on sufficient 
conditions for a given morphism to be a universal geometric quotient 
(For the definition of universal geometric quotient, see \cite{GIT}).

\begin{lemma}\label{lemma:univgeom} 
Let $G$ be an affine group scheme over an affine scheme $S$. 
Assume that the $S$-morphism $\sigma : G\times_{S} X \to X$ 
is a group action of $G$ on an $S$-scheme $X$ and that 
the action of $G$ on an $S$-scheme $Y$ is trivial.  
Let $\pi : X \to Y$ be an affine $G$-equivariant 
faithfully flat locally of finite presentation $S$-morphism. 
If the morphism $(\sigma, p_2) : G \times_{S} X \to X\times_{Y}X$ 
is faithfully flat, then 
$\pi : X \to Y$ is a universal geometric quotient by $G$.   
\end{lemma} 

\prf
From the assumption, we have 
$\pi \circ \sigma = \pi \circ p_2 : G\times_{S} X \to Y$. 
We also see that $\pi$ is surjective and that 
the image of the morphism $(\sigma, p_2) : 
G \times_{S} X \to X\times_{S} X$ is $X\times_{Y}X$. 
Since $\pi$ is faithfully flat locally of finite presentation, it is universally 
open ({\it cf}. \cite[Theorem~2.4.6]{EGA4}).   
For verifying $\pi_{\ast}({\mathcal O}_{X})^{G} ={\mathcal O}_Y$, 
we only need to check that ${\mathcal O}_X(\pi^{-1}(U))^{G} = 
{\mathcal O}_{Y}(U)$ for each affine open subscheme 
$U$ of $Y$. Set $U = {\rm Spec} A$ and $\pi^{-1}(U) = {\rm Spec} B$. 
Because $A \to B$ is faithfully flat, 
$0 \to A \to B \stackrel{\phi_1 - \phi_2}{\to} 
B\otimes_{A} B$ is exact, where $\phi_{1}(b) = b\otimes 1$ 
and $\phi_2(b) = 1\otimes b$ for $b \in B$ (for example, see 
\cite[Proposition~2.18]{Etale}).    
Put $S = {\rm Spec} \: C$ and $G = {\rm Spec} R$. 
The ring homomorphism $B \otimes_{A} B \to 
R \otimes_{C} B$ induced by $(\sigma, p_2) : 
G \times_{S} X \to X \times_{Y} X$ 
is faithfully flat.  
Since $B \otimes_{A} B \to 
R \otimes_{C} B$ is injective, 
$0 \to A \to B \stackrel{\sigma^{\ast} - p_2^{\ast}}{\to} 
R \otimes_{C} B$ is also exact. 
This implies that $B^{G} = A$ and that 
$\pi_{\ast}({\mathcal O}_{X})^{G} ={\mathcal O}_Y$. 
Hence we see that $\pi : X \to Y$ is a geometric 
quotient by $G$. 
For all morphisms $Y' \to Y$, $\pi' : X' := X\times_{Y} Y' 
\to Y'$ is also a geometric quotient by $G$ 
because $\pi'$ and $G\times_{S} X' \to X'\times_{Y'} X'$ 
are faithfully flat.  This completes the proof. 
\qed

\begin{lemma}\label{lemma:univgeom2} 
Let $X$ and $Y$ be schemes over a scheme $S$. 
Let $G$ be a group scheme over $S$. 
Assume that an $S$-morphism $\sigma : G\times_{S} X \to X$ 
is a group action of $G$ on $X$ and that 
the action of $G$ on $Y$ is trivial.  
Let $\pi : X \to Y$ be a $G$-equivariant $S$-morphism and 
$s : Y \to X$ be an $S$-morphism such that $\pi \circ s = 1_{Y}$.  
Suppose that 
\begin{enumerate}
\item\label{cond:2-1} $\pi$ is faithfully flat and locally of finite presentation, and 
\item\label{cond:2-2} $G \times_{S} Y \stackrel{(1_{G}, s)}{\longrightarrow} G \times_{S} X 
\stackrel{\sigma}{\to} X$ is faithfully flat. 
\end{enumerate}
Then $\pi : X \to Y$ is a universal geometric quotient by $G$.   
\end{lemma} 

\prf
From the assumption (\ref{cond:2-1}),  
$\pi : X \to Y$ is surjective and universally open.  
Let us show that the image of the morphism $(\sigma, p_2) : 
G \times_{S} X \to X\times_{S} X$ is $X\times_{Y}X$. 
Assume that $(f_1, f_2) \in X(\Omega)\times_{Y(\Omega)}X(\Omega)$ 
is a $\Omega$-valued point with a field $\Omega$. 
Set $h := \pi \circ f_1 = \pi \circ f_2 : {\rm Spec}\: \Omega 
\to Y$ and $f := s \circ h : {\rm Spec} \: \Omega \to X$. 
By the condition (\ref{cond:2-2}), 
there exist $(g_i, f) : {\rm Spec} \: \Omega \to 
G \times_{S} X$ such that 
$\sigma \circ (g_i, f) = f_i$ for each $i =1, 2$ 
(if necessary, take an extension field of $\Omega$). 
Then $(\sigma, p_2) \circ (g_1 g_2^{-1}, f_2) : {\rm Spec} \: \Omega 
\to G \times_{S} X \to X \times_{S} X$ is 
$(f_1, f_2)$. Hence the image of $(\sigma, p_2)$ is 
$X\times_{Y}X$. 

Let us show that $\pi_{\ast}({\mathcal O}_{X})^{G} ={\mathcal O}_Y$. 
Since $\pi$ is faithfully flat, ${\mathcal O}_Y \to \pi_{\ast}({\mathcal O}_X)^{G}$ 
is injective. 
For an open set $U$ of $Y$, put $V := \pi^{-1}(U)$.    
Assume that $\phi : V \to {\mathbb A}^{1}_{S}$ 
satisfies $\phi \circ (\sigma |_{G \times_{S} V})  = 
\phi \circ (p_2 |_{G \times_{S} V})$. 
Set $\psi := \phi \circ (s|_{\:U}) : U \to {\mathbb A}^{1}_{S}$. 
By the assumption (\ref{cond:2-2}), 
$G \times_{S} U \stackrel{(1_{G}, s\mid_{U})}{\longrightarrow} G \times_{S} V  
\stackrel{\sigma\mid_{G\times_{S} V}}{\to} V$ is faithfully flat. 
Put $\Phi := (\sigma\mid_{G\times_{S} V}) \circ (1_{G}, s\mid_{U})$. 
It is easy to verify that $\phi \circ \Phi = \psi \circ (\pi|_{V}) \circ \Phi$. 
Since ${\mathcal O}_V \to \Phi_{\ast}({\mathcal O}_{G \times_{S} V})$ is injective,  
$\phi = \psi \circ (\pi|_{V})$. 
This implies that ${\mathcal O}_Y \to \pi_{\ast}({\mathcal O}_X)^{G}$ 
is surjective, and hence $\pi_{\ast}({\mathcal O}_{X})^{G} ={\mathcal O}_Y$. 

The assumptions (\ref{cond:2-1}) and (\ref{cond:2-2}) 
are stable under any base change $Y' \to Y$. 
Therefore $\pi : X \to Y$ is a universal geometric quotient by $G$.   
\qed

\begin{remark}\rm 
We can extend Lemma \ref{lemma:univgeom2} in the following way: 
Let $S, G, X, Y$, and  $\pi : X \to Y$ be as above. 
Let $Y = \cup_{i \in I} Y_i$ be an open covering of $Y$. 
Set $X_i := \pi^{-1}(Y_i)$.  
Let $s_i : Y_i \to X_i$ be an $S$-morphism with $\pi\!\mid_{X_i} \circ s_i = 1_{Y_i}$ 
for each $i \in I$.  Suppose that (\ref{cond:2-1}) and the following (ii)' hold: 
\begin{enumerate}
\item[(ii)']\label{cond:2-3} $G \times_{S} Y_i \stackrel{(1_{G}, s_i)}{\longrightarrow} G \times_{S} X_i  
\stackrel{\sigma}{\to} X_i$ is faithfully flat for each $i \in I$. 
\end{enumerate}
Then $\pi : X \to Y$ is a universal geometric quotient by $G$. 
\end{remark} 

\smallskip 

\begin{theorem}\label{corssf1}
The morphism 
$\pi: {\rm Rep}_2(\Upsilon_1)_{{\rm rk} 2} \to {\rm Ch}_2(\Upsilon_1)$
is a universal geometric quotient by ${\rm PGL}_2$.
\end{theorem}

\prf
Obviously, $\pi$ is locally of finite presentation. 
By Propositions \ref{lemmassff} and \ref{prop:sectionff}, 
$\pi$ and  
${\rm Ch}_2(\Upsilon_1) \times {\rm PGL}_2 \stackrel{\sigma \circ (s, id)}{\to} 
{\rm Rep}_2(\Upsilon_1)_{{\rm rk} 2}$ are faithfully flat. 
Hence Lemma \ref{lemma:univgeom2} implies that 
$\pi: {\rm Rep}_2(\Upsilon_1)_{{\rm rk} 2} \to {\rm Ch}_2(\Upsilon_1)$
is a universal geometric quotient by ${\rm PGL}_2$.   
(Of course, Lemma \ref{lemma:univgeom2} also holds for right group actions.)  
\qed

\begin{corollary}\label{corssf2}
The prototype 
$\pi_{\Upsilon_1, {\rm s.s.}} : {\rm Rep}_2(\Upsilon_1)_{\rm s.s.} \to {\rm Ch}_2(\Upsilon_1)_{\rm s.s.}$
is a universal geometric 
quotient by ${\rm PGL}_2$.
\end{corollary}

\prf
The morphism $\pi_{\Upsilon_1, {\rm s.s.}} : {\rm Rep}_2(\Upsilon_1)_{\rm s.s.} 
\to {\rm Ch}_2(\Upsilon_1)_{\rm s.s.}$ 
is a base change of $\pi: {\rm Rep}_2(\Upsilon_1)_{{\rm rk} 2} \to {\rm Ch}_2(\Upsilon_1)$ 
by ${\rm Ch}_2(\Upsilon_1)_{\rm s.s.} \to {\rm Ch}_2(\Upsilon_1)$. 
The statement follows from Theorem \ref{corssf1}.  
\qed

\begin{corollary}\label{corssf3}
The prototype 
$\pi_{{\rm F}_1, {\rm s.s.}} : {\rm Rep}_2({\rm F}_1)_{\rm s.s.} \to {\rm Ch}_2({\rm F}_1)_{\rm s.s.}$ 
for group representations 
is a universal geometric 
quotient by ${\rm PGL}_2$.
\end{corollary}

\prf
The morphism $\pi_{{\rm F}_1, {\rm s.s.}} : {\rm Rep}_2({\rm F}_1)_{\rm s.s.} \to {\rm Ch}_2({\rm F}_1)_{\rm s.s.}$ 
is a base change of $\pi: {\rm Rep}_2(\Upsilon_1)_{\rm s.s.} \to {\rm Ch}_2(\Upsilon_1)_{\rm s.s.}$ 
by ${\rm Ch}_2({\rm F}_1)_{\rm s.s.} \to {\rm Ch}_2(\Upsilon_1)_{\rm s.s.}$. 
The statement follows from the previous corollary.  
\qed 

\begin{proposition}\label{propss}
Let ${\mathcal A} \subseteq {\rm M}_2(R)$ 
be a rank $2$ semi-simple mold over a commutative ring $R$.
Suppose that there exists $A \in {\mathcal A}$ such that 
$m(A) \in R^{\times}$.  
Then the following bilinear form is perfect:
\[
\begin{array}{ccccc}
  {\rm tr}(\cdot\cdot) & : & {\mathcal A} \times {\mathcal A} & \to & R \\
 & & (X, Y) & \mapsto & {\rm tr}(XY). 
\end{array}
\]
In other words, the $R$-linear map ${\mathcal A} \to {\rm Hom}_R ({\mathcal A}, R)$ 
defined by $X \mapsto (Y \mapsto {\rm tr}(XY))$ is an isomorphism.  
In particular, for each $X \in {\mathcal A}$, we have
\begin{eqnarray*}
X & = & ( I_2, A )
\left(
\begin{array}{cc}
{\rm tr}(I_2) & {\rm tr}(A) \\
{\rm tr}(A) & {\rm tr}(A^2) \\
\end{array}
\right)^{-1}
\left(
\begin{array}{c}
{\rm tr}(X) \\
{\rm tr}(AX)\\
\end{array}
\right) \\
& = & 
( I_2, A )\;
\frac{1}{m(A)}
\left(
\begin{array}{cc}
{\rm tr}(A^2) & -{\rm tr}(A) \\
-{\rm tr}(A) & 2 \\
\end{array}
\right)
\left(
\begin{array}{c}
{\rm tr}(X) \\
{\rm tr}(AX)\\
\end{array}
\right). 
\end{eqnarray*}
\end{proposition}

\prf
Remark that $\{ I_2, A \}$ forms a basis of ${\mathcal A}$ over $R$.
The determinant of the matrix
\[
\left(
\begin{array}{cc}
{\rm tr}(I_2) & {\rm tr}(A) \\
{\rm tr}(A) & {\rm tr}(A^2) \\
\end{array}
\right)
\]
is equal to $m(A) \in R^{\times}$, and hence the inverse matrix exists. 
 
For each $X =a I_2 + b A \in {\mathcal A}$ with $a, b \in R$, we have  
\[
\left(
\begin{array}{c} 
 {\rm tr}(X) \\
 {\rm tr}(AX)  \\ 
\end{array}
\right) 
=  
\left(
\begin{array}{cc}
{\rm tr}(I_2) & {\rm tr}(A) \\
{\rm tr}(A) & {\rm tr}(A^2) \\
\end{array} 
\right) 
\left(
\begin{array}{c}
a \\
b \\
\end{array}
\right) . 
\]
Since 
\[
\left(
\begin{array}{c}
a \\
b \\
\end{array}
\right) 
= 
\left(
\begin{array}{cc}
{\rm tr}(I_2) & {\rm tr}(A) \\
{\rm tr}(A) & {\rm tr}(A^2) \\
\end{array} 
\right)^{-1} 
\left(
\begin{array}{c} 
 {\rm tr}(X) \\
 {\rm tr}(AX)  \\ 
\end{array}
\right), 
\]
we have 
\begin{eqnarray*}
X & = & ( I_2, A )
\left(
\begin{array}{cc}
{\rm tr}(I_2) & {\rm tr}(A) \\
{\rm tr}(A) & {\rm tr}(A^2) \\
\end{array}
\right)^{-1}
\left(
\begin{array}{c}
{\rm tr}(X) \\
{\rm tr}(AX)\\
\end{array}
\right) \\
& = & 
( I_2, A )\;
\frac{1}{m(A)}
\left(
\begin{array}{cc}
{\rm tr}(A^2) & -{\rm tr}(A) \\
-{\rm tr}(A) & 2 \\
\end{array}
\right)
\left(
\begin{array}{c}
{\rm tr}(X) \\
{\rm tr}(AX)\\
\end{array}
\right). 
\end{eqnarray*}
Therefore we see that ${\rm tr}(\cdot\cdot)  :  {\mathcal A} \times {\mathcal A}  \to R$ 
is perfect. 
\qed


\begin{proposition}\label{lemmassgen}  
For each element $\gamma$ of a monoid $\Gamma$, 
$(A_2(\Gamma)_{{\rm rk} \le 2}^{\rm Ch})_{m(\sigma_{\Gamma}(\gamma))}$ is 
generated by $\{ {\rm tr}(\sigma_{\Gamma}(\delta)) \mid \delta \in \Gamma \}$ 
and $m(\sigma_{\Gamma}(\gamma)))^{-1}$ over ${\mathbb Z}$. 
\end{proposition} 

\prf 
By the definition of $A_2(\Gamma)_{{\rm rk} \le 2}^{\rm Ch}$, 
$(A_2(\Gamma)_{{\rm rk} \le 2}^{\rm Ch})_{m(\sigma_{\Gamma}(\gamma))}$ is generated by 
$\{ {\rm tr}(\sigma_{\Gamma}(\delta)), \det(\sigma_{\Gamma}(\delta)) 
\mid \delta \in \Gamma \}$ 
and $m(\sigma_{\Gamma}(\gamma)))^{-1}$ over ${\mathbb Z}$. 
Let $S$ be the subalgebra of $(A_2(\Gamma)_{{\rm rk} \le 2}^{\rm Ch})_{m(\sigma_{\Gamma}(\gamma))}$ 
generated by $\{ {\rm tr}(\sigma_{\Gamma}(\delta)) \mid \delta \in \Gamma \}$ 
and $m(\sigma_{\Gamma}(\gamma)))^{-1}$ 
over ${\mathbb Z}$. 
It only suffices to prove that  
$\det (\sigma_{\Gamma}(\delta)) \in S$ for each $\delta \in \Gamma$.  
Using Proposition \ref{propss}, we have 
\begin{eqnarray*} 
\sigma_{\Gamma}(\delta) = 
( \sigma_{\Gamma}(e), \sigma_{\Gamma}(\gamma) )
\left(
  \begin{array}{cc}
    {\rm tr}(\sigma_{\Gamma}(e)) & {\rm tr}(\sigma_{\Gamma}(\gamma)) \\
    {\rm tr}(\sigma_{\Gamma}(\gamma)) & 
    {\rm tr}(\sigma_{\Gamma}(\gamma^2))
  \end{array}
\right)^{-1}
\left(
  \begin{array}{c}
    {\rm tr}(\sigma_{\Gamma}(\delta)) \\
    {\rm tr}(\sigma_{\Gamma}(\gamma\delta)) 
  \end{array}
\right)
\end{eqnarray*}  
in ${\rm M}_2((A_2(\Gamma)_{{\rm rk} \le 2})_{m(\sigma_{\Gamma}(\gamma))})$ 
for each $\delta \in \Gamma$. 
Since the determinant of the matrix 
\[
T := \left(
  \begin{array}{cc}
    {\rm tr}(\sigma_{\Gamma}(e)) & {\rm tr}(\sigma_{\Gamma}(\gamma)) \\
    {\rm tr}(\sigma_{\Gamma}(\gamma)) & 
    {\rm tr}(\sigma_{\Gamma}(\gamma^2))
  \end{array}
\right)
\] 
is equal to $m(\sigma_{\Gamma}(\gamma)) = 2 {\rm tr}(\sigma_{\Gamma}(\gamma^2)) 
- ({\rm tr}(\sigma_{\Gamma}(\gamma)))^2$, $\sigma_{\Gamma}(\delta) = a I_2 + b 
\sigma_{\Gamma}(\gamma)$ with some $a, b \in S$ for each $\delta \in \Gamma$. 
By Lemma \ref{lemma:detformula}, the statement follows from the claim that 
$\det (\sigma_{\Gamma}(\gamma)) \in S$.   
Let us prove the claim. 
Putting $\delta = \gamma^2$, we have  
\begin{eqnarray*} 
\sigma_{\Gamma}(\gamma^2) 
 & = & 
( \sigma_{\Gamma}(e), \sigma_{\Gamma}(\gamma) ) 
\; T^{-1} 
\left(
  \begin{array}{c}
    {\rm tr}(\sigma_{\Gamma}(\gamma^2)) \\
    {\rm tr}(\sigma_{\Gamma}(\gamma^3)) 
  \end{array}
\right) \\ 
 & = & ( I_2, \sigma_{\Gamma}(\gamma) )
\displaystyle \frac{1}{m(\sigma_{\Gamma}(\gamma))} 
\left(
  \begin{array}{c}
    {\rm tr}(\sigma_{\Gamma}(\gamma^2))^2 - {\rm tr}(\sigma_{\Gamma}(\gamma)){\rm tr}(\sigma_{\Gamma}(\gamma^3)) \\
    -{\rm tr}(\sigma_{\Gamma}(\gamma)){\rm tr}(\sigma_{\Gamma}(\gamma^2))+2{\rm tr}(\sigma_{\Gamma}(\gamma^3)) 
  \end{array}
\right). 
\end{eqnarray*} 
We also obtain $\sigma_{\Gamma}(\gamma^2) = 
{\rm tr}(\sigma_{\Gamma}(\gamma)) \sigma_{\Gamma}(\gamma) - \det{\sigma_{\Gamma}(\gamma)} I_2$ 
by the Cayley-Hamilton Theorem. Comparing the coefficients of $I_2$, we have 
\[ \displaystyle 
\det(\sigma_{\Gamma}(\gamma)) = \frac{%
{\rm tr}(\sigma_{\Gamma}(\gamma)){\rm tr}(\sigma_{\Gamma}(\gamma^3)) - 
{\rm tr}(\sigma_{\Gamma}(\gamma^2))^2}{m(\sigma_{\Gamma}(\gamma))}. 
\] 
Hence we have proved the statement. 
\qed 

\bigskip 

Let $\Gamma_1, \Gamma_2$ be monoids.   
Let $\psi : \Gamma_1 \to \Gamma_2$ be a monoid homomorphism. 
Then $\psi$ induces canonical ring homomorphisms 
$\psi_{\ast} : 
A_2(\Gamma_1)_{{\rm rk} \le 2}
\to A_2(\Gamma_2)_{{\rm rk} \le 2}$ and 
$\psi_{\ast} : 
A_2(\Gamma_1)_{{\rm rk} \le 2}^{\rm Ch}
\to A_2(\Gamma_2)_{{\rm rk} \le 2}^{\rm Ch}$.  
Set $\gamma_2 := \psi(\gamma_1)$ for $\gamma_1 \in \Gamma_1$. 
We obtain the ring homomorphisms 
\begin{eqnarray*} 
\psi_{\ast} : (A_2(\Gamma_1)_{{\rm rk} \le 2})_{m(\sigma_{\Gamma_1}(\gamma_1))}
& \to & (A_2(\Gamma_2)_{{\rm rk} \le 2})_{m(\sigma_{\Gamma_2}(\gamma_2))},  \\   
\psi_{\ast} : (A_2(\Gamma_1)_{{\rm rk} \le 2}^{\rm Ch})_{m(\sigma_{\Gamma_1}(\gamma_1))}
& \to & (A_2(\Gamma_2)_{{\rm rk} \le 2}^{\rm Ch})_{m(\sigma_{\Gamma_2}(\gamma_2))}. 
\end{eqnarray*}   
Hence we have the morphisms 
\[
\begin{array}{ccc}
{\rm Rep}_2(\Gamma_2)_{{\rm rk} \le 2} & 
\stackrel{\psi^{\ast}}{\to} & {\rm Rep}_2(\Gamma_1)_{{\rm rk} \le 2} \\
\downarrow & & \downarrow \\ 
{\rm Ch}_2(\Gamma_2)_{{\rm rk} \le 2} & 
\stackrel{\psi^{\ast}}{\to} & {\rm Ch}_2(\Gamma_1)_{{\rm rk} \le 2} \\
\end{array}
\]
and 
\[
\begin{array}{ccc}
{\rm Rep}_2(\Gamma_2)_{{\rm s.s.}, \gamma_2} & 
\stackrel{\psi^{\ast}}{\to} & {\rm Rep}_2(\Gamma_1)_{{\rm s.s.}, \gamma_1} \\
\downarrow & & \downarrow \\ 
{\rm Ch}_2(\Gamma_2)_{{\rm s.s.}, \gamma_2} & 
\stackrel{\psi^{\ast}}{\to} & {\rm Ch}_2(\Gamma_1)_{{\rm s.s.}, \gamma_1}. \\
\end{array}
\]

\begin{lemma}\label{lemssror} 
Let $\gamma$ be an element of a monoid $\Gamma$.
Let $\psi : {\Upsilon}_1 \to \Gamma$ be the monoid homomorphism 
sending $\alpha$ to $\gamma$. 
Then $\psi$ induces the following diagram which is a fibre product:
\[
\begin{array}{ccc}
  {\rm Rep}_2(\Gamma)_{\rm s.s., \gamma} & \stackrel{\pi_{\Gamma, \rm s.s., \gamma}}{\to} & 
{\rm Ch}_2(\Gamma)_{\rm s.s., \gamma} \\
 \downarrow & & \downarrow \\
  {\rm Rep}_2(\Upsilon_1)_{\rm s.s., \alpha} & \to & 
{\rm Ch}_2(\Upsilon_1)_{\rm s.s., \alpha}. \\
\end{array}
\]
In particular, the morphism $\pi_{\Gamma, \rm s.s., \gamma}$ can 
be obtained by base change of the prototype for each $\gamma \in \Gamma$.
\end{lemma}

\prf
Put $X := {\rm Rep}_2(\Gamma)_{\rm s.s., \gamma}$ and 
$Y := {\rm Rep}_2(\Upsilon_1)_{\rm s.s., \alpha}\times_
{{\rm Ch}_2(\Upsilon_1)_{\rm s.s., \alpha}}
{\rm Ch}_2(\Gamma)_{\rm s.s., \gamma}$.
We shall show that the morphism   
$f : X \to Y$ induced by $\psi$ is an isomorphism.
Let $R$ be a commutative ring. 
For two $R$-valued points $\sigma_1, \sigma_2$ of $X$, 
assume that $f(\sigma_1) = f(\sigma_2)$ as $R$-valued points of 
$Y$.
Considering $\sigma_1$ and $\sigma_2$ as 
representations of degree $2$ for $\Gamma$ 
in $R$, we have 
\[
\sigma_i(\delta) = (I_2, \sigma_i(\gamma)) 
\left(
  \begin{array}{cc}
    {\rm tr}(I_2) & {\rm tr}(\sigma_i(\gamma)) \\
    {\rm tr}(\sigma_i(\gamma)) & {\rm tr}(\sigma_i(\gamma^2)) 
  \end{array}
\right)^{-1}
\left(
  \begin{array}{c}
    {\rm tr}(\sigma_i(\delta)) \\
    {\rm tr}(\sigma_i(\delta\gamma)) \\
  \end{array}
\right)
\]   
for each $\delta \in \Gamma$ and $i = 1, 2$.
Since $f(\sigma_1) = f(\sigma_2)$, we obtain 
$\sigma_1(\gamma) = f(\sigma_1)(\alpha) = f(\sigma_2)(\alpha) = \sigma_2(\gamma)$ and 
${\rm tr}(\sigma_1(\delta')) = {\rm tr}(\sigma_2(\delta'))$
for each $\delta' \in \Gamma$.
Hence we have $\sigma_1 = \sigma_2$. 

Let $y = (\rho, \chi)$ be an $R$-valued point of $Y$.
Here $\rho : \Upsilon_1 \to {\rm M}_2(R)$ and $\chi$ are $R$-valued points of 
${\rm Rep}_2(\Upsilon_1)_{\rm s.s., \alpha}$ 
and 
${\rm Ch}_2(\Gamma)_{\rm s.s., \gamma}$, respectively. 
Let us denote by $\chi(\delta)$ the image of ${\rm tr}(\sigma_{\Gamma}(\delta))$ 
by the ring homomorphism $\phi : (A_2(\Gamma)_{{\rm rk} \le 2}^{\rm Ch})_{m(\sigma_{\Gamma}(\gamma))} 
\to R$ associated to $\chi$.  
Then $\phi(m(\sigma_{\Gamma}(\gamma))) = m(\rho(\alpha))$ and 
$\chi(\gamma^m) = {\rm tr}(\rho(\alpha^m))$ for 
$m \in {\mathbb N}$.  
We define the map $\sigma : \Gamma \to {\rm M}_2(R)$ by
\[
\sigma(\delta) := (I_2, \rho(\alpha)) 
\left(
  \begin{array}{cc}
    {\rm tr}(I_2) & {\rm tr}(\rho(\alpha)) \\
    {\rm tr}(\rho(\alpha)) & {\rm tr}(\rho(\alpha^2)) 
  \end{array}
\right)^{-1}
\left(
  \begin{array}{c}
    \chi(\delta) \\
    \chi(\delta\gamma) \\
  \end{array}
\right)
\]
for $\delta \in \Gamma$. 
It is easy to see that 
$\sigma(e) = I_2$ and $\sigma(\gamma) = \rho(\alpha)$. 

Note that 
\[
\sigma_{\Gamma}(\delta) = (I_2, \sigma_{\Gamma}(\gamma)) 
\left(
  \begin{array}{cc}
    {\rm tr}(I_2) & {\rm tr}(\sigma_{\Gamma}(\gamma)) \\
    {\rm tr}(\sigma_{\Gamma}(\gamma)) & {\rm tr}(\sigma_{\Gamma}(\gamma^2)) 
  \end{array}
\right)^{-1}
\left(
  \begin{array}{c}
    {\rm tr}(\sigma_{\Gamma}(\delta)) \\
    {\rm tr}(\sigma_{\Gamma}(\gamma\delta)) \\
  \end{array}
\right)
\]
in ${\rm M}_2((A_2(\Gamma)_{{\rm rk} \le 2})_{m(\sigma_{\Gamma}(\gamma))})$ 
for each $\delta \in \Gamma$. 
By a similar discussion in the proof of \cite[Theorem~5.1]{Nkmt00}, 
we see that $\sigma$ is a representation, and hence 
that $\sigma$ can be regarded as an $R$-valued point of $X$.
Since ${\rm tr}(\sigma(\delta)) = \chi(\delta)$ for each $\delta \in \Gamma$, 
we also see that $f(\sigma) = y$ by using Proposition \ref{lemmassgen}. 
Therefore $f$ is an isomorphism.
\qed

\begin{remark}\rm 
We can also prove the group version of Lemma \ref{lemssror}: 
Let $\gamma$ be an element of a group $\Gamma$.
Let $\phi : {\rm F}_1 \to \Gamma$ be the group homomorphism 
sending $\alpha$ to $\gamma$.
Then $\phi$ induces the following diagram which is a fibre product:
\[
\begin{array}{ccc}
  {\rm Rep}_2(\Gamma)_{\rm s.s., \gamma} & \to & 
{\rm Ch}_2(\Gamma)_{\rm s.s., \gamma} \\
 \downarrow & & \downarrow \\
  {\rm Rep}_2({\rm F}_1)_{\rm s.s., \alpha} & \to & 
{\rm Ch}_2({\rm F}_1)_{\rm s.s., \alpha}. \\
\end{array}
\]
In particular,  
the morphism $\pi_{\Gamma, \rm s.s., \gamma}$ can 
be obtained by base change of the prototype for 
group representations.
\end{remark}

\begin{corollary}
The morphism 
$\pi_{\Gamma, {\rm s.s.}} : {\rm Rep}_2(\Gamma)_{\rm s.s.}
\to {\rm Ch}_2(\Gamma)_{\rm s.s.}$ is a universal 
geometric quotient 
by ${\rm PGL}_2$ for a group or a monoid $\Gamma$. 
\end{corollary}

\prf 
From Lemma \ref{lemssror} and Corollary \ref{corssf2} we 
see that $\pi_{\Gamma, {\rm s.s.}}$ gives a universal 
geometric quotient.
\qed 

\begin{remark}\rm
  The morphism $\pi_{\Gamma, {\rm s.s.}} : {\rm Rep}_2(\Gamma)_{\rm s.s.} \to 
{\rm Ch}_2(\Gamma)_{\rm s.s.}$ is smooth and surjective. 
Indeed, the prototype $\pi_{\Upsilon_1, {\rm s.s.}} : {\rm Rep}_2(\Upsilon_1)_{\rm s.s.} 
\to {\rm Ch}_2(\Upsilon_1)_{\rm s.s.}$ 
is smooth and surjective because it is obtained by base change of 
$\pi : {\rm Rep}_2(\Upsilon_{1})_{{\rm rk} 2} \to {\rm Ch}_2(\Upsilon_1)$ and 
$\pi$ is smooth and surjective by Proposition~\ref{lemmassff}.  
\end{remark}

\begin{remark}\label{remark:section}\rm 
For each point $x \in {\rm Ch}_2(\Gamma)_{\rm s.s.}$, there exists 
a local section $s_{x} : V_x \to {\rm Rep}_2(\Gamma)_{\rm s.s.}$ on 
a neighbourhood $V_x$ of $x$ such that 
$\pi_{\Gamma, {\rm s.s.}} \circ s_{x} = id_{V_x}$. 
Indeed, take $\gamma \in \Gamma$ such that 
$x \in {\rm Ch}_2(\Gamma)_{\rm s.s., \gamma}$. By Lemma \ref{lemssror}, 
$\pi_{\Gamma, {\rm s.s.}} : {\rm Rep}_2(\Gamma)_{{\rm s.s.}, \gamma} 
\to {\rm Ch}_2(\Gamma)_{{\rm s.s.}, \gamma}$ has a section $s_{\Gamma, \gamma}$ because  
${\rm Rep}_2(\Upsilon_1)_{{\rm s.s.}} \to {\rm Ch}_2(\Upsilon_1)_{\rm s.s.}$ 
has a section $s$.  
Hence we can take ${\rm Ch}_2(\Gamma)_{{\rm s.s.}, \gamma}$ 
as a neighbourhood $V_x$ of $x$. 
\end{remark}


\begin{lemma}\label{lemma:localeqss} 
Let $\rho_1, \rho_2$ be representations with semi-simple mold for a group (or a monoid) 
$\Gamma$ on a scheme $X$. 
Let $f_i : X \to {\rm Rep}_2(\Gamma)_{\rm s.s.}$ be the morphism 
associated to $\rho_i$ for $i = 1, 2$. 
If $\pi_{\Gamma, {\rm s.s.}} \circ f_1 = \pi_{\Gamma, {\rm s.s.}} \circ f_2 
: X \to {\rm Ch}_2(\Gamma)_{\rm s.s.}$, then for each $x \in X$ 
there exists $P_x \in {\rm GL}_2(\Gamma(V_x, {\mathcal O}_X))$ 
on a neighbourhood $V_x$ of $x$ such that 
$P_x^{-1}\rho_1 P_x = \rho_2$ on $V_x$.  
\end{lemma} 

\prf
For $x \in X$, take $\gamma \in \Gamma$ such that 
$(\pi_{\Gamma, {\rm s.s.}} \circ f_1) (x) = (\pi_{\Gamma, {\rm s.s.}} \circ f_2) (x) 
\in {\rm Ch}_2(\Gamma)_{\rm s.s., \gamma}$.  
We may assume that $f_i : X \to {\rm Rep}_2(\Gamma)_{{\rm s.s.}, \gamma}$ for 
$i =1, 2$ from the beginning. By Remark \ref{remark:section}, 
$\pi_{\Gamma, {\rm s.s.}} : {\rm Rep}_2(\Gamma)_{{\rm s.s.}, \gamma} 
\to {\rm Ch}_2(\Gamma)_{{\rm s.s.}, \gamma}$ has a section $s_{\Gamma, \gamma}$.  
Let $\rho_3$ be the representations with semi-simple mold on $X$ associated to 
$s_{\Gamma, \gamma} \circ \pi_{\Gamma, {\rm s.s.}} \circ f_1 = s_{\Gamma, \gamma} \circ 
\pi_{\Gamma, {\rm s.s.}} \circ f_2$. 
Note that 
\[
\rho_3(\gamma) = \left(
\begin{array}{cc}
0 & -\det(\rho_1(\gamma)) \\
1 & {\rm tr}(\rho_1(\gamma)) \\
\end{array}
\right) 
=
\left(
\begin{array}{cc}
0 & -\det(\rho_2(\gamma)) \\
1 & {\rm tr}(\rho_2(\gamma)) \\
\end{array}
\right) 
\] and that 
${\rm tr}(\rho_1(\delta)) = {\rm tr}(\rho_2(\delta)) = {\rm tr}(\rho_3(\delta))$ 
for each $\delta \in \Gamma$.  
There exist $Q_1, Q_2 \in {\rm GL}_2(\Gamma(V_x, {\mathcal O}_X))$ on a
neighbourhood $V_x$ of $x$ such that 
$Q_1^{-1}\rho_1(\gamma)Q_1 = \rho_3(\gamma)$ and 
$Q_2^{-1}\rho_2(\gamma)Q_2 = \rho_3(\gamma)$ by Lemma \ref{lemma-norm-nonsca}. 
Since 
\[
\rho_i(\delta) = (I_2, \rho_i(\gamma)) 
\left(
  \begin{array}{cc}
    {\rm tr}(I_2) & {\rm tr}(\rho_i(\gamma)) \\
    {\rm tr}(\rho_i(\gamma)) & {\rm tr}(\rho_i(\gamma^2)) 
  \end{array}
\right)^{-1}
\left(
  \begin{array}{c}
    {\rm tr}(\rho_i(\delta)) \\
    {\rm tr}(\rho_i(\gamma\delta)) \\
  \end{array}
\right)
\]
on $V_x$ for $\delta \in \Gamma$ and for $i=1, 2, 3$, we have 
 $Q_1^{-1}\rho_1(\delta)Q_1 = \rho_3(\delta)$ and 
$Q_2^{-1}\rho_2(\delta)Q_2 = \rho_3(\delta)$ for each $\delta \in \Gamma$. 
Hence $(Q_1 Q_2^{-1})^{-1} \rho_1 (Q_1 Q_2^{-1}) = \rho_2$ on $V_x$. 
This completes the proof. 
\qed

\begin{theorem}\label{th:eqtr}  
  Let $R$ be a local ring.
For two representations with semi-simple mold 
$\rho_1, \rho_2 : \Gamma \to {\rm GL}_2(R)$ for a group (or a monoid) $\Gamma$, 
$\rho_1$ and $\rho_2$ are equivalent to each other (in other words, there exists 
$P \in {\rm GL}_2(R)$ such that $P^{-1} \rho_1(\gamma) P = \rho_2(\gamma)$ 
for any $\gamma \in \Gamma$) if and only if 
${\rm tr}(\rho_1(\gamma)) = {\rm tr}(\rho_2(\gamma))$ for each $\gamma \in \Gamma$.
\end{theorem}

\prf 
Let $f_1, f_2$ be the $R$-valued points of ${\rm Rep}_2(\Gamma)_{\rm s.s.}$ associated to 
$\rho_1, \rho_2$, respectively. 
Using Proposition \ref{lemmassgen} and 
$m(\rho_i(\gamma)) = 2{\rm tr}(\rho_i(\gamma^2)) - ({\rm tr}(\rho_i(\gamma)))^2$ for $i=1, 2$, we see 
that ${\rm tr}(\rho_1(\gamma)) = {\rm tr}(\rho_2(\gamma))$ 
for each $\gamma \in \Gamma$ if and only if 
$\pi_{\Gamma, {\rm s.s.}} \circ f_1 = \pi_{\Gamma, {\rm s.s.}} \circ f_2$ as $R$-valued points of 
${\rm Ch}_2(\Gamma)_{\rm s.s.}$. 
The statement follows from Lemma \ref{lemma:localeqss}. 
\qed 

\bigskip 

Let us define 
${\mathcal E}q\mathcal{SS}_2(\Gamma)$ as the  
sheafification of the following contravariant functor with respect to Zariski topology: 
\[
\begin{array}{ccl}
 ({\bf Sch})^{op} & \to & ({\bf Sets}) \\
 X & \mapsto & \{ \rho \mid \mbox{ rep. 
with s.s. mold for $\Gamma$ on } X \}/\sim.
\end{array}
\]

By a {\it generalized representation with semi-simple mold} for $\Gamma$ on 
a scheme $X$, we understand pairs $\{ (U_i, \rho_i) \}_{i \in I}$ 
of an open set $U_i$ and a representation with semi-simple mold 
$\rho_i : \Gamma \to {\rm M}_2(\Gamma(U_i, {\mathcal O}_{X}))$ 
satisfying the following two conditions:
\begin{enumerate}
\item $\cup_{i \in I} U_ i = X$, 
\item for each $x \in U_i \cap U_j$, there exists $P_x \in 
{\rm GL}_2(\Gamma(V_x, {\mathcal O}_X))$ on 
a neighbourhood $V_x \subseteq U_i \cap U_j$ of $x$ such that 
$P_x^{-1} \rho_i P_x = \rho_j$ on $V_x$. 
\end{enumerate}

Generalized representations with semi-simple mold 
$\{ (U_i, \rho_i) \}_{i \in I}$ and $\{ (V_j, \sigma_j) \}_{j \in J}$ 
are called {\it equivalent} if $\{ (U_i, \rho_i) \}_{i \in I} \cup \{ (V_j, \sigma_j) \}_{j \in J}$ 
is a generalized representation with semi-simple mold again.  
We easily see that 
${\mathcal E}q\mathcal{SS}_2(\Gamma)(X)$ is the 
set of equivalence classes of 
generalized representations with semi-simple mold 
for $\Gamma$ on a scheme $X$. 

\begin{theorem}\label{th:moduliss}  
The scheme ${\rm Ch}_2(\Gamma)_{\rm s.s.}$ is a fine moduli scheme  
associated to the functor ${\mathcal E}q\mathcal{SS}_2(\Gamma)$ 
for a group or a monoid $\Gamma$:  
\[
\begin{array}{ccccl}
{\mathcal E}q\mathcal{SS}_2(\Gamma) & : & ({\bf Sch})^{op} & \to & ({\bf Sets}) \\
 & & X & \mapsto & 
\left\{ 
\begin{array}{r}  
\mbox{ gen. rep. 
with s.s. mold } \\ 
\mbox{ for $\Gamma$ on } X  
\end{array} 
\right\} 
\Big/\sim.
\end{array}
\] 
In other words, ${\rm Ch}_2(\Gamma)_{\rm s.s.}$ represents 
the functor ${\mathcal E}q\mathcal{SS}_2(\Gamma)$. 
The moduli ${\rm Ch}_2(\Gamma)_{\rm s.s.}$ is separated over ${\Bbb Z}$;
if $\Gamma$ is a finitely generated group or monoid, then 
${\rm Ch}_2(\Gamma)_{\rm s.s.}$ is of finite type over ${\Bbb Z}$. 
\end{theorem}

\prf 
It is easy to define a canonical morphism 
${\mathcal E}q\mathcal{SS}_2(\Gamma) \to 
h_{{\rm Ch}_2(\Gamma)_{\rm s.s.}} := {\rm Hom}(-, {\rm Ch}_2(\Gamma)_{\rm s.s.})$.  
Let us define a morphism $h_{{\rm Ch}_2(\Gamma)_{\rm s.s.}} 
\to {\mathcal E}q\mathcal{SS}_2(\Gamma)$.  
Let $g \in h_{{\rm Ch}_2(\Gamma)_{\rm s.s.}}(X)$ with a scheme $X$. 
For each $x \in X$, take $\gamma_x \in \Gamma$ such that 
$g(x) \in {\rm Ch}_2(\Gamma)_{{\rm s.s.}, \gamma_x}$. 
By using the section $s_{\Gamma, \gamma_x} : {\rm Ch}_2(\Gamma)_{{\rm s.s.}, \gamma_x} 
\to {\rm Rep}_2(\Gamma)_{{\rm s.s.}, \gamma_x}$ in Remark \ref{remark:section}, 
we can define a representation with semi-simple mold $\rho_x$ on a neighbourhood 
$U_x$ of $x$. By Lemma \ref{lemma:localeqss}, 
we see that $\{ (U_x, \rho_x) \}_{x \in X} \in 
{\mathcal E}q\mathcal{SS}_2(\Gamma)(X)$ and that the morphism 
$h_{{\rm Ch}_2(\Gamma)_{\rm s.s.}} 
\to {\mathcal E}q\mathcal{SS}_2(\Gamma)$ is well-defined.  
It is easy to see that ${\rm Ch}_2(\Gamma)_{\rm s.s.}$ represents 
the functor ${\mathcal E}q\mathcal{SS}_2(\Gamma)$. 

Since ${\rm Ch}_2(\Gamma)_{\rm s.s.}$ is an open subscheme of 
the affine scheme ${\rm Ch}_2(\Gamma)_{{\rm rk} \le 2}$, 
${\rm Ch}_2(\Gamma)_{\rm s.s.}$ is separated over ${\Bbb Z}$. 
Suppose that  $\Gamma$ is finitely generated. 
Then ${\rm Ch}_2(\Gamma)_{{\rm rk} \le 2}$ is of finite type over ${\Bbb Z}$
by Proposition \ref{prop:finitetype}. 
Hence ${\rm Ch}_2(\Gamma)_{\rm s.s.}$  
is also of finite type over ${\Bbb Z}$. 
\qed

\begin{remark}\label{remark:algss}\rm
Let $A$ be an associative algebra over a commutative ring $R$. 
For an $R$-scheme $X$, we say that an $R$-algebra homomorphism 
$\rho : A \to {\rm M}_2(\Gamma(X, {\mathcal O}_X))$ is a 
{\it $2$-dimensional representation} of $A$ on $X$. 
For a $2$-dimensional representation $\rho$ of $A$, 
$\rho$ is called a {\it representation with semi-simple mold} if 
the subalgebra $\rho(A)$ of ${\rm M}_2({\mathcal O}_X)$ generates  
a semi-simple mold on $X$. 
In a similar way as group or monoid cases, 
we can define generalized representations with semi-simple mold for $A$ 
on an $R$-scheme $X$. The contravariant functor 
${\mathcal E}q\mathcal{SS}_2(A)$ from the category of $R$-schemes to the 
category of sets is defined as 
\[
\begin{array}{ccccl}
{\mathcal E}q\mathcal{SS}_2(A) & : & ({\bf Sch}/R)^{op} & \to & ({\bf Sets}) \\
 & & X & \mapsto & 
\left\{ 
\begin{array}{r}  
\mbox{ gen. rep. 
with s.s. mold } \\ 
\mbox{ for $A$ on } X  
\end{array} 
\right\} 
\Big/\sim.
\end{array}
\] 
Then we can construct the fine moduli ${\rm Ch}_2(A)_{\rm s.s.}$ associated to 
${\mathcal E}q\mathcal{SS}_2(A)$ in the same way as Theorem \ref{th:moduliss}.  
The moduli ${\rm Ch}_2(A)_{\rm s.s.}$ is separated over $R$. 
If $A$ is a finitely generated algebra over $R$, then 
${\rm Ch}_2(A)_{\rm s.s.}$ is of finite type over $R$.  
For a local ring $S$ over $R$, we see that  
two representations with semi-simple mold 
$\rho_1, \rho_2 : A \to {\rm M}_2(S)$ are equivalent to each other 
(in other words, there exists 
$P \in {\rm GL}_2(S)$ such that $P^{-1} \rho_1(a) P = \rho_2(a)$ 
for any $a \in A$) if 
and only if ${\rm tr}(\rho_1(a)) = {\rm tr}(\rho_2(a))$ for each 
$a \in A$ (the associative algebra version of Theorem~\ref{th:eqtr}).  
\end{remark}

\begin{remark}\rm 
We have introduced the notion of 
generalized representations with semi-simple mold 
for describing the moduli functors ${\mathcal E}q\mathcal{SS}_2(\Gamma)$  
and ${\mathcal E}q\mathcal{SS}_2(A)$.  However, 
the moduli functors can also be described as 
${\mathcal E}q\mathcal{SS}'_2(\Gamma)$  
and ${\mathcal E}q\mathcal{SS}'_2(A)$ by using 
the notion of representations generating sheaves of algebras 
which define semi-simple molds.   
More precisely, see \S 8. 

\end{remark}

\section{Unipotent mold ($ch \neq 2$ case)}

Recall that a rank $2$ mold ${\mathcal A} \subseteq {\rm M}_2({\mathcal O}_X)$ 
over a ${\Bbb Z}[ 1/2 ]$-scheme $X$ is called {\it unipotent} if 
$m(s) := {\rm tr}(s)^2 - 4\det(s) = 0$ for each open subset $U \subseteq X$ 
and for each $s \in {\mathcal A}(U)$. 
In this section, all schemes are over ${\rm Spec} \: {\Bbb Z}[ 1/2 ]$ 
and all commutative rings are over ${\Bbb Z}[ 1/2 ]$. 
We construct the moduli of representations with unipotent mold
over ${\rm Spec} \: {\Bbb Z}[ 1/2 ]$.

As seen in Theorem~\ref{corssf1}, $\pi : {\rm Rep}_2(\Upsilon_1)_{{\rm rk} 2} \to {\rm Ch}_2(\Upsilon_1)$ 
is a universal geometric quotient by ${\rm PGL}_2$. Put $A := \sigma_{\Upsilon_1}(\alpha)$ 
for the universal representation $\sigma_{\Upsilon_1}$ of $\Upsilon_1 = \langle \alpha \rangle$. 
Let $Z$ be the closed subscheme of ${\rm Ch}_2(\Upsilon_1)$ 
defined by $m(A) = {\rm tr}(A)^2 - 4\det(A) = 0$. 
By base change, we obtain a universal geometric quotient 
$\pi' : {\rm Rep}_2(\Upsilon_1)_{{\rm rk} 2}\times_{{\rm Ch}_2(\Upsilon_1)} Z \to Z$ by 
${\rm PGL}_2$. 
However, this quotient $\pi'$ is not so good, because 
$Z$ has a singular fibre over ${\Bbb F}_2$ which is defined by ${\rm tr}(A)^2 = 0$.  
Therefore we assume that all schemes are over ${\rm Spec}{\Bbb Z}[ 1/2 ]$ in 
this section. The case of unipotent molds over ${\Bbb F}_2$ will be 
discussed in the next section. 

Assume that $R$ is a ${\Bbb Z}[ 1/2 ]$-algebra and that 
$A \subseteq {\rm M}_2(R)$ is a unipotent mold over $R$ through 
this section. 

\begin{notation}\rm
For $X \in {\rm M}_2(R)$, we define 
$\displaystyle \eta(X) := X- \frac{{\rm tr}(X)}{2}I_2$.  
\end{notation}

\begin{lemma}
Suppose that $X \in {\rm M}_2(R)$ satisfies 
$m(X) = {\rm tr}(X)^2 - 4\det(X) = 0$.
Then $\eta(X)^2 = 0$. 
\end{lemma}

{\it Proof.}
By the Cayley-Hamilton theorem, we have 
\begin{eqnarray*}
\eta(X)^2 & = & X^2 - {\rm tr}(X)X + \frac{({\rm tr}(X))^2}{4}I_2 \\
   & = & {\rm tr}(X)X -\det(X)I_2 - {\rm tr}(X)X + 
\frac{({\rm tr}(X))^2}{4}I_2 \\
 & = & 0, 
\end{eqnarray*} 
since ${\rm tr}(X)^2 = 4\det(X)$. 
\qed

\begin{lemma}\label{lemma-u-ra}
Let $R$ be a ${\Bbb Z}[ 1/2 ]$-algebra.
Let $A \subseteq {\rm M}_2(R)$ be a unipotent mold over $R$.
If $X, Y \in A$, then $2{\rm tr}(XY)={\rm tr}(X){\rm tr}(Y)$. 
\end{lemma}

{\it Proof.}
Since we have only to prove that the equality holds locally, 
we may assume that there exists $Z \in A$ such that 
$A = R\cdot I_2 + R\cdot Z$ and 
$m(Z) = {\rm tr}(Z)^2-4\det(Z)=0$. 
Put $X = a I_2 + b Z$ and $Y= c I_2 + d Z$. 
Then we have 
\begin{eqnarray*}
2{\rm tr}(XY) & = & 2ac\,{\rm tr}(I_2) + 2(ad+bc)\,{\rm tr}(Z) + 
2bd\,{\rm tr}(Z^2) \\
& = & 4ac + 2(ad+bc)\,{\rm tr}(Z) + 2bd\,({\rm tr}(Z))^2 - 
4bd\,\det(Z)  \\
& = & 4ac + 2(ad+bc)\,{\rm tr}(Z) + 
4bd\,\det(Z)
\end{eqnarray*}
and
\begin{eqnarray*}
{\rm tr}(X){\rm tr}(Y) & = & (2a + b\,{\rm tr}(Z))(2c + d\,{\rm tr}(Z)) \\
& = & 
4ac + 2(ad + bc)\,{\rm tr}(Z) + 
bd\,({\rm tr}(Z))^2 \\
& = & 4ac + 2(ad+bc)\,{\rm tr}(Z) + 
4bd\,\det(Z).
\end{eqnarray*}
This completes the proof. 
\qed

\begin{notation}\rm
Let $R$ and $A$ be as above. 
For $X \in {\rm M}_2(R)$, we denote ${\rm tr}(X)/2$ by $r(X)$. 
Note that $\eta(X) = X - r(X)I_2$. 
From Lemma \ref{lemma-u-ra}, we have 
$r(XY)=r(X)r(Y)$ for $X, Y \in A$.  
If $X = aI_2 + bZ$, then 
$\eta(X) = aI_2+bZ - r(aI_2+bZ)I_2 = b(Z -r(Z)I_2) = b\eta(Z)$. 
\end{notation}

\begin{lemma}\label{lemma-u-ra2}
Let $R$ and $A$ be as in Lemma \ref{lemma-u-ra}.
For $X, Y \in A$, 
\begin{eqnarray*}
\eta(XY) & = & r(X)\eta(Y) + \eta(X)r(Y). 
\end{eqnarray*}
\end{lemma}

{\it Proof.}
As in the proof of Lemma \ref{lemma-u-ra}, 
we may assume that there exists $Z \in A$ such that 
$A = R\cdot I_2 + R\cdot Z$ and 
$m(Z) = {\rm tr}(Z)^2-4\det(Z)=0$. 
For $X, Y \in A$, there exists $\lambda, \mu \in R$ such that 
$\eta(X) = \lambda \eta(Z)$ and $\eta(Y) = \mu \eta(Z)$. 
Since $\eta(X)\eta(Y) = \lambda\mu \eta(Z)^2 = 0$, we have 
\begin{eqnarray*}
\eta(XY) & = & XY - r(XY) I_2 \\
 & = & (X - r(X)I_2)(Y- r(Y)I_2) + r(Y)(X - r(X)I_2)  \\ 
 &   & \hspace*{10ex} + r(X)(Y -r(Y)I_2)  \\ 
 & = & \eta(X)\eta(Y) + r(X)\eta(Y) + \eta(X)r(Y) \\
 & = & r(X)\eta(Y) + \eta(X)r(Y). 
\end{eqnarray*}
This completes the proof. 
\qed 

\bigskip

\begin{notation}\rm
Let $\Gamma$ be a group or a monoid. 
Let $\rho : \Gamma \to {\rm M}_2(R)$ be 
a representation with the unipotent mold $A$. 
For each $\gamma \in 
\Gamma$, 
we denote $\eta(\rho(\gamma))$ and $r(\rho(\gamma))$ by $\eta(\gamma)$ and 
$r(\gamma)$, respectively.   
Assume that there exists $\alpha \in \Gamma$ such that 
$A = R\cdot I_2 + R\cdot \rho(\alpha)$.
Then $A = R\cdot I_2 + R\cdot \eta(\alpha)$ 
and for each $\gamma \in \Gamma$ 
there exists a unique $a_{\alpha}(\gamma) \in R$ such that 
$\eta(\gamma) = a_{\alpha}(\gamma)\eta(\alpha)$. 
\end{notation}

\begin{remark}\rm\label{remark-u-der}
Note that the map $r(\cdot) : \Gamma \to R$ is a 
character of $\Gamma$. In other words, 
$r(e) = 1$ and $r(\gamma\delta) = r(\gamma)r(\delta)$ for $\gamma, \delta \in \Gamma$. 
From Lemma \ref{lemma-u-ra2} 
we see that  
the map $a_{\alpha}(\cdot) : \Gamma \to R$ is a derivation 
with respect to $r$, that is, 
$a_{\alpha}$ satisfies the condition 
$a_{\alpha}(\gamma\delta) = r(\gamma)a_{\alpha}(\delta) + 
a_{\alpha}(\gamma)r(\delta)$ for each $\gamma, \delta \in \Gamma$. 
\end{remark}

\bigskip

For a representation $\rho : \Gamma \to {\rm M}_2(R)$ with 
the unipotent mold $A$ such that 
$A = R\cdot I_2 + R\cdot \rho(\alpha)$ 
for some $\alpha \in \Gamma$, 
we have a character $r : \Gamma \to R$ and 
a derivation $a_{\alpha} : \Gamma \to R$ with 
respect to $r$.  
Conversely, a character and a derivation give us a representation 
with unipotent mold.

\begin{lemma}\label{lemma-u-ror}     
Let $r : \Gamma \to R$ be a character  
and let $a : \Gamma \to R$ be a derivation 
with respect to $r$. 
Assume that there exists $\alpha \in \Gamma$ such that 
$a(\alpha) \in R^{\times}$.
Furthermore assume that there exists 
$Z \in {\rm M}_2(R)$ such that  
$A := R\cdot I_2 + R\cdot Z \subseteq {\rm M}_2(R)$ is a unipotent mold.
Then the map 
\[
\begin{array}{ccccl}
\rho & : & \Gamma & \to & {\rm M}_2(R) \\
 & & \gamma & \to &  r(\gamma)I_2 + a(\gamma)\eta(Z)
\end{array}
\] 
is a representation for $\Gamma$ with the unipotent mold $A$.
\end{lemma}

{\it Proof.}
For $\gamma, \delta \in \Gamma$, 
we have 
\begin{eqnarray*}
  \rho(\gamma)\rho(\delta) & = & (r(\gamma)I_2 + a(\gamma)\eta(Z))
(r(\delta)I_2 + a(\delta)\eta(Z)) \\
 & = & r(\gamma)r(\delta)I_2 + a(\gamma)r(\delta)\eta(Z) 
+ r(\gamma)a(\delta)\eta(Z) \\
 & = &  r(\gamma\delta)I_2 + a(\gamma\delta)\eta(Z) \\
 & = & \rho(\gamma\delta).  
\end{eqnarray*}
Since $\rho(e) = I_2$ and $\rho(\alpha) = (r(\alpha) - 
a(\alpha){\rm tr}(Z)/2)I_2 + a(\alpha)Z$, 
the map $\rho$ is a representation with the unipotent mold $A$.
\qed

\bigskip

\begin{definition}\rm
Let us denote ${\rm Rep}_2(\Gamma)\otimes_{\Bbb Z}{\Bbb Z}[ 1/2 ]$ by 
${\rm Rep}_2(\Gamma)[ 1/2 ]$. 
We define the subscheme  ${\rm Rep}_2(\Gamma)_u$ 
of ${\rm Rep}_2(\Gamma)[ 1/2 ]$ 
by 
\[
{\rm Rep}_2(\Gamma)_u := \{ \rho \in 
{\rm Rep}_2(\Gamma)[ 1/2 ] 
\mid \rho \mbox{ has a unipotent mold } \}. 
\]
We call ${\rm Rep}_2(\Gamma)_u$
the {\it unipotent part} of the representation variety of 
degree $2$ for $\Gamma$ over ${\Bbb Z}[ 1/2 ]$. 
Recall that ${\rm Rep}_2(\Gamma)_{{\rm rk} \le 2}$ 
is a closed subscheme of ${\rm Rep}_2(\Gamma)$ and that 
${\rm Rep}_2(\Gamma)_{{\rm rk} 2}$ is an open subscheme of ${\rm Rep}_2(\Gamma)_{{\rm rk} \le 2}$ 
(Definition \ref{def:reprkle2}).  Set 
${\rm Rep}_2(\Gamma)_{{\rm rk} \le 2}[1/2] := {\rm Rep}_2(\Gamma)_{{\rm rk} \le 2}\otimes_{{\Bbb Z}} {\Bbb Z}[1/2]$ 
and 
${\rm Rep}_2(\Gamma)_{{\rm rk} 2}[1/2] := {\rm Rep}_2(\Gamma)_{{\rm rk} 2}\otimes_{{\Bbb Z}} {\Bbb Z}[1/2]$. 
Then ${\rm Rep}_2(\Gamma)_u$ is a closed subscheme of ${\rm Rep}_2(\Gamma)_{{\rm rk} 2}[1/2]$ 
defined by $m(\sigma_{\Gamma}(\gamma)) = 0$ for all $\gamma \in \Gamma$. 
\end{definition}

\bigskip

\begin{definition}\rm
Let us ${\rm Rep}_1(\Gamma)[ 1/2 ]$ 
denote the representation variety  ${\rm Rep}_1(\Gamma)\otimes_{\Bbb Z}{\Bbb Z}[ 1/2 ]$ 
of degree $1$ for $\Gamma$ over ${\Bbb Z}[ 1/2 ]$. 
Let us denote by $A_{1}(\Gamma)[1/2]$ the coordinate ring 
of ${\rm Rep}_1(\Gamma)[ 1/2 ]$. 
For an $A_{1}(\Gamma)[1/2]$-module $M$, 
we define the {\it $A_{1}(\Gamma)[1/2]$-module 
of derivations} 
by 
\[
{\rm Der}(\Gamma, M) := 
\left\{     a : \Gamma \to M \; 
  \begin{array}{|c}
    a(\gamma\delta) = \chi_{\Gamma}(\gamma)a(\delta) + a(\gamma) 
   \chi_{\Gamma}(\delta)  \\
   \mbox{ for each } \gamma, \delta \in \Gamma 
  \end{array}
\right\}.
\]
Here we denote by $\chi_{\Gamma} : \Gamma \to 
A_{1}(\Gamma)[1/2]$ the universal representation of 
degree $1$ for $\Gamma$. 
\end{definition}

\bigskip

\begin{lemma}\label{lemma:omegamodule} 
  There exists a universal $A_1(\Gamma)[1/2]$-module $\Omega_{\Gamma}$ 
representing the covariant functor 
\[
\begin{array}{ccccc}
{\rm Der}(\Gamma, -) & : & (A_1(\Gamma)[1/2]\mbox{-}{\bf Mod}) & \to &  
(A_1(\Gamma)[1/2]\mbox{-}{\bf Mod}) \\ 
 & & M & \mapsto & {\rm Der}(\Gamma, M). 
\end{array} 
\]
In particular,  
\[
{\rm Der}(\Gamma, M) \stackrel{\sim}{\to} 
{\rm Hom}_{A_1(\Gamma)[1/2]}( \Omega_{\Gamma}, M)
\]
is an isomorphism for each $A_{1}(\Gamma)[1/2]$-module $M$. 
\end{lemma}

{\it Proof.} 
We define the $A_1(\Gamma)[1/2]$-module $\Omega_{\Gamma}$ 
by 
\[
  \Omega_{\Gamma} := (\oplus_{\gamma \in \Gamma} A_{1}(\Gamma)[1/2]\cdot d\gamma) 
/ N, 
\]
where $N$ is the $A_1(\Gamma)[1/2]$-submodule generated by 
$\{ \chi_{\Gamma}(\gamma) d\delta + \chi_{\Gamma}(\delta) d\gamma -
d(\gamma\delta) \mid \gamma, \delta \in \Gamma \}$ of the free 
$A_{1}(\Gamma)[1/2]$-module 
$\displaystyle \oplus_{\gamma \in \Gamma} A_{1}[1/2](\Gamma)\cdot d\gamma$. 
It is easy to check that $\Omega_{\Gamma}$ represents the above functor. 
\qed

\bigskip 

\begin{remark}\label{remark:fgu}\rm 
If $\Gamma$ is a finitely generated group or monoid, then  
$A_1(\Gamma)[1/2]$ is a finitely generated algebra over ${\Bbb Z}[1/2]$ and 
$\Omega_{\Gamma}$ is a finitely generated $A_1(\Gamma)[1/2]$-module. 
Indeed, let $S = \{ \alpha_1, \cdots, \alpha_n \}$ be a set of generators of $\Gamma$. 
We may assume that $\alpha_{i}^{-1}$ is also an element of $S$ for 
each $1 \le i \le n$ if $\Gamma$ is a group. 
Then $A_1(\Gamma)[1/2]$ is generated by 
$\{ \chi_{\Gamma}(\alpha_1), \ldots, \chi_{\Gamma}(\alpha_n) \}$ 
over ${\Bbb Z}[1/2]$ and $\Omega_{\Gamma}$ is generated by 
$\{ d(\alpha_1), \ldots, d(\alpha_n) \}$ over $A_1(\Gamma)[1/2]$.    
\end{remark}

\bigskip

\begin{definition}\rm
  We define the scheme ${\rm Ch}_2(\Gamma)_u$ over $A_{1}(\Gamma)[1/2]$ 
by 
\[
{\rm Ch}_2(\Gamma)_u := {\rm Proj}S(\Omega_{\Gamma}),  
\]
where $S(\Omega_{\Gamma})$ is the the symmetric algebra 
of $\Omega_{\Gamma}$ 
over $A_{1}(\Gamma)[1/2]$. 
\end{definition}

\bigskip

\begin{example}\label{example:omegau1}\rm
  Let ${\Upsilon}_1 = \langle \alpha_0 \rangle$ 
be the free monoid of rank $1$. 
The $A_{1}(\Upsilon_1)[1/2]$-module $\Omega_{\Upsilon_1}$ is 
isomorphic to $A_{1}(\Upsilon_1)[1/2]$. Indeed, 
the $A_{1}(\Upsilon_1)[1/2]$-module homomorphism   
\[
\begin{array}{ccc}
  A_{1}(\Upsilon_1)[1/2] & \to & \Omega_{\Upsilon_1} \\
  1  & \mapsto & d\alpha_{0}  
\end{array}
\]
gives an isomorphism. 
In particular, ${\rm Ch}_2(\Upsilon_1)_u
\cong {\rm Rep}_1(\Upsilon_1)[1/2]$.  
\end{example}

\bigskip

Let $\psi : X \to {\rm Rep}_{1}(\Gamma)[1/2]$ 
be a ${\mathbb Z}[1/2]$-morphism. 
Let us regard $\Omega_{\Gamma}$ as 
a quasi-coherent sheaf on ${\rm Rep}_1(\Gamma)[1/2]$.  
There exists a one-to-one correspondence  
\begin{multline*}
{\rm Hom}_{ \: {\rm Rep}_{1}(\Gamma)[1/2] \: }(X, {\rm Ch}_2(\Gamma)_{u})  \cong \\ 
\{ 
\psi^{\ast}(\Omega_{\Gamma}) \twoheadrightarrow {\mathcal L} \to 0  
\mid {\mathcal L} \mbox{ is a line bundle on } X  
\}/\sim.   
\end{multline*} 
Here we say that $\psi^{\ast}(\Omega_{\Gamma}) \stackrel{f_1}{\twoheadrightarrow} {\mathcal L}_1$ and 
$\psi^{\ast}(\Omega_{\Gamma}) \stackrel{f_2}{\twoheadrightarrow} {\mathcal L}_2$ are 
equivalent if there exists an isomorphism $g : {\mathcal L}_1 
\stackrel{\cong}{\to} {\mathcal L}_2$ such that $g \circ f_1 = f_2$.

\bigskip 

The group scheme ${\rm PGL}_2[1/2] := {\rm PGL}_2\otimes_{{\mathbb Z}} {\mathbb Z}[1/2]$ over ${\Bbb Z}[1/2]$ 
acts on ${\rm Rep}_2(\Gamma)_u$ by $\rho \mapsto P^{-1}\rho P$. 
We define the morphism $\lambda : 
{\rm Rep}_2(\Gamma)_{u} \to {\rm Rep}_1(\Gamma)[1/2]$ 
by $\rho \mapsto r = {\rm tr}(\rho)/2$. 
For $\alpha \in \Gamma$, we define the open subscheme 
${\rm Rep}_2(\Gamma)_{u, \alpha}$ of ${\rm Rep}_2(\Gamma)_u$ 
by 
\[
{\rm Rep}_2(\Gamma)_{u, \alpha} := 
\{ \rho \in {\rm Rep}_2(\Gamma)_u \mid \langle I_2, 
\rho(\alpha) \rangle \mbox{ generates a unipotent mold }  \}. 
\]
Then ${\rm Rep}_2(\Gamma)_{u, \alpha}$ is a ${\rm PGL}_2[1/2]$-invariant 
open subscheme of ${\rm Rep}_2(\Gamma)_{u}$. 
The derivation $a_{\alpha} : \Gamma \to 
\Gamma( {\rm Rep}_2(\Gamma)_{u, \alpha}, 
{\mathcal O}_{{\rm Rep}_2(\Gamma)_{u, \alpha}})$ in Remark \ref{remark-u-der}
induces  
the $A_{1}(\Gamma)[1/2]$-module homomorphism 
$\lambda^{\ast}(\Omega_{\Gamma}) \to {\mathcal O}_{{\rm Rep}_2(\Gamma)_{u, \alpha}}$, and 
hence we can define  
the morphism $\pi_{\alpha} : {\rm Rep}_2(\Gamma)_{u, \alpha} 
\to {\rm Ch}_2(\Gamma)_u$ over ${\rm Rep}_1(\Gamma)[1/2]$ 
associated to $\lambda^{\ast}(\Omega_{\Gamma}) \twoheadrightarrow 
{\mathcal O}_{{\rm Rep}_2(\Gamma)_{u, \alpha}}$. 
Gluing the morphisms $\{ \pi_{\alpha} \}_{\alpha \in \Gamma}$, 
we have
the morphism $\pi_{\Gamma, u} : {\rm Rep}_2(\Gamma)_{u} \to 
{\rm Ch}_2(\Gamma)_u$ over ${\rm Rep}_1(\Gamma)[1/2]$.  

\bigskip 

For $\alpha \in \Gamma$, 
we define the open subscheme ${\rm Ch}_2(\Gamma)_{u, \alpha}$ 
of ${\rm Ch}_2(\Gamma)_u$ by 
${\rm Ch}_2(\Gamma)_{u, \alpha} := D(d\alpha) = \{ d\alpha \neq 0 \}$. 
From the definition of $\pi_{\Gamma, u}$, we see that 
$\pi_{\Gamma, u}^{-1}({\rm Ch}_2(\Gamma)_{u, \alpha}) = {\rm Rep}_2(\Gamma)_{u, \alpha}$. 
For a ${\mathbb Z}[1/2]$-morphism $\psi : X \to {\rm Rep}_{1}(\Gamma)[1/2]$,   
there exists a one-to-one correspondence  
\begin{multline*}
{\rm Hom}_{ \: {\rm Rep}_{1}(\Gamma)[1/2] \: }(X, {\rm Ch}_2(\Gamma)_{u, \alpha})  \cong \\ 
\left\{ 
\psi^{\ast}(\Omega_{\Gamma}) \twoheadrightarrow {\mathcal L} \to 0  \; 
\begin{array}{|l}     
 {\mathcal L} \mbox{ is a line bundle on } X \mbox{ and }  \psi^{\ast}(d\alpha) \mbox{ is } \\ 
 \mbox{ nowhere vanishing as a section of } {\mathcal L}   
\end{array}
\right\}\bigg/\sim.   
\end{multline*} 
Since ${\mathcal L}$ is generated by $\psi^{\ast}(d\alpha)$, 
${\mathcal L}$ is isomorphic to ${\mathcal O}_X$.  
Let $r : \Gamma \to \Gamma(X, {\mathcal O}_X)$ be the character associated to 
$\psi : X \to {\rm Rep}_{1}(\Gamma)[1/2]$. 
Regarding $\psi^{\ast}(d\alpha)$ as $1$ of ${\mathcal O}_X$, 
we have the following: 
\begin{multline*}
{\rm Hom}_{ \: {\rm Rep}_{1}(\Gamma)[1/2] \: }(X, {\rm Ch}_2(\Gamma)_{u, \alpha})  \cong \\ 
\left\{ 
\; d \;\; 
\begin{array}{|l}     
 d  : \Gamma \to \Gamma(X, {\mathcal O}_X) \mbox{ is a derivation with respect to } r  \\ 
\mbox{ such that  } d(\alpha)=1   
\end{array}
\right\}.   
\end{multline*} 
Remark that $\pi_{\Gamma, u} : {\rm Rep}_2(\Gamma)_{u} \to 
{\rm Ch}_2(\Gamma)_u$ and 
$\pi_{\alpha} : {\rm Rep}_2(\Gamma)_{u, \alpha} \to 
{\rm Ch}_2(\Gamma)_{u, \alpha}$ are 
${\rm PGL}_2[1/2]$-equivariant morphisms, where 
the actions of ${\rm PGL}_2[1/2]$ on 
${\rm Ch}_2(\Gamma)_u$ and ${\rm Ch}_2(\Gamma)_{u, \alpha}$ 
are trivial. 

\bigskip

\begin{definition}\rm
  For the free monoid $\Upsilon_1 = \langle \alpha_0 \rangle$ of rank $1$, we say that  
the morphism $\pi_{\Upsilon_1, u} : {\rm Rep}_{2}(\Upsilon_1)_u \to 
{\rm Ch}_2(\Upsilon_1)_{u}$ is the {\it prototype} in the unipotent mold case.   
Remark that ${\rm Rep}_{2}(\Upsilon_1)_u = {\rm Rep}_{2}(\Upsilon_1)_{u, \alpha_0}$ 
and that ${\rm Ch}_2(\Upsilon_1)_{u} = {\rm Ch}_2(\Upsilon_1)_{u, \alpha_0}$. 
\end{definition}

\bigskip 

By Theorem~\ref{corssf1}, $\pi : {\rm Rep}_2(\Upsilon_1)_{{\rm rk} 2} \to 
{\rm Ch}_2(\Upsilon_1)$ is a universal geometric quotient by 
${\rm PGL}_2$. Taking the base change of $\pi$ by 
${\rm Spec} \: {\Bbb Z}[1/2] \to {\rm Spec} \: {\Bbb Z}$, 
we have ${\rm Rep}_2(\Upsilon_1)_{{\rm rk} 2}[1/2] \to 
{\rm Ch}_2(\Upsilon_1)[1/2]$.  
Here we denote $X\otimes_{{\Bbb Z}} {\Bbb Z}[1/2]$ by $X[1/2]$ 
for a ${\Bbb Z}$-scheme $X$.  
Let $Z$ be the closed subscheme of ${\rm Ch}_2(\Upsilon_1)[1/2]$ 
defined by $m(\sigma_{\Upsilon}(\alpha_0)) = 0$. 
Since ${\rm Ch}_2(\Upsilon_1) = {\rm Spec} \: {\Bbb Z}[T, D]$, 
the affine ring of $Z$ is isomorphic to ${\Bbb Z}[1/2, T]$.  Note that 
$r(\cdot) = {\rm tr}(\sigma_{\Upsilon_1}(\cdot))/2$ gives a character of $\Upsilon_1$ on 
$Z$ and that $r(\alpha_0) = T/2$. 
Hence $Z$ is isomorphic to ${\rm Ch}_2(\Upsilon_1)_u
\cong {\rm Rep}_1(\Upsilon_1)[1/2]$. 
Taking the base change of ${\rm Rep}_2(\Upsilon_1)_{{\rm rk} 2}[1/2] \to 
{\rm Ch}_2(\Upsilon_1)[1/2]$ by $Z \hookrightarrow {\rm Ch}_2(\Upsilon_1)[1/2]$, 
we have  $\pi_{\Upsilon_1, u} : {\rm Rep}_2(\Upsilon_1)_{u} \to {\rm Ch}_2(\Upsilon_1)_u$. 

Here we introduce the following lemma without proof: 

\begin{lemma}\label{lemma:basechangeGIT} 
Let $X \to Y$ be a (universal) geometric quotient by $G$ over $S$. 
For $S' \to S$, put $X_{S'} := X \times_{S} S'$, $Y_{S'} := Y \times_{S} S'$, and 
$G_{S'} := G \times_{S} S'$.  Then 
$X_{S'} \to Y_{S'}$ is a (resp. universal) geometric quotient by $G_{S'}$ over $S'$. 
\end{lemma}

\bigskip 

By the lemma above, we have: 

\begin{theorem}
The prototype $\pi_{\Upsilon_1, u} : {\rm Rep}_2(\Upsilon_1)_{u} \to {\rm Ch}_2(\Upsilon_1)_u$ 
is a universal geometric quotient by ${\rm PGL}_2[1/2]$. 
\end{theorem}

\bigskip

Let $\Gamma$ be a group or a monoid.  
For $\alpha \in \Gamma$, we define the monoid homomorphism 
$\phi : \Upsilon_1 =\langle \alpha_0 \rangle \to \Gamma$ by 
$\alpha_0 \mapsto \alpha$.  
By restricting representations and derivations of $\Gamma$ to 
those of $\Upsilon_1$ through $\phi$,
we can obtain the following 
commutative diagram:  
\[
\begin{array}{ccc}
  {\rm Rep}_2(\Gamma)_{u, \alpha} & \to & {\rm Ch}_2(\Gamma)_{u, \alpha} \\
  \downarrow & &\downarrow \\
  {\rm Rep}_2(\Upsilon_1)_u & \to &  {\rm Ch}_2(\Upsilon_1)_{u}. 
\end{array}
\]

Under this situation, we have the following lemma. 

\begin{lemma}\label{lemma:fibreproductu}
  The above diagram gives a fibre product. 
In particular, the morphism 
${\rm Rep}_2(\Gamma)_{u, \alpha}  \to  {\rm Ch}_2(\Gamma)_{u, \alpha}$ 
is obtained by base change of the prototype. 
\end{lemma}

{\it Proof.}
We claim that ${\rm Rep}_2(\Gamma)_{u, \alpha} \to 
{\rm Rep}_2(\Upsilon_1)_u\times_{{\rm Ch}_2(\Upsilon_1)_{u}}
{\rm Ch}_2(\Gamma)_{u, \alpha}$ is an isomorphism. 
Let $X$ be a ${\Bbb Z}[1/2]$-scheme. 
Assume that an $X$-valued point 
$\rho \in {\rm Rep}_2(\Gamma)_{u, \alpha}$
is sent to $(\rho', \sigma) \in 
{\rm Rep}_2(\Upsilon_1)_u\times_{{\rm Ch}_2(\Upsilon_1)_{u}}
{\rm Ch}_2(\Gamma)_{u, \alpha}$. 
We can regard the $X$-valued point $\sigma \in {\rm Ch}_2(\Gamma)_{u, \alpha}$ 
as a pair $(r, d)$ such that $d : \Gamma \to \Gamma(X, {\mathcal O}_X)$ is a derivation  
with respect to a character $r : \Gamma \to \Gamma(X, {\mathcal O}_X)$ and $d(\alpha)=1$. 
Since $\eta(\gamma) = d(\gamma) \eta(\alpha)$, 
\setcounter{equation}{5}
\begin{eqnarray}\label{eqn:u-1}
\rho(\gamma) =  r(\gamma)I_2 + d(\gamma)\eta(\rho'(\alpha_{0}))  
\end{eqnarray}
for each $\gamma \in \Gamma$.  
Hence $\rho$ is uniquely determined by $(\rho', \sigma)$.  

For an $X$-valued point $(\rho', \sigma) \in 
{\rm Rep}_2(\Upsilon_1)_u\times_{{\rm Ch}_2(\Upsilon_1)_{u}}
{\rm Ch}_2(\Gamma)_{u, \alpha}$, we define the map 
$\rho : \Gamma \to {\rm M}_2(\Gamma(X, {\mathcal O}_X))$ by (\ref{eqn:u-1}).
From Lemma \ref{lemma-u-ror}, we see that $\rho$ is an $X$-valued point 
of ${\rm Rep}_2(\Gamma)_{u, \alpha}$. 
Then the $X$-valued point $\rho$ 
is sent to $(\rho', \sigma) \in 
{\rm Rep}_2(\Upsilon_1)_u\times_{{\rm Ch}_2(\Upsilon_1)_{u}}
{\rm Ch}_2(\Gamma)_{u, \alpha}$. 
By these discussion, 
the diagram gives a fibre product. 
\qed

\bigskip 

\begin{theorem} 
The morphism $\pi_{\Gamma, u} : {\rm Rep}_2(\Gamma)_{u} \to 
{\rm Ch}_2(\Gamma)_u$ 
is a universal 
geometric quotient by ${\rm PGL}_2[1/2]$ 
for a group or a monoid $\Gamma$. 
\end{theorem}

{\it Proof.} 
For each $\alpha \in \Gamma$, 
$\pi_{\alpha} : {\rm Rep}_2(\Gamma)_{u, \alpha}
\to {\rm Ch}_2(\Gamma)_{u, \alpha}$ is 
a universal geometric quotient 
by ${\rm PGL}_2[1/2]$ because 
$\pi_{\alpha}$ is obtained by base change of 
the prototype. 
Hence this implies the statement.   
\qed 

\bigskip 

\begin{remark}\rm
  The morphism $\pi_{\Gamma, u} : {\rm Rep}_2(\Gamma)_{u} \to 
{\rm Ch}_2(\Gamma)_u$ is smooth and surjective. 
Indeed, the prototype $\pi_{\Upsilon_1, u} : {\rm Rep}_2(\Upsilon_1)_{u} \to {\rm Ch}_2(\Upsilon_1)_u$ 
is smooth and surjective because it is obtained by base change of 
$\pi : {\rm Rep}_2(\Upsilon_{1})_{{\rm rk} 2} \to {\rm Ch}_2(\Upsilon_1)$ and 
$\pi$ is smooth and surjective by Proposition~\ref{lemmassff}.  
\end{remark}

\begin{remark}\label{remark:sectionu}\rm 
For each point $x \in {\rm Ch}_2(\Gamma)_{u}$, there exists 
a local section $s_{x} : V_x \to {\rm Rep}_2(\Gamma)_{u}$ on 
a neighbourhood $V_x$ of $x$ such that 
$\pi_{\Gamma, u} \circ s_{x} = id_{V_x}$. 
Indeed, take $\alpha \in \Gamma$ such that 
$x \in {\rm Ch}_2(\Gamma)_{u, \alpha}$. 
The prototype ${\rm Rep}_2(\Upsilon_1)_{u} \to {\rm Ch}_2(\Upsilon_1)_{u}$ 
has a section $s$ since it is obtained by base change of 
${\rm Rep}_2(\Upsilon_1)_{{\rm rk} 2} \to 
{\rm Ch}_2(\Upsilon_1)$, which has a section 
(it has been defined just before Proposition~\ref{prop:sectionff}). 
By Lemma \ref{lemma:fibreproductu}, 
$\pi_{\alpha} : {\rm Rep}_2(\Gamma)_{u, \alpha} 
\to {\rm Ch}_2(\Gamma)_{u, \alpha}$ has a section $s_{\Gamma, \alpha}$.  
Hence we can take ${\rm Ch}_2(\Gamma)_{u, \alpha}$ 
as a neighbourhood $V_x$ of $x$. 
\end{remark}


\begin{lemma}\label{lemma:localequ} 
Let $\rho_1, \rho_2$ be representations with unipotent mold for a group (or a monoid) 
$\Gamma$ on a scheme $X$ over ${\Bbb Z}[1/2]$. 
Let $f_i : X \to {\rm Rep}_2(\Gamma)_{u}$ be the morphism 
associated to $\rho_i$ for $i = 1, 2$. 
If $\pi_{\Gamma, {u}} \circ f_1 = \pi_{\Gamma, {u}} \circ f_2 
: X \to {\rm Ch}_2(\Gamma)_{u}$, then for each $x \in X$ 
there exists $P_x \in {\rm GL}_2(\Gamma(V_x, {\mathcal O}_X))$ 
on a neighbourhood $V_x$ of $x$ such that 
$P_x^{-1}\rho_1 P_x = \rho_2$ on $V_x$.  
\end{lemma} 

\prf
For $x \in X$, take $\alpha \in \Gamma$ such that 
$(\pi_{\Gamma, u} \circ f_1) (x) = (\pi_{\Gamma, u} \circ f_2) (x) 
\in {\rm Ch}_2(\Gamma)_{u, \alpha}$.  
We may assume that $f_i : X \to {\rm Rep}_2(\Gamma)_{{u}, \alpha}$ for 
$i =1, 2$ from the beginning. By Remark \ref{remark:sectionu}, 
$\pi_{\alpha} : {\rm Rep}_2(\Gamma)_{u, \alpha} 
\to {\rm Ch}_2(\Gamma)_{u, \alpha}$ has a section $s_{\Gamma, \alpha}$.  
Let $\rho_3$ be the representations with unipotent mold on $X$ associated to 
$s_{\Gamma, \alpha} \circ \pi_{\Gamma, u} \circ f_1 = s_{\Gamma, \alpha} \circ 
\pi_{\Gamma, u} \circ f_2$. 
Then $\rho_i(\gamma) = r(\gamma) I_2 + d(\gamma)\eta(\rho_i(\alpha))$ for each 
$\gamma \in \Gamma$ and $i = 1, 2, 3$, where $d$ is the derivation with respect to 
the character $r$ associated to $\pi_{\Gamma, u} \circ f_1 = 
\pi_{\Gamma, u} \circ f_2$ such that 
$d(\alpha)=1$. 
Note that 
$\rho_3(\alpha) = 
\left(
\begin{array}{cc}
0 & -D \\
1 & T \\
\end{array}
\right)$  and that 
$T = {\rm tr}(\rho_1(\alpha)) = {\rm tr}(\rho_2(\alpha)) = {\rm tr}(\rho_3(\alpha))$ 
and $D = T^2/4$.  
There exist $Q_1, Q_2 \in {\rm GL}_2(\Gamma(V_x, {\mathcal O}_X))$ on a
neighbourhood $V_x$ of $x$ such that 
$Q_1^{-1}\rho_1(\alpha)Q_1 = \rho_3(\alpha)$ and 
$Q_2^{-1}\rho_2(\alpha)Q_2 = \rho_3(\alpha)$ by Lemma \ref{lemma-norm-nonsca}. 
Since 
$Q_1^{-1}\rho_1(\gamma)Q_1 = \rho_3(\gamma)$ and 
$Q_2^{-1}\rho_2(\gamma)Q_2 = \rho_3(\gamma)$ for each $\gamma \in \Gamma$, 
$(Q_1 Q_2^{-1})^{-1} \rho_1 (Q_1 Q_2^{-1}) = \rho_2$ on $V_x$. 
This completes the proof. 
\qed 

\bigskip 

Let us define 
${\mathcal E}q \:\mathcal{U}_2(\Gamma)$ as the  
sheafification of the following contravariant functor with respect to Zariski topology: 
\[
\begin{array}{ccl}
 ({\bf Sch}/{\Bbb Z}[1/2])^{op} & \to & ({\bf Sets}) \\
 X & \mapsto & \{ \rho \mid \mbox{ rep. 
with unipotent mold for $\Gamma$ on } X \}/\sim.
\end{array}
\]

By a {\it generalized representation with unipotent mold} for $\Gamma$ on 
a scheme $X$, we understand pairs $\{ (U_i, \rho_i) \}_{i \in I}$ 
of an open set $U_i$ and a representation with unipotent mold 
$\rho_i : \Gamma \to {\rm M}_2(\Gamma(U_i, {\mathcal O}_{X}))$ 
satisfying the following two conditions:
\begin{enumerate}
\item $\cup_{i \in I} U_ i = X$, 
\item for each $x \in U_i \cap U_j$, there exists $P_x \in 
{\rm GL}_2(\Gamma(V_x, {\mathcal O}_X))$ on 
a neighbourhood $V_x \subseteq U_i \cap U_j$ of $x$ such that 
$P_x^{-1} \rho_i P_x = \rho_j$ on $V_x$. 
\end{enumerate}

Generalized representations with unipotent mold 
$\{ (U_i, \rho_i) \}_{i \in I}$ and $\{ (V_j, \sigma_j) \}_{j \in J}$ 
are called {\it equivalent} if $\{ (U_i, \rho_i) \}_{i \in I} \cup \{ (V_j, \sigma_j) \}_{j \in J}$ 
is a generalized representation with unipotent mold again.  
We easily see that 
${\mathcal E}q \:\mathcal{U}_2(\Gamma)(X)$ is the 
set of equivalence classes of 
generalized representations with unipotent mold 
for $\Gamma$ on a scheme $X$. 

\begin{theorem}\label{th:moduliu}  
The scheme ${\rm Ch}_2(\Gamma)_{u}$ is a fine moduli scheme  
associated to the functor ${\mathcal E}q \:\mathcal{U}_2(\Gamma)$ 
for a group or a monoid $\Gamma$:  
\[
\begin{array}{ccccl}
{\mathcal E}q \:\mathcal{U}_2(\Gamma) & : & ({\bf Sch}/{\Bbb Z}[1/2])^{op} & \to & ({\bf Sets}) \\
 & & X & \mapsto & 
\left\{ 
\begin{array}{r}  
\mbox{ gen. rep. 
with unipotent  } \\ 
\mbox{ mold for $\Gamma$ on } X  
\end{array} 
\right\} 
\Big/\sim.
\end{array}
\] 
In other words, ${\rm Ch}_2(\Gamma)_{u}$ represents 
the functor ${\mathcal E}q \:\mathcal{U}_2(\Gamma)$. 
The moduli ${\rm Ch}_2(\Gamma)_{u}$ is separated over ${\Bbb Z}[1/2]$; 
if $\Gamma$ is a finitely generated group or monoid, then 
${\rm Ch}_2(\Gamma)_{u}$ is of finite type over ${\Bbb Z}[1/2]$. 
\end{theorem}

\prf 
Since $\pi_{\Gamma, u} : {\rm Rep}_2(\Gamma)_{u} 
\to {\rm Ch}_2(\Gamma)_{u}$ is a 
${\rm PGL}_2[1/2]$-equivariant morphism, we can define a canonical morphism 
${\mathcal E}q \: \mathcal{U}_2(\Gamma) \to 
h_{{\rm Ch}_2(\Gamma)_{u}}   := {\rm Hom}(-, {\rm Ch}_2(\Gamma)_{u})$.  
Let us define a morphism $h_{{\rm Ch}_2(\Gamma)_{u}} 
\to {\mathcal E}q \:\mathcal{U}_2(\Gamma)$.  
Let $g \in h_{{\rm Ch}_2(\Gamma)_{u}}(X)$ with a 
${\Bbb Z}[ 1/2 ]$-scheme $X$. 
For each $x \in X$, take $\alpha_x \in \Gamma$ such that 
$g(x) \in {\rm Ch}_2(\Gamma)_{u, \alpha_x}$. 
By using the section $s_{\Gamma, \alpha_x} : {\rm Ch}_2(\Gamma)_{u, \alpha_x} 
\to {\rm Rep}_2(\Gamma)_{u, \alpha_x}$ in Remark \ref{remark:sectionu}, 
we can define a representation with unipotent mold $\rho_x$ on a neighbourhood 
$U_x$ of $x$. By Lemma \ref{lemma:localequ}, 
we see that $\{ (U_x, \rho_x) \}_{x \in X} \in 
{\mathcal E}q \:\mathcal{U}_2(\Gamma)(X)$ and that the morphism 
$h_{{\rm Ch}_2(\Gamma)_{u}} 
\to {\mathcal E}q \:\mathcal{U}_2(\Gamma)$ is well-defined.  
It is easy to see that ${\rm Ch}_2(\Gamma)_{u}$ represents 
the functor ${\mathcal E}q \:\mathcal{U}_2(\Gamma)$. 

Since ${\rm Ch}_2(\Gamma)_{u}$ is 
defined as ${\rm Proj} S(\Omega_{\Gamma})$, 
it is separated over ${\Bbb Z}[1/2]$. 
If $\Gamma$ is finitely generated, then 
${\rm Ch}_2(\Gamma)_{u}$ is of finite type over ${\Bbb Z}[1/2]$
by Remark \ref{remark:fgu}. 
\qed 

\begin{remark}\label{remark:algu1}\rm
Let $A$ be an associative algebra over a commutative ring $R$ over ${\Bbb Z}[1/2]$. 
For a $2$-dimensional representation $\rho$ of $A$ on an $R$-scheme $X$, 
$\rho$ is called a {\it representation with unipotent mold} if 
the subalgebra $\rho(A)$ of ${\rm M}_2({\mathcal O}_X)$ generates  
a unipotent mold on $X$. 
In a similar way as group or monoid cases, 
we can define generalized representations with unipotent mold for $A$ 
on an $R$-scheme $X$. The contravariant functor 
${\mathcal E}q \;\mathcal{U}_2(A)$ from the category of $R$-schemes to the 
category of sets is defined as 
\[
\begin{array}{ccccl}
{\mathcal E}q \;\mathcal{U}_2(A) & : & ({\bf Sch}/R)^{op} & \to & ({\bf Sets}) \\
 & & X & \mapsto & 
\left\{ 
\begin{array}{r}  
\mbox{ gen. rep. with unipotent 
 } \\ 
\mbox{  mold for $A$ on } X  
\end{array} 
\right\} 
\Big/\sim.
\end{array}
\] 
Then we can construct the fine moduli ${\rm Ch}_2(A)_{u}$ associated to 
${\mathcal E}q \;\mathcal{U}_2(A)$ in the same way as Theorem \ref{th:moduliu} 
(for details, see Remark~\ref{remark:constructu}). 
The moduli ${\rm Ch}_2(A)_{u}$ is separated over $R$. 
If $A$ is a finitely generated algebra over $R$, then 
${\rm Ch}_2(A)_{u}$ is of finite type over $R$.  
\end{remark}

\begin{remark}\label{remark:constructu}\rm 
For an associative algebra $A$ over a commutative ring $R$ over 
${\Bbb Z}[1/2]$, we can construct ${\rm Ch}_2(A)_{u}$ in the following way. 
We define the contravariant functor 
${\rm Rep}_1(A)$ from the category of $R$-schemes 
to the category of sets by $X \mapsto \{ 
\varphi : A \to \Gamma(X, {\mathcal O}_X) \mid R\mbox{-algebra hom. } 
\}$.  The functor ${\rm Rep}_1(A)$ is representable by an affine scheme, and 
let us denote its coordinate ring by $A_1(A)$.  
Let $d : A \to A_1(A)$ be the universal $R$-algebra homomorphism. 
For an $A_1(A)$-module $M$, 
put 
\[ {\rm Der}(A, M) := \left\{ 
\delta : A \to M  \; 
\begin{array}{|c} 
 \delta \mbox{ : R-linear and for } a, b \in A,  \\
 \delta(ab) = d(a)\delta(b)+\delta(a)d(b)  
\end{array} 
\right\}.
\]   
The functor ${\rm Der}(A, -) : (A_1(A)\mbox{-}{\bf Mod}) \to (A_1(A)\mbox{-}{\bf Mod})$ 
defined by $M \mapsto {\rm Der}(A, M)$ is representable by 
some $A_1(A)$-module $\Omega_{A/R}$.  
Let ${\rm Rep}_2(A)$ be the representation variety of degree $2$ for $A$ over 
$R$, that is, the affine scheme representing the contravariant 
functor from the category of $R$-schemes to the category of sets 
which is 
defined by $X \mapsto \{ 2\mbox{-dim. rep. of $A$ on $X$} \}$. 
Let ${\rm Rep}_2(A)_{u}$ be the subscheme of 
${\rm Rep}_2(A)$ consisting of representations 
with unipotent mold, and   
$\sigma : A \to 
{\rm M}_2(\Gamma({\rm Rep}_2(A)_{u}, {\mathcal O}_{{\rm Rep}_2(A)_{u}}))$ 
the universal representation with unipotent mold.  
Then ${\rm tr}(\sigma(\cdot))/2 : A \to 
\Gamma({\rm Rep}_2(A)_{u}, {\mathcal O}_{{\rm Rep}_2(A)_{u}})$ 
is an $R$-algebra homomorphism, and it defines 
a morphism ${\rm Rep}_2(A)_{u} \to {\rm Rep}_1(A)$.  In a similar way as 
group or monoid cases, 
we can define a ${\rm Rep}_1(A)$-morphism $\pi : {\rm Rep}_2(A)_{u} 
\to {\rm Ch}_2(A)_{u} := {\rm Proj} S(\Omega_{A/R})$, where 
$S(\Omega_{A/R})$ is the symmetric algebra of 
$\Omega_{A/R}$ over $A_1(A)$.  
We can verify that $\pi$ is a universal geometric quotient by 
${\rm PGL}_2\otimes_{\Bbb Z} R$ and that ${\rm Ch}_2(A)_{u}$ represents 
${\mathcal E}q \;\mathcal{U}_2(A)$. 
\end{remark}

\begin{remark}\rm 
We have introduced the notion of 
generalized representations with unipotent mold 
for describing the moduli functors ${\mathcal E}q\mathcal{U}_2(\Gamma)$  
and ${\mathcal E}q\mathcal{U}_2(A)$.  However, 
the moduli functors can also be described as 
${\mathcal E}q\mathcal{U}'_2(\Gamma)$  
and ${\mathcal E}q\mathcal{U}'_2(A)$ by using 
the notion of representations generating sheaves of algebras 
which define unipotent molds.   
More precisely, see \S 8. 

\end{remark}

\section{Unipotent mold over ${\mathbb F}_2$}

In this section, all schemes are over ${\rm Spec} \: {\Bbb F}_2$ 
and all commutative rings are over ${\Bbb F}_2$. 
Recall that a rank $2$ mold ${\mathcal A} \subseteq {\rm M}_2({\mathcal O}_X)$ 
over an ${\Bbb F}_2$-scheme $X$ 
is called {\it unipotent over ${\Bbb F}_2$} if 
${\rm tr}(s) = 0$ for each open subset $U \subseteq X$ 
and for each $s \in {\mathcal A}(U)$. 
We construct the moduli of representations with unipotent mold
over ${\Bbb F}_2$. 

\begin{definition}\label{def:abcoeff}\rm 
Let $X$ be an ${\Bbb F}_2$-scheme. Let ${\mathcal A} 
\subseteq {\rm M}_2({\mathcal O}_X)$ be a unipotent 
mold over ${\Bbb F}_2$ on $X$. 
Let $U \subseteq X$ be an open subset of $X$. 
Suppose that $Z \in {\mathcal A}(U)$ 
satisfies ${\mathcal A}|_{U} = {\mathcal O}_{U} \cdot I_2 \oplus 
{\mathcal O}_{U} \cdot Z$.  
For each $Y \in {\mathcal A}(U)$, 
we set $Y = a_{Z}(Y) I_2 + b_{Z}(Y) Z$. 
We call $a_Z(Y), b_Z(Y) \in {\mathcal O}_X(U)$ 
the {\it $(a, b)$-coefficients} of $Y$ with respect to $Z$.  
\end{definition} 

\begin{lemma}\label{lemma:abproperty} 
Let $U$ and $Z$ be as in Definition~\ref{def:abcoeff}.  
Assume that $\rho : \Gamma \to {\rm M}_2(\Gamma(X, {\mathcal O}_X))$ 
is a representation with unipotent mold ${\mathcal A}$ on $X$ for a
group or a monoid $\Gamma$. 
For each $\gamma \in \Gamma$, let us 
denote $a_Z(\rho(\gamma))$ and $b_Z(\rho(\gamma))$ by $a(\gamma)$ and $b(\gamma)$, 
respectively.  Then for $\gamma, \delta \in \Gamma$, we have 
\begin{eqnarray*}
a(e) & = & 1, \\
b(e) & = & 0, \\
a(\gamma\delta) & = & a(\gamma)a(\delta) + b(\gamma)b(\delta) \det Z,  \\
b(\gamma\delta) & = & a(\gamma)b(\delta) + b(\gamma)a(\delta).  
\end{eqnarray*} 
\end{lemma} 

\prf
Since $\rho(e) = I_2 = 1 \cdot I_2+ 0 \cdot Z$, 
$a(e) =  1$ and  
$b(e)  =  0$. 
By the Cayley-Hamilton theorem, 
$Z^2 - {\rm tr}(Z) Z + \det(Z)I_2=0$. Hence we have 
$Z^2 = -\det(Z)I_2 = \det(Z)I_2$ by 
${\rm tr}(Z)=0$.  We see that 
\begin{eqnarray*}
\rho(\gamma\delta) & = & \rho(\gamma)\rho(\delta) \\
 & = & (a(\gamma)I_2+b(\gamma)Z)(a(\delta)I_2+b(\delta)Z) \\
 & = & \{ a(\gamma)a(\delta) + b(\gamma)b(\delta) \det(Z) \} I_2 + 
\{ a(\gamma)b(\delta) + b(\gamma)a(\delta) \} Z.  
\end{eqnarray*} 
Comparing the coefficients of $\rho(\gamma\delta) = a(\gamma\delta) I_2 + 
b(\gamma\delta)Z$, we obtain $a(\gamma\delta) = a(\gamma)a(\delta) + b(\gamma)b(\delta) \det Z$ and 
$b(\gamma\delta) = a(\gamma)b(\delta) + b(\gamma)a(\delta)$.  
\qed 

\bigskip 

Let $R$ be an algebra over ${\Bbb F}_2$. Let $\rho : \Gamma 
\to {\rm M}_2(R)$ be a representation with unipotent mold ${\mathcal A}$ over ${\Bbb F}_2$ 
such that ${\mathcal A} = R\cdot I_2 \oplus R\cdot\rho(\alpha)$. For each $\gamma \in \Gamma$, 
we denote $a_{\rho(\alpha)}(\gamma), b_{\rho(\alpha)}(\gamma)$ 
by $a(\gamma), b(\gamma)$, respectively.  
Then $a(\cdot)$ and $b(\cdot)$ satisfy the formulas in Lemma~\ref{lemma:abproperty}, where 
$Z = \rho(\alpha)$. 
Furthermore, $a(\alpha)=0$ and $b(\alpha)=1$. 
Conversely, a character and $(a, b)$-coefficients give a representation with 
unipotent mold over ${\Bbb F}_2$:

\begin{lemma}\label{lemma:uf2rep} 
Let $d : \Gamma \to R$ be a character. 
Let ${\mathcal A} = R\cdot I_2 \oplus R\cdot Z \subseteq {\rm M}_2(R)$ be a 
unipotent mold over ${\Bbb F}_2$ such that 
${\rm tr}(Z)=0$ and $\det (Z) = d(\alpha)$ for some $\alpha \in \Gamma$. 
Assume that $a : \Gamma \to R$, $b : \Gamma \to R$, and $\alpha \in \Gamma$ satisfy the equalities 
\begin{eqnarray*}
a(\gamma\delta) & = & a(\gamma)a(\delta) + b(\gamma)b(\delta)d(\alpha), \\
b(\gamma\delta) & = & a(\gamma)b(\delta) + b(\gamma)a(\delta), \\ 
d(\gamma) &=& a(\gamma)^2 + b(\gamma)^2 d(\alpha)  
\end{eqnarray*} 
for each $\gamma, \delta \in \Gamma$. Furthermore, 
assume that $a(\alpha)=0$ and that $b(\alpha)=1$. 
Then the map 
\[
\begin{array}{ccccl}
\rho & : & \Gamma & \to & {\rm M}_2(R) \\
 & & \gamma & \mapsto & a(\gamma) I_2 + b(\gamma) Z
\end{array}
\]
is a representation with unipotent mold ${\mathcal A}$ over ${\Bbb F}_2$ such that 
$\det (\rho(\gamma)) = d(\gamma)$ for 
each $\gamma \in \Gamma$.  
\end{lemma}

\prf
First we show that $a(e)=1$ and $b(e)=0$. By the assumption, 
\begin{eqnarray*}
b(e) = b(e\cdot e) & = & a(e)b(e) + b(e)a(e) \\
  & = & 2a(e)b(e) = 0 
\end{eqnarray*}
and 
\begin{eqnarray*}
a(e) = a(e\cdot e) & = & a(e)a(e) + b(e)b(e)d(\alpha) \\
  & = & d(e) =1.  
\end{eqnarray*} 
Hence $\rho(e) = a(e)I_2 + b(e)Z = I_2$. 

Next we show that $\rho(\gamma\delta) = \rho(\gamma)\rho(\delta)$ for each $\gamma, \delta \in \Gamma$. 
Since $Z^2 = \det(Z) I_2 = d(\alpha) I_2$, 
\begin{eqnarray*}
\rho(\gamma)\rho(\delta) & = & (a(\gamma) I_2 + b(\gamma) Z)(a(\delta) I_2 + b(\delta) Z) \\
 & = & \{ a(\gamma)a(\delta) + b(\gamma)b(\delta) d(\alpha) \} I_2 +  
 \{ a(\gamma)b(\delta) + b(\gamma)a(\delta) \} Z \\
  & = & a(\gamma\delta) I_2 + b(\gamma\delta) Z \\
  & = & \rho(\gamma\delta). 
\end{eqnarray*}
By Lemma~\ref{lemma:detformula}, $\det \rho(\gamma) = \det (a(\gamma) I_2 + b(\gamma)Z) 
= a(\gamma)^2\det I_2 + b(\gamma)^2 \det Z = a(\gamma)^2 + b(\gamma)^2 d(\alpha) = d(\gamma)$. 

Finally, $\rho(\alpha) = 0\cdot I_2 + 1\cdot Z = Z$ implies that $\rho(\Gamma)$ generates 
${\mathcal A}$. Hence $\rho$ is a representation with unipotent mold ${\mathcal A}$ over ${\Bbb F}_2$. 
\qed   


\begin{definition}\rm 
Let $d : \Gamma \to R$ be a character in an ${\Bbb F}_2$-algebra $R$. 
For $a : \Gamma \to R$, $b : \Gamma \to R$ and $\alpha \in \Gamma$, 
we say that $a$ and $b$ are {\it $(a, b)$-coefficients} with respect to $(d, \alpha)$ 
if  
$$a(e)=1, b(e)=0, a(\alpha)=0, b(\alpha)=1,$$ \vspace*{-6ex} 
\begin{eqnarray*}
a(\gamma\delta) & = & a(\gamma)a(\delta) + b(\gamma)b(\delta)d(\alpha), \\
b(\gamma\delta) & = & a(\gamma)b(\delta) + b(\gamma)a(\delta), \mbox{ and}  \\ 
d(\gamma) & = & a(\gamma)^2 + b(\gamma)^2 d(\alpha)     
\end{eqnarray*}                                                            
hold for each $\gamma, \delta \in \Gamma$. 
\end{definition} 

\begin{definition}\rm
Set ${\rm Rep}_2(\Gamma)_{{\Bbb F}_2} := {\rm Rep}_2(\Gamma)\otimes_{{\Bbb Z}} {{\Bbb F}_2}$ and 
${\rm Rep}_2(\Gamma)_{{\rm rk} 2/{\Bbb F}_2} := {\rm Rep}_2(\Gamma)_{{\rm rk} 2}\otimes_{{\Bbb Z}} {{\Bbb F}_2}$. 
Let us define ${\rm Rep}_2(\Gamma)_{u/{\Bbb F}_2}$ as a closed subscheme 
of ${\rm Rep}_2(\Gamma)_{{\rm rk} 2/{\Bbb F}_2}$ by 
\[
{\rm Rep}_2(\Gamma)_{u/{\Bbb F}_2} := 
\{ \rho \in {\rm Rep}_2(\Gamma)_{{\rm rk} 2/{\Bbb F}_2} 
\mid {\rm tr}(\rho(\gamma)) = 0 
\mbox{ for each $\gamma \in \Gamma$ } \}. 
\]
For $\gamma \in \Gamma$, we define 
${\rm Rep}_2(\Gamma)_{u/{\Bbb F}_2, \gamma}$ 
as an open subscheme of ${\rm Rep}_2(\Gamma)_{u/{\Bbb F}_2}$ 
by 
\[
{\rm Rep}_2(\Gamma)_{u/{\Bbb F}_2, \gamma} 
:= \left\{ 
\rho \in {\rm Rep}_2(\Gamma)_{u/{\Bbb F}_2}  \; 
\begin{array}{|l} 
 I_2 \mbox{ and } \rho(\gamma) \mbox{ generate } \\ 
 \mbox{ a unipotent mold over } {\Bbb F}_2 
\end{array} 
\right\}.  
\]
Note that 
\[
\displaystyle 
{\rm Rep}_2(\Gamma)_{u/{\Bbb F}_2}
= \bigcup_{\gamma \in \Gamma} 
{\rm Rep}_2(\Gamma)_{u/{\Bbb F}_2, \gamma}.  
\]
\end{definition} 

\begin{definition}\label{def:aringchuf2}\rm 
Set ${\rm Rep}_1(\Gamma)_{{\Bbb F}_2} := {\rm Rep}_1(\Gamma)\otimes_{{\Bbb Z}} {\Bbb F}_2$.  
Let $A_1(\Gamma)_{{\Bbb F}_2}$ be the coordinate ring of 
${\rm Rep}_1(\Gamma)_{{\Bbb F}_2}$, and let $d : \Gamma \to A_1(\Gamma)_{{\Bbb F}_2}$ be 
the universal character of $\Gamma$.  
For $\alpha \in \Gamma$, we define the $A_1(\Gamma)_{{\Bbb F}_2}$-algebra 
$A_2(\Gamma)_{u/{\Bbb F}_2, \alpha}^{\rm Ch}$ by $A_1(\Gamma)_{{\Bbb F}_2}[a(\gamma), b(\gamma) \mid \gamma \in \Gamma]/I$, 
where $a(\gamma), b(\gamma)$ are indeterminates for each $\gamma \in \Gamma$ and 
$I$ is generated by 
\begin{eqnarray*}
a(e)-1, b(e), a(\alpha), b(\alpha)-1 \\
a(\gamma\delta) - a(\gamma)a(\delta) - b(\gamma)b(\delta)d(\alpha) \\
b(\gamma\delta) - a(\gamma)b(\delta) - b(\gamma)a(\delta) \\
a(\gamma)^2 - b(\gamma)^2 d(\alpha) -d(\gamma) 
\end{eqnarray*} 
for all $\gamma, \delta \in \Gamma$. 
We set ${\rm Ch}_2(\Gamma)_{u/{\Bbb F}_2, \alpha} := 
{\rm Spec} A_2(\Gamma)_{u/{\Bbb F}_2, \alpha}^{\rm Ch}$.  
\end{definition} 

\bigskip 

For $\alpha \in \Gamma$, ${\rm Ch}_2(\Gamma)_{u/{\Bbb F}_2, \alpha}$ 
is a ${\rm Rep}_1(\Gamma)_{{\Bbb F}_2}$-scheme. 
For a ${\rm Rep}_1(\Gamma)_{{\Bbb F}_2}$-scheme $X$, denote by  
$\chi : \Gamma \to \Gamma(X, {\mathcal O}_X)$ the character of $\Gamma$
associated to $X \to {\rm Rep}_1(\Gamma)_{{\Bbb F}_2}$. There exists a $1$-$1$ 
correspondence 
\[
{\rm Hom_{{\rm Rep}_1(\Gamma)_{{\Bbb F}_2}}}( X, {\rm Ch}_2(\Gamma)_{u/{\Bbb F}_2, \alpha}) 
\cong \left\{  (a, b) \; 
\begin{array}{|l} 
a, b : \Gamma \to \Gamma(X, {\mathcal O}_X) \\  
\mbox{ are } (a, b)\mbox{-coefficients } \\
\mbox{ with respect to }  (\chi, \alpha)  
\end{array} 
\right\}. 
\] 

Let $\sigma_{\Gamma, u/{\Bbb F}_2}$  and $\sigma_{\Gamma, u/{\Bbb F}_2, \alpha}$
be the universal representation with unipotent mold over ${\Bbb F}_2$ on 
${\rm Rep}_2(\Gamma)_{u/{\Bbb F}_2}$ and ${\rm Rep}_2(\Gamma)_{u/{\Bbb F}_2, \alpha}$, 
respectively. 
Put $d(\gamma) :=\det( \sigma_{\Gamma, u/{\Bbb F}_2, \alpha} (\gamma) )$ for $\gamma 
\in \Gamma$.  
By Lemma~\ref{lemma:abproperty}, 
$a_{\sigma_{\Gamma, u/{\Bbb F}_2, \alpha}(\alpha)}(\cdot)$ and $b_{\sigma_{\Gamma, u/{\Bbb F}_2, \alpha}(\alpha)}(\cdot)$ 
are $(a, b)$-coefficients with respect to $(d, \alpha)$ on 
${\rm Rep}_2(\Gamma)_{u/{\Bbb F}_2, \alpha}$. 
Hence we have the morphism 
$\pi_{\Gamma, u/{\Bbb F}_2, \alpha} : {\rm Rep}_2(\Gamma)_{u/{\Bbb F}_2, \alpha} 
\to {\rm Ch}_2(\Gamma)_{u/{\Bbb F}_2, \alpha}$ 
associated to 
$a_{\sigma_{\Gamma, u/{\Bbb F}_2, \alpha}(\alpha)}(\cdot)$ and $b_{\sigma_{\Gamma, u/{\Bbb F}_2, \alpha}(\alpha)}(\cdot)$. 
Note that ${\rm Rep}_2(\Gamma)_{u/{\Bbb F}_2, \alpha} 
\to {\rm Rep}_1(\Gamma)_{{\Bbb F}_2}$ is given by 
$\rho \mapsto \det(\rho)$. 
The group scheme ${\rm PGL}_2\otimes_{\Bbb Z} {\Bbb F}_2$ over 
${\Bbb F}_2$ acts on ${\rm Rep}_2(\Gamma)_{u/{\Bbb F}_2}$ by 
$\rho \mapsto P^{-1}\rho P$. 
The open subscheme ${\rm Rep}_2(\Gamma)_{u/{\Bbb F}_2, \alpha}$ of 
${\rm Rep}_2(\Gamma)_{u/{\Bbb F}_2}$ is ${\rm PGL}_2\otimes_{\Bbb Z} {\Bbb F}_2$-invariant for 
each $\alpha \in \Gamma$.  
Let ${\rm PGL}_2\otimes_{\Bbb Z} {\Bbb F}_2$ act on 
${\rm Ch}_2(\Gamma)_{u/{\Bbb F}_2, \alpha}$ trivially. Then we have:  

\begin{proposition}\label{prop:abinv} 
$\pi_{\Gamma, u/{\Bbb F}_2, \alpha} : {\rm Rep}_2(\Gamma)_{u/{\Bbb F}_2, \alpha} 
\to {\rm Ch}_2(\Gamma)_{u/{\Bbb F}_2, \alpha}$ 
is ${\rm PGL}_2\otimes_{\Bbb Z} {\Bbb F}_2$-equivariant. 
\end{proposition} 

{\it Proof.} 
Let $(\rho, P)$ be an $X$-valued point of 
${\rm Rep}_2(\Gamma)_{u/{\Bbb F}_2, \alpha} \times ({\rm PGL}_2
\otimes_{\Bbb Z} {\Bbb F}_2)$ for an ${\Bbb F}_2$-scheme $X$.   
The representations $\rho$ and $P^{-1} \rho P$ on $X$ have 
the same determinants. For proving that 
$\rho$ and $P^{-1} \rho P$ induce the same 
$X$-valued point of ${\rm Ch}_2(\Gamma)_{u/{\Bbb F}_2, \alpha}$, 
it only suffices to 
show that $\rho$ and $P^{-1} \rho P$ have the same 
$(a, b)$-coefficients. We may assume that $X$ is an affine scheme 
${\rm Spec} R$ and that $P \in {\rm GL}_2(R)$. 
Note that $R[\rho(\Gamma)] = 
R\cdot I_2 + R\cdot \rho(\alpha) \subset {\rm M}_2(R)$ 
is a unipotent mold over ${\Bbb F}_2$. 
For each $\gamma \in \Gamma$, $\rho(\gamma) = a(\gamma)I_2 + 
b(\gamma)\rho(\alpha)$, where 
$a(\cdot)$ and $b(\cdot)$ are 
the $(a, b)$-coefficients of $\rho$. 
Multiplying the both sides by $P^{-1}$ from the left and by 
$P$ from the right, we have  $P^{-1} \rho(\gamma) P = a(\gamma) I_2 
+ b(\gamma) P^{-1} \rho(\alpha)P$.   
Hence the $(a, b)$-coefficients of $P\rho P^{-1}$ with 
respect to $P^{-1}\rho(\alpha)P$ coincide 
with the $(a, b)$-coefficients of $\rho$ 
with respect to $\rho(\alpha)$.
This implies that 
$(a, b)$-coefficients are ${\rm PGL}_2\otimes_{\Bbb Z} {\Bbb F}_2$-invariant, 
which completes the proof. 
\qed

\begin{example}\label{ex:upsilon2}\rm
Let $\Upsilon_2 = \langle \alpha_1, \alpha_2 \rangle$ be the 
free monoid of rank $2$. 
Let $\sigma := \sigma_{\Upsilon_2, u/{\Bbb F}_2, \alpha_1}$ denote 
the universal representation of $\Upsilon_2$ on 
${\rm Rep}_2(\Upsilon_2)_{u/{\Bbb F}_2, \alpha_1}$. 
Then we can write 
$\sigma(\alpha_1) = \left( 
\begin{array}{cc} 
a & b \\
c & a \\
\end{array} \right)$,  
$\sigma(\alpha_2) = \left( 
\begin{array}{cc} 
d & e \\
f & d \\
\end{array} \right)$,  and 
\[
{\rm Rep}_2(\Upsilon_2)_{u/{\Bbb F}_2, \alpha_1} 
= D(b)\cup D(c) \subset {\rm Spec} \: {\Bbb F}_2[a, b, c, d, e, f]/(bf-ce).
\] 
The $(a, b)$-coefficients of $\sigma(\alpha_2)$ with respect to 
$\sigma(\alpha_1)$ are given by 
$\displaystyle a(\alpha_2) = \frac{bd-ae}{b}$, 
$\displaystyle b(\alpha_2) = \frac{e}{b}$ on $D(b)$, 
and $\displaystyle a(\alpha_2) = \frac{cd-af}{c}$, 
$\displaystyle b(\alpha_2) = \frac{f}{c}$ on $D(c)$. 
For the universal matrix 
$P = \left( 
\begin{array}{cc} 
p & q \\
r & s \\
\end{array} \right) \in {\rm PGL}_2\otimes_{{\Bbb Z}} {\Bbb F}_2$,  
\[
\displaystyle 
P^{-1} \sigma(\alpha_1) P = \frac{1}{\Delta} \left( 
\begin{array}{cc} 
a\Delta+brs-cpq  & bs^2-cq^2 \\
-br^2+cp^2 & a\Delta-brs+cpq \\
\end{array} \right)  
\] 
and 
\[
\displaystyle 
P^{-1} \sigma(\alpha_2) P = \frac{1}{\Delta} \left( 
\begin{array}{cc} 
d\Delta+ers-fpq  & es^2-fq^2 \\
-er^2+fp^2 & d\Delta-ers+fpq \\
\end{array} \right), 
\]
where $\Delta = ps-qr$. 
By direct calculation, we can verify that $a(\alpha_2)$ and $b(\alpha_2)$ are 
${\rm PGL}_2\otimes_{{\Bbb Z}} {\Bbb F}_2$-invariant functions on 
${\rm Rep}_2(\Upsilon_2)_{u/{\Bbb F}_2, \alpha_1}$. 
Hence $\pi_{\Upsilon_2, u/{\Bbb F}_2, \alpha_1} : 
{\rm Rep}_2(\Upsilon_2)_{u/{\Bbb F}_2, \alpha_1} 
\to {\rm Ch}_2(\Upsilon_2)_{u/{\Bbb F}_2, \alpha_1}$ 
is ${\rm PGL}_2\otimes_{\Bbb Z} {\Bbb F}_2$-equivariant. 
\end{example} 

\begin{remark}\rm 
For any group $\Gamma$ and for any $\alpha, \gamma \in \Gamma$, 
let us define the group homomorphism $\phi : \Upsilon_2 = \langle 
\alpha_1, \alpha_2 \rangle \to \Gamma$ by 
$\alpha_1 \mapsto \alpha$ and $\alpha_2 \mapsto \gamma$. 
Then $\phi$ induces a  
morphism $\phi^{\ast} :  {\rm Rep}_2(\Gamma)_{u/{\Bbb F}_2, \alpha} \to  
{\rm Rep}_2(\Upsilon_2)_{u/{\Bbb F}_2, \alpha_1}$ 
by $\rho \mapsto \rho\circ \phi$. 
By Example \ref{ex:upsilon2}, the $(a, b)$-coefficients 
$a(\alpha_2)$ and $b(\alpha_2)$ are 
${\rm PGL}_2\otimes_{{\Bbb Z}} {\Bbb F}_2$-invariant functions on 
${\rm Rep}_2(\Upsilon_2)_{u/{\Bbb F}_2, \alpha_1}$.  
Let $\sigma_{\Gamma}$ be the universal representation of 
$\Gamma$ on ${\rm Rep}_2(\Gamma)_{u/{\Bbb F}_2, \alpha}$. Let 
$a_{\alpha}(\gamma)$ and $b_{\alpha}(\gamma)$ denote the $(a, b)$-coefficients of $\sigma_{\Gamma}(\gamma)$ 
with respect to $\sigma_{\Gamma}(\alpha)$. 
Note that $\phi_{\ast}(a(\alpha_2))=a_{\alpha}(\gamma)$ and 
$\phi_{\ast}(b(\alpha_2))=b_{\alpha}(\gamma)$, where 
$\phi_{\ast}$ denotes the ring homomorphism 
$\Gamma({\rm Rep}_2(\Upsilon_2)_{u/{\Bbb F}_2, \alpha_1}, 
{\mathcal O}_{{\rm Rep}_2(\Upsilon_2)_{u/{\Bbb F}_2, \alpha_1}}) 
\to \Gamma({\rm Rep}_2(\Gamma)_{u/{\Bbb F}_2, \alpha}, 
{\mathcal O}_{{\rm Rep}_2(\Gamma)_{u/{\Bbb F}_2, \alpha}})$ 
associated with $\phi^{\ast} :  {\rm Rep}_2(\Gamma)_{u/{\Bbb F}_2, \alpha} \to  
{\rm Rep}_2(\Upsilon_2)_{u/{\Bbb F}_2, \alpha_1}$. 
Since $\phi^{\ast}$ is ${\rm PGL}_2\otimes_{{\Bbb Z}} {\Bbb F}_2$-equivariant, 
$a_{\alpha}(\gamma)$ and $b_{\alpha}(\gamma)$ 
are ${\rm PGL}_2\otimes_{{\Bbb Z}} {\Bbb F}_2$-invariant 
functions on ${\rm Rep}_2(\Gamma)_{u/{\Bbb F}_2, \alpha}$.  
Hence any $(a, b)$-coefficients are 
${\rm PGL}_2\otimes_{{\Bbb Z}} {\Bbb F}_2$-invariant and 
$\pi_{\Gamma, u/{\Bbb F}_2, \alpha} : 
{\rm Rep}_2(\Gamma)_{u/{\Bbb F}_2, \alpha} 
\to {\rm Ch}_2(\Gamma)_{u/{\Bbb F}_2, \alpha}$ 
is ${\rm PGL}_2\otimes_{\Bbb Z} {\Bbb F}_2$-equivariant. 
This is another proof of Proposition \ref{prop:abinv}.  
\end{remark}

\begin{definition}\label{def:prototypeuf2}\rm 
  For the free monoid $\Upsilon_1 = \langle \alpha_0 \rangle$ of rank $1$, we say that  
the morphism $\pi_{\Upsilon_1, u/{\Bbb F}_2, \alpha_0} : {\rm Rep}_{2}(\Upsilon_1)_{u/{\Bbb F}_2, \alpha_0} \to 
{\rm Ch}_2(\Upsilon_1)_{u/{\Bbb F}_2, \alpha_0}$ is the {\it prototype} in the unipotent mold over ${\Bbb F}_2$ case.   
Remark that ${\rm Rep}_{2}(\Upsilon_1)_{u/{\Bbb F}_2} = {\rm Rep}_{2}(\Upsilon_1)_{u/{\Bbb F}_2, \alpha_0}$  
and that ${\rm Ch}_2(\Upsilon_1)_{u/{\Bbb F}_2} = 
{\rm Ch}_2(\Upsilon_1)_{u/{\Bbb F}_2, \alpha_0}$.   
\end{definition}

\bigskip 

The coordinate ring $A_2(\Upsilon_1)_{u/{\Bbb F}_2, \alpha_0}^{\rm Ch}$ of ${\rm Ch}_2(\Upsilon_1)_{u/{\Bbb F}_2, \alpha_0}$ 
is isomorphic to $A_{1}(\Upsilon_1)_{{\Bbb F}_2}$. 
Indeed, $a(\alpha_0)=0$ and $b(\alpha_0)=1$. 
By induction on $n$, we can verify that 
$a(\alpha_0^n) = 0, b(\alpha_0^n) = d(\alpha_0)^{(n-1)/2}$ for each positive odd integer $n$ 
and that $a(\alpha_0^n) = d(\alpha_0)^{n/2}, b(\alpha_0^n) = 0$ for each positive even integer $n$.  
Hence $A_2(\Upsilon_1)_{u/{\Bbb F}_2, \alpha_0}^{\rm Ch}  
\cong A_{1}(\Upsilon_1)_{{\Bbb F}_2}$ and ${\rm Ch}_2(\Upsilon_1)_{u/{\Bbb F}_2, \alpha_0} 
\cong {\rm Rep}_1(\Upsilon_1)_{{\Bbb F}_2}$. 

By Theorem~\ref{corssf1}, $\pi : {\rm Rep}_2(\Upsilon_1)_{{\rm rk} 2} \to 
{\rm Ch}_2(\Upsilon_1)$ is a universal geometric quotient by 
${\rm PGL}_2$. Taking the base change of $\pi$ by 
${\rm Spec} \; {\Bbb F}_2\to {\rm Spec} \: {\Bbb Z}$, 
we have ${\rm Rep}_2(\Upsilon_1)_{{\rm rk} 2/{\Bbb F}_2} \to 
{\rm Ch}_2(\Upsilon_1)_{{\Bbb F}_2}$, where   
${\rm Rep}_2(\Upsilon_1)_{{\rm rk} 2/{\Bbb F}_2} 
:= {\rm Rep}_2(\Upsilon_1)_{{\rm rk} 2}\otimes_{\Bbb Z} {\Bbb F}_2$ 
and ${\rm Ch}_2(\Upsilon_1)_{{\Bbb F}_2} := 
{\rm Ch}_2(\Upsilon_1)\otimes_{\Bbb Z} {\Bbb F}_2$. 
Let $Z$ be the closed subscheme of ${\rm Ch}_2(\Upsilon_1)_{{\Bbb F}_2}$ 
defined by ${\rm tr}(\sigma_{\Upsilon_1}(\alpha_0)) = 0$. 
Since ${\rm Ch}_2(\Upsilon_1) = {\rm Spec} \: {\Bbb Z}[T, D]$, 
the affine ring of $Z$ is isomorphic to ${\Bbb F}_2[D]$.  
Hence $Z$ is isomorphic to ${\rm Ch}_2(\Upsilon_1)_{u/{\Bbb F}_2, \alpha_0}
\cong {\rm Rep}_1(\Upsilon_1)_{{\Bbb F}_2}$. 
Taking the base change of 
${\rm Rep}_2(\Upsilon_1)_{{\rm rk} 2/{\Bbb F}_2} \to 
{\rm Ch}_2(\Upsilon_1)_{{\Bbb F}_2}$ by $Z \hookrightarrow {\rm Ch}_2(\Upsilon_1)_{{\Bbb F}_2}$, 
we have the prototype $\pi_{\Upsilon_1, u/{\Bbb F}_2, \alpha_0} : 
{\rm Rep}_2(\Upsilon_1)_{u/{\Bbb F}_2, \alpha_0} \to 
{\rm Ch}_2(\Upsilon_1)_{u/{\Bbb F}_2, \alpha_0}$.  

By Lemma~\ref{lemma:basechangeGIT}, we have: 

\begin{theorem}\label{th:prototypeuf2geom}
The prototype $\pi_{\Upsilon_1, u/{\Bbb F}_2, \alpha_0} : 
{\rm Rep}_2(\Upsilon_1)_{u/{\Bbb F}_2, \alpha_0} \to 
{\rm Ch}_2(\Upsilon_1)_{u/{\Bbb F}_2, \alpha_0}$ is a universal 
geometric quotient by ${{\rm PGL}_2}\otimes_{{\Bbb Z}} {{\Bbb F}_2}$.   
\end{theorem}

\bigskip

Let $\Gamma$ be a group or a monoid.  
For $\alpha \in \Gamma$, we define the monoid homomorphism 
$\phi : \Upsilon_1 =\langle \alpha_0 \rangle \to \Gamma$ by 
$\alpha_0 \mapsto \alpha$.  
By restricting representations and $(a, b)$-coefficients of $\Gamma$ to 
those of $\Upsilon_1$ through $\phi$,
we can obtain the following 
commutative diagram:  
\[
\begin{array}{ccc}
  {\rm Rep}_2(\Gamma)_{u/{\Bbb F}_2, \alpha} & \to & {\rm Ch}_2(\Gamma)_{u/{\Bbb F}_2, \alpha} \\
  \downarrow & &\downarrow \\
  {\rm Rep}_2(\Upsilon_1)_{u/{\Bbb F}_2, \alpha_0} & \to &  {\rm Ch}_2(\Upsilon_1)_{u/{\Bbb F}_2, \alpha_0}. 
\end{array}
\]

Under this situation, we have the following lemma. 
\begin{lemma}\label{lemma:fibreproductuf2}
  The above diagram gives a fibre product. 
In particular, the morphism 
${\rm Rep}_2(\Gamma)_{u/{\Bbb F}_2, \alpha}  \to  {\rm Ch}_2(\Gamma)_{u/{\Bbb F}_2, \alpha}$ 
is obtained by base change of the prototype. 
\end{lemma}

{\it Proof.} 
Set $Z := {\rm Rep}_2(\Upsilon_1)_{u/{\Bbb F}_2, \alpha_0}\times_{{\rm Ch}_2(\Upsilon_1)_{u/{\Bbb F}_2, \alpha_0}}
{\rm Ch}_2(\Gamma)_{u/{\Bbb F}_2, \alpha}$. 
We claim that ${\rm Rep}_2(\Gamma)_{u/{\Bbb F}_2, \alpha} \to Z$ is an isomorphism. 
Let $X$ be an ${\Bbb F}_2$-scheme. 
Assume that an 
$X$-valued point $\rho \in {\rm Rep}_2(\Gamma)_{u/{\Bbb F}_2, \alpha}$
is sent to $(\rho', \sigma) \in Z$. 
We can regard the $X$-valued point $\sigma \in {\rm Ch}_2(\Gamma)_{u/{\Bbb F}_2, \alpha}$ 
as a pair $(d, (a, b))$ such that $a, b : \Gamma \to \Gamma(X, {\mathcal O}_X)$ are $(a, b)$-coefficients with 
respect to $(d, \alpha)$, where $d(\cdot) := \det(\rho(\cdot)) : \Gamma \to \Gamma(X, {\mathcal O}_X)$.   
Since 
\begin{eqnarray}\label{eqn-u-1}
\rho(\gamma) =  a(\gamma)I_2 + b(\gamma)\rho'(\alpha_{0}) 
\end{eqnarray}
for each $\gamma \in \Gamma$,   
$\rho$ is uniquely determined by $(\rho', \sigma)$.  

For an $X$-valued point $(\rho', \sigma) \in Z$, we define the map 
$\rho : \Gamma \to {\rm M}_2(\Gamma(X, {\mathcal O}_X))$ by (\ref{eqn-u-1}).
From Lemma \ref{lemma:uf2rep}, we see that $\rho$ is an $X$-valued point 
of ${\rm Rep}_2(\Gamma)_{u/{\Bbb F}_2, \alpha}$. 
Then the $X$-valued point $\rho$ 
is sent to $(\rho', \sigma) \in Z$. 
By these discussion, we see that 
${\rm Rep}_2(\Gamma)_{u/{\Bbb F}_2, \alpha} \to Z$ 
is an isomorphism, and hence that  
the diagram gives a fibre product. 
\qed

\bigskip 

\begin{theorem}\label{th:univgeomuf2alpha}  
The morphism $\pi_{\Gamma, u/{\Bbb F}_2, \alpha} : {\rm Rep}_2(\Gamma)_{u/{\Bbb F}_2, \alpha} \to 
{\rm Ch}_2(\Gamma)_{u/{\Bbb F}_2, \alpha}$ 
is a universal 
geometric quotient by ${\rm PGL}_2\otimes_{{\Bbb Z}} {\Bbb F}_2$ 
for each $\alpha \in \Gamma$. 
\end{theorem}

{\it Proof.} 
The statement follows from that $\pi_{\Gamma, u/{\Bbb F}_2, \alpha}$ 
is obtained by base change of 
the prototype. 
\qed 

\bigskip 

Let $\alpha, \beta \in \Gamma$.  
Let $U_{\alpha, \beta} \subseteq {\rm Ch}_2(\Gamma)_{u/{\Bbb F}_2, \alpha}$ 
be the open subscheme defined by $\{ b(\beta) \neq 0 \}$. 
The inverse image $\pi_{\Gamma, u/{\Bbb F}_2, \alpha}^{-1}(U_{\alpha, \beta})$ by  
$\pi_{\Gamma, u/{\Bbb F}_2, \alpha} : {\rm Rep}_2(\Gamma)_{u/{\Bbb F}_2, \alpha} \to 
{\rm Ch}_2(\Gamma)_{u/{\Bbb F}_2, \alpha}$ 
coincides with ${\rm Rep}_2(\Gamma)_{u/{\Bbb F}_2, \alpha} 
\cap {\rm Rep}_2(\Gamma)_{u/{\Bbb F}_2, \beta}$. 
Then 
$\pi_{\Gamma, u/{\Bbb F}_2, \alpha}^{-1}(U_{\alpha, \beta}) = 
{\rm Rep}_2(\Gamma)_{u/{\Bbb F}_2, \alpha} 
\cap {\rm Rep}_2(\Gamma)_{u/{\Bbb F}_2, \beta} 
\to U_{\alpha, \beta}$ is a universal geometric quotient 
by ${\rm PGL}_2\otimes_{{\Bbb Z}} {\Bbb F}_2$. 
Hence $U_{\alpha, \beta} \cong U_{\beta, \alpha}$, and 
let us denote the canonical isomorphism by 
$\varphi_{\alpha, \beta} : U_{\alpha, \beta} \to U_{\beta, \alpha}$.  
Note that 
$\varphi_{\beta, \alpha} = \varphi_{\alpha, \beta}^{-1}$. 
For $\alpha, \beta, \gamma \in \Gamma$, 
$\varphi_{\alpha, \beta}(U_{\alpha, \beta} \cap U_{\alpha, \gamma}) = U_{\beta, \alpha} \cap U_{\beta, \gamma}$ 
and $\varphi_{\alpha, \gamma} = \varphi_{\beta, \gamma} \circ \varphi_{\alpha, \beta}$ 
on $U_{\alpha, \beta} \cap U_{\alpha, \gamma}$.   
Gluing the schemes $\{ {\rm Ch}_2(\Gamma)_{u/{\Bbb F}_2, \alpha} \}_{\alpha \in \Gamma}$, 
we obtain a scheme, which we call ${\rm Ch}_2(\Gamma)_{u/{\Bbb F}_2}$  
 (for example, see \cite[Chap. II, Ex. 2.12]{Hartshorne}).  
Gluing $\{ \pi_{\Gamma, u/{\Bbb F}_2, \alpha} \}_{\alpha \in \Gamma}$, we also 
obtain $\pi_{\Gamma, u/{\Bbb F}_2} : {\rm Rep}_2(\Gamma)_{u/{\Bbb F}_2} 
\to {\rm Ch}_2(\Gamma)_{u/{\Bbb F}_2}$. 


\begin{corollary}\label{cor:univgeomuf2-total}
The morphism $\pi_{\Gamma, u/{\Bbb F}_2} : {\rm Rep}_2(\Gamma)_{u/{\Bbb F}_2} \to 
{\rm Ch}_2(\Gamma)_{u/{\Bbb F}_2}$ 
is a universal 
geometric quotient by ${\rm PGL}_2\otimes_{{\Bbb Z}} {\Bbb F}_2$.  
\end{corollary} 

\begin{remark}\label{remark:dpidet}\rm
By Definition \ref{def:aringchuf2}, ${\rm Ch}_2(\Gamma)_{u/{\Bbb F}_2, \alpha}$ 
is a ${\rm Rep}_1(\Gamma)_{{\Bbb F}_2}$-scheme for each $\alpha \in \Gamma$. 
Let $d_{\alpha} : {\rm Ch}_2(\Gamma)_{u/{\Bbb F}_2, \alpha} \to 
{\rm Rep}_1(\Gamma)_{{\Bbb F}_2}$ be the canonical morphism for each $\alpha  
\in \Gamma$. We can obtain a morphism $d : {\rm Ch}_2(\Gamma)_{u/{\Bbb F}_2} \to {\rm Rep}_1(\Gamma)_{{\Bbb F}_2}$ by gluing  the canonical morphisms 
$\{ d_{\alpha} \}_{\alpha \in \Gamma}$.  Hence 
${\rm Ch}_2(\Gamma)_{u/{\Bbb F}_2}$ is also 
a ${\rm Rep}_1(\Gamma)_{{\Bbb F}_2}$-scheme. 
Let us denote by $\det : {\rm Rep}_2(\Gamma)_{u/{\Bbb F}_2} 
\to {\rm Rep}_1(\Gamma)_{{\Bbb F}_2}$ the morphism corresponding to 
the character $\det(\sigma_{\Gamma, u/{\Bbb F}_2}(\cdot)))$. 
Then $d \circ \pi_{\Gamma, u/{\Bbb F}_2} = \det$. 
\end{remark} 

\bigskip  

The open subscheme $U_{\alpha, \beta} \subseteq {\rm Ch}_2(\Gamma)_{u/{\Bbb F}_2, \alpha}$ 
is affine and its coordinate ring is isomorphic to 
the localization $A_2(\Gamma)_{u/{\Bbb F}_2, \alpha}^{\rm Ch}[b_{\alpha}(\beta)^{-1}]$ 
of the coordinate ring $A_2(\Gamma)_{u/{\Bbb F}_2, \alpha}^{\rm Ch}$ of ${\rm Ch}_2(\Gamma)_{u/{\Bbb F}_2, \alpha}$ 
by $b_{\alpha}(\beta)^{-1}$.  
Here we denote by $a_{\alpha}, b_{\alpha}$ the $(a, b)$-coefficients 
of the universal representations $\sigma_{\Gamma, u/{\Bbb F}_2, \alpha}$ with respect to 
$(\det(\sigma_{\Gamma, u/{\Bbb F}_2, \alpha}), \alpha)$. 
Let $\varphi_{\alpha, \beta}^{\ast} : 
A_2(\Gamma)_{u/{\Bbb F}_2, \beta}^{\rm Ch}[b_{\beta}(\alpha)^{-1}] \to 
A_2(\Gamma)_{u/{\Bbb F}_2, \alpha}^{\rm Ch}[b_{\alpha}(\beta)^{-1}]$  
denote by the ring isomorphism associated to 
$\varphi_{\alpha, \beta} : U_{\alpha, \beta} \to 
U_{\beta, \alpha}$.   
Then $\varphi_{\alpha, \beta}^{\ast}(a_{\beta}(\gamma)) 
= a_{\alpha}(\gamma) + a_{\beta}(\alpha)b_{\alpha}(\gamma)$ 
and $\varphi_{\alpha, \beta}^{\ast}(b_{\beta}(\gamma)) 
= b_{\alpha}(\gamma)b_{\beta}(\alpha)$ for each $\gamma \in \Gamma$. 
Note that $b_{\beta}(\alpha) = b_{\alpha}(\beta)^{-1}$ 
and $a_{\beta}(\alpha)=-a_{\alpha}(\beta)b_{\alpha}(\beta)^{-1}$ on 
$U_{\alpha, \beta} \cong U_{\beta, \alpha}$. 
Since $A_2(\Gamma)_{u/{\Bbb F}_2, \alpha}^{\rm Ch}\otimes_{{\Bbb F}_2} A_2(\Gamma)_{u/{\Bbb F}_2, \beta}^{\rm Ch}  
\to A_2(\Gamma)_{u/{\Bbb F}_2, \alpha}^{\rm Ch}[b_{\alpha}(\beta)^{-1}]$ is a surjective 
ring homomorphism, 
the diagonal morphism $U_{\alpha, \beta} \stackrel{\Delta}{\to} 
{\rm Ch}_2(\Gamma)_{u/{\Bbb F}_2, \alpha}\times_{{\Bbb F}_2} {\rm Ch}_2(\Gamma)_{u/{\Bbb F}_2, \beta}$ 
is a closed immersion. 
Hence we have: 

\begin{proposition}\label{prop:separateduf2} 
For a group or a monoid $\Gamma$, 
${\rm Ch}_2(\Gamma)_{u/{\Bbb F}_2}$
is separated over ${\Bbb F}_2$. 
\end{proposition}

\begin{remark}\rm  
  The morphism $\pi_{\Gamma, u/{\Bbb F}_2} : {\rm Rep}_2(\Gamma)_{u/{\Bbb F}_2} \to 
{\rm Ch}_2(\Gamma)_{u/{\Bbb F}_2}$  is smooth and surjective for each group or monoid $\Gamma$. 
Indeed, the prototype $\pi_{\Upsilon_1, u/{\Bbb F}_2, \alpha_0} : {\rm Rep}_2(\Upsilon_1)_{u/{\Bbb F}_2, \alpha_0} 
\to {\rm Ch}_2(\Upsilon_1)_{u/{\Bbb F}_2, \alpha_0}$ 
is smooth and surjective because it is obtained by base change of 
$\pi : {\rm Rep}_2(\Upsilon_{1})_{{\rm rk} 2} \to {\rm Ch}_2(\Upsilon_1)$ and 
$\pi$ is smooth and surjective by Proposition~\ref{lemmassff}. 
Hence ${\rm Rep}_2(\Gamma)_{u/{\Bbb F}_2, \alpha} = \pi_{\Gamma, u/{\Bbb F}_2}^{-1}({\rm Ch}_2(\Gamma)_{u/{\Bbb F}_2, \alpha}) 
\to {\rm Ch}_2(\Gamma)_{u/{\Bbb F}_2, \alpha}$ is smooth and surjective for each $\alpha \in \Gamma$.  
Therefore, so is $\pi_{\Gamma, u/{\Bbb F}_2}$. 
\end{remark}

\begin{remark}\label{remark:sectionuf2}\rm 
For each point $x \in {\rm Ch}_2(\Gamma)_{u/{\Bbb F}_2}$, there exists 
a local section $s_{x} : V_x \to {\rm Rep}_2(\Gamma)_{u/{\Bbb F}_2}$ on 
a neighbourhood $V_x$ of $x$ such that 
$\pi_{\Gamma, u/{\Bbb F}_2} \circ s_{x} = id_{V_x}$. 
Indeed, take $\alpha \in \Gamma$ such that 
$x \in {\rm Ch}_2(\Gamma)_{u/{\Bbb F}_2, \alpha}$. 
The prototype ${\rm Rep}_2(\Upsilon_1)_{u/{\Bbb F}_2, \alpha_0} \to {\rm Ch}_2(\Upsilon_1)_{u/{\Bbb F}_2, \alpha_0}$ 
has a section $s$ since it is obtained by base change of 
${\rm Rep}_2(\Upsilon_1)_{{\rm rk} 2} \to 
{\rm Ch}_2(\Upsilon_1)$, which has a section 
(it has been defined just before Proposition~\ref{prop:sectionff}). 
By Lemma \ref{lemma:fibreproductuf2}, we see that 
${\rm Rep}_2(\Gamma)_{u/{\Bbb F}_2, \alpha} 
\to {\rm Ch}_2(\Gamma)_{u/{\Bbb F}_2, \alpha}$ has a section $s_{\Gamma, \alpha}$.  
Hence we can take ${\rm Ch}_2(\Gamma)_{u/{\Bbb F}_2, \alpha}$ 
as a neighbourhood $V_x$ of $x$. 
\end{remark}


\begin{lemma}\label{lemma:localequf2} 
Let $\rho_1, \rho_2$ be representations with unipotent mold over ${\Bbb F}_2$ for a group (or a monoid) 
$\Gamma$ on a scheme $X$ over ${\Bbb F}_2$. 
Let $f_i : X \to {\rm Rep}_2(\Gamma)_{u/{\Bbb F}_2}$ be the morphism 
associated to $\rho_i$ for $i = 1, 2$. 
If $\pi_{\Gamma, u/{\Bbb F}_2} \circ f_1 = \pi_{\Gamma, u/{\Bbb F}_2} \circ f_2 
: X \to {\rm Ch}_2(\Gamma)_{u/{\Bbb F}_2}$, then for each $x \in X$ 
there exists $P_x \in {\rm GL}_2(\Gamma(V_x, {\mathcal O}_X))$ 
on a neighbourhood $V_x$ of $x$ such that 
$P_x^{-1}\rho_1 P_x = \rho_2$ on $V_x$.  
\end{lemma} 

\prf
For $x \in X$, take $\alpha \in \Gamma$ such that 
$(\pi_{\Gamma, u/{\Bbb F}_2} \circ f_1) (x) = (\pi_{\Gamma, u/{\Bbb F}_2} \circ f_2) (x) 
\in {\rm Ch}_2(\Gamma)_{u/{\Bbb F}_2, \alpha}$.  
We may assume that $f_i : X \to {\rm Rep}_2(\Gamma)_{u/{\Bbb F}_2, \alpha}$ for 
$i =1, 2$ from the beginning. By Remark \ref{remark:sectionuf2}, 
$\pi_{\Gamma, u/{\Bbb F}_2, \alpha} : {\rm Rep}_2(\Gamma)_{u/{\Bbb F}_2, \alpha} 
\to {\rm Ch}_2(\Gamma)_{u/{\Bbb F}_2, \alpha}$ has a section $s_{\Gamma, \alpha}$.  
Let $\rho_3$ be the representations with unipotent mold over ${\Bbb F}_2$ on $X$ associated to 
$s_{\Gamma, \alpha} \circ \pi_{\Gamma, u/{\Bbb F}_2} \circ f_1 = s_{\Gamma, \alpha} \circ 
\pi_{\Gamma, u/{\Bbb F}_2} \circ f_2$. 
Then $\rho_i(\gamma) = a(\gamma) I_2 + b(\gamma) \rho_i(\alpha)$ for each 
$\gamma \in \Gamma$ and $i = 1, 2, 3$, where $(a, b)$ is the $(a, b)$-coefficients with respect to 
$(d, \alpha)$ and $d(\cdot) := \det(\rho_1(\cdot)) = \det(\rho_2(\cdot)) = \det(\rho_3(\cdot))$. 

Note that 
$\rho_3(\alpha) = 
\left(
\begin{array}{cc}
0 & -D \\
1 & 0 \\
\end{array}
\right)$ and that 
$D = d(\alpha)$.  
There exist $Q_1, Q_2 \in {\rm GL}_2(\Gamma(V_x, {\mathcal O}_X))$ on a
neighbourhood $V_x$ of $x$ such that 
$Q_1^{-1}\rho_1(\alpha)Q_1 = \rho_3(\alpha)$ and 
$Q_2^{-1}\rho_2(\alpha)Q_2 = \rho_3(\alpha)$ by Lemma \ref{lemma-norm-nonsca}. 
Then 
$Q_1^{-1}\rho_1(\gamma)Q_1 = \rho_3(\gamma)$ and 
$Q_2^{-1}\rho_2(\gamma)Q_2 = \rho_3(\gamma)$ for each $\gamma \in \Gamma$, 
and hence 
$(Q_1 Q_2^{-1})^{-1} \rho_1 (Q_1 Q_2^{-1}) = \rho_2$ on $V_x$. 
This completes the proof. 
\qed 

\bigskip 

Let us define 
${\mathcal E}q \:\mathcal{U}_2(\Gamma)_{{\Bbb F}_2}$ as the  
sheafification of the following contravariant functor with respect to Zariski topology:   
\[
\begin{array}{ccl} 
({\bf Sch}/{\Bbb F}_2)^{op} & \to & ({\bf Sets}) \\
 X & \mapsto & 
\left\{ 
\rho \; 
\begin{array}{|l} 
 \mbox{ rep. 
with unip. mold } \\
\mbox{ over ${\Bbb F}_2$ for $\Gamma$ on } X 
\end{array} 
\right\}\Big/\sim.
\end{array}
\]

By a {\it generalized representation with unipotent mold over ${\Bbb F}_2$} for $\Gamma$ on 
an ${\Bbb F}_2$-scheme $X$, we understand pairs $\{ (U_i, \rho_i) \}_{i \in I}$ 
of an open set $U_i$ and a representation $\rho_i : \Gamma \to {\rm M}_2(\Gamma(U_i, {\mathcal O}_{X}))$ 
with unipotent mold over ${\Bbb F}_2$ 
satisfying the following two conditions:
\begin{enumerate}
\item $\cup_{i \in I} U_ i = X$, 
\item for each $x \in U_i \cap U_j$, there exists $P_x \in 
{\rm GL}_2(\Gamma(V_x, {\mathcal O}_X))$ on 
a neighbourhood $V_x \subseteq U_i \cap U_j$ of $x$ such that 
$P_x^{-1} \rho_i P_x = \rho_j$ on $V_x$. 
\end{enumerate}

Generalized representations $\{ (U_i, \rho_i) \}_{i \in I}$ and $\{ (V_j, \sigma_j) \}_{j \in J}$ 
with unipotent mold over ${\Bbb F}_2$ 
are called {\it equivalent} if $\{ (U_i, \rho_i) \}_{i \in I} \cup \{ (V_j, \sigma_j) \}_{j \in J}$ 
is a generalized representation with unipotent mold over ${\Bbb F}_2$ again.  
We easily see that 
${\mathcal E}q \:\mathcal{U}_2(\Gamma)_{{\Bbb F}_2}(X)$ is the 
set of equivalence classes of 
generalized representations with unipotent mold over ${\Bbb F}_2$ 
for $\Gamma$ on an ${\Bbb F}_2$-scheme $X$.

\begin{theorem}\label{th:moduliuf2}  
The scheme ${\rm Ch}_2(\Gamma)_{u/{\Bbb F}_2}$ is a fine moduli scheme  
associated to the functor ${\mathcal E}q \:\mathcal{U}_2(\Gamma)_{{\Bbb F}_2}$ 
for a group or a monoid $\Gamma$:  
\[
\begin{array}{ccccl}
{\mathcal E}q \:\mathcal{U}_2(\Gamma)_{{\Bbb F}_2} & : & ({\bf Sch}/{\Bbb F}_2)^{op} & \to & ({\bf Sets}) \\
 & & X & \mapsto & 
\left\{ 
\begin{array}{r}  
\mbox{ gen. rep. 
with unipotent  } \\ 
\mbox{ mold over ${\Bbb F}_2$ for $\Gamma$ on } X  
\end{array} 
\right\} 
\Big/\sim.
\end{array}
\] 
In other words, ${\rm Ch}_2(\Gamma)_{u/{\Bbb F}_2}$ represents 
the functor ${\mathcal E}q \:\mathcal{U}_2(\Gamma)_{{\Bbb F}_2}$. 
The moduli ${\rm Ch}_2(\Gamma)_{u/{\Bbb F}_2}$ is separated over ${\Bbb F}_2$; 
if $\Gamma$ is a finitely generated group or monoid, then 
${\rm Ch}_2(\Gamma)_{u/{\Bbb F}_2}$ is of finite type over ${\Bbb F}_2$. 
\end{theorem}

\prf 
Since $\pi_{\Gamma, u/{\Bbb F}_2} : {\rm Rep}_2(\Gamma)_{u/{\Bbb F}_2} 
\to {\rm Ch}_2(\Gamma)_{u/{\Bbb F}_2}$ is a 
${\rm PGL}_2\otimes_{{\Bbb Z}} {\Bbb F}_2$-equivariant morphism, we can define a canonical morphism 
${\mathcal E}q \: \mathcal{U}_2(\Gamma)_{{\Bbb F}_2} \to 
h_{{\rm Ch}_2(\Gamma)_{u/{\Bbb F}_2}}   := {\rm Hom}(-, {\rm Ch}_2(\Gamma)_{u/{\Bbb F}_2})$.  
We define a morphism $h_{{\rm Ch}_2(\Gamma)_{u/{\Bbb F}_2}} 
\to {\mathcal E}q \:\mathcal{U}_2(\Gamma)_{{\Bbb F}_2}$ as follows.  
Let $g \in h_{{\rm Ch}_2(\Gamma)_{u/{\Bbb F}_2}}(X)$ with an ${\Bbb F}_2$-scheme $X$. 
For each $x \in X$, take $\alpha_x \in \Gamma$ such that 
$g(x) \in {\rm Ch}_2(\Gamma)_{u/{\Bbb F}_2, \alpha_x}$. 
By using the section $s_{\Gamma, \alpha_x} : {\rm Ch}_2(\Gamma)_{u/{\Bbb F}_2, \alpha_x} 
\to {\rm Rep}_2(\Gamma)_{u/{\Bbb F}_2, \alpha_x}$ in Remark \ref{remark:sectionuf2}, 
we can define a representation $\rho_x$ with unipotent mold over ${\Bbb F}_2$ 
on a neighbourhood 
$U_x$ of $x$. By Lemma \ref{lemma:localequf2}, 
we see that $\{ (U_x, \rho_x) \}_{x \in X} \in 
{\mathcal E}q \:\mathcal{U}_2(\Gamma)_{{\Bbb F}_2}(X)$ and that the morphism 
$h_{{\rm Ch}_2(\Gamma)_{u/{\Bbb F}_2}}  
\to {\mathcal E}q \:\mathcal{U}_2(\Gamma)_{{\Bbb F}_2}$ is well-defined.  
It is easy to see that ${\rm Ch}_2(\Gamma)_{u/{\Bbb F}_2}$ represents 
the functor ${\mathcal E}q \:\mathcal{U}_2(\Gamma)_{{\Bbb F}_2}$. 

By Proposition~\ref{prop:separateduf2}, 
${\rm Ch}_2(\Gamma)_{u/{\Bbb F}_2}$ is separated over ${\Bbb F}_2$. 
If $\Gamma$ is finitely generated, 
then we can verify that 
$A_2(\Gamma)_{u/{\Bbb F}_2, \alpha}^{\rm Ch}$ is 
a finitely generated algebra over ${\Bbb F}_2$ 
in a similar way as Remark \ref{remark:fgu}.  
Let $S = \{ \alpha_1, \ldots, \alpha_n \}$ be a 
set of generators of $\Gamma$.  
Then ${\rm Ch}_2(\Gamma)_{u/{\Bbb F}_2}$ 
is covered by finitely many affine open subschemes 
${\rm Ch}_2(\Gamma)_{u/{\Bbb F}_2, \alpha_i} (1 \le i \le n)$. Hence  
${\rm Ch}_2(\Gamma)_{u/{\Bbb F}_2}$ is of finite type over ${\Bbb F}_2$.  
\qed 

\bigskip 

In the following Example \ref{ex:describe-chuf2m}, 
we describe ${\rm Ch}_2(\Upsilon_m)_{u/{\Bbb F}_2}$ 
for the free monoid $\Upsilon_m = \langle \alpha_1, \ldots, \alpha_m \rangle$ 
of rank $m$. This description has been inspired by the referee. 

\begin{example}\rm\label{ex:describe-chuf2m} 
Let us describe ${\rm Ch}_2(\Upsilon_m)_{u/{\Bbb F}_2}$ 
for the free monoid $\Upsilon_m = \langle \alpha_1, \ldots, \alpha_m \rangle$ 
of rank $m$. 
Put $C(m) := {\rm Ch}_2(\Upsilon_m)_{u/{\Bbb F}_2}$ and 
$C(m)_i := {\rm Ch}_2(\Upsilon_m)_{u/{\Bbb F}_2, \alpha_i}$ for $1 \le i \le m$. 
Let us denote by $A(m)_{i}$ the $A_1(\Upsilon_m)_{{\Bbb F}_2}$-algebra 
$A_2(\Upsilon_m)^{\rm Ch}_{u/{\Bbb F}_2, \alpha_i}$ for 
$1 \le i \le m$ in Definition \ref{def:aringchuf2}.  
We can write 
$A_1(\Upsilon_m)_{{\Bbb F}_2} = {\Bbb F}_2[d(\alpha_1), \ldots, 
d(\alpha_m)]$ and  
\[ 
A(m)_i = {\Bbb F}_2[d(\alpha_j), a_i(\alpha_j), b_i(\alpha_j) \mid 1 \le j  \le m]/I(m)_i,  
\] 
where $I(m)_i$ is the ideal of ${\Bbb F}_2[d(\alpha_j), a_i(\alpha_j), b_i(\alpha_j) \mid 1 \le j  \le m]$ generated by $a_i(\alpha_i)$, $b_i(\alpha_i)-1$, and 
$d(\alpha_j) -a_i(\alpha_j)^2 - b_i(\alpha_j)^2d(\alpha_i)$ for $1 \le j  \le m$. 
Note that 
\begin{multline*} 
C(m)_{i}  \cong {\Bbb A}_{{\Bbb F}_2}^{2m-1}  =  \{ (a_{i}(\alpha_1), 
\ldots, a_{i}(\alpha_{i-1}), 
a_{i}(\alpha_{i+1}), \ldots, a_{i}(\alpha_{m}), \\ 
 b_{i}(\alpha_1), \ldots, b_{i}(\alpha_{i-1}), 
b_{i}(\alpha_{i+1}), \ldots, b_{i}(\alpha_{m}),  
d(\alpha_{i})) \in {\Bbb A}_{{\Bbb F}_2}^{2m-1} 
\}.
\end{multline*}    
We set $U_{ij} := \{ b_i(\alpha_j) \neq 0 \} \subset 
C(m)_i = {\rm Spec} A(m)_i$. 
For $1 \le i \neq j \le m$, the isomorphism $\varphi_{ij} : U_{ij} \to U_{ji}$ 
is given by the ${\Bbb F}_2$-algebra isomorphism 
$\varphi^{\ast}_{ij} : A(m)_{j}[b_j(\alpha_i)^{-1}] \to 
A(m)_{i}[b_i(\alpha_j)^{-1}]$ which is defined by 
$\varphi^{\ast}_{ij}(a_j(\alpha_k)) = a_i(\alpha_k) + b_i(\alpha_k) a_i(\alpha_j)/b_i(\alpha_j)$, 
$\varphi^{\ast}_{ij}(b_j(\alpha_k)) = b_i(\alpha_k)/b_i(\alpha_j)$, and 
$\varphi^{\ast}_{ij}(d(\alpha_k)) 
 = d(\alpha_k)$ for $1 \le k \le m$. 

On the other hand, let us define the closed subvariety $D(m)$ of 
${\Bbb P}_{{\Bbb F}_2}^{m^2+m-1} \times {\Bbb A}_{{\Bbb F}_2}^{m}$ over 
${\Bbb F}_2$ in the following way: 
\begin{multline*} 
D(m) := \{ ([a_{ij} : b_1: \cdots : b_m ]_{1 \le i, j \le m}, (d_1, \ldots, d_m) ) \in {\Bbb P}_{{\Bbb F}_2}^{m^2+m-1} \times {\Bbb A}_{{\Bbb F}_2}^{m}
\mid \\ 
a_{ji} = a_{ij}, \; a_{ii}=0 \; (1 \le i, j \le m),  \\ 
a_{ij}^2+b_{i}^2d_j + b_{j}^2d_{i} = 0  \; (1 \le i, j \le m),  \\
a_{ij}b_{k} + a_{jk}b_{i} + a_{ki}b_{j}=0 \; (1 \le i, j, k \le m) 
\}. 
\end{multline*} 
(Note that $D(m)$ can also be defined as a closed subvariety of 
${\Bbb P}_{{\Bbb F}_2}^{\frac{m(m+1)}{2}-1} \times {\Bbb A}_{{\Bbb F}_2}^{m}$ 
by using homogeneous coordinates $\{ a_{ij} \}_{1 \le i < j \le m}$ and 
$\{ b_i \}_{1 \le i \le m}$.)
Put $D(m)_i := \{ b_i \neq 0 \} \subset D(m)$ for 
$1 \le i  \le m$. 
By using inhomogeneous coordinates 
$\overline{a}_{jk} = a_{jk}/b_{i}$ and $\overline{b}_{j} = b_{j}/b_{i}$ 
for $1 \le j, k \le m$, we easily see that 
\begin{multline*} 
D(m)_i \cong {\Bbb A}_{{\Bbb F}_2}^{2m-1} \\ 
= \{ (\overline{a}_{i1}, \ldots, \overline{a}_{i, i-1}, 
\overline{a}_{i, i+1}, \ldots, \overline{a}_{im}, \overline{b}_{1}, \ldots, \overline{b}_{i-1}, 
\overline{b}_{i+1}, \ldots, \overline{b}_{m}, d_{i}) \in {\Bbb A}_{{\Bbb F}_2}^{2m-1} 
\}.
\end{multline*}     
Note that $\displaystyle D(m) = \cup_{i=1}^m D(m)_i$. 

Let us define the ${\Bbb F}_2$-isomorphism 
$\psi_i : C(m)_i \to D(m)_i \subset D(m)$ by the 
${\Bbb F}_2$-algebra isomorphism 
$\psi^{\ast}_i : A(D(m)_i) \to A(m)_i$ which is defined by  
$\psi^{\ast}_i(\overline{a}_{ij}) = a_{i}(\alpha_j)$, 
$\psi^{\ast}_i(\overline{b}_{j}) = b_{i}(\alpha_j)$ for $j \neq i$, and 
$\psi^{\ast}_i(d_{i}) = d(\alpha_i)$, where 
$A(D(m)_i)$ is the coordinate ring of $D(m)_i$. 
We can glue $\{ \psi_i : C(m)_i \to D(m)_i \}_{ 1 \le i  \le m }$ 
(remark that $-1$ equals to $1$ in 
characteristic $2$ for the verification), and 
hence we have an isomorphism $\psi : C(m) \to D(m)$.   
Note that ${\rm Ch}_2(\Upsilon_m)_{u/{\Bbb F}_2} = C(m) \cong D(m)$ 
is a $(2m-1)$-dimensional smooth irreducible variety over ${\Bbb F}_2$.

Let $\rho_i$ be the representation of $\Upsilon_m$ 
on $D(m)_i \cong C(m)_i$ 
defined by  
\[
\rho_i(\alpha_j) = \left(
\begin{array}{cc}
\overline{a}_{ij} & \overline{b}_{j} d_i \\ 
\overline{b}_{j} & \overline{a}_{ij} 
\end{array} 
\right) 
\] 
for $1 \le j \le m$. Then $\{ (C(m)_i, \rho_i) \}_{1 \le i \le m}$ 
is the universal equivalence class of generalized representations with unipotent 
mold over ${\Bbb F}_2$ of $\Upsilon_m$ on ${\rm Ch}_2(\Upsilon_m)_{u/{\Bbb F}_2}$. 
\end{example}

\begin{remark}\label{remark:alguf2}\rm
Let $A$ be an associative algebra over a commutative ring $R$ over ${\Bbb F}_2$. 
For a $2$-dimensional representation $\rho$ of $A$ on an $R$-scheme $X$, 
$\rho$ is called a {\it representation with unipotent mold over ${\Bbb F}_2$} if 
the subalgebra $\rho(A)$ of ${\rm M}_2({\mathcal O}_X)$ generates  
a unipotent mold over ${\Bbb F}_2$ on $X$. 
In a similar way as group or monoid cases, 
we can define generalized representations with unipotent mold 
over ${\Bbb F}_2$ for $A$ 
on an $R$-scheme $X$. The contravariant functor 
${\mathcal E}q \;\mathcal{U}_2(A)_{{\Bbb F}_2}$ from the category of $R$-schemes to the 
category of sets is defined as 
\[
\begin{array}{ccccl}
{\mathcal E}q \;\mathcal{U}_2(A)_{{\Bbb F}_2} & : & ({\bf Sch}/R)^{op} & \to & ({\bf Sets}) \\
 & & X & \mapsto & 
\left\{ 
\begin{array}{r}  
\mbox{ gen. rep. with unipotent 
 } \\ 
\mbox{  mold over ${\Bbb F}_2$ for $A$ on } X  
\end{array} 
\right\} 
\Big/\sim.
\end{array}
\] 
We can construct the fine moduli ${\rm Ch}_2(A)_{u/{\Bbb F}_2}$ associated to 
${\mathcal E}q \;\mathcal{U}_2(A)_{{\Bbb F}_2}$ in the same way as Theorem \ref{th:moduliuf2} (for details, see Remark~\ref{remark:constructuf2}).  
The moduli ${\rm Ch}_2(A)_{u/{\Bbb F}_2}$ is separated over $R$. 
If $A$ is a finitely generated algebra over $R$, then 
${\rm Ch}_2(A)_{u/{\Bbb F}_2}$ is of finite type over $R$.  
\end{remark} 

\begin{remark}\label{remark:constructuf2}\rm 
For an associative algebra $A$ over a commutative ring $R$ over 
${\Bbb F}_2$, we can construct ${\rm Ch}_2(A)_{u/{\Bbb F}_2}$ in the following way. 
We define the contravariant functor 
${\rm Rep}_1'(A)$ from the category of $R$-schemes 
to the category of sets by 
\[ 
{\rm Rep}_1'(A)(X) := \left\{ 
d : A \to \Gamma(X, {\mathcal O}_X) 
\begin{array}{|c} 
d \mbox{ is a ring homomorphism and } \\  
d(\alpha a) = \alpha^2 d(a) \mbox{ for } a \in A, \alpha \in R  
\end{array} 
\right\} 
\] 
for each $R$-scheme $X$.  
The functor ${\rm Rep}_1'(A)$ is representable by an affine scheme, and 
let us denote its coordinate ring by $A_1'(A)$.  
Let $d : A \to A_1'(A)$ be the universal ring homomorphism. 
Let ${\rm Rep}_2(A)$ be the representation variety of degree $2$ for 
$A$ over $R$ introduced in Remark \ref{remark:constructu}.  
Let ${\rm Rep}_2(A)_{u/{\Bbb F}_2}$ be the subscheme of 
${\rm Rep}_2(A)$ consisting of representations 
with unipotent mold over ${\Bbb F}_2$, and 
$\sigma_{A, u/{\Bbb F}_2} : A \to 
{\rm M}_2(\Gamma({\rm Rep}_2(A)_{u/{\Bbb F}_2}, {\mathcal O}_{{\rm Rep}_2(A)_{u/{\Bbb F}_2}}))$ 
the universal representation with unipotent mold over ${\Bbb F}_2$.  
Note that $\det \sigma_{A, u/{\Bbb F}_2} : A \to \Gamma({\rm Rep}_2(A)_{u/{\Bbb F}_2}, 
{\mathcal O}_{{\rm Rep}_2(A)_{u/{\Bbb F}_2}})$ gives a morphism 
${\rm Rep}_2(A)_{u/{\Bbb F}_2} \to {\rm Rep}_1'(A)$. 
For $c \in A$, define the open subscheme 
${\rm Rep}_2(A)_{u/{\Bbb F}_2, c} := \{ \langle I_2, 
\sigma_{A, u/{\Bbb F}_2}(c) \rangle \mbox{ is a unipotent mold over ${\Bbb F}_2$ } \}$ of ${\rm Rep}_2(A)_{u/{\Bbb F}_2}$.  
For a scheme $X$ over $A_1'(A)$ and for $c \in A$, we say that 
$a, b : A \to \Gamma(X, {\mathcal O}_X)$ are {\it $(a, b)$-coefficients} with 
respect to $(d, c)$ on $X$ if $a, b$ are $R$-linear maps satisfying 
\begin{eqnarray*}
a(1)=1, a(c)=0, b(1)=0, b(c)=1, \\ 
a(c_1c_2) = a(c_1)a(c_2)+b(c_1)b(c_2)d(c),  \\
b(c_1c_2)=a(c_1)b(c_2)+b(c_1)a(c_2), \\
a(c_1)^2+b(c_1)^2d(c)=d(c_1), 
\end{eqnarray*} 
for all $c_1, c_2 \in A$. Here $d : A \to \Gamma(X, {\mathcal O}_X)$ 
denotes the ring homomorphism associated to $X \to {\rm Rep}_1'(A)$. 
There exists a commutative ring $A_2(A)_{u/{\Bbb F}_2, c}^{\rm Ch}$ over 
$A_1'(A)$ such that ${\rm Ch}_2(A)_{u/{\Bbb F}_2, c} := 
{\rm Spec} A_2(A)_{u/{\Bbb F}_2, c}^{\rm Ch}$ represents 
the functor corresponding $X$ to the set of 
$(a, b)$-coefficients with respect to $(d, c)$ on $X$ for each scheme 
$X$ over $A_1'(A)$.  
In a similar way as 
group or monoid cases, 
we can define a ${\rm Rep}_1'(A)$-morphism $\pi_{c} : {\rm Rep}_2(A)_{u/{\Bbb F}_2, c} 
\to {\rm Ch}_2(A)_{u/{\Bbb F}_2, c}$. 
We see that $\pi_{c}$ is a universal geometric quotient by 
${\rm PGL}_2\otimes_{\Bbb Z} R$. 
Gluing schemes $\{ {\rm Ch}_2(A)_{u/{\Bbb F}_2, c} \}_{c \in A}$, we have 
a scheme ${\rm Ch}_2(A)_{u/{\Bbb F}_2}$ over ${\rm Rep}_1'(A)$.   
Gluing $\{ \pi_{c} \}_{c \in A}$, we also have 
a morphism $\pi : {\rm Rep}_2(A)_{u/{\Bbb F}_2} \to {\rm Ch}_2(A)_{u/{\Bbb F}_2}$ 
over ${\rm Rep}_1'(A)$. 
Hence $\pi$ is a universal geometric quotient by 
${\rm PGL}_2\otimes_{\Bbb Z} R$ and ${\rm Ch}_2(A)_{u/{\Bbb F}_2}$ 
represents ${\mathcal E}q \;\mathcal{U}_2(A)_{{\Bbb F}_2}$. 
\end{remark}

\begin{remark}\rm 
We have introduced the notion of 
generalized representations with unipotent mold over ${\Bbb F}_2$ 
for describing the moduli functors ${\mathcal E}q\mathcal{U}_2(\Gamma)_{{\Bbb F}_2}$  
and ${\mathcal E}q\mathcal{U}_2(A)_{{\Bbb F}_2}$.  However, 
the moduli functors can also be described as 
${\mathcal E}q\mathcal{U}'_2(\Gamma)_{{\Bbb F}_2}$  
and ${\mathcal E}q\mathcal{U}'_2(A)_{{\Bbb F}_2}$ by using 
the notion of representations generating sheaves of algebras 
which define unipotent molds over ${\Bbb F}_2$.    
More precisely, see \S 8. 
\end{remark}

\section{Another approach for unipotent molds over ${\Bbb F}_2$ } 
In this section, we construct the moduli scheme 
${\rm Ch}_2(\Gamma)_{u/{\Bbb F}_2}$ in a different way 
from \S 6. 
When we take the quotient of ${\rm Rep}_2(\Gamma)_{u/{\Bbb F}_2}$ by 
${\rm PGL}_2\otimes_{{\Bbb Z}} {\Bbb F}_2$, 
we need to introduce the notion of $(a, b)$-coefficients 
because there exist no eigenvalue of $\sigma_{\Gamma, u/{\Bbb F}_2}(\gamma)$ 
on ${\rm Rep}_2(\Gamma)_{u/{\Bbb F}_2}$ in general, 
where $\sigma_{\Gamma, u/{\Bbb F}_2}$ is the universal 
representation on  ${\rm Rep}_2(\Gamma)_{u/{\Bbb F}_2}$. 
However, we can obtain eigenvalues of $\sigma_{\Gamma, u/{\Bbb F}_2}(\gamma)$ 
by taking the pull-back of $\sigma_{\Gamma, u/{\Bbb F}_2}$ by 
a faithfully flat finite morphism $p
 : \widetilde{{\rm Rep}_2(\Gamma)}_{u/{\Bbb F}_2}    
\to {\rm Rep}_2(\Gamma)_{u/{\Bbb F}_2}$. 
Then by discussing derivations we can construct a universal geometric quotient 
$\tilde{\pi}_{\Gamma, u/{\Bbb F}_2} :   
\widetilde{{\rm Rep}_2(\Gamma)}_{u/{\Bbb F}_2} 
\to \widetilde{{\rm Ch}_2(\Gamma)}_{u/{\Bbb F}_2}$ 
by  ${\rm PGL}_2\otimes_{{\Bbb Z}} {\Bbb F}_2$ 
in the same way as the unipotent mold 
($ch \neq 2$) case in \S 5. 
Considering the ``descent'' of $\tilde{\pi}_{\Gamma, u/{\Bbb F}_2} :   
\widetilde{{\rm Rep}_2(\Gamma)}_{u/{\Bbb F}_2} 
\to \widetilde{{\rm Ch}_2(\Gamma)}_{u/{\Bbb F}_2}$, 
we have a universal geometric quotient 
$\pi_{\Gamma, u/{\Bbb F}_2} :  
{\rm Rep}_2(\Gamma)_{u/{\Bbb F}_2}   
\to {\rm Ch}_2(\Gamma)_{u/{\Bbb F}_2}$. 
In this section, we will use the same notation as \S 6. 
Without Lemma \ref{lemma:fibreproductuf2}, we will 
prove Theorem \ref{th:univgeomuf2alpha}. 
It should be pointed out that 
this section was inspired by the referee.

\bigskip 

Let $\Gamma$ be a group or a monoid. 
For $\alpha \in \Gamma$, let us consider 
the scheme ${\rm Rep}_2(\Gamma)_{u/{\Bbb F}_2, \alpha}$ 
over ${\Bbb F}_2$. 
Recall that  
\[
{\rm Rep}_2(\Gamma)_{u/{\Bbb F}_2, \alpha} 
= \left\{ 
\rho \in {\rm Rep}_2(\Gamma)_{u/{\Bbb F}_2}  \; 
\begin{array}{|l} 
 I_2 \mbox{ and } \rho(\alpha) \mbox{ generate } \\ 
 \mbox{ a unipotent mold over } {\Bbb F}_2 
\end{array} 
\right\}. 
\] 
Denote by $\sigma_{\Gamma, u/{\Bbb F}_2, \alpha}$ the universal 
representation with unipotent mold over ${\Bbb F}_2$ on 
${\rm Rep}_2(\Gamma)_{u/{\Bbb F}_2, \alpha}$.  
There exists no eigenvalue of $\sigma_{\Gamma, u/{\Bbb F}_2, \alpha}(\alpha)$ 
on ${\rm Rep}_2(\Gamma)_{u/{\Bbb F}_2, \alpha}$ in general,  
and hence we will construct a 
faithfully flat finite morphism $p_{\alpha} :   
\widetilde{{\rm Rep}_2(\Gamma)}_{u/{\Bbb F}_2, \alpha}   
\to {\rm Rep}_2(\Gamma)_{u/{\Bbb F}_2, \alpha}$ such that 
there exist eigenvalues of $p_{\alpha}^{\ast}(\sigma_{\Gamma, u/{\Bbb F}_2, \alpha}(\alpha))$ on 
$\widetilde{{\rm Rep}_2(\Gamma)}_{u/{\Bbb F}_2, \alpha}$. 

\begin{definition}\label{def:widetilde}\rm 
Let us define 
a quasi-coherent sheaf ${\mathcal A}_{\alpha}$ 
of ${\mathcal O}_{{\rm Rep}_2(\Gamma)_{u/{\Bbb F}_2, \alpha}}$-algebras 
on ${\rm Rep}_2(\Gamma)_{u/{\Bbb F}_2, \alpha}$  
by 
\[
{\mathcal A}_{\alpha} := {\mathcal O}_{{\rm Rep}_2(\Gamma)_{u/{\Bbb F}_2, \alpha}}[ X_{\alpha} ]/( X_{\alpha}^2 - \det \sigma_{\Gamma, u/{\Bbb F}_2, \alpha}(\alpha) ). 
\] 
Then set $\widetilde{{\rm Rep}_2(\Gamma)}_{u/{\Bbb F}_2, \alpha}    
:= {\mathcal Spec} {\mathcal A}_{\alpha}$. 
\end{definition}

\begin{remark}\rm 
The canonical morphism $p_{\alpha} :   
\widetilde{{\rm Rep}_2(\Gamma)}_{u/{\Bbb F}_2, \alpha}   
\to {\rm Rep}_2(\Gamma)_{u/{\Bbb F}_2, \alpha}$ 
is faithfully flat and finite. 
Since ${\rm tr} (\sigma_{\Gamma, u/{\Bbb F}_2, \alpha}(\alpha)) = 0$, $X_{\alpha}$ is an 
eigenvalue of $p_{\alpha}^{\ast}(\sigma_{\Gamma, u/{\Bbb F}_2, \alpha}(\alpha))$ on 
$\widetilde{{\rm Rep}_2(\Gamma)}_{u/{\Bbb F}_2, \alpha}$. 
We see that 
\[  
\widetilde{{\rm Rep}_2(\Gamma)}_{u/{\Bbb F}_2, \alpha}    
= \left\{ 
(\rho, X)   \; 
\begin{array}{|l} 
\rho \in {\rm Rep}_2(\Gamma)_{u/{\Bbb F}_2, \alpha}  
 \mbox{ and } 
X^2 = \det \rho(\alpha) 
\end{array} 
\right\}. 
\] 
\end{remark} 

\bigskip 
 
For simplicity, we put 
$\widetilde{R}_{\alpha} :=  
\widetilde{{\rm Rep}_2(\Gamma)}_{u/{\Bbb F}_2, \alpha}$ and 
$R_{\alpha} := {\rm Rep}_2(\Gamma)_{u/{\Bbb F}_2, \alpha}$. 
Remark that 
${\mathcal O}_{R_{\alpha}}[\sigma_{\Gamma, u/{\Bbb F}_2, \alpha}(\Gamma)] 
= {\mathcal O}_{R_{\alpha}}\cdot  I_2 + 
{\mathcal O}_{R_{\alpha}}\cdot \sigma_{\Gamma, u/{\Bbb F}_2, \alpha}(\alpha)$ 
is a unipotent mold over ${\Bbb F}_2$ on $R_{\alpha}$.   
For each $\gamma \in \Gamma$, we can write 
$\sigma_{\Gamma, u/{\Bbb F}_2, \alpha}(\gamma) = a_{\alpha}(\gamma) I_2 + 
b_{\alpha}(\gamma) \sigma_{\Gamma, u/{\Bbb F}_2, \alpha}(\alpha)$. 
Note that $a_{\alpha}(e) = 1, b_{\alpha}(e)=0$ and that 
$a_{\alpha}(\alpha) = 0, b_{\alpha}(\alpha)=1$. 
 
\begin{definition}\rm 
For each $\gamma$, we define $r_{\alpha}(\gamma) 
\in {\mathcal O}_{\widetilde{R}_{\alpha}}(\widetilde{R}_{\alpha})$ by 
$r_{\alpha}(\gamma) = a_{\alpha}(\gamma) + b_{\alpha}(\gamma) X_{\alpha}$. 
(In the sequel, we will omit $p_{\alpha}^{\ast}$.) 
\end{definition} 

\begin{proposition}\label{prop:reigen} 
For each $\gamma  \in \Gamma$, $r_{\alpha}(\gamma)$ is an eigenvalue of  
$\sigma_{\Gamma, u/{\Bbb F}_2, \alpha}(\gamma)$ on $\widetilde{R}_{\alpha}$. 
In other words, $r_{\alpha}(\gamma)$ is a root of the characteristic polynomial of  $\sigma_{\Gamma, u/{\Bbb F}_2, \alpha}(\gamma)$. 
\end{proposition} 

{\it Proof}.  
By using $X_{\alpha}^2 = \det \sigma_{\Gamma, u/{\Bbb F}_2, \alpha}(\alpha)$ 
and Lemma \ref{lemma:detformula}, we have 
\begin{eqnarray*} 
r_{\alpha}(\gamma)^2 & = & (a_{\alpha}(\gamma) + b_{\alpha}(\gamma) X_{\alpha})^2 \\ 
& = & a_{\alpha}(\gamma) ^2+ b_{\alpha}(\gamma) ^2X_{\alpha}^2 \\ 
& = & a_{\alpha}(\gamma) ^2+ b_{\alpha}(\gamma) ^2\det \sigma_{\Gamma, u/{\Bbb F}_2, \alpha}(\alpha) 
= \det \sigma_{\Gamma, u/{\Bbb F}_2, \alpha}(\gamma).  
\end{eqnarray*} 
Since ${\rm tr}(\sigma_{\Gamma, u/{\Bbb F}_2, \alpha}(\gamma)) = 0$, the characteristic polynomial of $\sigma_{\Gamma, u/{\Bbb F}_2, \alpha}(\gamma)$ is 
$x^2 - \det \sigma_{\Gamma, u/{\Bbb F}_2, \alpha}(\gamma)$. Hence $r_{\alpha}(\gamma)$ is an eigenvalue of  
$\sigma_{\Gamma, u/{\Bbb F}_2, \alpha}(\gamma)$. 
\qed 

\begin{proposition}\label{prop:chuf2}
For each $\alpha \in \Gamma$, 
$r_{\alpha} : \Gamma \to {\mathcal O}_{\widetilde{R}_{\alpha}}(\widetilde{R}_{\alpha})$ 
is a character.  
\end{proposition} 

{\it Proof}. 
Note that $a_{\alpha}(\gamma)$ and $b_{\alpha}(\gamma)$ are 
the $(a, b)$-coefficients  
of  $\sigma_{\Gamma, u/{\Bbb F}_2, \alpha}(\gamma)$ with respect to 
$\sigma_{\Gamma, u/{\Bbb F}_2, \alpha}(\alpha)$. By Lemma \ref{lemma:abproperty},   
$r_{\alpha}(e)=a_{\alpha}(e)+b_{\alpha}(e)X_{\alpha} = 1$ and 
\begin{multline*}
r_{\alpha}(\gamma)r_{\alpha}(\delta)  =  (a_{\alpha}(\gamma)+b_{\alpha}(\gamma)X_{\alpha})(a_{\alpha}(\delta)+b_{\alpha}(\delta)X_{\alpha}) \\ 
 =  (a_{\alpha}(\gamma)a_{\alpha}(\delta) + b_{\alpha}(\gamma)b_{\alpha}(\delta) X_{\alpha}^2) + (a_{\alpha}(\gamma)b_{\alpha}(\delta) + b_{\alpha}(\gamma)a_{\alpha}(\delta)) X_{\alpha} \\
= a_{\alpha}(\gamma\delta) + b_{\alpha}(\gamma\delta)X_{\alpha}  = 
r_{\alpha}(\gamma\delta) 
\end{multline*} 
for $\gamma, \delta \in \Gamma$.  Here we used $X_{\alpha}^2 = \det \sigma_{\Gamma, u/{\Bbb F}_2, \alpha}(\alpha)$. Hence $r_{\alpha}$ is a character 
of $\Gamma$. 
\qed 

\bigskip 

For $\alpha, \beta \in \Gamma$, let us 
consider the pull-backs of 
the intersection $R_{\alpha} \cap R_{\beta} \subseteq 
{\rm Rep}_2(\Gamma)_{u/{\Bbb F}_2}$ by 
$p_{\alpha} : \widetilde{R}_{\alpha} \to R_{\alpha}$ 
and $p_{\beta} : \widetilde{R}_{\beta} \to R_{\beta}$. 
Set $\widetilde{R}_{\alpha\beta} = p_{\alpha}^{-1}(R_{\alpha}\cap R_{\beta})$  
and $\widetilde{R}_{\beta\alpha} = p_{\beta}^{-1}(R_{\beta}\cap R_{\alpha})$.   
We define the morphism $\phi_{\beta\alpha} : \widetilde{R}_{\alpha\beta} 
\to \widetilde{R}_{\beta\alpha}$ over $R_{\alpha} \cap R_{\beta}$ by 
$X_{\beta} \mapsto a_{\alpha}(\beta) + b_{\alpha}(\beta) X_{\alpha}$. 
Since $(a_{\alpha}(\beta) + b_{\alpha}(\beta) X_{\alpha})^2 
= a_{\alpha}(\beta)^2 + b_{\alpha}(\beta)^2 X_{\alpha}^2 = 
\det \sigma_{\Gamma, u/{\Bbb F}_2, \alpha}(\beta)$ (see the proof of Proposition \ref{prop:reigen}), 
we can define the morphism $\phi_{\beta\alpha}$.

It is easy to check that $\phi_{\alpha\alpha} = 1$, 
$\phi_{\alpha\beta} \circ \phi_{\beta\alpha} = 1$ and 
$\phi_{\gamma\beta} \circ \phi_{\beta\alpha} = \phi_{\gamma\alpha}$ 
over $R_{\alpha} \cap R_{\beta} \cap R_{\gamma}$  
for $\alpha, \beta, \gamma \in \Gamma$.   
Gluing $\{ \widetilde{R}_{\alpha} \}_{\alpha \in \Gamma}$, we have 
a scheme $\widetilde{{\rm Rep}_2(\Gamma)}_{u/{\Bbb F}_2}$ 
over ${\rm Rep}_2(\Gamma)_{u/{\Bbb F}_2}$ by \cite[Chap. II, Ex. 2.12]{Hartshorne}.  
The canonical morphism $p :  \widetilde{{\rm Rep}_2(\Gamma)}_{u/{\Bbb F}_2} 
\to {\rm Rep}_2(\Gamma)_{u/{\Bbb F}_2}$ is a faithfully flat finite 
morphism. 

Let us define a ${\rm PGL}_2\otimes_{{\Bbb Z}} {\Bbb F}_2$-action on $\widetilde{{\rm Rep}_2(\Gamma)}_{u/{\Bbb F}_2}$. 
First, we define the action $\sigma_{\alpha}$ of 
${\rm PGL}_2\otimes_{{\Bbb Z}} {\Bbb F}_2$ on $\widetilde{R}_{\alpha}$ as follows: 
For a $Z$-valued point $((\rho, X), P)$ of $\widetilde{R}_{\alpha} 
\times {\rm PGL}_2\otimes_{{\Bbb Z}} {\Bbb F}_2$ with an 
${\Bbb F}_2$-scheme $Z$, 
we set $\sigma_{\alpha}((\rho, X), P) := (P^{-1}\rho P, X)$ as 
a $Z$-valued point of $\widetilde{R}_{\alpha}$. 
Since $X^2 = \det(\rho(\alpha)) = \det(P^{-1}\rho(\alpha)P)$, 
the morphism $\sigma_{\alpha} :  \widetilde{R}_{\alpha} 
\times {\rm PGL}_2\otimes_{{\Bbb Z}} {\Bbb F}_2 \to 
\widetilde{R}_{\alpha}$ can be defined. 
It is easy to see that $\sigma_{\alpha}$ gives a group action. 

Next, let us glue the actions $\{ \sigma_{\alpha} \}_{\alpha \in \Gamma}$ 
of ${\rm PGL}_2\otimes_{{\Bbb Z}} {\Bbb F}_2$. 
Recall that  $\phi_{\beta\alpha} : \widetilde{R}_{\alpha\beta} 
\to \widetilde{R}_{\beta\alpha}$ over $R_{\alpha} \cap R_{\beta}$ 
is given by $X_{\beta} \mapsto a_{\alpha}(\beta) + b_{\alpha}(\beta) X_{\alpha}$.  
The proof of Proposition~\ref{prop:abinv} shows that  
$a_{\alpha}(\beta)$ and 
$b_{\alpha}(\beta)$ are ${\rm PGL}_2\otimes_{{\Bbb Z}} {\Bbb F}_2$-invariant 
on $R_{\alpha}$.    
Thereby, the actions $\sigma_{\alpha}$ and $\sigma_{\beta}$ are 
compatible over  $R_{\alpha} \cap R_{\beta}$, and hence 
we obtain the action $\sigma$ of ${\rm PGL}_2\otimes_{{\Bbb Z}} {\Bbb F}_2$ 
on $\widetilde{{\rm Rep}_2(\Gamma)}_{u/{\Bbb F}_2}$ by gluing  
$\{ \sigma_{\alpha} \}_{\alpha \in \Gamma}$. 
Finally, remark that the canonical morphism 
$p :  \widetilde{{\rm Rep}_2(\Gamma)}_{u/{\Bbb F}_2} 
\to {\rm Rep}_2(\Gamma)_{u/{\Bbb F}_2}$
is ${\rm PGL}_2\otimes_{{\Bbb Z}} {\Bbb F}_2$-equivariant. 

\bigskip 
 
Let $A_1(\Gamma)$ be the coordinate ring of 
the affine scheme ${\rm Rep}_1(\Gamma)$. 
Set $A_1(\Gamma)_{{\Bbb F}_2} = 
A_1(\Gamma) \otimes_{{\Bbb Z}} {\Bbb F}_2$ 
and ${\rm Rep}_1(\Gamma)_{{\Bbb F}_2} 
= {\rm Spec} A_1(\Gamma)_{{\Bbb F}_2}$.  
Let $\chi_{\Gamma} : \Gamma \to A_1(\Gamma)_{{\Bbb F}_2}$ 
be the universal character of $\Gamma$ over ${\Bbb F}_2$. 

\begin{definition}\rm 
For a $A_1(\Gamma)_{{\Bbb F}_2}$-module $M$, 
we define 
\[
{\rm Der}(\Gamma, M) = 
\left\{ \delta : \Gamma \to M \:  
\begin{array}{|c} 
  \delta(\alpha\beta) = \chi_{\Gamma}(\alpha)\delta(\beta)
+ \delta(\alpha)\chi_{\Gamma}(\beta) \\ 
 \mbox{ for } \alpha, \beta \in \Gamma 
\end{array} 
\right\}.  
\]
\end{definition} 

We can prove the following lemma in the same way as Lemma \ref{lemma:omegamodule}. 

\begin{lemma}
  There exists a universal $A_1(\Gamma)_{{\Bbb F}_2}$-module  $\Omega_{\Gamma/{\Bbb F}_2}$ 
representing the covariant functor 
\[
\begin{array}{ccccc}
{\rm Der}(\Gamma, -) & : & (A_1(\Gamma)_{{\Bbb F}_2}\mbox{-}{\bf Mod}) & \to &  
(A_1(\Gamma)_{{\Bbb F}_2} \mbox{-}{\bf Mod}) \\ 
 & & M & \mapsto & {\rm Der}(\Gamma, M). 
\end{array} 
\]
In particular,  
\[
{\rm Der}(\Gamma, M) \stackrel{\sim}{\to} 
{\rm Hom}_{A_1(\Gamma)_{{\Bbb F}_2}}( \Omega_{\Gamma/{\Bbb F}_2}, M)
\]
is an isomorphism for each $A_{1}(\Gamma)_{{\Bbb F}_2}$-module $M$. 
\end{lemma}

\begin{remark}\label{remark:omegaisfgmoduletilde}\rm 
Let $d : \Gamma \to \Omega_{\Gamma/{\Bbb F}_2}$ be the 
universal derivation of $\Gamma$. We see that 
$\Omega_{\Gamma/{\Bbb F}_2}$ is generated by 
$\{ d\gamma \mid \gamma \in \Gamma \}$ as an 
$A_1(\Gamma)_{{\Bbb F}_2}$-module.   
As in Remark \ref{remark:fgu}, we see that 
if $\Gamma$ is finitely generated, then 
$\Omega_{\Gamma/{\Bbb F}_2}$ is a finitely generated $A_1(\Gamma)_{{\Bbb F}_2}$-module.  
\end{remark} 

\begin{definition}\label{def:chtilde}\rm 
We define the scheme $\widetilde{{\rm Ch}_2(\Gamma)}_{u/{\Bbb F}_2}$ over 
${\rm Rep}_1(\Gamma)_{{\Bbb F}_2}$ by 
\[
\widetilde{{\rm Ch}_2(\Gamma)}_{u/{\Bbb F}_2} = 
{\rm Proj} \: S({\Omega}_{\Gamma/{\Bbb F}_2}), 
\] 
where $S(\Omega_{\Gamma/{\Bbb F}_2})$ is the symmetric 
algebra of $\Omega_{\Gamma/{\Bbb F}_2}$ over $A_1(\Gamma)_{{\Bbb F}_2}$. 
\end{definition}

\begin{example}[{\it cf.} Example \ref{example:omegau1}]\rm
  Let ${\Upsilon}_1 = \langle \alpha_0 \rangle$ 
be the free monoid of rank $1$. 
The $A_{1}(\Upsilon_1)_{{\Bbb F}_2}$-module $\Omega_{\Upsilon_1/{\Bbb F}_2}$ is 
isomorphic to $A_{1}(\Upsilon_1)_{{\Bbb F}_2}$ by 
\[
\begin{array}{ccc}
  A_{1}(\Upsilon_1)_{{\Bbb F}_2} & \to & \Omega_{\Upsilon_1/{\Bbb F}_2} \\
  1  & \mapsto & d\alpha_{0} .  
\end{array}
\]
In particular, $\widetilde{{\rm Ch}_2(\Upsilon_1)}_{u/{\Bbb F}_2} 
\cong {\rm Rep}_1(\Upsilon_1)_{{\Bbb F}_2}$.  
\end{example}

\bigskip

Let $\psi : X \to {\rm Rep}_{1}(\Gamma)_{{\Bbb F}_2}$ 
be an ${\mathbb F}_2$-morphism. 
Let us regard $\Omega_{\Gamma/{\Bbb F}_2}$ as 
a quasi-coherent sheaf on ${\rm Rep}_1(\Gamma)_{{\Bbb F}_2}$.  
There exists a one-to-one correspondence  
\begin{multline*}
{\rm Hom}_{ \: {\rm Rep}_{1}(\Gamma)_{{\Bbb F}_2} \: }
(X, \widetilde{{\rm Ch}_2(\Gamma)}_{u/{\Bbb F}_2})  \cong \\ 
\{ 
\psi^{\ast}(\Omega_{\Gamma/{\Bbb F}_2}) \twoheadrightarrow {\mathcal L} \to 0  
\mid {\mathcal L} \mbox{ is a line bundle on } X  
\}/\sim.   
\end{multline*} 
Here we say that $\psi^{\ast}(\Omega_{\Gamma/{\Bbb F}_2}) \stackrel{f_1}{\twoheadrightarrow} {\mathcal L}_1$ and 
$\psi^{\ast}(\Omega_{\Gamma/{\Bbb F}_2}) 
\stackrel{f_2}{\twoheadrightarrow} {\mathcal L}_2$ are 
equivalent if there exists an isomorphism $g : {\mathcal L}_1 
\stackrel{\cong}{\to} {\mathcal L}_2$ such that $g \circ f_1 = f_2$. 

For $\alpha \in \Gamma$, we define the open subscheme 
$\widetilde{{\rm Ch}_2(\Gamma)}_{u/{\Bbb F}_2, \alpha}$ of 
$\widetilde{{\rm Ch}_2(\Gamma)}_{u/{\Bbb F}_2}$ by 
$\widetilde{{\rm Ch}_2(\Gamma)}_{u/{\Bbb F}_2, \alpha} := D(d\alpha) = 
\{ d\alpha  \neq 0 \}$. 
For an ${\mathbb F}_2$-morphism $\psi : X \to {\rm Rep}_{1}(\Gamma)_{{\Bbb F}_2}$, 
there exists a one-to-one correspondence  
\begin{multline*}
{\rm Hom}_{ \: {\rm Rep}_{1}(\Gamma)_{{\Bbb F}_2} \: }(X, 
\widetilde{{\rm Ch}_2(\Gamma)}_{u/{\Bbb F}_2, \alpha})  \cong \\ 
\left\{ 
\psi^{\ast}(\Omega_{\Gamma/{\Bbb F}_2}) \twoheadrightarrow {\mathcal L} \to 0  \; 
\begin{array}{|l}     
 {\mathcal L} \mbox{ is a line bundle on } X \mbox{ and }  \psi^{\ast}(d\alpha) \mbox{ is } \\ 
 \mbox{ nowhere vanishing as a section of } {\mathcal L}   
\end{array}
\right\}\bigg/\sim.   
\end{multline*} 
When ${\mathcal L}$ is generated by $\psi^{\ast}(d\alpha)$, 
${\mathcal L}$ is isomorphic to ${\mathcal O}_X$.  
Let $r : \Gamma \to \Gamma(X, {\mathcal O}_X)$ be the character associated to 
$\psi : X \to {\rm Rep}_{1}(\Gamma)_{{\Bbb F}_2}$. 
Regarding $\psi^{\ast}(d\alpha)$ as $1$ of ${\mathcal O}_X$, 
we have the following: 
\begin{multline*}
{\rm Hom}_{ \: {\rm Rep}_{1}(\Gamma)_{{\Bbb F}_2} \: }(X, 
\widetilde{{\rm Ch}_2(\Gamma)}_{u/{\Bbb F}_2, \alpha})  \cong \\ 
\left\{ 
\; d \;\; 
\begin{array}{|l}     
 d  : \Gamma \to \Gamma(X, {\mathcal O}_X) \mbox{ is a derivation with respect to } r  \\ 
\mbox{ such that  } d(\alpha)=1   
\end{array}
\right\},    
\end{multline*} 
where we say that $d : \Gamma \to \Gamma(X, {\mathcal O}_X)$ 
is a {\it derivation with respect to $r$} if 
$d(\gamma\delta) = r(\gamma)d(\delta) + d(\gamma)r(\delta)$ holds for
each $\gamma, \delta \in \Gamma$. 

\bigskip 
 
We construct a morphism $\lambda :  
\widetilde{{\rm Rep}_2(\Gamma)}_{u/{\Bbb F}_2}  
\to {\rm Rep}_1(\Gamma)_{{\Bbb F}_2}$. 
By Proposition \ref{prop:chuf2}, $r_{\alpha} : \Gamma \to {\mathcal O}_{\widetilde{R}_{\alpha}}(\widetilde{R}_{\alpha})$ is a character for each 
$\alpha \in \Gamma$. 
It gives us 
a morphism $\lambda_{\alpha} : \widetilde{R}_{\alpha} \to {\rm Rep}_1(\Gamma)_{{\Bbb F}_2}$. 
Through the isomorphism $\phi_{\beta\alpha} : 
\widetilde{R}_{\alpha\beta} \to 
\widetilde{R}_{\beta\alpha}$,  $\lambda_{\alpha}$ and $\lambda_{\beta}$ coincide 
on $\widetilde{R}_{\alpha} \cap \widetilde{R}_{\beta}$ 
for each $\alpha, \beta \in \Gamma$. 
Indeed, $\phi_{\beta\alpha}$ is given by $X_{\beta} 
\mapsto a_{\alpha}(\beta) + b_{\alpha}(\beta) X_{\alpha}$. 
By comparing 
\begin{multline*} 
\sigma_{\Gamma, u/{\Bbb F}_2, \beta}(\gamma) = 
a_{\beta}(\gamma)I_2 + 
b_{\beta}(\gamma)\sigma_{\Gamma, u/{\Bbb F}_2, \beta}(\beta) \\  
= a_{\beta}(\gamma)I_2 + b_{\beta}(\gamma)(a_{\alpha}(\beta)I_2 + b_{\alpha}(\beta)\sigma_{\Gamma, u/{\Bbb F}_2, \alpha}(\alpha))  \\ 
= (a_{\beta}(\gamma) + b_{\beta}(\gamma)a_{\alpha}(\beta))I_2 + b_{\beta}(\gamma)b_{\alpha}(\beta)\sigma_{\Gamma, u/{\Bbb F}_2, \alpha}(\alpha) 
\end{multline*} 
with  
\[
\sigma_{\Gamma, u/{\Bbb F}_2, \alpha}(\gamma) =  
a_{\alpha}(\gamma)I_2 + b_{\alpha}(\gamma)
\sigma_{\Gamma, u/{\Bbb F}_2, \alpha}(\alpha),  
\] 
we have $a_{\alpha}(\gamma) = a_{\beta}(\gamma)+b_{\beta}(\gamma)a_{\alpha}(\beta)$ 
and $b_{\alpha}(\gamma) = b_{\beta}(\gamma)b_{\alpha}(\beta)$ on 
$R_{\alpha} \cap R_{\beta}$ for each $\gamma \in \Gamma$. 
The isomorphism $\phi_{\beta\alpha}$ induces   
$\lambda_{\beta}(\gamma) = 
a_{\beta}(\gamma) + b_{\beta}(\gamma)X_{\beta} 
\mapsto a_{\beta}(\gamma) + 
b_{\beta}(\gamma)(a_{\alpha}(\beta) + b_{\alpha}(\beta) X_{\alpha}) 
= a_{\alpha}(\gamma) + b_{\alpha}(\gamma)X_{\alpha} 
=\lambda_{\alpha}(\gamma)$. 
Hence $\lambda_{\alpha}$ and $\lambda_{\beta}$ coincide 
on $\widetilde{R}_{\alpha} \cap \widetilde{R}_{\beta}$. 
By gluing $\{ \lambda_{\alpha} \}_{\alpha \in \Gamma}$, 
we obtain a morphism $\lambda : \widetilde{{\rm Rep}_2(\Gamma)}_{u/{\Bbb F}_2}    
\to {\rm Rep}_1(\Gamma)_{{\Bbb F}_2}$.  
We regard $\widetilde{{\rm Rep}_2(\Gamma)}_{u/{\Bbb F}_2}$ as 
a ${\rm Rep}_1(\Gamma)_{{\Bbb F}_2}$-scheme by $\lambda$. 

We construct a ${\rm Rep}_1(\Gamma)_{{\Bbb F}_2}$-morphism $\widetilde{\pi}_{\Gamma, u/{\Bbb F}_2} : 
\widetilde{{\rm Rep}_2(\Gamma)}_{u/{\Bbb F}_2} \to 
\widetilde{{\rm Ch}_2(\Gamma)}_{u/{\Bbb F}_2}$. 
First, let us define $\widetilde{\pi}_{\Gamma, u/{\Bbb F}_2, \alpha} : 
\widetilde{R}_{\alpha} = \widetilde{{\rm Rep}_2(\Gamma)}_{u/{\Bbb F}_2, \alpha} \to \widetilde{{\rm Ch}_2(\Gamma)}_{u/{\Bbb F}_2, \alpha}$ 
for each $\alpha \in \Gamma$. 
Put $\widetilde{C}_{\alpha} = \widetilde{{\rm Ch}_2(\Gamma)}_{u/{\Bbb F}_2, \alpha}$ 
and $\widetilde{\pi}_{\alpha} := \widetilde{\pi}_{\Gamma, u/{\Bbb F}_2, \alpha}$ 
for simplicity. 
Set $\eta_{\alpha}(\gamma) := \sigma_{\Gamma, u/{\Bbb F}_2, \alpha}(\gamma) - 
r_{\alpha}(\gamma) I_2 \in {\rm M}_2({\mathcal O}_{\widetilde{R}_{\alpha}}({\mathcal O}_{\widetilde{R}_{\alpha}}))$ for $\gamma \in \Gamma$. 
By Proposition \ref{prop:reigen}, 
$\eta_{\alpha}(\gamma)^2 = \sigma_{\Gamma, u/{\Bbb F}_2, \alpha}(\gamma)^2 + 
r_{\alpha}(\gamma)^2 I_2 = \sigma_{\Gamma, u/{\Bbb F}_2, \alpha}(\gamma)^2 + 
\det(\sigma_{\Gamma, u/{\Bbb F}_2, \alpha}(\gamma)) I_2 = 0$. 
Note that 
${\mathcal O}_{\widetilde{R}_{\alpha}}[\sigma_{\Gamma, u/{\Bbb F}_2, \alpha}(\Gamma)] = {\mathcal O}_{\widetilde{R}_{\alpha}}\cdot I_2 + {\mathcal O}_{\widetilde{R}_{\alpha}} \cdot \sigma_{\Gamma, u/{\Bbb F}_2, \alpha}(\alpha) = {\mathcal O}_{\widetilde{R}_{\alpha}}\cdot I_2 + {\mathcal O}_{\widetilde{R}_{\alpha}} \cdot \eta_{\alpha}(\alpha)$. 
For each $\gamma \in \Gamma$,  
\begin{eqnarray*} 
\eta_{\alpha}(\gamma) & = &  
(a_{\alpha}(\gamma) - r_{\alpha}(\gamma)) I_2 + b_{\alpha}(\gamma) 
\sigma_{\Gamma, u/{\Bbb F}_2, \alpha}(\alpha) \\ 
& = & - b_{\alpha}(\gamma) X_{\alpha} I_2 + b_{\alpha}(\gamma) 
\sigma_{\Gamma, u/{\Bbb F}_2, \alpha}(\alpha) \\ 
& = & b_{\alpha}(\gamma) \eta_{\alpha}(\alpha), 
\end{eqnarray*} 
since $r_{\alpha}(\alpha) = X_{\alpha}$. 

\begin{proposition} 
For each $\alpha \in \Gamma$, 
$b_{\alpha}(\cdot)  :  \Gamma \to 
{\mathcal O}_{\widetilde{R}_{\alpha}}(\widetilde{R}_{\alpha})$ 
is a derivation with respect to $r_{\alpha}$. 
\end{proposition} 

{\it Proof}. 
By calculating $\sigma_{\Gamma, u/{\Bbb F}_2, \alpha}(\gamma\delta)$ and 
$\sigma_{\Gamma, u/{\Bbb F}_2, \alpha}(\gamma)
\sigma_{\Gamma, u/{\Bbb F}_2, \alpha}(\delta)$, 
we have $b_{\alpha}(\gamma\delta) = a_{\alpha}(\gamma)b_{\alpha}(\delta) 
+ b_{\alpha}(\gamma)a_{\alpha}(\delta)$ for each $\gamma, \delta \in \Gamma$. 
It follows that 
$b_{\alpha}(\gamma\delta) = r_{\alpha}(\gamma)b_{\alpha}(\delta) 
+ b_{\alpha}(\gamma)r_{\alpha}(\delta)$. 
\qed 

\bigskip 

Hence we have a morphism $\widetilde{\pi}_{\alpha} : \widetilde{R}_{\alpha} 
\to \widetilde{C}_{\alpha}$ by 
the derivation $b_{\alpha}(\cdot)  :  \Gamma \to 
{\mathcal O}_{\widetilde{R}_{\alpha}}(\widetilde{R}_{\alpha})$ 
with $b_{\alpha}(\alpha)=1$. 

Next, let us glue the morphisms $\{ \widetilde{\pi}_{\alpha} : \widetilde{R}_{\alpha} 
\to \widetilde{C}_{\alpha} \}_{\alpha \in \Gamma}$.  
Because $b_{\alpha}(\gamma) = b_{\beta}(\gamma)b_{\alpha}(\beta)$ on 
$R_{\alpha} \cap R_{\beta}$ 
for each $\gamma \in \Gamma$ 
and $b_{\alpha}(\beta) \in ({\mathcal O}_{{R}_{\alpha} 
\cap {R}_{\beta}})^{\times}$, we have the following 
commutative diagram: 
\[
\begin{array}{ccccc}
\lambda_{\alpha}^{\ast} (\Omega_{u/{\Bbb F}_2})\mid_{\widetilde{R}_{\alpha}\cap \widetilde{R}_{\beta}} & \twoheadrightarrow 
& {\mathcal O}_{\widetilde{R}_{\alpha}\cap \widetilde{R}_{\beta}} & \to 0  \\ 
|| & & \uparrow_{b_{\alpha}(\beta)\cdot} & & \\
\lambda_{\beta}^{\ast} (\Omega_{u/{\Bbb F}_2})\mid_{\widetilde{R}_{\alpha}\cap \widetilde{R}_{\beta}} & \twoheadrightarrow 
& {\mathcal O}_{\widetilde{R}_{\alpha}\cap \widetilde{R}_{\beta}} & \to 0 ,  
\end{array}
\]
where the morphism $b_{\alpha}(\beta)\cdot : {\mathcal O}_{\widetilde{R}_{\alpha}\cap \widetilde{R}_{\beta}} \to {\mathcal O}_{\widetilde{R}_{\alpha}\cap \widetilde{R}_{\beta}}$ 
defined by $\varphi \mapsto b_{\alpha}(\beta) \varphi$ 
is an isomorphism. 
It follows that $\widetilde{\pi}_{\alpha}\mid_{\widetilde{R}_{\alpha}\cap \widetilde{R}_{\beta}} = \widetilde{\pi}_{\beta}\mid_{\widetilde{R}_{\alpha}\cap \widetilde{R}_{\beta}}$ for each $\alpha, \beta \in \Gamma$. 
Therefore we have a ${\rm Rep}_1(\Gamma)_{{\Bbb F}_2}$-morphism 
$\widetilde{\pi}_{\Gamma, u/{\Bbb F}_2} : 
\widetilde{{\rm Rep}_2(\Gamma)}_{u/{\Bbb F}_2} \to 
\widetilde{{\rm Ch}_2(\Gamma)}_{u/{\Bbb F}_2}$. 

Let ${\rm PGL}_2\otimes_{{\Bbb Z}} {\Bbb F}_2$ act on 
$\widetilde{{\rm Ch}_2(\Gamma)}_{u/{\Bbb F}_2}$ trivially. 
Then we have the following proposition: 

\begin{proposition}
The morphism 
$\widetilde{\pi}_{\Gamma, u/{\Bbb F}_2} : 
\widetilde{{\rm Rep}_2(\Gamma)}_{u/{\Bbb F}_2} \to 
\widetilde{{\rm Ch}_2(\Gamma)}_{u/{\Bbb F}_2}$ is 
${\rm PGL}_2\otimes_{{\Bbb Z}} {\Bbb F}_2$-equivariant. 
\end{proposition} 

{\it Proof}. 
It suffices to show that 
$\widetilde{\pi}_{\alpha} : 
\widetilde{R}_{\alpha} \to 
\widetilde{C}_{\alpha}$ is 
${\rm PGL}_2\otimes_{{\Bbb Z}} {\Bbb F}_2$-equivariant 
for each $\alpha \in \Gamma$. 
The character $r_{\alpha}$ on $\widetilde{R}_{\alpha}$ is 
given by $r_{\alpha}(\gamma) = 
a_{\alpha}(\gamma)+b_{\alpha}(\gamma)X_{\alpha}$ for 
$\gamma \in \Gamma$.  
The proof of Proposition~\ref{prop:abinv} shows that 
$a_{\alpha}(\gamma)$ and $b_{\alpha}(\gamma)$ are 
${\rm PGL}_2\otimes_{{\Bbb Z}} {\Bbb F}_2$-invariant on 
$R_{\alpha}$. 
From the definition of the action $\sigma_{\alpha}$ on 
$\widetilde{R}_{\alpha}$, $X_{\alpha}$ is 
${\rm PGL}_2\otimes_{{\Bbb Z}} {\Bbb F}_2$-invariant.    
Hence $r_{\alpha}(\gamma)$ is also 
${\rm PGL}_2\otimes_{{\Bbb Z}} {\Bbb F}_2$-invariant for each 
$\gamma \in \Gamma$.  
The morphism $\widetilde{\pi}_{\alpha}$ 
is given by the derivation $b_{\alpha}(\cdot)$ with respect to $r_{\alpha}$. 
Since $b_{\alpha}(\gamma)$ is ${\rm PGL}_2\otimes_{{\Bbb Z}} {\Bbb F}_2$-invariant, 
the morphism $\widetilde{\pi}_{\alpha}$ is 
${\rm PGL}_2\otimes_{{\Bbb Z}} {\Bbb F}_2$-equivariant. 
This completes the proof. 
\qed 

\bigskip 

Let us define a morphism $q_{\alpha} : 
\widetilde{{\rm Ch}_2(\Gamma)}_{u/{\Bbb F}_2, \alpha} 
\to {\rm Ch}_2(\Gamma)_{u/{\Bbb F}_2, \alpha}$ for each $\alpha \in \Gamma$. 
By the definitions, $\widetilde{{\rm Ch}_2(\Gamma)}_{u/{\Bbb F}_2, \alpha}$ 
and ${\rm Ch}_2(\Gamma)_{u/{\Bbb F}_2, \alpha}$ are  
${\rm Rep}_1(\Gamma)_{{\Bbb F}_2}$-schemes.  
Let $r$ and $d$ be the universal characters on 
$\widetilde{{\rm Ch}_2(\Gamma)}_{u/{\Bbb F}_2, \alpha}$ 
and ${\rm Ch}_2(\Gamma)_{u/{\Bbb F}_2, \alpha}$, respectively. 
Consider the character $r^2$ on 
$\widetilde{{\rm Ch}_2(\Gamma)}_{u/{\Bbb F}_2, \alpha}$ instead of $r$. 
For constructing $q_{\alpha}$, it suffices to define 
$(a, b)$-coefficients with respect to $(r^2, \alpha)$ 
on $\widetilde{{\rm Ch}_2(\Gamma)}_{u/{\Bbb F}_2, \alpha}$. 
Denote by $\delta$ the universal derivation with respect to 
$r$ on $\widetilde{{\rm Ch}_2(\Gamma)}_{u/{\Bbb F}_2, \alpha}$ such that 
$\delta(\alpha)=1$. 
Then we define $a(\gamma) = r(\gamma) - r(\alpha)\delta(\gamma)$ and 
$b(\gamma) =\delta(\gamma)$ for $\gamma \in \Gamma$. 

\begin{proposition}\label{prop:abab}  
Let $a$ and $b$ be as above. 
Then $a$ and $b$ are $(a, b)$-coefficients with respect to $(r^2, \alpha)$ 
on $\widetilde{{\rm Ch}_2(\Gamma)}_{u/{\Bbb F}_2, \alpha}$.   
\end{proposition} 

{\it Proof}. 
It is easy to check that 
$a(e) = 1, b(e) = 0$ and that 
$a(\alpha) = 0, b(\alpha)=1$. 
By direct calculations, we have 
\begin{eqnarray*} 
a(\gamma_1\gamma_2) & = & r(\gamma_1\gamma_2)-r(\alpha)\delta(\gamma_1\gamma_2) \\  
 & = & r(\gamma_1)r(\gamma_2) - r(\alpha)r(\gamma_1)\delta(\gamma_2) 
-r(\alpha)\delta(\gamma_1)r(\gamma_2) \\
 & = & (r(\gamma_1)-r(\alpha)\delta(\gamma_1))(r(\gamma_2)-r(\alpha)\delta(\gamma_2)) + 
r(\alpha)^2\delta(\gamma_1)\delta(\gamma_2) \\ 
 & = & a(\gamma_1)a(\gamma_2) + r^2(\alpha) b(\gamma_1)b(\gamma_2),  
\end{eqnarray*} 
\begin{eqnarray*} 
b(\gamma_1\gamma_2) & = & \delta(\gamma_1\gamma_2) \\ 
 & = & r(\gamma_1)\delta(\gamma_2)+\delta(\gamma_1)r(\gamma_2) \\ 
 & = & (r(\gamma_1) - r(\alpha)\delta(\gamma_1))\delta(\gamma_2) 
+ \delta(\gamma_1)(r(\gamma_2)-r(\alpha)\delta(\gamma_2)) \\ 
 & = & a(\gamma_1)b(\gamma_2)+b(\gamma_1)a(\gamma_2),  
\end{eqnarray*} 
and 
\begin{eqnarray*} 
a(\gamma)^2 + b(\gamma)^2 r^2(\alpha) & = & 
(a(\gamma)+b(\gamma)r(\alpha))^2 \\ 
 & = & (r(\gamma) - r(\alpha)\delta(\gamma) + \delta(\gamma)r(\alpha))^2    \\ 
 & = & r^2(\gamma). 
\end{eqnarray*} 
Hence we have proved the statement. 
\qed 

\bigskip 

By the $(a, b)$-coefficients on $\widetilde{{\rm Ch}_2(\Gamma)}_{u/{\Bbb F}_2, \alpha}$, 
we obtain a morphism $q_{\alpha}  : 
\widetilde{{\rm Ch}_2(\Gamma)}_{u/{\Bbb F}_2, \alpha} 
\to {\rm Ch}_2(\Gamma)_{u/{\Bbb F}_2, \alpha}$. 
Denote by $r^2 : \widetilde{{\rm Ch}_2(\Gamma)}_{u/{\Bbb F}_2, \alpha} 
\to {\rm Rep}_1(\Gamma)_{{\Bbb F}_2}$ and 
$d : {\rm Ch}_2(\Gamma)_{u/{\Bbb F}_2, \alpha} 
\to {\rm Rep}_1(\Gamma)_{{\Bbb F}_2}$ 
the morphisms induced by the characters 
$r^2$ and $d$, respectively. 
Then $d \circ q_{\alpha} = r^2$. 

Thus, we have the following commutative diagram for each $\alpha \in \Gamma$: 
\[
\begin{array}{ccc}
\widetilde{{\rm Rep}_2(\Gamma)}_{u/{\Bbb F}_2, \alpha} 
& \stackrel{p_{\alpha}}{\to} & {\rm Rep}_2(\Gamma)_{u/{\Bbb F}_2, \alpha} \\
\widetilde{\pi}_{\alpha} \downarrow & & \downarrow \pi_{\alpha} \\
\widetilde{{\rm Ch}_2(\Gamma)}_{u/{\Bbb F}_2, \alpha} 
& \stackrel{q_{\alpha}}{\to} & {\rm Ch}_2(\Gamma)_{u/{\Bbb F}_2, \alpha},   
\end{array}  
\]  
where we denote by $\pi_{\alpha}$ the morphism $\pi_{\Gamma, u/{\Bbb F}_2, \alpha} : 
{\rm Rep}_2(\Gamma)_{u/{\Bbb F}_2, \alpha} \to 
{\rm Ch}_2(\Gamma)_{u/{\Bbb F}_2, \alpha}$ which was defined in \S 6.  
These morphisms are ${\rm PGL}_2 \otimes_{{\Bbb Z}} {\Bbb F}_2$-equivariant. 
Let $\lambda_{\alpha}^2 : \widetilde{{\rm Rep}_{2}(\Gamma)}_{u/{\Bbb F}_2, \alpha}  
\to {\rm Rep}_1(\Gamma)_{{\Bbb F}_2}$ and 
$\det : {\rm Rep}_{2}(\Gamma)_{u/{\Bbb F}_2, \alpha}  
\to {\rm Rep}_1(\Gamma)_{{\Bbb F}_2}$ 
be the morphisms corresponding to the characters 
$r_{\alpha}^2$ and $\det(\sigma_{\Gamma, u/{\Bbb F}_2, \alpha}(\cdot))$, respectively. 
Then $\det \circ \: p_{\alpha} = \lambda_{\alpha}^2$. 

\begin{proposition}\label{prop:uf2fibre}
The commutative diagram 
\[
\begin{array}{ccc}
\widetilde{{\rm Rep}_2(\Gamma)}_{u/{\Bbb F}_2, \alpha} 
& \stackrel{p_{\alpha}}{\to} & {\rm Rep}_2(\Gamma)_{u/{\Bbb F}_2, \alpha} \\
\widetilde{\pi}_{\alpha} \downarrow & & \downarrow \pi_{\alpha} \\
\widetilde{{\rm Ch}_2(\Gamma)}_{u/{\Bbb F}_2, \alpha} 
& \stackrel{q_{\alpha}}{\to} & {\rm Ch}_2(\Gamma)_{u/{\Bbb F}_2, \alpha}.  
\end{array}  
\]  gives  
a fibre product for each $\alpha \in \Gamma$. 
\end{proposition} 

{\it Proof}. 
Put $R_{\alpha} := {\rm Rep}_2(\Gamma)_{u/{\Bbb F}_2, \alpha}$ 
and $C_{\alpha} := {\rm Ch}_2(\Gamma)_{u/{\Bbb F}_2, \alpha}$.  
We show that $(p_{\alpha}, \widetilde{\pi}_{\alpha}) : \widetilde{R}_{\alpha} \to R_{\alpha}\times_{C_{\alpha}} \widetilde{C}_{\alpha}$ is an isomorphism.  
It suffices to prove that $(p_{\alpha}, \widetilde{\pi}_{\alpha})$ 
induces a bijective map between the sets of $Z$-valued points 
for any ${\Bbb F}_2$-scheme $Z$. 
Let $(\rho_1, X_1)$ and $(\rho_2, X_2)$ be $Z$-valued points 
of $\widetilde{R}_{\alpha}$ such that the images by  
$(p_{\alpha}, \widetilde{\pi}_{\alpha})$ coincide. 
Obviously, $\rho_1 = \rho_2$. 
By the assumption, $(\rho_1, X_1)$ and $(\rho_2, X_2)$ induce 
the same character $r$ on $\widetilde{R}_{\alpha}$. 
Since $X_1 = r(\alpha) = X_2$ on $Z$, 
$(\rho_1, X_1) = (\rho_2, X_2)$. Hence 
we have proved the injectivity. 

Let $(\rho, (r, \delta))$ be a $Z$-valued point of 
$R_{\alpha}\times_{C_{\alpha}} \widetilde{C}_{\alpha}$, 
where $r : \Gamma \to {\mathcal O}_Z(Z)$ is a character and 
$\delta : \Gamma \to {\mathcal O}_Z(Z)$ is a derivation 
with respect to $r$ such that $\delta(\alpha)=1$.  
The $Z$-valued point $(a, b)$ of 
$C_{\alpha}$ induced by $(r, \delta)$ is 
given by 
$a(\gamma) = r(\gamma) - r(\alpha)\delta(\gamma)$ and 
$b(\gamma) = \delta(\gamma)$ for each $\gamma \in 
\Gamma$. 
Note that $a$ and $b$ are $(a, b)$-coefficients with respect to 
$(r^2, \alpha)$ by Proposition~\ref{prop:abab}. 
Because $\rho$ and $(r, \delta)$ induce the same 
$Z$-valued point of $C_{\alpha}$, 
$r^2(\gamma) = \det \rho(\gamma)$ and 
$\rho(\gamma) = a(\gamma)I_2 + b(\gamma) \rho(\alpha) 
= (r(\gamma) - r(\alpha)\delta(\gamma))I_2 + \delta(\gamma)\rho(\alpha)$ 
for each $\gamma \in \Gamma$. 
Set $X := r(\alpha)$. Then $(\rho, X)$ is  a $Z$-valued point of 
$\widetilde{R}_{\alpha}$ by $X^2 = r^2(\alpha) = \det \rho(\alpha)$. 
It is obvious that $p_{\alpha} (\rho, X) = \rho$. 
Denote by $(r', \delta')$ the image of 
$(\rho, X)$ by $\widetilde{\pi}_{\alpha}$.  
For each $\gamma \in \Gamma$, 
$\delta'(\gamma) = b(\gamma) = \delta(\gamma)$ and 
\begin{eqnarray*}
r'(\gamma) & = & a(\gamma) + b(\gamma)X \\
 & = & (r(\gamma) - r(\alpha)\delta(\gamma)) + \delta(\gamma) X \\
 & = & r(\gamma). 
\end{eqnarray*} 
Thus $\widetilde{\pi}_{\alpha}(\rho, X) = (r, \delta)$.  
Hence $(p_{\alpha}, \widetilde{\pi}_{\alpha})(\rho, X) = (\rho, (r, \delta))$, 
which implies the surjectivity.  
Therefore we have proved the statement. 
\qed 

\begin{definition}\rm 
Let $\Upsilon_1 = \langle \alpha_0 \rangle$ be the free monoid 
of rank $1$. 
As in Definition \ref{def:prototypeuf2}, we call the morphism  
$\widetilde{\pi}_{\Upsilon_1, u/{\Bbb F}_2, \alpha_0} : 
\widetilde{{\rm Rep}_{2}(\Upsilon_1)}_{u/{\Bbb F}_2, \alpha_0} \to 
\widetilde{{\rm Ch}_2(\Upsilon_1)}_{u/{\Bbb F}_2, \alpha_0}$ 
the {\it prototype} in the unipotent mold over ${\Bbb F}_2$ case.     
Remark that $\widetilde{{\rm Rep}_{2}(\Upsilon_1)}_{u/{\Bbb F}_2} = 
\widetilde{{\rm Rep}_{2}(\Upsilon_1)}_{u/{\Bbb F}_2, \alpha_0}$ 
and that $\widetilde{{\rm Ch}_2(\Upsilon_1)}_{u/{\Bbb F}_2} 
= \widetilde{{\rm Ch}_2(\Upsilon_1)}_{u/{\Bbb F}_2, \alpha_0}$.  
\end{definition} 

\begin{remark}\label{remark:prototypeuf2tilde}\rm 
Recall that $\pi_{\Upsilon_1, u/{\Bbb F}_2, \alpha_0} : {\rm Rep}_{2}(\Upsilon_1)_{u/{\Bbb F}_2, \alpha_0} \to 
{\rm Ch}_2(\Upsilon_1)_{u/{\Bbb F}_2, \alpha_0}$ is 
a universal geometric quotient by ${\rm PGL}_2\otimes_{{\Bbb Z}} {\Bbb F}_2$ 
(Theorem \ref{th:prototypeuf2geom}). 
The prototype $\pi_{\Upsilon_1, u/{\Bbb F}_2, \alpha_0}$ is described by 
${\rm Spec} ({\Bbb F}_2[a, b, c, d]/(a+d))  \supset D(b) \cup D(c)  
\to {\rm Spec} {\Bbb F}_2[D]$, where 
$D$ is mapped to $ad-bc$.  
Then $\widetilde{{\rm Rep}_{2}(\Upsilon_1)}_{u/{\Bbb F}_2, \alpha_0}$ 
is isomorphic to 
$D(b) \cup D(c) \subset {\rm Spec} ({\Bbb F}_2[a, b, c, d][X]/(a+d, X^2 -ad+bc))$.  
Since $\Omega_{\Upsilon_1, u/{\Bbb F}_2}$ is isomorphic to 
$A_1(\Upsilon_1)_{{\Bbb F}_2}$, $\widetilde{{\rm Ch}_2(\Upsilon_1)}_{u/{\Bbb F}_2, \alpha_0}$ is isomorphic to ${\rm Rep}_1(\Upsilon_{1})_{{\Bbb F}_2} 
\cong {\rm Spec} {\Bbb F}_2[\chi]$. 
Here the indeterminate $\chi$ corresponds to the value at $\alpha_0$ of 
the universal character on ${\rm Rep}_1(\Upsilon_{1})_{{\Bbb F}_2}$.  
Therefore $\widetilde{\pi}_{\Upsilon_1, u/{\Bbb F}_2, \alpha_0} : 
\widetilde{{\rm Rep}_{2}(\Upsilon_1)}_{u/{\Bbb F}_2, \alpha_0} \to 
\widetilde{{\rm Ch}_2(\Upsilon_1)}_{u/{\Bbb F}_2, \alpha_0}$ 
is described by ${\rm Spec} ({\Bbb F}_2[a, b, c, d][X]/(a+d, X^2 -ad+bc)) 
\supset D(b) \cup D(c)  \to {\rm Spec} {\Bbb F}_2[\chi]$, where 
$\chi$ is mapped to $X$. 
\end{remark} 

By Proposition \ref{prop:uf2fibre}, the prototype 
$\widetilde{\pi}_{\Upsilon_1, u/{\Bbb F}_2, \alpha_0} : 
\widetilde{{\rm Rep}_{2}(\Upsilon_1)}_{u/{\Bbb F}_2, \alpha_0} \to 
\widetilde{{\rm Ch}_2(\Upsilon_1)}_{u/{\Bbb F}_2, \alpha_0}$ 
is obtained by base 
change of $\pi_{\Upsilon_1, u/{\Bbb F}_2, \alpha_0} : 
{\rm Rep}_{2}(\Upsilon_1)_{u/{\Bbb F}_2, \alpha_0} 
\to {\rm Ch}_2(\Upsilon_1)_{u/{\Bbb F}_2, \alpha_0}$. 
Theorem \ref{th:prototypeuf2geom} implies the following:  

\begin{theorem}\label{th:prototypeuf2tilde}  
The prototype 
\[ 
\widetilde{\pi}_{\Upsilon_1, u/{\Bbb F}_2, \alpha_0} : 
\widetilde{{\rm Rep}_{2}(\Upsilon_1)}_{u/{\Bbb F}_2, \alpha_0} \to 
\widetilde{{\rm Ch}_2(\Upsilon_1)}_{u/{\Bbb F}_2, \alpha_0}
\] 
is a universal geometric quotient by 
${\rm PGL}_2\otimes_{{\Bbb Z}} {\Bbb F}_2$. 
\end{theorem} 

\bigskip

Let $\Gamma$ be a group or a monoid.  
For $\alpha \in \Gamma$, we define the monoid homomorphism 
$\phi : \Upsilon_1 =\langle \alpha_0 \rangle \to \Gamma$ by 
$\alpha_0 \mapsto \alpha$.  
By restricting representations, characters, and derivations of $\Gamma$ to 
those of $\Upsilon_1$ through $\phi$,
we can obtain the following 
commutative diagram:  
\[
\begin{array}{ccc}
  \widetilde{{\rm Rep}_2(\Gamma)}_{u/{\Bbb F}_2, \alpha} & \to & \widetilde{{\rm Ch}_2(\Gamma)}_{u/{\Bbb F}_2, \alpha} \\
  \downarrow & &\downarrow \\
  \widetilde{{\rm Rep}_2(\Upsilon_1)}_{u/{\Bbb F}_2, \alpha_0} & \to &  
\widetilde{{\rm Ch}_2(\Upsilon_1)}_{u/{\Bbb F}_2, \alpha_0}. 
\end{array}
\]

Under this situation, we have the following lemma. 
\begin{lemma}\label{lemma:fibreproductuf2-2}
  The above diagram gives a fibre product. 
In particular, the morphism 
$\widetilde{{\rm Rep}_2(\Gamma)}_{u/{\Bbb F}_2, \alpha}  \to  
\widetilde{{\rm Ch}_2(\Gamma)}_{u/{\Bbb F}_2, \alpha}$ 
is obtained by base change of the prototype. 
\end{lemma}

{\it Proof}.
Here we prove the statement without using Lemma \ref{lemma:fibreproductuf2}. 
Put $\widetilde{R}_{\alpha_0} := \widetilde{{\rm Rep}_2(\Upsilon_1)}_{u/{\Bbb F}_2, \alpha_0}$ and $\widetilde{C}_{\alpha_0} :=  \widetilde{{\rm Ch}_2(\Upsilon_1)}_{u/{\Bbb F}_2, \alpha_0}$. It suffices to prove that 
$\widetilde{R}_{\alpha} \to \widetilde{R}_{\alpha_0} \times_{\widetilde{C}_{\alpha_0}} 
\widetilde{C}_{\alpha}$ 
induces a bijective map between the sets of $Z$-valued points 
for any ${\Bbb F}_2$-scheme $Z$. 
Let $(\rho_1, X_1)$ and $(\rho_2, X_2)$ be $Z$-valued points of 
$\widetilde{R}_{\alpha}$ whose images coincide. 
By the assumption, $X_1 = X_2$ and $\rho_1(\alpha) = \rho_1(\phi(\alpha_0)) 
= \rho_2(\phi(\alpha_0)) = \rho_2(\alpha)$. 
Since $(\rho_1, X_1)$ and $(\rho_2, X_2)$ induce the same $Z$-valued point 
$(r, \delta)$ of $\widetilde{C}_{\alpha}$, $a(\gamma) = r(\gamma)-r(\alpha)\delta(\gamma)$ and $b(\gamma) = \delta(\gamma)$ 
have the same values for $(\rho_1, X_1)$ and $(\rho_2, X_2)$.   
It follows that $\rho_1(\gamma) = a(\gamma) I_2 + b(\gamma) \rho_1(\alpha) 
= a(\gamma) I_2 + b(\gamma) \rho_2(\alpha) = \rho_2(\gamma)$ 
for each $\gamma \in \Gamma$. Hence $(\rho_1, X_1) = (\rho_2, X_2)$, 
which implies the injectivity. 

Let $((\rho_{0}, X_{0}), (r, \delta))$ be a $Z$-valued point 
of $\widetilde{R}_{\alpha_0} \times_{\widetilde{C}_{\alpha_0}} 
\widetilde{C}_{\alpha}$, where $r : \Gamma \to {\mathcal O}_Z(Z)$ is 
a character and $\delta : \Gamma \to {\mathcal O}_Z(Z)$ is a derivation 
with respect to $r$ such that $\delta(\alpha)=1$. 
Put $a(\gamma) = r(\gamma)-r(\alpha)\delta(\gamma)$ and $b(\gamma) = \delta(\gamma)$ for each $\gamma \in \Gamma$. 
By Proposition \ref{prop:abab}, 
$a$ and $b$ are $(a, b)$-coefficients 
with respect to $(r^2, \alpha)$.  
Set $X := X_0$ and 
$\rho(\gamma) := a(\gamma)I_2 + b(\gamma) \rho_{0}(\alpha_0)$ 
for $\gamma \in \Gamma$. 
Since $a(\alpha) = 0$ and $b(\alpha)=1$, 
$\rho(\alpha) = \rho_{0}(\alpha_0)$.  Note that 
$r^2(\alpha) = r^2(\phi(\alpha_0)) = X_{0}^2 = \det \rho_0(\alpha_0)$. 
It follows from 
Lemma \ref{lemma:uf2rep} that $\rho$ is a representation with 
unipotent mold over ${\Bbb F}_2$ such that 
$\det \rho(\gamma) = r^2(\gamma)$ for each 
$\gamma \in \Gamma$.  
Then $(\rho, X)$ is a $Z$-valued point of $\widetilde{R}_{\alpha}$. 
It is easy to check that $(\rho, X)$ is mapped to 
$((\rho_{0}, X_{0}), (r, \delta))$. 
This implies the surjectivity. Hence we have proved the statement. 
\qed 
 
\begin{theorem}\label{th:univgeomuf2tilde} 
The morphism $\widetilde{\pi}_{\Gamma, u/{\Bbb F}_2, \alpha} : 
\widetilde{{\rm Rep}_2(\Gamma)}_{u/{\Bbb F}_2, \alpha} \to 
\widetilde{{\rm Ch}_2(\Gamma)}_{u/{\Bbb F}_2, \alpha}$ 
is a universal 
geometric quotient by ${\rm PGL}_2\otimes_{{\Bbb Z}} {\Bbb F}_2$ 
for each $\alpha \in \Gamma$. 
Furthermore, $\widetilde{\pi}_{\Gamma, u/{\Bbb F}_2} : 
\widetilde{{\rm Rep}_2(\Gamma)}_{u/{\Bbb F}_2} \to 
\widetilde{{\rm Ch}_2(\Gamma)}_{u/{\Bbb F}_2}$ 
is a universal 
geometric quotient by ${\rm PGL}_2\otimes_{{\Bbb Z}} {\Bbb F}_2$. 
\end{theorem}

{\it Proof.} 
The statement follows from Theorem \ref{th:prototypeuf2tilde} and 
Lemma \ref{lemma:fibreproductuf2-2}.  
\qed 

\bigskip

The following Lemma states the ``descent'' of universal geometric 
quotients. The proof was suggested by Michiaki Inaba. 

\begin{lemma}\label{lemma:univgeomdescent}
Let $G$ be a group scheme separated of finite type over a scheme $S$. 
Let $\phi : X \to Y$ be a $G$-equivariant 
separated morphism of finite type over $S$, 
where $G$ acts on $Y$ trivially. 
For a faithfully flat and quasi-compact 
morphism $f : Y' \to Y$, put $X' := X\times_{Y} Y'$ and $\phi' : X' \to Y'$. 
If $\phi'$ is a (universal) geometric quotient by $G$, then 
$\phi$ is also a (resp. universal) geometric quotient by $G$. 
\end{lemma} 

{\it Proof}. 
It suffices to prove that if $\phi'$ is a geometric quotient, 
so is $\phi$. 
It is easy to see that $\phi$ is surjective and that 
the image of $G \times X \to X\times_{S} X$ is equal to 
$X\times_{Y} X$. 
If $\phi'$ is (universally) submersive, then so is $\phi$ by 
\cite[Lemma 15.7.11.1]{EGA4-3}. 
Let $\sigma : G\times_{S} X \to X$ and $\sigma' :  G\times_{S} X' \to X'$ 
be the groups actions of $G$ on $X$ and $X'$, respectively.  
Denote the second projections by 
$p_2 : G\times_{S} X \to X$ and $p_2' :  G\times_{S} X' \to X'$. 
Put $\tau := \phi \circ \sigma = \phi \circ p_2$,  
$\tau' := \phi' \circ \sigma' = \phi' \circ p_2'$, and 
$f' : X' \to X$.  
For proving that $\phi_{\ast}({\mathcal O}_X)^G = {\mathcal O}_Y$, 
we show that 
$0 \to {\mathcal O}_Y \to \phi_{\ast}({\mathcal O}_X) 
\stackrel{\sigma_{\ast}-p_{2 \ast}}{\to} \tau_{\ast}({\mathcal O}_{G\times_{S} X})$ 
is exact. 
Taking the pull-back by $f$, 
we have  
$0 \to f^{\ast}{\mathcal O}_Y \to f^{\ast}\phi_{\ast}({\mathcal O}_X) 
\stackrel{f^{\ast}(\sigma_{\ast}-p_{2 \ast})}{\to} 
f^{\ast}\tau_{\ast}({\mathcal O}_{G\times_{S} X})$. 
By \cite[Proposition 1.4.15]{EGA3}, 
$f^{\ast}\phi_{\ast}({\mathcal O}_X) \cong \phi'_{\ast}f'^{\ast}({\mathcal O}_X) 
\cong \phi'_{\ast}({\mathcal O}_{X'})$ 
and $f^{\ast}\tau_{\ast}({\mathcal O}_{G\times_{S} X}) 
\cong \tau'_{\ast}(1_{G}\times f')^{\ast}({\mathcal O}_{G\times_{S} X}) 
\cong \tau'_{\ast}({\mathcal O}_{G\times_{S} X'})$. 
Then we obtain the following commutative diagram:  
\[
\begin{array}{ccccccc}
0 & \to & f^{\ast}({\mathcal O}_Y) & \to & f^{\ast}\phi_{\ast}({\mathcal O}_X) &  
\stackrel{f^{\ast}(\sigma_{\ast}-p_{2 \ast})}{\to} & 
f^{\ast}\tau_{\ast}({\mathcal O}_{G\times_{S} X}) \\ 
 & &|| & & \downarrow \cong & & \downarrow \cong \\
0 & \to & {\mathcal O}_{Y'} & \to  & \phi'_{\ast}({\mathcal O}_X') &  
\stackrel{\sigma'_{\ast}-p'_{2 \ast}}{\to} & \tau'_{\ast}({\mathcal O}_{G\times_{S} X'}). 
\end{array} 
\]
Since $\phi'$ is a geometric quotient, 
the complex $0 \to {\mathcal O}_{Y'} \to  \phi'_{\ast}({\mathcal O}_X')   
\stackrel{\sigma'_{\ast}-p'_{2 \ast}}{\to}  \tau'_{\ast}({\mathcal O}_{G\times_{S} X'})$ is exact, and hence so is 
$0 \to f^{\ast}{\mathcal O}_Y \to f^{\ast}\phi_{\ast}({\mathcal O}_X) 
\stackrel{f^{\ast}(\sigma_{\ast}-p_{2 \ast})}{\to} 
f^{\ast}\tau_{\ast}({\mathcal O}_{G\times_{S} X})$.  
Because $f$ is faithfully flat, 
$0 \to {\mathcal O}_Y \to \phi_{\ast}({\mathcal O}_X) 
\stackrel{\sigma_{\ast}-p_{2 \ast}}{\to} \tau_{\ast}({\mathcal O}_{G\times_{S} X})$ 
is also exact. 
Thus we have proved the statement. 
\qed 

\bigskip 

Now we can prove Theorem~\ref{th:univgeomuf2alpha} 
by another approach. 

\begin{theorem}[Theorem~\ref{th:univgeomuf2alpha}]
The morphism $\pi_{\Gamma, u/{\Bbb F}_2, \alpha} : {\rm Rep}_2(\Gamma)_{u/{\Bbb F}_2, \alpha} 
\to {\rm Ch}_2(\Gamma)_{u/{\Bbb F}_2, \alpha}$ 
is a universal 
geometric quotient by ${\rm PGL}_2\otimes_{{\Bbb Z}} {\Bbb F}_2$ 
for each $\alpha \in \Gamma$. 
\end{theorem} 

{\it Proof}. 
The statement follows from 
Proposition~\ref{prop:uf2fibre}, Theorem~\ref{th:univgeomuf2tilde}, 
and Lemma~\ref{lemma:univgeomdescent}. 
\qed 

\bigskip 

By the same discussion in \S 6, 
we can construct ${\rm Ch}_2(\Gamma)_{u/{\Bbb F}_2}$ by gluing 
$\{ {\rm Ch}_2(\Gamma)_{u/{\Bbb F}_2, \alpha} \}_{\alpha \in \Gamma}$. 
Then we have Corollary~\ref{cor:univgeomuf2-total}, which states that  
$\pi_{\Gamma, u/{\Bbb F}_2} : {\rm Rep}_2(\Gamma)_{u/{\Bbb F}_2} 
\to {\rm Ch}_2(\Gamma)_{u/{\Bbb F}_2}$ 
is a universal geometric quotient by ${\rm PGL}_2\otimes_{{\Bbb Z}} {\Bbb F}_2$. 
Similarly, 
we have a morphism 
$q :  \widetilde{{\rm Ch}_2(\Gamma)}_{u/{\Bbb F}_2} 
\to {\rm Ch}_2(\Gamma)_{u/{\Bbb F}_2}$ by gluing 
$\{ q_{\alpha} : \widetilde{{\rm Ch}_2(\Gamma)}_{u/{\Bbb F}_2, \alpha} 
\to {\rm Ch}_2(\Gamma)_{u/{\Bbb F}_2, \alpha} \}_{\alpha \in \Gamma}$.

\begin{remark}\rm
By Remark~\ref{remark:prototypeuf2tilde}, we see that 
$q : \widetilde{{\rm Ch}_2(\Upsilon_1)}_{u/{\Bbb F}_2, \alpha_0} \to 
{\rm Ch}_2(\Upsilon_1)_{u/{\Bbb F}_2, \alpha_0}$ is described as 
${\rm Spec} {\Bbb F}_2[\chi] \to {\rm Spec} {\Bbb F}_2[D]$, where 
$D$ is mapped to $\chi^2$.   
We also see that $p : 
\widetilde{{\rm Rep}_2(\Upsilon_1)}_{u/{\Bbb F}_2, \alpha_0} \to 
{\rm Rep}_2(\Upsilon_1)_{u/{\Bbb F}_2, \alpha_0}$ is described as 
${\rm Spec} \:{\Bbb F}_2 [a, b, c, X]/(X^2 + a^2 + bc) \supset 
D(b) \cup D(c) \to 
D(b) \cup D(c) \subset 
{\rm Spec} {\Bbb F}_2 [a, b, c]$. 
Hence $p$ and $q$ are faithfully flat finite morphisms 
of finite presentation, but not smooth morphisms. 
\end{remark}

For $\alpha \in \Gamma$, the monoid homomorphism 
$\phi : \Upsilon_{1} \to \Gamma$ by 
$\alpha_0 \mapsto \alpha$ induces the following 
commutative diagrams: 
\[
\begin{array}{ccc}
 \widetilde{{\rm Rep}_2(\Gamma)}_{u/{\Bbb F}_2, \alpha} 
& {\to} & {\rm Rep}_2(\Gamma)_{u/{\Bbb F}_2, \alpha} \\
 \downarrow & & \downarrow  \\
\widetilde{{\rm Rep}_2(\Upsilon_1)}_{u/{\Bbb F}_2, \alpha_0} 
& {\to} & {\rm Rep}_2(\Upsilon_1)_{u/{\Bbb F}_2, \alpha_0}  
\end{array}  
\] 
and 
\[
\begin{array}{ccc}
 \widetilde{{\rm Ch}_2(\Gamma)}_{u/{\Bbb F}_2, \alpha} 
& {\to} & {\rm Ch}_2(\Gamma)_{u/{\Bbb F}_2, \alpha} \\
 \downarrow & & \downarrow  \\
\widetilde{{\rm Ch}_2(\Upsilon_1)}_{u/{\Bbb F}_2, \alpha_0} 
& {\to} & {\rm Ch}_2(\Upsilon_1)_{u/{\Bbb F}_2, \alpha_0}.  
\end{array}  
\] 
We can show that the diagrams above give fibre products 
in the same way as Lemma \ref{lemma:fibreproductuf2-2}.

From the discussions above, we obtain the following commutative 
diagram which gives a fibre product: 
\[
\begin{array}{ccr}
\phantom{AAAAA} \widetilde{{\rm Rep}_2(\Gamma)}_{u/{\Bbb F}_2} 
& \stackrel{p_{}}{\to} & {\rm Rep}_2(\Gamma)_{u/{\Bbb F}_2} \\
\widetilde{\pi}_{\Gamma, u/{\Bbb F}_2} \downarrow & & \downarrow 
\pi_{\Gamma, u/{\Bbb F}_2} \\
\phantom{AAAAA} \widetilde{{\rm Ch}_2(\Gamma)}_{u/{\Bbb F}_2} 
& \stackrel{q_{}}{\to} & {\rm Ch}_2(\Gamma)_{u/{\Bbb F}_2}.  
\end{array}  
\] 
All morphisms are ${\rm PGL}_2\otimes_{\Bbb Z} {\Bbb F}_2$-equivariant, and 
$p$ and $q$ are faithfully flat finite morphisms of finite presentation. 
Note that $p$ and $q$ are not smooth in general. 

Recall that $\widetilde{{\rm Rep}_2(\Gamma)}_{u/{\Bbb F}_2}$ can be regarded 
as a ${\rm Rep}_1(\Gamma)_{{\Bbb F}_2}$-scheme by $\lambda : 
\widetilde{{\rm Rep}_2(\Gamma)}_{u/{\Bbb F}_2} \to 
{\rm Rep}_1(\Gamma)_{{\Bbb F}_2}$. 
Let $r$ be the corresponding character on 
$\widetilde{{\rm Rep}_2(\Gamma)}_{u/{\Bbb F}_2}$ to $\lambda$, and let 
$\lambda^2 : \widetilde{{\rm Rep}_2(\Gamma)}_{u/{\Bbb F}_2} \to 
{\rm Rep}_1(\Gamma)_{{\Bbb F}_2}$ be the morphism induced by the character $r^2$. 
Denote by $\det : {\rm Rep}_2(\Gamma)_{u/{\Bbb F}_2} 
\to {\rm Rep}_1(\Gamma)_{{\Bbb F}_2}$ the morphism corresponding to 
the character $\det(\sigma_{\Gamma, u/{\Bbb F}_2}(\cdot)))$. 
Then $\det \circ\; p = \lambda^2$. 

By Definition \ref{def:chtilde}, $\widetilde{{\rm Ch}_2(\Gamma)}_{u/{\Bbb F}_2}$ is 
a ${\rm Rep}_1(\Gamma)_{{\Bbb F}_2}$-scheme. 
Let us denote by $r : \widetilde{{\rm Ch}_2(\Gamma)}_{u/{\Bbb F}_2} \to 
{\rm Rep}_1(\Gamma)_{{\Bbb F}_2}$ the canonical morphism. 
We also denote by the same symbol $r$ the corresponding character on 
$\widetilde{{\rm Ch}_2(\Gamma)}_{u/{\Bbb F}_2}$ to $r$.  
We define $r^2 : \widetilde{{\rm Ch}_2(\Gamma)}_{u/{\Bbb F}_2} \to 
{\rm Rep}_1(\Gamma)_{{\Bbb F}_2}$ as the morphism induced by 
the character $r^2$ on $\widetilde{{\rm Ch}_2(\Gamma)}_{u/{\Bbb F}_2}$. 
By Remark \ref{remark:dpidet}, 
${\rm Ch}_2(\Gamma)_{u/{\Bbb F}_2}$ is also a   
${\rm Rep}_1(\Gamma)_{{\Bbb F}_2}$-scheme by 
$d : {\rm Ch}_2(\Gamma)_{u/{\Bbb F}_2} \to 
{\rm Rep}_1(\Gamma)_{{\Bbb F}_2}$.  
Then $d \circ q=r^2$.

\begin{remark}\rm  
  The morphism $\widetilde{\pi}_{\Gamma, u/{\Bbb F}_2} : 
\widetilde{{\rm Rep}_2(\Gamma)}_{u/{\Bbb F}_2} \to 
\widetilde{{\rm Ch}_2(\Gamma)}_{u/{\Bbb F}_2}$  is smooth and surjective for each group or monoid $\Gamma$. 
Indeed, $\widetilde{\pi}_{\Gamma, u/{\Bbb F}_2, \alpha} : \widetilde{{\rm Rep}_2(\Gamma)}_{u/{\Bbb F}_2, \alpha}  \to  
\widetilde{{\rm Ch}_2(\Gamma)}_{u/{\Bbb F}_2, \alpha}$ 
is obtained by base change of the prototype by Lemma \ref{lemma:fibreproductuf2-2}. 
The prototype $\widetilde{\pi}_{\Upsilon_1, u/{\Bbb F}_2, \alpha_0} : 
\widetilde{{\rm Rep}_2(\Upsilon_1)}_{u/{\Bbb F}_2, \alpha_0} 
\to \widetilde{{\rm Ch}_2(\Upsilon_1)}_{u/{\Bbb F}_2, \alpha_0}$ 
is smooth and surjective because it is obtained by base change of 
$\pi : {\rm Rep}_2(\Upsilon_{1})_{{\rm rk} 2} \to {\rm Ch}_2(\Upsilon_1)$ and 
$\pi$ is smooth and surjective by Proposition~\ref{lemmassff}. 
\end{remark}

\begin{example}\label{example:freemonoidtildemoduli}\rm 
Let us describe $\widetilde{{\rm Ch}_2(\Upsilon_m)}_{u/{\Bbb F}_2}$ 
for the free monoid $\Upsilon_m = \langle \alpha_1, \ldots, \alpha_m \rangle$ 
of rank $m$. 
Put $\widetilde{C(m)} := \widetilde{{\rm Ch}_2(\Upsilon_m)}_{u/{\Bbb F}_2}$. 
Let $A_1(m)_{{\Bbb F}_2}$ denote the coordinate ring $A_1(\Upsilon_m)_{{\Bbb F}_2}$ 
of ${\rm Rep}_1(\Upsilon_m)_{{\Bbb F}_2}$. We can write 
$A_1(m)_{{\Bbb F}_2} = {\Bbb F}_2[\chi(\alpha_1), \ldots, \chi(\alpha_m)]$, 
where $\chi(\alpha_1), \ldots, \chi(\alpha_m)$ are indeterminates. 
It is easy to see that the $A_1(m)_{{\Bbb F}_2}$-module 
$\Omega_{{\Upsilon_m}/{\Bbb F}_2}$ 
is isomorphic to the free module $\oplus_{i=1}^m A_1(m)_{{\Bbb F}_2}\cdot d\alpha_i$. 
Hence $\widetilde{C(m)} = {\rm Proj} \: S(\Omega_{{\Upsilon_m}/{\Bbb F}_2})$ 
is isomorphic to ${\Bbb A}_{{\Bbb F}_2}^m \times {\Bbb P}_{{\Bbb F}_2}^{m-1}$.  
The projection $\widetilde{C(m)} \to {\rm Rep}_1(\Upsilon_m)_{{\Bbb F}_2}$ 
can be described by the first projection $p_1 : {\Bbb A}_{{\Bbb F}_2}^m 
\times {\Bbb P}_{{\Bbb F}_2}^{m-1} \to {\Bbb A}_{{\Bbb F}_2}^m = {\rm Rep}_1(\Upsilon_m)_{{\Bbb F}_2}$. 

Put $\widetilde{C(m)}_i := \{ d\alpha_i \neq 0 \} \subset \widetilde{C(m)}$ 
for $1 \le i \le m$. Note that $\widetilde{C(m)}_i \cong {\Bbb A}_{{\Bbb F}_2}^m 
\times {\Bbb A}_{{\Bbb F}_2}^{m-1} \subset {\Bbb A}_{{\Bbb F}_2}^m 
\times {\Bbb P}_{{\Bbb F}_2}^{m-1}$.  
In Example \ref{ex:describe-chuf2m}, we have described 
$C(m) = {\rm Ch}_2(\Upsilon_m)_{u/{\Bbb F}_2}$ and 
$C(m)_i = {\rm Ch}_2(\Upsilon_m)_{u/{\Bbb F}_2, \alpha_i}$ for  
$1 \le i \le m$. 
The morphism $q : \widetilde{C(m)} \to C(m)$ can be 
described as follows: Let $q_i : \widetilde{C(m)}_i \to C(m)_i$ be 
the restriction of $q$ to $\widetilde{C(m)}_i$ for $1 \le i \le m$. 
For $(r, \delta) \in \widetilde{C(m)}_i$, $q_i(r, \delta) = 
(a, b) \in C(m)_i$ is given by   
$a(\gamma)= r(\gamma)-r(\alpha_i)\delta(\gamma)$ and 
$b(\gamma)=\delta(\gamma)$ for $\gamma \in \Upsilon_m$. 
Recall that the isomorphisms $\widetilde{C(m)}_i \cong {\Bbb A}_{{\Bbb F}_2}^{2m-1}$ 
and $C(m)_i \cong {\Bbb A}_{{\Bbb F}_2}^{2m-1}$ are 
given by 
\begin{multline*}
(r, \delta) \mapsto (r(\alpha_1), \ldots, r(\alpha_m), 
\delta(\alpha_1)/\delta(\alpha_i), \ldots, \delta(\alpha_{i-1})/\delta(\alpha_i), \\ 
\delta(\alpha_{i+1})/\delta(\alpha_i), \ldots, \delta(\alpha_{m})/\delta(\alpha_i))
\end{multline*} and 
\begin{multline*}
(a_i, b_i, d) \mapsto (a_i(\alpha_1), \ldots, a_i(\alpha_{i-1}), a_i(\alpha_{i+1}), 
\ldots, a_i(\alpha_m), 
b_i(\alpha_1), \ldots, \\ b_i(\alpha_{i-1}), 
b_i(\alpha_{i+1}), \ldots, b_i(\alpha_{m}), d(\alpha_i)), 
\end{multline*} 
respectively. 
By these isomorphisms, 
$q_i : {\Bbb A}_{{\Bbb F}_2}^{2m-1} 
\to {\Bbb A}_{{\Bbb F}_2}^{2m-1}$ is described by 
\begin{multline*} 
q_i(r_1, \cdots, r_m, \overline{\delta}_1, \ldots, \overline{\delta}_{i-1}, 
\overline{\delta}_{i+1}, \ldots, \overline{\delta}_m) =  \\ 
(r_1-r_i \overline{\delta}_1, \ldots, r_{i-1}-r_i \overline{\delta}_{i-1}, r_{i+1}-r_i \overline{\delta}_{i+1}, \ldots, r_m - r_i\overline{\delta}_{m},  \\ 
\overline{\delta}_1, \ldots, \overline{\delta}_{i-1}, 
\overline{\delta}_{i+1}, \ldots, \overline{\delta}_m, r_i^2). 
\end{multline*} 

Set $\widetilde{R(m)} := \widetilde{{\rm Rep}_2(\Upsilon_m)}_{{\Bbb F}_2}$ 
and $R(m) := {\rm Rep}_2(\Upsilon_m)_{{\Bbb F}_2}$. 
For $1 \le i \le m$, put 
$\widetilde{R(m)}_i := \widetilde{{\rm Rep}_2(\Upsilon_m)}_{{\Bbb F}_2, \alpha_i}$ 
and $R(m)_i := {\rm Rep}_2(\Upsilon_m)_{{\Bbb F}_2, \alpha_i}$. 
Let $p_i : \widetilde{R(m)}_i \to R(m)_i$ be 
the restriction of $p : \widetilde{R(m)} \to R(m)$ 
to $\widetilde{R(m)}_i$ for $1 \le i \le m$. We can describe  
$p_i : \widetilde{R(m)}_i = \{ (A_1, \ldots, A_m, X_i) \mid (A_1, \ldots, A_m) 
\in R(m)_i \mbox{ and } X_i^2 = \det A_i  \} 
\to R(m)_i = \{ (A_1, \ldots, A_m) \mid \langle A_1, \ldots, A_m \rangle = 
\langle A_i \rangle 
\mbox{  is a unipotent 
mold over } {\Bbb F}_2 \}$ by $(A_1, \ldots, A_m, X_i) \allowbreak
\mapsto (A_1, \ldots, A_m)$, where 
$A_j := \rho(\alpha_j)$ for $1 \le j \le m$ and for each representation 
$\rho$. Let $\widetilde{\pi}_i$ and $\pi_i$ denote the restrictions of 
$\widetilde{\pi}_{\Upsilon_m, u/{\Bbb F}_2} : 
\widetilde{R(m)} \to \widetilde{C(m)}$ and $\pi_{\Upsilon_m, u/{\Bbb F}_2} : 
R(m) \to C(m)$ to $\widetilde{R(m)}_i$ and $R(m)_i$, respectively. 
For $(A_1, \ldots, A_m) \in R(m)_i$, we can write 
$A_j = \overline{a}_{ij}I_2 + \overline{b}_{j}A_i$ for $1 \le j \le m$.  
Then $\widetilde{\pi}_i : \widetilde{R(m)}_i \to \widetilde{C(m)}_i$ is described by 
$(A_1, \ldots, A_m, X_i) \mapsto (\overline{a}_{i1}+\overline{b}_{1}X_i, 
\ldots, \overline{a}_{im}+\overline{b}_{m}X_i, \overline{b}_1, \ldots, \overline{b}_{i-1}, 
\overline{b}_{i+1}, \ldots, \overline{b}_{m})$ and 
$\pi_i : R(m)_i \to C(m)_i$ is described by $(A_1, \ldots, A_m) \mapsto 
\allowbreak 
(\overline{a}_{i1}, 
\ldots, \overline{a}_{i, i-1}, \overline{a}_{i, i+1}, 
\ldots, \overline{a}_{im}, \overline{b}_1, \ldots, \overline{b}_{i-1}, \\ 
\overline{b}_{i+1}, \ldots, \overline{b}_{m}, \det A_i)$. 
Remark that 
\[
\begin{array}{ccr}
\phantom{AA} \widetilde{R(m)}_{i} 
& \stackrel{p_{i}}{\to} & R(m)_{i} \\
\widetilde{\pi}_{i} \downarrow & & \downarrow 
\pi_{i} \\
\phantom{AA} \widetilde{C(m)}_{i} 
& \stackrel{q_{i}}{\to} & C(m)_{i}.  
\end{array}  
\] 
gives a fibre product. 
\end{example}

\bigskip 

\begin{definition}\label{def:tilderep}\rm  
Let $X$ be an ${\Bbb F}_2$-scheme. 
By a {\it tilde representation with unipotent mold over ${\Bbb F}_2$} for $\Gamma$ 
on $X$, we understand a pair $(\rho, \lambda)$ of 
a representation $\rho$ of with unipotent mold over ${\Bbb F}_2$ 
for $\Gamma$ on $X$ and a character $\lambda : \Gamma 
\to {\mathcal O}_X(X)$ satisfying the following conditions: 
\begin{enumerate} 
\item $\det(\rho(\gamma)) = \lambda(\gamma)^2$ for each $\gamma \in \Gamma$. 
\item\label{cond:sublinebundle-tilde} 
$\{ \rho(\gamma) - \lambda(\gamma) I_2 \mid \gamma \in \Gamma \}$ 
spans a sub-line bundle of 
${\mathcal O}_X[\rho(\Gamma)]$. 
\end{enumerate} 
\end{definition} 

\begin{remark}\label{remark:tilderep}\rm
Let $(\rho, \lambda)$ be a tilde representation  with unipotent mold over 
${\Bbb F}_2$ for $\Gamma$ 
on an ${\Bbb F}_2$-scheme $X$.  For each point $x \in X$, 
choose $\alpha_x \in \Gamma$ and a neighbourhood $U_x$ of $x$ such that 
${\mathcal O}_{U_x}[\rho(\Gamma)]={\mathcal O}_{U_x} \cdot I_2 \oplus 
{\mathcal O}_{U_x} \cdot \rho(\alpha_x)$.  
The condition (\ref{cond:sublinebundle-tilde}) in Definition \ref{def:tilderep} means that 
for each $\gamma \in \Gamma$ there exists $c \in {\mathcal O}_{U_x}(U_{x})$ such that 
$\rho(\gamma) - \lambda(\gamma) I_2 = c(\rho(\alpha_x) - \lambda(\alpha_x) I_2)$. 
Since $\rho(\gamma) = (\lambda(\gamma) - c\lambda(\alpha_x))I_2 + c\rho(\alpha_x)$, 
the $(a, b)$-coefficients of $\rho(\gamma)$ with respect to 
$\rho(\alpha_x)$ are given by $a_{\alpha_x}(\gamma) = \lambda(\gamma) - c\lambda(\alpha_x)$ and $b_{\alpha_x}(\gamma) = c$. 
Then $\lambda(\gamma) = a_{\alpha_x}(\gamma) + 
b_{\alpha_x}(\gamma)\lambda(\alpha_x)$ for each $\gamma \in \Gamma$. Note that 
$(\rho\!\!\mid_{U_x}, \lambda(\alpha_x))$ gives a $U_x$-valued point of 
$\widetilde{{\rm Rep}_2(\Gamma)}_{u/{\Bbb F}_2, \alpha_x}$. 
Considering the definition of $\widetilde{{\rm Rep}_2(\Gamma)}_{u/{\Bbb F}_2}$, 
we see that we can obtain an $X$-valued point of 
$\widetilde{{\rm Rep}_2(\Gamma)}_{u/{\Bbb F}_2}$ by gluing $\{ (\rho\!\!\mid_{U_x}, \lambda(\alpha_x))\}_{x\in X}$.
\end{remark} 

\bigskip 

By Remark \ref{remark:tilderep}, we have:  

\begin{proposition}\label{prop:tilderepfunctor} 
The following functor is representable by 
$\widetilde{{\rm Rep}_2(\Gamma)}_{u/{\Bbb F}_2}$:  
\[
\begin{array}{ccccl}
\widetilde{{\rm Rep}_2(\Gamma)}_{u/{\Bbb F}_2} 
 & : & ({\bf Sch}/{\Bbb F}_2)^{op} & \to & ({\bf Sets}) \\
 & & X & \mapsto & 
\left\{ 
\begin{array}{r}  
\mbox{ tilde rep. 
with unipotent  } \\ 
\mbox{ mold over ${\Bbb F}_2$ for $\Gamma$ on } X  
\end{array} 
\right\} .
\end{array}
\] 
\end{proposition} 

\begin{remark}\rm 
The condition (\ref{cond:sublinebundle-tilde}) in Definition \ref{def:tilderep} 
is necessary for Proposition \ref{prop:tilderepfunctor}. 
Indeed, in the case of the free monoid 
$\Upsilon_2 = \langle \alpha, \beta \rangle$ of rank $2$, let $\rho : \Upsilon_2 \to 
{\rm M}_2(k[\epsilon]/(\epsilon^2))$ be the representation defined by 
\[
\rho(\alpha) = 
\left( 
\begin{array}{cc} 
0 & 0 \\
1 & 0 \\
\end{array} 
\right), \; 
\rho(\beta) = 
\left( 
\begin{array}{cc} 
1 & 0 \\
1 & 1 \\
\end{array} 
\right), 
\]
where $k$ is a field.  
Let $\lambda : \Upsilon_2 \to k[\epsilon]/(\epsilon^2)$ be the character 
defined by $\lambda(\alpha)=\epsilon$ and $\lambda(\beta) = 1$. 
Then $\rho$ is a representation with unipotent mold over ${\Bbb F}_2$ for $\Upsilon_2$. 
For each $\gamma \in \Upsilon_2$, $\det \rho(\gamma) = \lambda(\gamma)^2$. 
However, the condition (\ref{cond:sublinebundle-tilde}) in Definition \ref{def:tilderep} fails.  
The $(a, b)$-coefficients of $\rho(\beta)$ with respect to 
$\rho(\alpha)$ is given by $a_{\alpha}(\beta)=1$ and $b_{\alpha}(\beta)=1$. 
The equality $\lambda(\beta) = a_{\alpha}(\beta)+ b_{\alpha}(\beta)\lambda(\alpha)$ 
does not hold. 
Hence $(\rho, \lambda(\alpha)) \in \widetilde{{\rm Rep}_2(\Gamma)}_{u/{\Bbb F}_2, \alpha}(k[\epsilon]/(\epsilon^2))$ and $(\rho, \lambda(\beta)) \in \widetilde{{\rm Rep}_2(\Gamma)}_{u/{\Bbb F}_2, \beta}(k[\epsilon]/(\epsilon^2))$ induce different 
morphisms from $k[\epsilon]/(\epsilon^2)$ to $\widetilde{{\rm Rep}_2(\Gamma)}_{u/{\Bbb F}_2}$.  
This means that $(\rho, \lambda)$ does not canonically  
induce a morphism to $\widetilde{{\rm Rep}_2(\Gamma)}_{u/{\Bbb F}_2}$ 
without the condition  (\ref{cond:sublinebundle-tilde}).  
\end{remark}

\begin{remark}\label{remark:sectionuf2tilde}\rm 
For each point $x \in \widetilde{{\rm Ch}_2(\Gamma)}_{u/{\Bbb F}_2}$, there exists 
a local section $\widetilde{s_{x}} : V_x \to 
\widetilde{{\rm Rep}_2(\Gamma)}_{u/{\Bbb F}_2}$ on 
a neighbourhood $V_x$ of $x$ such that 
$\widetilde{\pi}_{\Gamma, u/{\Bbb F}_2} \circ \widetilde{s_{x}} = id_{V_x}$. 
Indeed, take $\alpha \in \Gamma$ such that 
$x \in \widetilde{{\rm Ch}_2(\Gamma)}_{u/{\Bbb F}_2, \alpha}$. 
By Proposition \ref{prop:uf2fibre}, $\widetilde{\pi_{\Gamma, u/{\Bbb F}_2, \alpha}} : 
\widetilde{{\rm Rep}_2(\Gamma)}_{u/{\Bbb F}_2, \alpha} \to 
\widetilde{{\rm Ch}_2(\Gamma)}_{u/{\Bbb F}_2, \alpha}$ is 
obtained by base change of $\pi_{\Gamma, u/{\Bbb F}_2, \alpha} : 
{\rm Rep}_2(\Gamma)_{u/{\Bbb F}_2, \alpha} \to 
{\rm Ch}_2(\Gamma)_{u/{\Bbb F}_2, \alpha}$. 
Remark \ref{remark:sectionuf2} follows that  $\pi_{\Gamma, u/{\Bbb F}_2, \alpha}$ 
has a section $s_{\Gamma, \alpha}$. Hence 
$\widetilde{\pi}_{\Gamma, u/{\Bbb F}_2, \alpha}$ has a section 
$\widetilde{s_{\Gamma, \alpha}}$. 
We can take $\widetilde{{\rm Ch}_2(\Gamma)}_{u/{\Bbb F}_2, \alpha}$ 
as a neighbourhood $V_x$ of $x$. It is easy to see that 
$(\rho, \lambda) = \widetilde{s_{\Gamma, \alpha}}(r, \delta)$ is described by 
$\rho(\gamma) = a(\gamma)I_2+b(\gamma)
\left( 
\begin{array}{cc}
0 & -r(\alpha)^2 \\
1 & 0 \\
\end{array}
\right)$ and $\lambda(\gamma)  
=r(\gamma)$ for $\gamma \in \Gamma$, where 
$a(\gamma)=r(\gamma)-r(\alpha)\delta(\gamma)$ and 
$b(\gamma) = \delta(\gamma)$. 
\end{remark}


\begin{lemma}\label{lemma:localequf2tilde} 
Let $(\rho_1, \lambda_1),  (\rho_2, \lambda_2)$ be tilde 
representations with unipotent mold over ${\Bbb F}_2$ for a group (or a monoid) 
$\Gamma$ on a scheme $X$ over ${\Bbb F}_2$. 
Let $f_i : X \to \widetilde{{\rm Rep}_2(\Gamma)}_{u/{\Bbb F}_2}$ be the morphism 
associated to $(\rho_i, \lambda_i)$ for $i = 1, 2$. 
If $\widetilde{\pi}_{\Gamma, u/{\Bbb F}_2} \circ f_1 = 
\widetilde{\pi}_{\Gamma, u/{\Bbb F}_2} \circ f_2 
: X \to \widetilde{{\rm Ch}_2(\Gamma)}_{u/{\Bbb F}_2}$, then for each $x \in X$ 
there exists $P_x \in {\rm GL}_2(\Gamma(V_x, {\mathcal O}_X))$ 
on a neighbourhood $V_x$ of $x$ such that 
$P_x^{-1}\rho_1 P_x = \rho_2$ and $\lambda_1 = \lambda_2$ on $V_x$.  
\end{lemma} 

\prf In the same way as Lemma \ref{lemma:localequf2}, we can prove 
the statement.  
\qed

\bigskip

By a {\it generalized tilde representation with unipotent mold over ${\Bbb F}_2$} for $\Gamma$ on $X$, we understand triples $\{ (U_i, \rho_i, \lambda_i) \}_{i \in I}$ 
of an open set $U_i$ and a tilde representation $(\rho_i, \lambda_i)$  
with unipotent mold over ${\Bbb F}_2$ for $\Gamma$ on $U_i$  
satisfying the following three conditions:
\begin{enumerate}
\item $\cup_{i \in I} U_ i = X$, 
\item for each $x \in U_i \cap U_j$, there exists $P_x \in {\rm GL}_2(\Gamma(V_x, {\mathcal O}_X))$ on 
a neighbourhood $V_x \subseteq U_i \cap U_j$ of $x$ such that 
$P_x^{-1} \rho_i P_x = \rho_j$ on $V_x$. 
\item $\lambda_i = \lambda_j$ on $U_i \cap U_j$ for each $i, j$. 
\end{enumerate}

Generalized tilde representations 
$\{ (U_i, \rho_i, \lambda_i) \}_{i \in I}$ and $\{ (V_j, \sigma_j, \mu_j) \}_{j \in J}$ 
with unipotent mold over ${\Bbb F}_2$ 
are called {\it equivalent} if $\{ (U_i, \rho_i, \lambda_i) \}_{i \in I} 
\cup \{ (V_j, \sigma_j, \mu_j) \}_{j \in J}$ 
is a generalized tilde representation with unipotent mold over ${\Bbb F}_2$ again.  
Let us define the contravariant functor $\widetilde{{\mathcal E}q \:\mathcal{U}_2(\Gamma)}_{{\Bbb F}_2}$:  
\[
\begin{array}{ccccl}
\widetilde{{\mathcal E}q \:\mathcal{U}_2(\Gamma)}_{{\Bbb F}_2} & : & ({\bf Sch}/{\Bbb F}_2)^{op} & \to & ({\bf Sets}) \\
 & & X & \mapsto & 
\left\{ 
\begin{array}{r}  
\mbox{ gen. tilde rep. 
with unip.  } \\ 
\mbox{  mold over ${\Bbb F}_2$ for $\Gamma$ on } X  
\end{array} 
\right\} 
\Big/\sim.
\end{array}
\]
 

\begin{theorem}\label{th:moduliuf2tilde}  
The scheme $\widetilde{{\rm Ch}_2(\Gamma)}_{u/{\Bbb F}_2}$ is a fine moduli scheme  
associated to the functor 
$\widetilde{{\mathcal E}q \:\mathcal{U}_2(\Gamma)}_{{\Bbb F}_2}$ 
for a group or a monoid $\Gamma$.   
In other words, $\widetilde{{\rm Ch}_2(\Gamma)}_{u/{\Bbb F}_2}$ represents 
the functor $\widetilde{{\mathcal E}q \:\mathcal{U}_2(\Gamma)}_{{\Bbb F}_2}$. 
The moduli $\widetilde{{\rm Ch}_2(\Gamma)}_{u/{\Bbb F}_2}$ 
is separated over ${\Bbb F}_2$; 
if $\Gamma$ is a finitely generated group or monoid, then 
$\widetilde{{\rm Ch}_2(\Gamma)}_{u/{\Bbb F}_2}$ is of finite type over ${\Bbb F}_2$. 
\end{theorem}

\prf
In the same way as Theorem \ref{th:moduliuf2}, we can prove that $\widetilde{{\rm Ch}_2(\Gamma)}_{u/{\Bbb F}_2}$ represents 
the functor $\widetilde{{\mathcal E}q \:\mathcal{U}_2(\Gamma)}_{{\Bbb F}_2}$ by 
using Lemma \ref{lemma:localequf2tilde}. It follows from 
Definition \ref{def:chtilde} that $\widetilde{{\rm Ch}_2(\Gamma)}_{u/{\Bbb F}_2}$ 
is separated over ${\Bbb F}_2$. If $\Gamma$ is finitely generated, then 
$\Omega_{\Gamma/{\Bbb F}_2}$ is a finitely generated module over 
$A_1(\Gamma)_{{\Bbb F}_2}$ by
Remark \ref{remark:omegaisfgmoduletilde}. Hence 
$\widetilde{{\rm Ch}_2(\Gamma)}_{u/{\Bbb F}_2}$ is of finite type over ${\Bbb F}_2$. 
\qed

\begin{remark}\label{remark:widetildeuf2algebra}\rm 
For an associative algebra $A$ over a commutative ring $R$ 
over ${\Bbb F}_2$, we can construct 
$\widetilde{\pi}_{A, u/{\Bbb F}_2} : 
\widetilde{{\rm Rep}_2(A)}_{u/{\Bbb F}_2} 
\to \widetilde{{\rm Ch}_2(A)}_{u/{\Bbb F}_2}$ 
in the same way as group or monoid cases.  
Indeed, for $c \in A$, $\widetilde{{\rm Rep}_2(A)}_{u/{\Bbb F}_2, c}$ is 
defined as in Definition \ref{def:widetilde}. 
By gluing $\{ \widetilde{{\rm Rep}_2(A)}_{u/{\Bbb F}_2, c} \}_{c \in A}$, 
we have an $R$-scheme $\widetilde{{\rm Rep}_2(A)}_{u/{\Bbb F}_2}$  such that 
$p : 
\widetilde{{\rm Rep}_2(A)}_{u/{\Bbb F}_2} \to 
{\rm Rep}_2(A)_{u/{\Bbb F}_2}$ is a faithfully flat finite morphism. 
Let $A_1(A)$ be the coordinate ring of ${\rm Rep}_1(A)$. 
As in Remark~\ref{remark:constructu}, we can  
construct $A_1(A)$-module $\Omega_{A/R}$ such that 
${\rm Der}(A, M) \cong {\rm Hom}_{A_1(A)}(\Omega_{A/R}, M)$ 
for any $A_1(A)$-module $M$. Set $\widetilde{{\rm Ch}_2(A)}_{u/{\Bbb F}_2} 
:= {\rm Proj} S(\Omega_{A/R})$. Then we can construct 
a ${\rm Rep}_1(A)$-morphism 
$\widetilde{\pi}_{A, u/{\Bbb F}_2} : 
\widetilde{{\rm Rep}_2(A)}_{u/{\Bbb F}_2} 
\to \widetilde{{\rm Ch}_2(A)}_{u/{\Bbb F}_2}$, 
which is a universal geometric quotient by 
${\rm PGL}_2\otimes_{{\Bbb Z}} R$. 
By a {\it tilde representation with unipotent mold over ${\Bbb F}_2$} for $A$ 
on an $R$-scheme $X$, we understand a pair $(\rho, \lambda)$ of 
a representation $\rho$ of with unipotent mold over ${\Bbb F}_2$ 
for $A$ on $X$ and an $R$-homomorphism $\lambda : A 
\to {\mathcal O}_X(X)$ satisfying the following conditions: 
\begin{enumerate} 
\item $\det(\rho(c)) = \lambda(c)^2$ for each $c \in A$. 
\item\label{cond:sublinebundle-tilde-alg} 
$\{ \rho(c) - \lambda(c) I_2 \mid c \in A \}$ 
spans a sub-line bundle of 
${\mathcal O}_X[\rho(A)]$. 
\end{enumerate} 
As in Proposition \ref{prop:tilderepfunctor},  
we see that $\widetilde{{\rm Rep}_2(A)}_{u/{\Bbb F}_2}$ represents the following 
contravariant functor: 
\[
\begin{array}{ccl}
 ({\bf Sch}/R)^{op} & \to & ({\bf Sets}) \\
 X & \mapsto & 
\left\{ 
\begin{array}{r}  
\mbox{ tilde rep. 
with unipotent  } \\ 
\mbox{ mold over ${\Bbb F}_2$ for $A$ on } X  
\end{array} 
\right\} .
\end{array}
\]  
In a similar way as group or monoid cases, 
we can define generalized tilde representations with unipotent mold 
over ${\Bbb F}_2$ for $A$ 
on an $R$-scheme $X$. The contravariant functor 
$\widetilde{{\mathcal E}q \;\mathcal{U}_2(A)}_{{\Bbb F}_2}$ from the category of $R$-schemes to the 
category of sets is defined as 
\[
\begin{array}{ccccl}
\widetilde{{\mathcal E}q \;\mathcal{U}_2(A)}_{{\Bbb F}_2} & : & ({\bf Sch}/R)^{op} & \to & ({\bf Sets}) \\
 & & X & \mapsto & 
\left\{ 
\begin{array}{r}  
\mbox{ gen. tilde rep. with unip.  
 } \\ 
\mbox{  mold over ${\Bbb F}_2$ for $A$ on } X  
\end{array} 
\right\} 
\Big/\sim.
\end{array}
\] 
We can prove that 
$\widetilde{{\rm Ch}_2(A)}_{u/{\Bbb F}_2}$ is the fine moduli associated to 
$\widetilde{{\mathcal E}q \;\mathcal{U}_2(A)}_{{\Bbb F}_2}$ in the same way as Theorem \ref{th:moduliuf2tilde}.  
The moduli $\widetilde{{\rm Ch}_2(A)}_{u/{\Bbb F}_2}$ is separated over $R$;  
if $A$ is a finitely generated algebra over $R$, then 
$\widetilde{{\rm Ch}_2(A)}_{u/{\Bbb F}_2}$ is of finite type over $R$.  
\end{remark} 

\begin{example}\label{example:purelyinseparablealgebra}\rm 
Let $k$ be a field of characteristic $2$. 
Let $K = k(\alpha)$ be a purely inseparable extension of $k$ of 
degree $2$ with $\beta = \alpha^2 \in k$. 
Regarding $K$ as a $k$-vector space of dimension $2$,  we have 
a $k$-algebra homomorphism $\rho : K \to {\rm End}_k(K) \cong {\rm M}_2(k)$ 
by $c \mapsto (c' \mapsto cc')$, 
which is a representation with unipotent mold over ${\Bbb F}_2$. 
The matrix $\rho(\alpha) \in {\rm M}_2(k)$ has no 
eigenvalue in $k$, but has an eigenvalue after base change $k \to K$. 
This example is an interesting case of representations with unipotent 
mold over ${\Bbb F}_2$. 

The universal representation ${\sigma}_{K, u/{\Bbb F}_2}$ 
with unipotent mold over ${\Bbb F}_2$ on ${\rm Rep}_2(K)_{u/{\Bbb F}_2}$ 
is characterized by $\sigma_{K, u/{\Bbb F}_2}(\alpha)$ because 
it is a $k$-algebra homomorphism and $K=k(\alpha)$. 
Since ${\rm tr}(\sigma_{K, u/{\Bbb F}_2}(\alpha)) = 0$ and 
$\sigma_{K, u/{\Bbb F}_2}(\alpha)^2 = \beta I_2$, 
we can write ${\rm Rep}_2(K)_{u/{\Bbb F}_2} = D(b) \cup D(c) 
\subset {\rm Spec} \: k[a, b, c]/(a^2+bc+\beta)$ 
and 
$\sigma_{K, u/{\Bbb F}_2}(\alpha) = 
\left( 
\begin{array}{cc}
a & b \\
c & a \\
\end{array}
\right)$. 
Then $\widetilde{{\rm Rep}_2(K)}_{u/{\Bbb F}_2} = D(b) \cup D(c) 
\subset {\rm Spec} \: k[a, b, c, X]/(a^2+bc+\beta, X^2 - \beta)$. 
Identifying $K$ with $k[X]/(X^2 - \beta)$, we have 
$\widetilde{{\rm Rep}_2(K)}_{u/{\Bbb F}_2} = D(b) \cup D(c) 
\subset {\rm Spec} \: K[a, b, c]/(a^2+bc+\beta)$. 
In particular, $\widetilde{{\rm Rep}_2(K)}_{u/{\Bbb F}_2} 
= {\rm Rep}_2(K)_{u/{\Bbb F}_2}\otimes_{k} K$. 
Remark that ${\rm Rep}_2(K)_{u/{\Bbb F}_2} = 
{\rm Rep}_2(K)_{u/{\Bbb F}_2, \alpha}$ and that 
$\widetilde{{\rm Rep}_2(K)}_{u/{\Bbb F}_2} = 
\widetilde{{\rm Rep}_2(K)}_{u/{\Bbb F}_2, \alpha}$. 

Let us use the same notation in Remark~\ref{remark:constructuf2}.  
The universal character $d'_{K}$ on ${\rm Rep}'_1(K)$ is 
characterized by $d'_{K}(\alpha)$, and it satisfies $d'_{K}(\alpha)^2 = \beta^2$. 
Hence we can write ${\rm Rep}'_1(K) = {\rm Spec} k[x]/(x^2 - \beta^2)$ and 
$d'_{K}(\alpha)=x$. 
The universal $(a, b)$-coefficients with respect to $(\alpha, d'_{K})$ on 
${\rm Ch}_2(K)_{u/{\Bbb F}_2, \alpha}$ satisfies 
$a(1)=1, b(1)=0, a(\alpha)=0$, and $b(\alpha)=1$. 
By the condition $\beta = 
a(\beta) = a(\alpha^2) = a(\alpha)^2+b(\alpha)^2 d'_{K}(\alpha)$,    
we have $x = \beta$. 
Thus, we see that ${\rm Ch}_2(K)_{u/{\Bbb F}_2, \alpha} = {\rm Spec} \: k[x]/(x-\beta) = {\rm Spec} \: k$. 
On the other hand, 
${\rm Rep}_1(K) = {\rm Spec} \: k[x]/(x^2 -\beta) = 
{\rm Spec} K$ because the universal character $d_{K}$ on 
${\rm Rep}_1(K)$ satisfies $d_{K}(\alpha)^2 = \beta$. 
The $A_1(K)$-module $\Omega_{K/k}$ introduced in  Remark~\ref{remark:widetildeuf2algebra} is isomorphic to the free module 
$A_{1}(K) d\alpha = K d\alpha$.  
Hence $\widetilde{{\rm Ch}_2(K)}_{u/{\Bbb F}_2} = {\rm Proj} \: S(A_1(K) d\alpha)  
= {\rm Rep}_1(K) = {\rm Spec} \: K$. 
Therefore, the commutative diagram 
\[
\begin{array}{ccr}
\phantom{AAAAA} \widetilde{{\rm Rep}_2(K)}_{u/{\Bbb F}_2} 
& \stackrel{p_{}}{\to} & {\rm Rep}_2(K)_{u/{\Bbb F}_2} \\
\widetilde{\pi}_{K, u/{\Bbb F}_2} \downarrow & & \downarrow 
\pi_{K, u/{\Bbb F}_2} \\
\phantom{AAAAA} \widetilde{{\rm Ch}_2(K)}_{u/{\Bbb F}_2} 
& \stackrel{q_{}}{\to} & {\rm Ch}_2(K)_{u/{\Bbb F}_2}   
\end{array}  
\] 
is identified with 
\[
\begin{array}{ccc}
{\rm Rep}_2(K)_{u/{\Bbb F}_2}\otimes_{k} K 
& \to & {\rm Rep}_2(K)_{u/{\Bbb F}_2} \\
 \downarrow & & \downarrow  \\
{\rm Spec} \: K & \to & {\rm Spec} \: k,   
\end{array}  
\] 
which gives a fibre product. 
The representation $\rho : K \to {\rm End}_k(K)$ gives the only  
equivalence class of $2$-dimensional 
representations of $K$ over $k$ which have   
unipotent molds over ${\Bbb F}_2$.  
\end{example}

\begin{remark}\rm
For understanding the difference between the moduli schemes 
$\widetilde{{\rm Ch}_2(\Gamma)}_{u/{\Bbb F}_2}$ 
and ${\rm Ch}_2(\Gamma)_{u/{\Bbb F}_2}$, let us pay attention to  
the morphism $q : \widetilde{{\rm Ch}_2(\Gamma)}_{u/{\Bbb F}_2} 
\to {\rm Ch}_2(\Gamma)_{u/{\Bbb F}_2}$. 
Assume that a point $\widetilde{x}$ of 
$\widetilde{{\rm Ch}_2(\Gamma)}_{u/{\Bbb F}_2}$ is mapped to a point 
$x$ of ${\rm Ch}_2(\Gamma)_{u/{\Bbb F}_2}$. We can write $\widetilde{x} = [ (\rho, \lambda) ]$ and $x = [\rho]$, 
where $\rho : \Gamma \to {\rm M}_2(k(x))$ is a representation 
with unipotent mold over ${\Bbb F}_2$ on the residue field $k(x)$ of $x$ 
and $\lambda : \Gamma \to k(\widetilde{x})$ 
is a character on the residue field $k(\widetilde{x})$ of $\widetilde{x}$ 
such that $(\rho, \lambda)$ is a tilde representation 
with unipotent mold over 
${\Bbb F}_2$. 
It is easy to see that $\widetilde{{\rm Ch}_2(\Gamma)}_{u/{\Bbb F}_2}(K)  
\to {\rm Ch}_2(\Gamma)_{u/{\Bbb F}_2}(K)$ is injective for any field $K$. 
Hence $q$ is universally injective (or radical) (see \cite[Definition~3.5.4]{EGA1}). 
Then $k(\widetilde{x})$ is a purely inseparable extension of $k(x)$. 
In that meaning, 
$q$ is a generalization of purely inseparable extension of fields and it is globally defined. 

Although $q : \widetilde{{\rm Ch}_2(\Gamma)}_{u/{\Bbb F}_2} 
\to {\rm Ch}_2(\Gamma)_{u/{\Bbb F}_2}$ is surjective, 
$\widetilde{{\rm Ch}_2(\Gamma)}_{u/{\Bbb F}_2}(K)  
\to {\rm Ch}_2(\Gamma)_{u/{\Bbb F}_2}(K)$ is not surjective in general. 
In the free monoid case, $q : \widetilde{C(m)} \to C(m)$ is described in 
Example \ref{example:freemonoidtildemoduli}.  When $m =1$,  
$q : {\Bbb A}^1_{{\Bbb F}_2} \to {\Bbb A}^1_{{\Bbb F}_2}$ is given by $r_1 \mapsto r_1^2$, 
where $r_1 = r(\alpha_1)$ and $\alpha_1$ is 
the generator of the free monoid $\Upsilon_1 = \langle \alpha_1 \rangle$.  
Let $\beta$ be an element of 
a field $k$ of characteristic $2$ such that 
$\alpha = \sqrt{\beta} \notin k$. 
Let $x$ be the $k$-rational point of $C(1) \cong {\Bbb A}^1_{{\Bbb F}_2}$ given by  $r_1^2 = r(\alpha_1)^2= \beta \in k$, 
and let $\widetilde{x}$ be the $k(\alpha)$-rational point 
of $\widetilde{C(1)} \cong {\Bbb A}^1_{{\Bbb F}_2}$ given by  $r_1 = r(\alpha_1) = 
\alpha \in k(\alpha)$. Then $\widetilde{x}$ corresponds to $x$ and 
$k(\alpha)$ is a purely inseparable extension of degree $2$ over $k$ 
({\it cf.} Example \ref{example:purelyinseparablealgebra}).  
In particular, $\widetilde{C(1)}(k) \to C(1)(k)$ is not surjective, since 
$x$ is not contained in the image. 
Remark that if $\Gamma$ is finitely generated, then the residue field 
$k(x)$ of a closed point $x$ of ${\rm Ch}_2(\Gamma)_{u/{\Bbb F}_2}$ is a finite field. 
In this case, $k(\widetilde{x}) = k(x)$ for the unique point $\widetilde{x}$ lying over $x$. 
Note that $q : \widetilde{C(m)} \to C(m)$ induces a bijection of sets 
$\widetilde{C(m)}(K) \cong C(m)(K)$ 
if $K$ is an algebraically closed field of characteristic 
$2$ and that $q$ induces a purely inseparable 
extension of the function fields of degree $2$ (see also \cite[Remark 3.3]{topos2}). 
\end{remark}

\begin{remark}\rm
We have introduced the notion of 
generalized tilde representations with unipotent mold over ${\Bbb F}_2$ 
for describing the moduli functors 
$\widetilde{{\mathcal E}q\mathcal{U}_2(\Gamma)}_{{\Bbb F}_2}$  
and $\widetilde{{\mathcal E}q\mathcal{U}_2(A)}_{{\Bbb F}_2}$.  However, 
the moduli functors can also be described as 
$\widetilde{{\mathcal E}q\mathcal{U}'_2(\Gamma)}_{{\Bbb F}_2}$  
and $\widetilde{{\mathcal E}q\mathcal{U}'_2(A)}_{{\Bbb F}_2}$ by using 
the notion of tilde representations generating sheaves of algebras 
which define unipotent molds over ${\Bbb F}_2$.    
More precisely, see \S 8. 
\end{remark}

\section{Representations in sheaves of algebras} 

For describing the moduli functor ${\mathcal E}q\mathcal{SS}_2(\Gamma)$ 
(${\mathcal E}q\mathcal{U}_2(\Gamma)$, or 
${\mathcal E}q\mathcal{U}_2(\Gamma)_{{\Bbb F}_2}$), 
we introduced the notion of generalized representations with 
semi-simple mold (unipotent mold, unipotent mold over ${\Bbb F}_2$, respectively) 
in \S 4-\S 6. However, we can also formulate these moduli functors  
by using representations generating sheaves of 
${\mathcal O}_X$-algebras which define molds of rank $2$ on 
a scheme $X$. In this section, we 
discuss this formulation for describing the moduli schemes.  
We also reformulate $\widetilde{{\mathcal E}q\mathcal{U}_2(\Gamma)}_{{\Bbb F}_2}$ 
in \S 7 by using tilde representations generating sheaves 
of ${\mathcal O}_X$-algebras which 
define unipotent molds over ${\Bbb F}_2$.  
The author needs to say that this section was inspired by the referee.  

\begin{definition}\rm 
Let $\Gamma$ be a group or a monoid. 
Let ${\mathcal A}$ be a sheaf of ${\mathcal O}_X$-algebras 
on a scheme $X$.  
We say that a homomorphism $\rho : \Gamma \to {\mathcal A}(X)$ is 
a {\it representation in} ${\mathcal A}$ of $\Gamma$. 
For two representations $\rho_1 : \Gamma \to {\mathcal A}_1(X)$ and 
$\rho_2 : \Gamma \to {\mathcal A}_2(X)$, we say that 
$\rho_1$ and $\rho_2$ are {\it equivalent} if there exists an 
isomorphism $\phi : {\mathcal A}_1 \to {\mathcal A}_2$ 
as sheaves of ${\mathcal O}_X$-algebras such that 
$\phi \circ \rho_1 = \rho_2$. 
We call a representation $\rho : \Gamma \to {\mathcal A}(X)$ 
a {\it representation generating ${\mathcal A}$} if 
${\mathcal O}_X[\rho(\Gamma)] = {\mathcal A}$.

\end{definition}

Let ${\mathcal A}$ be a sheaf of ${\mathcal O}_X$-algebras 
on a scheme $X$ which is locally free of rank $2$. 
We define $\Phi_{{\mathcal A}} : 
{\mathcal A} \to {\rm End}_{{\mathcal O}_X} ({\mathcal A})$ 
by $a \mapsto (b \mapsto ab)$ for each open subset $U$ of $X$ and for each 
$a, b \in {\mathcal A}(U)$. Then $\Phi_{{\mathcal A}}$ is injective. 
Remark that ${\mathcal A}(U)$ is a commutative ring since ${\mathcal A}$ is 
locally free of rank $2$ and $1 \in {\mathcal A}(X)$.  
For each $x \in X$, choose a neighbourhood $U_x$ of $x$ such 
that ${\mathcal A}\!\mid_{U_{x}} \cong {\mathcal O}_{U_{x}}^{\oplus 2}$.  
By considering the inclusion  
${\mathcal A}\!\mid_{U_x} 
\cong \Phi_{{\mathcal A}}\!\mid_{U_x}\!\!({\mathcal A}\!\mid_{U_x}) 
\subset {\rm End}_{{\mathcal O}_{U_x}}({\mathcal A}\!\mid_{U_x})   
\cong {\rm M}_2({\mathcal O}_{U_x})$, 
we obtain a mold of rank $2$ on $U_x$. 

\begin{definition}\rm 
If $\Phi_{{\mathcal A}}\!\mid_{U_x}\!\!({\mathcal A}\!\mid_{U_x})$ is 
a semi-simple mold (unipotent mold, or unipotent mold over ${\Bbb F}_2$) 
for each $x \in X$, we say that ${\mathcal A}$ {\it defines a 
semi-simple mold (unipotent mold, or unipotent mold over ${\Bbb F}_2$, 
respectively)}. 
This definition does not depend  
on choices of neighbourhoods  $U_x$ of $x$ and 
isomorphisms ${\mathcal A}\!\mid_{U_{x}} \cong {\mathcal O}_{U_{x}}^{\oplus 2}$.  
\end{definition} 

For a generalized representation $\{ (U_i, \rho_i) \}_{i \in I}$ 
with semi-simple mold (unipotent mold, unipotent mold over ${\Bbb F}_2$, respectively) 
of $\Gamma$ on a scheme $X$, we define a sheaf ${\mathcal A}$ 
of ${\mathcal O}_X$-algebras 
which is a locally free sheaf of rank $2$ as follows:  
Set ${\mathcal A}_i := {\mathcal O}_{U_i}[\rho_i(\Gamma)]$.  
Let us define an isomorphism $\varphi_{ij} : {\mathcal A}_{i} \!\mid_{U_i \cap U_j} 
\to {\mathcal A}_{j} \!\mid_{U_i \cap U_j}$ by 
$\rho_i(\gamma) \mapsto \rho_j(\gamma)$ for each $\gamma \in \Gamma$. 
It is easy to check that $\varphi_{ik} = \varphi_{jk} \circ \varphi_{ij} : 
{\mathcal A}_{i} \!\mid_{U_i \cap U_j \cap U_k} 
\to {\mathcal A}_{k} \!\mid_{U_i \cap U_j \cap U_k}$ and 
that $\varphi_{ij} = \varphi_{ji}^{-1}$ and $\varphi_{ii} = {\rm id}$.  
Hence by using \cite[Chap.~II, Ex.1.22]{Hartshorne},  
we obtain a unique sheaf ${\mathcal A}$ 
of  ${\mathcal O}_X$-algebras on $X$ (up to isomorphism), together with isomorphisms 
$\psi_i : {\mathcal A}\!\mid_{U_i} \stackrel{\sim}{\to} {\mathcal A}_i$ 
such that $\psi_{j} = \varphi_{ij} \circ \psi_{i}$ on $U_i \cap U_j$ for each 
$i, j$.  Obviously, ${\mathcal A}$ is locally free of rank $2$. 
If $\{ (U_i, \rho_i) \}_{i \in I}$ is a generalized representation with 
semi-simple mold (unipotent mold, or unipotent mold over ${\Bbb F}_2$), 
then ${\mathcal A}$ defines  a semi-simple mold 
(unipotent mold, or unipotent mold over ${\Bbb F}_2$, 
respectively). 
By gluing $\{ \rho_i : \Gamma \to {\rm M}_2({\mathcal O}_{U_i}(U_i)) \}_{i \in I}$,  
we have a representation $\rho : \Gamma \to {\mathcal A}(X)$ in ${\mathcal A}$ 
of $\Gamma$ such that  $\Gamma \stackrel{\rho}{\to} {\mathcal A}(X) 
\stackrel{{\rm res}}{\to} {\mathcal A}(U_i) \stackrel{\psi_i}{\to} {\mathcal A}_i(U_i)$ 
coincides with $\rho_i$ for each $i \in I$.  
Then $\rho$ is a representation generating ${\mathcal A}$. 
Thus each generalized representation $\{ (U_i, \rho_i) \}_{i \in I}$ 
corresponds to a representation $\rho : \Gamma 
\to {\mathcal A}(X)$ of $\Gamma$ generating ${\mathcal A}$  
on $X$. 

\begin{definition}\rm 
Let us define the contravariant functor ${\mathcal E}q\mathcal{SS}_2'(\Gamma)$ 
from the category of schemes to the category of sets as follows: 
\[
\begin{array}{cccl}
{\mathcal E}q\mathcal{SS}_2'(\Gamma)  : & ({\bf Sch})^{op} & \to & ({\bf Sets}) \\
  & X & \mapsto & 
\left\{ 
\begin{array}{l}  
\mbox{ a representation of $\Gamma$ generating  } \\ 
\mbox{ a sheaf of ${\mathcal O}_X$-algebras   
${\mathcal A}$ which   }  \\
\mbox{ is locally free of rank $2$ and   }   \\
\mbox{ defines a semi-simple mold on } X
\end{array} 
\right\} 
\Big/\sim.
\end{array}
\] 
We also define the contravariant functors ${\mathcal E}q\mathcal{U}_2'(\Gamma) : ({\bf Sch}/{\Bbb Z}[1/2])^{op}  \to  ({\bf Sets})$ and 
${\mathcal E}q\mathcal{U}_2'(\Gamma)_{{\Bbb F}_2} : 
({\bf Sch}/{\Bbb F}_2)^{op}  \to  ({\bf Sets})$ in the same way. 
\end{definition} 
  
By the correspondence above, we obtain natural transformations  
$\sigma_{s.s.} :  {\mathcal E}q\mathcal{SS}_2(\Gamma) \to {\mathcal E}q\mathcal{SS}_2'(\Gamma)$, 
$\sigma_{u} :  {\mathcal E}q\mathcal{U}_2(\Gamma) \to {\mathcal E}q\mathcal{U}_2'(\Gamma)$, and 
$\sigma_{u/{\Bbb F}_2} :  {\mathcal E}q\mathcal{U}_2(\Gamma)_{{\Bbb F}_2} \to {\mathcal E}q\mathcal{U}_2'(\Gamma)_{{\Bbb F}_2}$.  

Let us define natural transformations 
$\tau_{s.s.} :  {\mathcal E}q\mathcal{SS}_2'(\Gamma) \to {\mathcal E}q\mathcal{SS}_2(\Gamma)$, 
$\tau_{u} :  {\mathcal E}q\mathcal{U}_2'(\Gamma) \to {\mathcal E}q\mathcal{U}_2(\Gamma)$, and 
$\tau_{u/{\Bbb F}_2} :  {\mathcal E}q\mathcal{U}_2'(\Gamma)_{{\Bbb F}_2} \to {\mathcal E}q\mathcal{U}_2(\Gamma)_{{\Bbb F}_2}$ in the following way. 
Let $\rho$ be a representation of $\Gamma$ generating ${\mathcal A}$ on a scheme. 
Assume that ${\mathcal A}$ defines a semi-simple mold, 
unipotent mold, or unipotent mold over ${\Bbb F}_2$ on $X$.  
For each $x \in X$, choose a neighbourhood $U_x$ such that 
${\mathcal A}\!\!\mid_{U_x} \cong {\mathcal O}_{U_{x}}^{\oplus 2}$.  
Then by considering ${\mathcal A}\!\!\mid_{U_x} 
\cong \Phi_{{\mathcal A}}\!\mid_{U_x}\!\!({\mathcal A}\!\mid_{U_x}) \subset 
{\rm End}_{{\mathcal O}_{U_x}}({\mathcal A}\!\mid_{U_x}) 
\cong {\rm M}_2({\mathcal O}_{U_x})$, 
we have a representation $\rho_{x} : \Gamma \to 
{\rm M}_2({\mathcal O}_{U_x})$ with the corresponding mold on $U_x$. 
It is easy to check that $\{ (U_{x}, \rho_x ) \}_{x \in X}$ is 
a generalized representation  
with the corresponding mold on $X$ 
and that the equivalence class of $\{ (U_{x}, \rho_x ) \}_{x \in X}$ is 
well-defined.  
The equivalence class of $\{ (U_{x}, \rho_x ) \}_{x \in X}$ 
does not depend on choosing a representative of the equivalence 
class of  $\rho : X \to {\mathcal A}(X)$.     
This correspondence defines $\tau_{s.s.}$, $\tau_{u}$, and 
$\tau_{u/{\Bbb F}_2}$. 

It is straightforward to verify that 
$\tau_{s.s.} \circ \sigma_{s.s.} = 1_{{\mathcal E}q\mathcal{SS}_2(\Gamma)}$ 
and $\sigma_{s.s.} \circ \tau_{s.s.}= 1_{{\mathcal E}q\mathcal{SS}_2'(\Gamma)}$ and 
so on.   Hence we can obtain the following:  
\begin{proposition}  
There are canonical isomorphisms: 
\begin{eqnarray*}
{\mathcal E}q\mathcal{SS}_2(\Gamma) & \cong & {\mathcal E}q\mathcal{SS}_2'(\Gamma), \\
{\mathcal E}q\mathcal{U}_2(\Gamma) & \cong & 
{\mathcal E}q\mathcal{U}_2'(\Gamma), \\
{\mathcal E}q\mathcal{U}_2(\Gamma)_{{\Bbb F}_2} & \cong & {\mathcal E}q\mathcal{U}_2'(\Gamma)_{{\Bbb F}_2}. 
\end{eqnarray*} 
In particular, 
${\rm Ch}_2(\Gamma)_{s.s.}$, ${\rm Ch}_2(\Gamma)_{u}$, and 
${\rm Ch}_2(\Gamma)_{u/{\Bbb F}_2}$ represent 
${\mathcal E}q\mathcal{SS}_2'(\Gamma)$, 
${\mathcal E}q\mathcal{U}_2'(\Gamma)$, and 
${\mathcal E}q\mathcal{U}_2'(\Gamma)_{{\Bbb F}_2}$, respectively.
\end{proposition} 

In the case of  unipotent molds over ${\Bbb F}_2$, we can define 
another functor ${\mathcal E}q\mathcal{U}_2''(\Gamma)_{{\Bbb F}_2}$:    

\begin{definition}\rm 
Let ${\mathcal A}$ be a sheaf of ${\mathcal O}_X$-algebras 
which is locally free of rank $2$ on a scheme $X$ over ${\Bbb F}_2$. 
We say that $a \in {\mathcal A}(X)$ is {\it scalar} if  
there exists $f \in {\mathcal O}_X(X)$ such that 
$a = f \cdot 1_{{\mathcal A}}$. We define ${\mathcal E}q\mathcal{U}_2''(\Gamma)_{{\Bbb F}_2}$ by 
\[
\begin{array}{cccl}
{\mathcal E}q\mathcal{U}_2''(\Gamma)_{{\Bbb F}_2}  : & 
({\bf Sch}/{\Bbb F}_2)^{op} & \to & ({\bf Sets}) \\
  & X & \mapsto & 
\left\{ 
\begin{array}{l}  
\mbox{ a representation $\rho$ of $\Gamma$ } \\ 
\mbox{ generating  a sheaf of  }  \\ 
\mbox{ ${\mathcal O}_X$-algebras 
${\mathcal A}$ which is } \\
\mbox{ locally free of rank $2$ on } X   \\
\mbox{ such that $\rho(\gamma)^2$ is scalar } \\
\mbox{ for each $\gamma \in \Gamma$  } 
\end{array} 
\right\} 
\Big/\sim.
\end{array}
\] 
\end{definition} 

\begin{proposition}
There are canonical isomorphisms  
\[
{\mathcal E}q\mathcal{U}_2(\Gamma)_{{\Bbb F}_2}  \cong  {\mathcal E}q\mathcal{U}_2'(\Gamma)_{{\Bbb F}_2} 
\cong {\mathcal E}q\mathcal{U}_2''(\Gamma)_{{\Bbb F}_2}. 
\] 
\end{proposition}

{\it Proof}.
It suffices to prove that ${\mathcal E}q\mathcal{U}_2'(\Gamma)_{{\Bbb F}_2} 
\cong {\mathcal E}q\mathcal{U}_2''(\Gamma)_{{\Bbb F}_2}$. 
For $[\rho : \Gamma \to {\mathcal A}(X)] \in {\mathcal E}q\mathcal{U}_2'(\Gamma)_{{\Bbb F}_2}(X)$ 
with a scheme $X$ over ${\Bbb F}_2$, 
$\rho(\gamma)^2 = \det(\rho(\gamma)) \cdot 1_{{\mathcal A}}$ is scalar 
for each $\gamma \in \Gamma$. Hence $[\rho : \Gamma \to {\mathcal A}(X)] \in {\mathcal E}q\mathcal{U}_2''(\Gamma)_{{\Bbb F}_2}(X)$. 
Conversely, let $[\rho : \Gamma \to {\mathcal A}(X)] \in {\mathcal E}q\mathcal{U}_2''(\Gamma)_{{\Bbb F}_2}(X)$.  
For $x \in X$, there exist $\gamma \in \Gamma$ and a neighbourhood 
$U$ of $x$ such that ${\mathcal A}\!\mid_{U} \cong 
{\mathcal O}_U\cdot 1_{U} \oplus {\mathcal O}_U\cdot \rho(\gamma)$. 
Since $\rho(\gamma)^2$ is scalar, $\rho(\gamma)^2 = c \cdot 1_{{\mathcal A}}$ 
for some $c \in {\mathcal O}_U(U)$. By the Cayley-Hamilton theorem, 
$\rho(\gamma)^2 - {\rm tr}(\rho(\gamma)) \rho(\gamma) + \det(\rho(\gamma)) I_2
=0$ on $U$. Thus we have ${\rm tr}(\rho(\gamma)) = 0$ and 
$\det(\rho(\gamma)) = c$ on $U$. 
For any $\gamma' \in \Gamma$, $\rho(\gamma') = a I_2 + b \rho(\gamma)$ 
on $U$ for some $a, b \in {\mathcal O}_U(U)$. 
This implies that ${\rm tr}(\rho(\gamma')) = 
a {\rm tr}(I_2) + b {\rm tr}(\rho(\gamma)) = 0$.  Hence  
${\mathcal A}$ defines a unipotent mold over ${\Bbb F}_2$ and that 
$[\rho : \Gamma \to {\mathcal A}(X)] \in {\mathcal E}q\mathcal{U}_2'(\Gamma)_{{\Bbb F}_2}(X)$. Therefore we have proved that 
${\mathcal E}q\mathcal{U}_2'(\Gamma)_{{\Bbb F}_2} 
\cong {\mathcal E}q\mathcal{U}_2''(\Gamma)_{{\Bbb F}_2}$. 
\qed 

\bigskip 

Let ${\mathcal A}$ be a sheaf of ${\mathcal O}_X$-algebras 
which is locally free of rank $2$ on a scheme $X$ over ${\Bbb F}_2$. 
Let $\rho : \Gamma \to {\mathcal A}(X)$ be a representation of 
$\Gamma$ generating ${\mathcal A}$, and let $\chi : \Gamma 
\to {\mathcal O}_X(X)$ be a character. 
We say that a pair $(\rho, \chi)$ is a {\it tilde representation with 
unipotent mold over ${\Bbb F}_2$ for $\Gamma$ 
generating ${\mathcal A}$ on $X$} if  
$\{ \rho(\gamma) - \chi(\gamma) \cdot 1_{{\mathcal A}} \}_{\gamma \in \Gamma}$ 
spans a sub-line bundle 
of ${\mathcal A}$ and $(\rho(\gamma) - \chi(\gamma) \cdot 1_{{\mathcal A}})^2 = 0$ 
for any $\gamma \in \Gamma$. 
(Then we can prove that ${\mathcal A}$ defines a unipotent mold over 
${\Bbb F}_2$ as in the proof of Proposition \ref{prop:eqtildeprime}.) 
For two tilde representations $(\rho_1, \chi_1)$ and $(\rho_2, \chi_2)$ 
with unipotent mold over ${\Bbb F}_2$ for $\Gamma$ on $X$, we say that 
they are {\it equivalent} if $\chi_1 = \chi_2$ and 
there exists an isomorphism $\phi : {\mathcal A}_1 \to 
{\mathcal A}_2$ as sheaves of  ${\mathcal O}_X$-algebras such that 
$\phi \circ \rho_1 = \rho_2$, where $\rho_i$ is a homomorphism $\rho_i : \Gamma \to 
{\mathcal A}_i(X)$ for $i=1, 2$.

\begin{definition}\rm 
\[
\begin{array}{cccl}
\widetilde{{\mathcal E}q\mathcal{U}_2'(\Gamma)_{{\Bbb F}_2}}  : 
& ({\bf Sch}/{\Bbb F}_2)^{op} & \to & ({\bf Sets}) \\
  & X & \mapsto & 
\left\{ (\rho, \chi) \: 
\begin{array}{|l}  
\mbox{ $(\rho, \chi)$ is 
a tilde representation   } \\ 
\mbox{ with unipotent mold over ${\Bbb F}_2$ } \\ 
\mbox{ for $\Gamma$ generating a sheaf of } \\ 
\mbox{  ${\mathcal O}_X$-algebras  
${\mathcal A}$    }  \\
\end{array} 
\right\} 
\Big/\sim.
\end{array}
\] 
\end{definition} 

\begin{proposition}\label{prop:eqtildeprime}
There is a canonical isomorphism  
\[
\widetilde{{\mathcal E}q\mathcal{U}_2(\Gamma)}_{{\Bbb F}_2}  \cong  
\widetilde{{\mathcal E}q\mathcal{U}_2'(\Gamma)}_{{\Bbb F}_2}. 
\] 
\end{proposition}

{\it Proof}.
Let $\{ (U_i, \rho_i, \lambda_i) \}_{i \in I} 
\in \widetilde{{\mathcal E}q\mathcal{U}_2(\Gamma)}_{{\Bbb F}_2}(X)$  
be a generalized tilde representation  
with unipotent mold over ${\Bbb F}_2$ for $\Gamma$ on 
a scheme $X$ over ${\Bbb F}_2$.  
Since $\{ (U_i, \rho_i) \}_{i \in I} 
\in {\mathcal E}q\mathcal{U}_2(\Gamma)_{{\Bbb F}_2}(X)$, 
we have 
$[\rho : \Gamma \to {\mathcal A}(X)] \in 
{\mathcal E}q\mathcal{U}_2'(\Gamma)_{{\Bbb F}_2}(X)$ 
by ${\mathcal E}q\mathcal{U}_2(\Gamma)_{{\Bbb F}_2}(X)  \cong  
{\mathcal E}q\mathcal{U}_2'(\Gamma)_{{\Bbb F}_2}(X)$. 
By gluing $\{ (U_i, \lambda_i) \}_{i \in I}$,  
we can define a character $\chi : \Gamma \to 
{\mathcal O}_X(X)$ such that $\chi\!\!\mid_{U_i} = \lambda_i$ for $i \in I$. 
Note that $\det \rho_i(\gamma) = \chi(\gamma)^2$ on $U_i$. 
It is easy to see that $(\rho, \chi) \in \widetilde{{\mathcal E}q\mathcal{U}_2'(\Gamma)}_{{\Bbb F}_2}(X)$. 
This correspondence induces a natural transformation 
$\widetilde{\sigma_{u/{\Bbb F}_2}} : \widetilde{{\mathcal E}q\mathcal{U}_2(\Gamma)}_{{\Bbb F}_2}  \to  
\widetilde{{\mathcal E}q\mathcal{U}_2'(\Gamma)}_{{\Bbb F}_2}$. 

Conversely, let $(\rho, \chi) \in 
\widetilde{{\mathcal E}q\mathcal{U}_2'(\Gamma)}_{{\Bbb F}_2}(X)$. 
For each point $x \in X$, there exist $\alpha_x \in \Gamma$ and 
a neighbourhood $U_x$ of $x$ such that 
${\mathcal A}\!\mid_{U_x} 
\cong {\mathcal O}_{U_x}\cdot 1_{{\mathcal A}} 
\oplus {\mathcal O}_{U_x}\cdot \rho(\alpha_x)$. 
Denote  $\Gamma \stackrel{\rho}{\to} {\mathcal A}\!\!\mid_{U_x} 
\cong \Phi_{{\mathcal A}}\!\mid_{U_x}\!\!({\mathcal A}\!\mid_{U_x}) \subset 
{\rm End}_{{\mathcal O}_{U_x}}({\mathcal A}\!\mid_{U_x}) 
\cong {\rm M}_2({\mathcal O}_{U_x})$ by $\rho_x$. 
By the assumption, $(\rho_x(\alpha_x)-\chi(\alpha_x)I_2)^2=0$. 
Then $\rho_x(\alpha_x)^2-\chi(\alpha_x)^2I_2 = {\rm tr}(\rho_x(\alpha_x))\rho_x(\alpha_x) 
-\det(\rho_x(\alpha_x))I_2-\chi(\alpha_x)^2I_2 = 0$. Hence  
${\rm tr}(\rho_x(\alpha_x)) = 0$ and $\det(\rho_x(\alpha_x))=\chi(\alpha_x)^2$.  
Since $\{ \rho(\gamma) - \chi(\gamma) \cdot 1_{{\mathcal A}} \}_{\gamma \in \Gamma}$ spans a sub-line bundle of ${\mathcal A}$, for each $\gamma \in \Gamma$  
there exists $c \in {\mathcal O}_X(U_x)$ such that 
$\rho(\gamma) - \chi(\gamma)\cdot 1_{{\mathcal A}} 
= c(\rho(\alpha_x) - \chi(\alpha_x)\cdot 1_{{\mathcal A}})$ on $U_x$. 
We have $\rho_x(\gamma) = (\chi(\gamma)-c\chi(\alpha_x))I_2 +c\rho_x(\alpha_x)$. 
Putting $a(\gamma) = \chi(\gamma)-c\chi(\alpha_x)$ and 
$b(\gamma)= c$, we obtain $\rho_x(\gamma) = a(\gamma)I_2+b(\gamma)\rho_x(\alpha_x)$ and $\chi(\gamma)=a(\gamma)+b(\gamma)\chi(\alpha_x)$. 
Thereby ${\rm tr}(\rho_{x}(\gamma))=0$ and $\det \rho_x(\gamma) = \chi(\gamma)^2$ 
for each $\gamma \in \Gamma$. It is easy to check that ${\mathcal A}$ defines 
a unipotent mold over ${\Bbb F}_2$ and that 
$\{ (U_x, \rho_x, \chi\!\!\mid_{U_x}) \}_{x \in X} \in \widetilde{{\mathcal E}q\mathcal{U}_2(\Gamma)}_{{\Bbb F}_2}(X)$. 
Therefore this correspondence induces a natural transformation 
$\widetilde{\tau_{u/{\Bbb F}_2}} : 
\widetilde{{\mathcal E}q\mathcal{U}_2'(\Gamma)}_{{\Bbb F}_2}  \to  
\widetilde{{\mathcal E}q\mathcal{U}_2(\Gamma)}_{{\Bbb F}_2}$. 

It is easy to see that $\widetilde{\tau_{u/{\Bbb F}_2}} \circ 
\widetilde{\sigma_{u/{\Bbb F}_2}} = 1_{\widetilde{{\mathcal E}q\mathcal{U}_2(\Gamma)}_{{\Bbb F}_2}}$ and that 
$\widetilde{\sigma_{u/{\Bbb F}_2}} \circ 
\widetilde{\tau_{u/{\Bbb F}_2}} = 1_{\widetilde{{\mathcal E}q\mathcal{U}_2'(\Gamma)}_{{\Bbb F}_2}}$. This completes the proof. 
\qed

\bigskip 

In the case of representations of an associative algebra over a 
commutative ring, we have similar results as the group or monoid cases. 

\begin{definition}\rm 
Let $A$ be an associative algebra over a commutative ring $R$. 
Let ${\mathcal A}$ be a sheaf of ${\mathcal O}_X$-algebras on 
an $R$-scheme $X$. 
We say that an $R$-homomorphism $\rho : A \to {\mathcal A}(X)$ is 
a {\it representation in} ${\mathcal A}$ of $A$. 
For two representations $\rho_1 : A \to {\mathcal A}_1(X)$ and 
$\rho_2 : A \to {\mathcal A}_2(X)$, we say that 
$\rho_1$ and $\rho_2$ are {\it equivalent} if there exists an 
isomorphism $\phi : {\mathcal A}_1 \to {\mathcal A}_2$
as sheaves of ${\mathcal O}_X$-algebras such that 
$\phi \circ \rho_1 = \rho_2$. 
We call a representation $\rho : A \to {\mathcal A}(X)$ 
a {\it representation generating ${\mathcal A}$} if 
${\mathcal O}_X[\rho(A)] = {\mathcal A}$. 
\end{definition} 

In the same way as group or monoid cases, 
we define ${\mathcal E}q\mathcal{SS}_2'(A)$, 
${\mathcal E}q\mathcal{U}_2'(A)$, and 
${\mathcal E}q\mathcal{U}_2'(A)_{{\Bbb F}_2}$. 
Similarly, we have 

\begin{proposition} 
There are canonical isomorphisms: 
\begin{eqnarray*}
{\mathcal E}q\mathcal{SS}_2(A) & \cong & {\mathcal E}q\mathcal{SS}_2'(A), \\
{\mathcal E}q\mathcal{U}_2(A) & \cong & {\mathcal E}q\mathcal{U}_2'(A), \\
{\mathcal E}q\mathcal{U}_2(A)_{{\Bbb F}_2} & \cong & {\mathcal E}q\mathcal{U}_2'(A)_{{\Bbb F}_2}. 
\end{eqnarray*} 
\end{proposition} 

Hence we can conclude that 
${\rm Ch}_2(A)_{s.s.}$, ${\rm Ch}_2(A)_{u}$, and 
${\rm Ch}_2(A)_{u/{\Bbb F}_2}$ represent 
${\mathcal E}q\mathcal{SS}_2'(A)$, 
${\mathcal E}q\mathcal{U}_2'(A)$, and 
${\mathcal E}q\mathcal{U}_2'(A)_{{\Bbb F}_2}$, respectively.

\begin{definition}\rm 
Let ${\mathcal A}$ be a sheaf of ${\mathcal O}_X$-algebras which is 
locally free of rank $2$ on 
an $R$-scheme $X$. 
Let $\rho : A \to {\mathcal A}(X)$ be a representation 
generating ${\mathcal A}$ on $X$, and let 
$\chi : A \to {\mathcal O}_X(X)$ be an $R$-homomorphism.   
We say that $(\rho, \chi)$ is 
a {\it tilde representation with unipotent mold over ${\Bbb F}_2$ for $A$ generating}  
${\mathcal A}$ on $X$ if $\{ \rho(c) - \chi(c) \cdot 1_{{\mathcal A}} \}_{c \in A}$ 
spans a sub-line bundle 
of ${\mathcal A}$ and $(\rho(c) - \chi(c) \cdot 1_{{\mathcal A}})^2 = 0$ 
for any $c \in A$. 
\end{definition} 

We can also define $\widetilde{{\mathcal E}q\mathcal{U}_2'(A)}_{{\Bbb F}_2}$. 
Similarly, we have 

\begin{proposition} 
There are a canonical isomorphism: 
\begin{eqnarray*}
\widetilde{{\mathcal E}q\mathcal{U}_2(A)}_{{\Bbb F}_2} & \cong & 
\widetilde{{\mathcal E}q\mathcal{U}_2'(A)}_{{\Bbb F}_2}. 
\end{eqnarray*} 
In particular, $\widetilde{{\rm Ch}_2(A)}_{u/{\Bbb F}_2}$ represents 
$\widetilde{{\mathcal E}q\mathcal{U}_2'(A)}_{{\Bbb F}_2}$. 
\end{proposition}

\section{Appendix: Discriminants}

In this section we deal with the discriminant locus 
of the representation variety of degree $2$.
The discriminant locus is exactly the subset
consisting of representations which are not
absolutely irreducible.
We describe the absolutely irreducible representation part of 
the representation variety of degree $2$ explicitly ({\it cf.} \cite{kj-Saito93} or 
\cite{kj-Saito96}).

\begin{definition}[\cite{kj-Saito93}, \cite{kj-Saito96}]\rm
Let $R$ be a commutative ring.
For $A, B \in \mat{2}{R}$ we define the {\it discriminant} 
$\Delta(A, B)$ by 
\[\begin{array}{ccl}
\Delta(A, B) & := & \tr(A)^2\det(B) + \tr(B)^2\det(A) 
+ \tr(AB)^2 \\ 
& & \hspace{10ex} - \tr(A)\tr(B)\tr(AB) -4\det(A)\det(B). 
\end{array}
\]
From the definition we see that 
$\Delta(A, B) = \Delta(B, A)$. 
If $A, B \in {\rm GL}_2(R)$, then 
$\Delta(A, B) = \det(A) \det(B) {\rm tr}(ABA^{-1}B^{-1} - I_2)$. 
\end{definition}

\begin{remark}\label{remarktwodeltas}\rm
The discriminant $\Delta(A, B)$ above 
is closely related to the discriminant in 
\cite{Nkmt00}.
For $A_1, A_2, A_3, A_4 \in \mat{2}{R}$ we define the {\it discriminant} 
of degree $2$ in \cite{Nkmt00}
by 
\[
\Delta(A_1, A_2, A_3, A_4) := \det \left(
\begin{array}{cccc}
a(1)_{11} & a(1)_{12} & a(1)_{21} & a(1)_{22} \\
a(2)_{11} & a(2)_{12} & a(2)_{21} & a(2)_{22} \\
a(3)_{11} & a(3)_{12} & a(3)_{21} & a(3)_{22} \\
a(4)_{11} & a(4)_{12} & a(4)_{21} & a(4)_{22} \\
\end{array}
\right),
\]
where 
$A_i=\left(
\begin{array}{cc}
a(i)_{11} & a(i)_{12} \\
a(i)_{21} & a(i)_{22}
\end{array}
\right)$ for $i=1, 2, 3, 4$. 
Then we have 
\[
\Delta(A, B)=-\Delta(I_2, A, B, AB).
\]
Note that $\Delta(A_1, A_2, A_3, A_4) \in R^{\times}$ 
if and only if $\{ A_1, A_2, A_3, A_4 \}$ is an $R$-basis of $\mat{2}{R}$.
\end{remark}

\begin{lemma}\label{lemmaabel}
Let $k$ be a field. Assume that ${\cal A} \subseteq \mat{2}{k}$ is 
a subalgebra over $k$. Then ${\cal A} \neq \mat{2}{k}$ if and only if 
${\cal A}$ is commutative
or there exists a $1$-dimensional 
${\cal A}$-invariant subspace of $k^{2}$.  
\end{lemma}
\prf
The ``if'' part is easy.
We only need to prove the ``only if'' part.
Suppose that ${\cal A} \neq \mat{2}{k}$.
If ${\cal A}$ has no nontrivial invariant subspace, then 
$k^2$ is a simple ${\cal A}$-module.
Since the Jacobson radical 
${\rm Jac}{\cal A}$ is equal to $\cap_{M : \,{\rm simple}}
{\rm Ann}M = 0$, the algebra ${\cal A}$ is semi-simple.
From the Wedderburn Theorem we see that 
${\cal A}$ is a product of the full matrix rings 
over division algebras over $k$.
Because ${\cal A}$ has a faithful $2$-dimensional simple module 
and $\dim {\cal A} \le 3$, the algebra ${\cal A}$
is isomorphic to a quadratic extension 
of $k$. Hence ${\cal A}$ is commutative.
\qedd

\begin{prop}\label{propdelta2}
Let $k$ be a field.
Suppose that $A,B \in \mat{2}{k}$ and ${\cal A}$ is the $k$-subalgebra 
of $\mat{2}{k}$ generated by $A$ and $B$. 
Then the following conditions  are equivalent: 
\begin{enumerate}
\item\label{propitem1} $\Delta(A, B)=0$. 
\item\label{propitem2} ${\cal A} \neq \mat{2}{k}$.
\item\label{propitem3} $A$ and $B$ are commutative or 
$A$ and $B$ have a common invariant subspace 
of dimension $1$.
\end{enumerate}
\end{prop}

\prf 
The part (\ref{propitem2}) $\Leftrightarrow$ (\ref{propitem3}) follows from
Lemma \ref{lemmaabel}.
The part (\ref{propitem2}) $\Rightarrow$ (\ref{propitem1})
follows from that  $\{ I_2, A, B, AB \}$ is not a basis of $\mat{2}{k}$
and that 
$\Delta(A, B)=-\Delta(I_2, A, B, AB)=0$ by Remark \ref{remarktwodeltas}. 
We only have to prove (\ref{propitem1}) $\Rightarrow$ (\ref{propitem2}).

We assume that $\Delta(A, B)=0$.
First, we show that $AB$ and $BA$ can be expressed as 
linear combinations of $\{ I_2, A, B\}$.
Here let us consider the case $AB$.
From $\Delta(A,B)=0$, we see that $\{I_2, A, B, AB\}$
is linearly dependent.
Hence there exist $c_i \in k$ ($1 \le i \le 4$) such that $c_1I_2 + c_2A + c_3B + c_4AB=0$
and some $c_i \neq 0$.
If $c_4 \neq 0$, then the claim is true.
If $c_4 = 0$, then either $c_2 \neq 0$ or $c_3 \neq 0$ holds.
When $c_2 \neq 0$, the matrix 
$A$ is expressed as a linear combination of $I_2$ and
$B$, and hence $AB$ can be expressed as a polynomial of $B$.
By the Cayley-Hamilton Theorem, $AB$ can be expressed as a linear combination 
of $I_2$ and  $B$. 
We can also prove the claim for the $c_3 \neq 0$ case.
Thus we have shown that 
$AB$ can be expressed as a linear combination of $\{ I_2, A, B\}$.
We can also prove the $BA$ case in the same way.

Next, we show that any monomial of $A$ and $B$ can be expressed 
as a linear combination of $\{I_2, A, B\}$. 
This implies that ${\cal A} \neq \mat{2}{k}$.  
We prove the claim by induction on the length of monomials.
The length $0$, $1$ and $2$ cases are true. 
Suppose that the length $n-1$ case is true for $n \ge 3$. 
Let $X$ be a monomial whose length is $n$. 
If $X$ has a subsequence $AB$ or $BA$, then 
$X$ can be reduced to the the length $n-1$ case
from the above claim. 
If $X$ has a subsequence $AA$ or $BB$, then from the Cayley-Hamilton Theorem
we also see that $X$ can be reduced to the length $n-1$ case.
This completes the proof.
\qedd 
 
\begin{corollary}\label{coreq}
Let $k$ be a field. Suppose that ${\cal A}$ is a $k$-subalgebra of 
$\mat{2}{k}$.
Then the following conditions are equivalent: 
\begin{enumerate}
\item\label{item1-1} $\Delta(A, B)=0$ for each $A, B \in {\cal A}$.
\item\label{item1-2} ${\cal A} \neq \mat{2}{k}$, or equivalently the
${\cal A}$-module $k^{2}$ is not absolutely irreducible.
\item\label{item1-3} ${\cal A}$ is commutative or 
${\cal A}$ has an invariant subspace 
of dimension $1$.
\end{enumerate}
\end{corollary}

\prf
The part
(\ref{item1-2}) $\Leftrightarrow$ (\ref{item1-3}) follows from Lemma
\ref{lemmaabel}.
Suppose that (\ref{item1-2}) holds.
For any $A, B \in {\cal A}$, the matrices 
$A$ and $B$ do not generate the full matrix ring
$\mat{2}{k}$, so $\Delta(A, B)=0$ from the above proposition.
This shows that (\ref{item1-2}) $\Rightarrow$ (\ref{item1-1}) holds.

We now prove (\ref{item1-1}) $\Rightarrow$ (\ref{item1-2}).
Assume that (\ref{item1-1}) holds.
Suppose that ${\cal A}=\mat{2}{k}$.
For  
\[A=\left(
\begin{array}{cc}
 1 & 1 \\
 0 & 1
\end{array}
\right),
B=\left(
\begin{array}{cc}
 1 & 0 \\
 1 & 1
\end{array}
\right) \in {\mathcal A}, \] 
we have $\Delta(A, B) = 1 \neq 0$. 
This is a contradiction. 
Therefore we have proved (\ref{item1-1}) $\Rightarrow$ (\ref{item1-2}).
\qedd

\bigskip

Here we introduce another invariant.

\begin{definition}\rm
Let $R$ be a commutative ring. For $A, B, C \in \mat{2}{R}$ we define 
$\tau(A, B, C)$ by
\[\tau(A, B, C):=\tr(ABC)-\tr(ACB) \]
or equivalently, 
\[\tau(A, B, C)=2\tr(ABC)-\tr(A)\tr(BC)-\tr(B)\tr(CA)-\tr(C)\tr(AB)\]
\[ +\tr(A)\tr(B)\tr(C).\]
\end{definition}

\begin{remark}\rm
The above $\tau$ is closely 
related to the discriminant defined in \cite{Nkmt00}.
Indeed, 
$\tau(A,B,C) = \Delta(A, B, C, I_2)$ holds 
for $A, B, C \in \mat{2}{R}$. 
\end{remark}

\begin{definition}\rm
Let $k$ be a field.
Pick $x \in k^{2} - \{ 0 \}$.
We denote 
by $\overline{x}$ the equivalence class containing $x$ in 
${\mathbb P}^{1}_{k} := (k^{2} - \{ 0 \}) / k^{\times}$.
For $A \in \mat{2}{k}$ we say that $\overline{x}$ is an 
{\it $A$-fixed point} 
if $x$ is an eigenvector of $A$.
In particular if $A \in \gl{2}{R}$, then $A$ can be regarded as 
an automorphism of ${\mathbb P}^{1}_{k}$, and so 
$\overline{x}$ is an $A$-fixed point if and only if 
$\overline{x}$ is fixed by $A$ as a point of 
${\mathbb P}^{1}_{k}$.
\end{definition}

\begin{remark}\rm
If $A \in \mat{2}{k}$ is not a scalar matrix, then $A$ has at most 
two fixed points in ${\mathbb P}^{1}_{k}$.
\end{remark}

\begin{prop}\label{prop:deltacondition}
Let $k$ be a field. Suppose that $A, B, C \in \mat2{k}$ and that
${\cal A}$ is a $k$-subalgebra of 
$\mat{2}{k}$ generated by $A, B, C$.
Then the following conditions are equivalent.
\begin{enumerate}
\item\label{item2-1} $\Delta(A, B)=\Delta(B, C)=\Delta(C, A)=\tau(A, B, C)=0$. 
\item\label{item2-2} ${\cal A} \neq \mat{2}{k}$, or equivalently
the ${\cal A}$-module $k^{2}$ is not absolutely irreducible.
\end{enumerate}
Furthermore, if $A, B, C \in {\rm GL}_2(k)$, then 
the following condition is also equivalent to the above two conditions.

\begin{enumerate} 
\setcounter{enumi}{2}
\item\label{item2-4} $\Delta(A, B)=\Delta(B, C)=\Delta(C, A)=
\Delta(AB, C)=0$.  
\end{enumerate}
\end{prop}

\prf
Since we may replace $k$ with an algebraic closure $\overline{k}$ of $k$,
we assume that $k$ is an algebraic closed field from the beginning. 
Note that ${\mathcal A} \neq {\rm M}_2(k)$ 
if and only if the ${\mathcal A}$-module $k^2$ is not irreducible 
when $k = \overline{k}$. 

(\ref{item2-2}) $\Rightarrow$ (\ref{item2-1})
The ${\cal A}$-module $k^{2}$ is not irreducible, so there exists
$P \in \gl{2}{k}$ such that 
\[P^{-1}AP=\left(
\begin{array}{cc}
 \ast & \ast \\
  0   & \ast
\end{array}
\right),
P^{-1}BP=\left(
\begin{array}{cc}
 \ast & \ast \\
  0   & \ast
\end{array}
\right),
P^{-1}CP=\left(
\begin{array}{cc}
 \ast & \ast \\
  0   & \ast
\end{array}
\right).\]
This immediately implies that $\tr(ABC)=\tr(ACB)$.
Hence we have $\tau(A, B, C)=0$. 
By Corollary \ref{coreq}, 
$\Delta(A, B)=\Delta(B, C)=\Delta(C, A)=0$. 

(\ref{item2-1}) $\Rightarrow$ (\ref{item2-2})
If one of $A, B, C$ is a scalar matrix, then 
(\ref{item2-1}) $\Rightarrow$ (\ref{item2-2})
follows from Proposition \ref{propdelta2}. 
Hence we may assume that none of $A, B, C$ are scalar matrices.
Note that if $X \in \mat{2}{k}$ is not a scalar matrix, then 
$X$ has at most two fixed points in ${\mathbb P}^1_{k}$. 
Suppose that one of $A, B, C$ has exactly one fixed point in ${\mathbb P}^1_{k}$.  
Then the others have the same fixed point in ${\mathbb P}^1_k$ 
since it follows from $\Delta(A, B)=\Delta(B, C)=\Delta(C, A)=0$ and  
Proposition \ref{propdelta2} that $k^2$ is not an 
irreducible module over the subalgebras generated by any two of $A, B, C$.
Hence the ${\mathcal A}$-module $k^2$ is not irreducible. 
 
Now let us consider the case that each of $A, B, C$ has exactly two 
fixed points in ${\mathbb P}^1_k$. 
If $A, B, C$ have a common fixed point, then 
we see that (\ref{item2-1}) $\Rightarrow$ (\ref{item2-2}). 
Suppose that $A, B, C$ have no common fixed point. 
By $\Delta(A, B)=\Delta(B, C)=\Delta(C, A)=0$, we may assume that 
$A$ has eigenvectors $u$ and $v$, $B$ has 
$v$ and $w$, and $C$ has $w$ and $u$.
With respect to the basis $\{u, v\}$, the matrices 
$A, B, C$ have the following forms: 
\[A=\left(
\begin{array}{cc}
 a_1 & 0 \\
 0   & a_4 
\end{array}
\right),
B=\left(
\begin{array}{cc}
 b_1 & 0 \\
 b_3 & b_4 
\end{array}
\right),
C=\left(
\begin{array}{cc}
 c_1 & c_2 \\
 0   & c_4 
\end{array}
\right).\] 
Since $A, B, C$ have no common eigenvector, 
$b_3 \neq 0$ and $c_2 \neq 0$.
Hence we have $\tr(ABC)=a_1b_1c_1+a_4b_3c_2+a_4b_4c_4$ and
$\tr(ACB)=a_1b_1c_1+a_1b_3c_2+a_4b_4c_4$.
Therefore
$\tau(A, B, C)=(a_4-a_1)b_3c_2 \neq 0$, since $A$ is not a scalar matrix.
This is a contradiction. Therefore we have shown that 
(\ref{item2-1}) $\Rightarrow$ (\ref{item2-2}).

(\ref{item2-2}) $\Rightarrow$ (\ref{item2-4}) This follows from 
Corollary \ref{coreq}.

(\ref{item2-4}) $\Rightarrow$ (\ref{item2-2})
In the same way as the discussion above in 
the (\ref{item2-1}) $\Rightarrow$ (\ref{item2-2}) part,
we only need to consider the case that $A, B, C$ have exactly 
two fixed points in ${\mathbb P}^1_k$. 
Suppose that $A, B, C$ have no fixed point. 
From the assumption that 
$\Delta(A, B)=\Delta(B, C)=\Delta(C, A)=0$ we may assume that 
$A$ has eigenvectors $u$ and $v$, B has $v$ and $w$, and 
$C$ has $w$ and $u$, where $u$, $v$, and $w$ are distinct 
up to scalar multiplication. 
The assumption that $\Delta(AB, C)=0$ implies that 
$AB$ and $C$ have a common eigenvector. 
It is $w$ or $u$.
If $w$ is an eigenvector of $AB$, then it is also an eigenvector of $A$ 
because $B \in {\rm GL}_2(k)$. 
This is a contradiction.
If $u$ is an eigenvector of $AB$, then it is also an 
eigenvector of $B$ because $A \in {\rm GL}_2(k)$. 
This is also a contradiction.
Hence the matrices $A, B, C$ have a common eigenvector, which implies that 
${\mathcal A} \neq {\rm M}_2(k)$. 
\qed

\bigskip

From the above proposition we obtain the 
following proposition.

\begin{proposition}
Let $\Gamma$ be a group with a subset $G = \{ \alpha_i \}_{i \in I}$
generating $\Gamma$.
Assume that the index set $I$ is a totally ordered set.
The air part ${\rm Rep}_2(\Gamma)_{\rm air}$ 
of the representation variety of degree $2$ for 
$\Gamma$ is equal to 
\[
\bigcup_{i < j} D(\Delta(\sigma_{\Gamma}(\alpha_i), \sigma_{\Gamma}(\alpha_j)))
\cup \bigcup_{i < j < k}D(\tau(\sigma_{\Gamma}(\alpha_i), 
\sigma_{\Gamma}(\alpha_j), \sigma_{\Gamma}(\alpha_k)))
\] 
or 
\[ 
\bigcup_{i < j} D(\Delta(\sigma_{\Gamma}(\alpha_i), \sigma_{\Gamma}(\alpha_j)))
\cup \bigcup_{i < j < k}D(\Delta(\sigma_{\Gamma}(\alpha_i\alpha_j), 
\sigma_{\Gamma}(\alpha_k))).  
\]
Here $\sigma_{\Gamma}$ is the universal representation and 
$D(\ast)$ is the open subset where $\ast$ does not vanish.
\end{proposition}

\prf
Let $k(x)$ be the residue field of a point $x$ of ${\rm Rep}_2(\Gamma)$.  
Note that $x \in {\rm Rep}_2(\Gamma)_{\rm air}$ if and only if 
$k(x)[\sigma_{\Gamma}(\Gamma)] = {\rm M}_2(k(x))$.  
By the Cayley-Hamilton theorem, we have $\rho(\alpha^{-1}) 
\in k(x)[\rho(\alpha)]$ for $\alpha \in \Gamma$. 
Since $\dim_{k(x)} {\rm M}_2(k(x)) = 4$, 
we see that $x \in {\rm Rep}_2(\Gamma)_{\rm air}$ if and only if 
there exist $\alpha_1, \alpha_2, \alpha_3 \in G$ which 
generate ${\rm M}_2(k(x))$ as a $k(x)$-algebra.  
Hence we can verify the statement by Proposition \ref{prop:deltacondition}.  
\qed 

\bigskip

In a similar way, we obtain the following proposition.

\begin{proposition}
Let $\Gamma$ be a monoid with a subset $G = \{ \alpha_i \}_{i \in I}$
generating $\Gamma$.
Assume that the index set $I$ is a totally ordered set.
The air part ${\rm Rep}_2(\Gamma)_{\rm air}$ 
of the representation variety of degree $2$ for 
$\Gamma$ is equal to 
\[
\bigcup_{i < j} D(\Delta(\sigma_{\Gamma}(\alpha_i), \sigma_{\Gamma}(\alpha_j)))
\cup \bigcup_{i < j < k}D(\tau(\sigma_{\Gamma}(\alpha_i), 
\sigma_{\Gamma}(\alpha_j), \sigma_{\Gamma}(\alpha_k))). 
\] 
\end{proposition}

\bibliographystyle{amsplain}


\end{document}